%% file: pap.tex
\documentclass{elsarticle}

\input{preamble.pap.tex}

\input{preamble.macros.tex}

\input{preamble.tikz.tex}

\usepackage{setspace}

\newbool{fastcompile}
\setbool{fastcompile}{false}

\newcommand{\mapped}[1]{{#1}_X}
\newcommand{\mappedTwo}[2]{{#1}_{X_{#2}}}
\newcommand{\map}{\Gcal}
\newcommand{\defgrad}{G}
\newcommand{\defjac}{g}
\newcommand{\dom}{\Omega}
\newcommand{\refdom}{\Omega_0}
\newcommand{\stvc}{U}
\newcommand{\flux}{F}
\newcommand{\src}{S}

\newcommand{\numflux}{\Hcal}

\newcommand{\msh}{\mathrm{msh}}
\newcommand{\err}{\mathrm{err}}
\newcommand{\adv}{\mathrm{adv}}
\newcommand{\up}{\mathrm{up}}

\begin{document}

\title{Implicit shock tracking using an optimization-based
       high-order discontinuous Galerkin method}

\author[rvt1]{M.~J.~Zahr\fnref{fn1}\corref{cor1}}
\ead{mzahr@nd.edu}

\author[rvt2]{A.~Shi\fnref{fn2}}
\ead{andrewshi94@berkeley.edu}

\author[rvt2,rvt3]{P.-O.~Persson\fnref{fn3}}
\ead{persson@berkeley.edu}

\address[rvt1]{Department of Aerospace and Mechanical Engineering, University
               of Notre Dame, Notre Dame, IN 46556, United States}
\address[rvt2]{Department of Mathematics, University of California, Berkeley,
                    Berkeley, CA 94720, United States}
\address[rvt3]{Mathematics Group, Lawrence Berkeley National Laboratory,
               1 Cyclotron Road, Berkeley, CA 94720, United States}
\cortext[cor1]{Corresponding author}
\fntext[fn1]{Assistant Professor, Department of Aerospace and Mechanical
             Engineering, University of Notre Dame}
\fntext[fn2]{Graduate student, Department of Mathematics, University of
             California, Berkeley}
\fntext[fn3]{Professor, Department of Mathematics, University of
             California, Berkeley}

\begin{keyword}
 shock tracking, %
 shock fitting, %
 $r$-adaptivity, %
 high-order methods, %
 discontinuous Galerkin, %
 high-speed flows
\end{keyword}


\begin{abstract}
A novel framework for resolving discontinuous solutions of conservation laws,
e.g., contact lines, shock waves, and interfaces, using
\textit{implicit tracking} and a high-order discontinuous
Galerkin (DG) discretization was introduced in \cite{zahr2018optimization}.
Central to the framework is an optimization problem whose solution is a
discontinuity-aligned mesh and the corresponding high-order approximation
to the flow that does not require explicit meshing of the unknown discontinuity
surface. The method was shown to deliver highly accurate solutions on coarse,
high-order discretizations without nonlinear stabilization and recover
optimal convergence rates  $\Ocal(h^{p+1})$ even for problems with
discontinuous solutions. This work extends the implicit tracking framework
such that robustness is improved and convergence accelerated. In particular,
we introduce an improved formulation of the central optimization problem and
an associated sequential quadratic programming (SQP) solver. The new
error-based objective function penalizes violation of the DG residual in
an \textit{enriched test space} and is shown to have excellent tracking
properties. The SQP solver simultaneously converges the nodal coordinates
of the mesh and DG solution to their optimal values and is equipped with a
number of features to ensure robust, fast convergence:
Levenberg-Marquardt approximation of the Hessian with weighted
elliptic regularization, backtracking line search based on the $\ell_1$ merit
function, and rigorous convergence criteria. We use the proposed method
to solve a range of inviscid conservation laws of varying difficulty.
We show the method is able to deliver accurate solutions on coarse,
high-order meshes and the SQP solver is robust and usually able to drive the
first-order optimality system to tight tolerances.
\end{abstract}

\maketitle

\section{Introduction}
\label{sec:intro}



High-order methods, such as the discontinuous Galerkin (DG) method
\cite{cockburn01rkdg,hesthaven2007nodal}, are widely believed to be superior
to traditional low-order schemes for simulation of turbulent flow
problems. However, in the presence of shocks and other discontinuities such as
contact lines or interfaces, the lack of nonlinear stability proves to be a
fundamental challenge. These features are ubiquitous in engineering and
science applications, particularly in high-speed flow problems. Many solutions
have been proposed, but most require excessively refined meshes to resolve all
the features, which in practices makes it difficult to accurately predict high
Reynolds, high Mach flows that feature shocks, boundary layers, and
interactions between them. Therefore, new advances are required to make
high-order schemes sufficiently robust and competitive for real-world problems.

Most of the techniques for addressing shocks are based on so-called
\emph{shock capturing}, that is, the numerical discretization somehow
incorporates the discontinuities independently of the computational grid.  One
simple method is to use a sensor that identifies the mesh elements that
contain shocks, and reduce their polynomial degrees
\cite{BauOden,burbeau01limiter}. For the DG method, this essentially leads to
a standard cell-centered finite volume scheme locally, which is well-known to
handle shocks robustly. Related, more sophisticated approaches include
limiting, such as the weighted essentially non-oscillatory (WENO) schemes
\cite{eno1,weno1,weno2}. For high-order methods, artificial viscosity has also
proven to be highly competitive, since it can smoothly resolve the jumps in
the solution without introducing additional discontinuities between the
elements \cite{persson06shock}. The main problem with all these approaches is
that they reduce to first order accuracy in the affected elements, which
translates into a globally first order accurate scheme. This can be remedied
by using local mesh refinement around the shock ($h$-adaptivity)
\cite{dervieux03adaptation}, although the anisotropic elements that are
required for efficiency are difficult to generate and excessively fine
elements are needed around the shock.

An alternative approach is \emph{shock tracking} or \emph{shock fitting},
where the computational mesh is moved such that their faces are aligned with
the discontinuities in the solution \cite{shubin1981steady, shubin1982steady,
  bell1982fully, harten1983self, rosendale1994floating,
  trepanier1996conservative, zhong98tracking, taghaddosi1999adaptive,
  baines2002multidimensional, roe2002adaptive, glimm2003conservative,
  palaniappan2008sub}. This is very natural in the setting of a DG method
since the numerical scheme already incorporates jumps between the elements
and the approximate Riemann solvers employed on the element faces handle the
discontinuities correctly. However, it is a difficult meshing problem since it
essentially requires generating a fitted mesh to the (unknown) shock
surface. Also, in the early approaches to shock fitting, it was applied to
low-order schemes where the relative advantage over shock capturing is smaller
than for high-order methods. For these reasons, shock tracking is largely not
used in practical CFD today.

In \cite{zahr2018optimization}, we introduced a novel approach to
shock tracking that does not require explicitly generating a mesh of the
unknown discontinuity surface. Rather, the conservation law is discretized
on a mesh without knowledge of the discontinuity surface and an optimization
problem is formulated such that its solution is a mesh that aligns with
discontinuities in the flow and the corresponding solution of the discrete
conservation law. That is, tracking of the discontinuities is
implicitly defined through the solution of the optimization problem and
will be referred to as \textit{implicit shock tracking}.
While this approach works with any discretization that allows for
inter-element discontinuities, we focus on high-order DG
methods due to the high degree of accuracy attainable on coarse meshes,
proper treatment of discontinuities with approximate Riemann solvers, and
the ability to used curved elements to track discontinuities with curvature.
The optimization problem is solved by simultaneously converging the mesh
and solution to their optimal values, which never requires the fully
converged DG solution on non-aligned meshes and does not require nonlinear
stabilization. The combination of implicit tracking with a DG discretization
is truly high-order accurate, since the solution is smooth within each
element, and very accurate solutions can be obtained on coarse meshes.

This paper extends our prior work with a new error-based objective function,
a sequential quadratic programming solver for the optimization problem, and
a number of practical considerations. The proposed objective function
penalizes violation of the DG residual in an enriched test space, which
is a surrogate for violation of the weak formulation of the conservation
law. Even though a traditional DG solution will oscillate about
discontinuities in the solution on a non-aligned mesh, this violates
the true conservation law; as a result, the proposed objective function
promotes alignment of the mesh with discontinuities. This formulation has
the added benefit of $r$-adaptive behavior; even in smooth regions of the
flow, nodes will adjust to improve the approximation of the conservation law.
The other main contribution of this work is an SQP solver for the optimization
problem that leverages its structure. Due to the minimum-residual
structure of the objective function, we employ a Levenberg-Marquardt
approximation of the Hessian. We propose to use the stiffness matrix of
a linear elliptic partial differential equation (PDE) with the coefficient
chosen inversely proportional to the local element size as the regularization
matrix. This tends to smooth out the search directions for the mesh coordinates
and is particularly important for problems with elements of significantly
different size. The SQP method is globalized with a line search based on
the $\ell_1$ merit function and equipped with termination conditions based
on the first-order optimality criteria. For the method to be practical
for difficult problems, we initialize the solve with the $p = 0$ DG solution
and use continuation in the polynomial degree (solution and mesh) for
high-order ($p > 1$) discretizations. Finally, we identified the critical role
of smoothness of the DG numerical flux function with respect to variations
in the element normal in the context of implicit shock tracking, mainly
with respect to solver convergence, and use smoothed versions of traditional
numerical fluxes throughout.


To our knowledge, the only other approach to implicit shock tracking was
proposed in \cite{corrigan2018application, corrigan2019convergence,
                  corrigan2019moving, corrigan2019unsteady, kercher2020moving}.
where the authors enforce a DG
discretization with unconventional numerical fluxes and the Rankine-Hugoniot
interface conditions in a minimum-residual sense. Interestingly, enforcement
of the interface condition circumvented traditional stability requirements
for the DG numerical fluxes, allowing them to solely rely on fluxes interior
to an element. Their method was shown to successfully track even complex
discontinuity surfaces and provide accurate approximations to the conservation
law on traditionally coarse, high-order meshes. The present work incorporates
some aspects of their method, in particular topological mesh operations
and some aspects of the Hessian approximation and regularization. However, our
approach that directly enforces a stable DG discretization (with zero residual)
inherits many attractive features of DG methods such as guaranteed
conservation and a rigorous framework for high-order convergence. Furthermore,
by choosing an error-based objective function rather than a physics-based one,
the extension to viscous problems only requires treatment of second-order
terms in the DG setting, which has been well-established
\cite{arnold2002unified}.

 

The remainder of the paper is organized as follows. Section~\ref{sec:disc}
introduces the governing system of inviscid conservation laws and 
its discretization using a discontinuous Galerkin method.
Section~\ref{sec:optim} recalls the implicit tracking framework originally
proposed in \cite{zahr2018optimization}, introduces the new error-based
objective function, and discusses a parametrization of the mesh deformation
that ensures the boundaries of the computational domain remain on the
boundaries of the actual domain. Section~\ref{sec:solver} introduces the
proposed SQP solver for the central optimization problem that incorporates
a Levenberg-Marquardt approximation of the Hessian with a novel
regularization matrix, a line search based on a $\ell_1$ merit function,
and termination criteria based on the first-order optimality conditions.
Section~\ref{sec:practical} discusses two important details required to
make the proposed tracking framework work in practice: initialization
of the SQP solver and topological mesh operations to remove small
elements. Finally, Section~\ref{sec:num-exp} presents a number of numerical
experiments that demonstrate the method is able to accurately approximate
complex flows using coarse, high-order meshes and the SQP solver is able
to quickly converge to a mesh that tracks all discontinuities and
exhibits deep convergence to the first-order optimality conditions.

\section{Governing equations and high-order numerical discretization}
\label{sec:disc}
Consider a general system of $M$ inviscid conservation laws, defined on the
physical domain $\dom \subset \Rbb^d$ and subject to appropriate boundary
conditions,
\begin{equation} \label{eqn:claw-phys}
 \nabla\cdot \flux(\stvc) = \src(\stvc) \quad \text{in}~~\dom
\end{equation}
where $\func{\stvc}{\dom}{\Rbb^M}$ is the solution of the system of
conservation laws, $\func{\flux}{\Rbb^M}{\Rbb^{M\times d}}$ is the physical
flux, $\func{\src}{\Rbb^M}{\Rbb^M}$ is the source term, and
$\ds{\nabla \coloneqq (\partial_{x_1},\dots,\partial_{x_d})}$
is the gradient operator in the physical domain such that
$\nabla w(x) =
 \begin{bmatrix} \partial_{x_1} w(x) & \cdots &
                 \partial_{x_d} w(x)
 \end{bmatrix} \in \Rbb^{N \times d}$
for any $N$ vector-valued function $w$ over $\Omega$
($w(x) \in \Rbb^N$ for $x \in \Omega$). The boundary of
the domain is $\partial\Omega$ with outward unit normal
$\func{n}{\partial\dom}{\Rbb^d}$.
The formulation of the conservation law in
(\ref{eqn:claw-phys}) is sufficiently general to encapsulate
steady conservation laws in a $d$-dimensional spatial domain or
time-dependent conservation laws in a $(d-1)$-dimensional domain, i.e.,
a $d$-dimensional space-time domain. In general, the solution $\stvc(x)$ may
contain discontinuities, in which case, the conservation law
(\ref{eqn:claw-phys}) holds away from the discontinuities
and the Rankine-Hugoniot conditions \cite{majda2012compressible}
\begin{equation}
 F(U^+)n=F(U^-)n
\end{equation}
hold for $x\in\Gamma_s$, where $\Gamma_s\subset\Omega$ is a surface along
which $U$ is discontinuous, $U^+(x), U^-(x)\in\Rbb^M$ are the values of
$U(x)$ on either side of the discontinuity, and $n$ is a normal vector
to the surface $\Gamma_s$.

Similar to our previous work \cite{zahr2018optimization},
we will construct a numerical method that directly tracks discontinuities
with the computational grid, which places three requirements on the
discretization:
\begin{inparaenum}[1)]
 \item represents a stable and convergent discretization of the conservation
       law in (\ref{eqn:claw-phys}),
 \item allows for deformation of the computational domain, and
 \item employs a solution basis that supports discontinuities between
       computational cells or elements.
\end{inparaenum}
To achieve high-order accuracy, the tracking framework introduced in
\cite{zahr2018optimization} is built upon a standard high-order DG method
given their proven ability \cite{zahr2018optimization, corrigan2019moving} to
deliver accurate solutions on very coarse discretizations provided
discontinuities are tracked. The extension to other discretizations that
support inter-element discontinuities such as finite volumes, flux
reconstruction, hybridizable DG, and the DG spectral element method are
possible, but beyond the scope of this work.

The remainder of this section will detail the discretization of the
conservation law (\ref{eqn:claw-phys}) using DG such that it reduces
to the discrete form: given $\xbm \in \Rbb^{N_\xbm}$, find
$\ubm \in \Rbb^{N_\ubm}$
\begin{equation} \label{eqn:claw-disc}
 \rbm(\ubm,\,\xbm) = \zerobold
\end{equation}
where $\ubm$ is the discrete representation of the conservation law state
$U$, $\xbm$ is the discrete representation of the conservation law domain
$\Omega$ (nodal coordinates of mesh nodes), and
$\func{\rbm}{\Rbb^{N_\ubm}\times\Rbb^{N_\xbm}}{\Rbb^{N_\ubm}}$ is
the discretized conservation law. The same discretization will be used
to define a residual function based on an enriched \textit{test} space
$\func{\Rbm}{\Rbb^{N_\ubm}\times\Rbb^{N_\xbm}}{\Rbb^{N_\ubm}}$,
which will be used in the new definition of the proposed objective function.

\subsection{Transformed conservation law from deformation of physical domain}
\label{sec:disc:transf}
Before introducing a discretization of (\ref{eqn:claw-phys}) it is convenient
to explicitly treat deformations to the domain of the conservation law $\Omega$,
which will eventually be induced by deformation to the mesh as nodal
coordinates are moved to track discontinuities,
by transforming to a fixed reference domain $\Omega_0 \subset \Rbb^d$.
Suppose the physical domain can be taken as the result of a
diffeomorphism applied to a reference domain (Figure~\ref{fig:dom-map})
\begin{equation} \label{eqn:dom-map}
 \dom = \map(\refdom),
\end{equation}
where $\refdom \subset \Rbb^d$ is a fixed reference domain and
$\func{\map}{\Rbb^d}{\Rbb^d}$ is the diffeomorphism defining the
domain mapping.
\ifbool{fastcompile}{}{
\begin{figure}[H]
  \centering
  \includegraphics[width=2.5in]{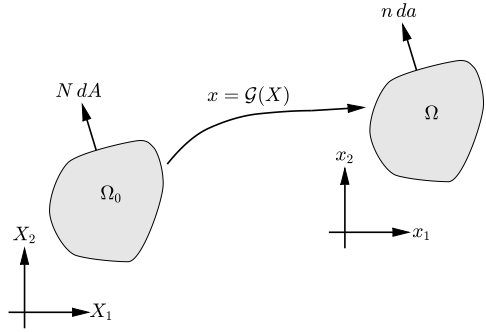}
  \caption{Mapping between reference and physical domains.}
  \label{fig:dom-map}
\end{figure}
}
For convenience, the conservation law on the physical domain $\dom$ is
transformed to a conservation law on the reference domain
\begin{equation} \label{eqn:claw-ref}
 \nabla_X \cdot \mapped\flux(\mapped\stvc, \map) =
 \mapped\src(\mapped\stvc,\map) \quad \text{in}~~\refdom
\end{equation}
where $\nabla_X \coloneqq (\partial_{X_1},\dots,\partial_{X_d})$ denotes
spatial derivatives with respect to the reference domain $\Omega_0$ with
coordinates $X$, $\func{\mapped\stvc}{\refdom}{\Rbb^M}$ is the mapped
state vector we define as
\begin{equation}
 \mapped\stvc = \stvc \circ \map,
\end{equation}
$\mapped\flux(\mapped\stvc,\map) \in \Rbb^{M\times d}$ is the
transformed flux function, and
$\mapped\src(\mapped\stvc,\map) \in \Rbb^M$ is the transformed source term.
The unit normal in the reference and physical domain are related by
\begin{equation} \label{eqn:transf-normal}
 n = \frac{\defjac\defgrad^{-T}N}{\norm{\defjac\defgrad^{-T}N}},
\end{equation}
where $G(X) = \pder{}{X}\Gcal(X)$ is the deformation gradient
of the domain mapping and $g(X) = \det G(X)$ is the
Jacobian.

For the transformed conservation
law (\ref{eqn:claw-ref}) and original conservation law (\ref{eqn:claw-phys})
to be equivalent, we require
\begin{equation} \label{eqn:claw-transf-req}
 \int_V (\nabla_X \cdot \mapped\flux-\mapped\src) \, dV =
 \int_{\map(V)} (\nabla \cdot \flux-\src) \, dv
\end{equation}
holds over an arbitrary volume $V \subset \refdom$ and arguments have been
dropped for convenience. Considering the source term and flux individually,
\begin{equation} \label{eqn:claw-transf-src}
 \int_{\map(V)} S \, dv = \int_V \defjac S \, dV
\end{equation}
follows directly from a change of variables in the integral and
\begin{equation} \label{eqn:claw-transf-flux}
 \int_{\map(V)} \nabla\cdot\flux \, dv =
 \int_{\partial \map(V)}  \flux \cdot n \, ds =
 \int_{\partial V} \defjac \flux \cdot \defgrad^{-T} N \, dS =
 \int_V \mapped\nabla \cdot (\defjac \flux \cdot \defgrad^{-T}) \, dS.
\end{equation}
follows from the divergence theorem in the physical domain,
change of variables for surface integrals (Nanson's formula),
and the divergence theorem in the reference domain.
Combining (\ref{eqn:claw-transf-req})-(\ref{eqn:claw-transf-flux}) and
invoking arbitrariness of the volume $V$, the transformed flux and source
term take the following form
\begin{equation} \label{eqn:claw-transf-srcflux}
 \mapped\flux(\mapped\stvc,\map) = \defjac\flux(\stvc)\defgrad^{-T}, \qquad
 \mapped\src(\mapped\stvc,\map) = \defjac\src(\stvc).
\end{equation}

\subsection{Discontinuous Galerkin discretization of
            transformed conservation law}
\label{sec:disc:dg}
We use a standard nodal discontinuous Galerkin method
\cite{cockburn01rkdg,hesthaven08dgbook} to discretize
the transformed conservation law (\ref{eqn:claw-ref}).
Let $\Ecal_{h,q}$ represent a discretization of the reference domain
$\refdom$ into non-overlapping, potentially curved, computational elements,
where $h$ is a mesh element size parameter and $q$ is the polynomial order
associated with the curved elements. The DG construction begins with the
elementwise weak form of the conservation
law (\ref{eqn:claw-ref}) that results from multiplying each equation by a test
function $\mapped\psi$, integrating over a single element
$K \in \Ecal_{h,q}$, and applying the divergence theorem
\begin{equation} \label{eqn:claw-weak-elem}
 \int_{\partial K} \mapped\psi^+ \cdot \mapped\flux(\mapped\stvc,\map) N \, dS -
 \int_K F(\mapped\stvc,\map):\mapped\nabla \mapped\psi \, dV = 0,
\end{equation}
where $N$ is the outward normal to the surface $\partial K$ and
$\mapped\psi^+$ denotes the trace of $\psi$ \textit{interior} to
element $K$. To ensure the face integrals are single-valued, we replace
$\mapped\flux(\mapped\stvc,\map) N$ in the first term with a numerical flux
function $\mapped\numflux(\mapped\stvc^+,\mapped\stvc^-,N,\map)$
\begin{equation} \label{eqn:claw-weak-numflux}
 \int_{\partial K} \mapped\psi^+ \cdot
   \mapped\numflux(\mapped\stvc^+,\mapped\stvc^-,N,\map) \, dS -
 \int_K F(\mapped\stvc,\map):\mapped\nabla \mapped\psi \, dV = 0,
\end{equation}
where $\mapped\stvc^+$ denotes the interior trace of $\mapped\stvc$
and $\mapped\stvc^-$ denotes the exterior trace of $\mapped\stvc$
if $\partial K$ is an interior face, i.e.,
$\partial K \cap \partial\refdom = \emptyset$, otherwise
$\mapped\stvc^-$ is a boundary state $\mapped\stvc^\partial(\mapped\stvc,n)$
constructed to enforce the appropriate boundary condition.
We defer a detailed discussion of numerical flux functions to
Section~\ref{sec:disc:numflux} and boundary conditions to the
specific conservation laws considered in Section~\ref{sec:num-exp}.

To establish the finite-dimensional form of (\ref{eqn:claw-weak-numflux}),
we introduce the mapped finite element space of piecewise polynomial
functions associated with the mesh $\Ecal_{h,q}$:
\begin{equation*}
 \Vcal_{h,p} = \left\{v \in [L^2(\refdom)]^M \mid
         \left.v\right|_K \circ \Tcal_K \in [\Pcal_p(K_0)]^M
         ~\forall K \in \Ecal_{h,q}\right\}
\end{equation*}
where $\Pcal_p(K_0)$ is the space of polynomial functions of degree at most
$p \geq 1$ on the parent element $K_0$ and $K = \Tcal_K(K_0)$ defines a
mapping from the parent element to element $K \in \Ecal_{h,q}$. For notational
brevity, we assume all elements map from a single parent element. The
finite-dimensional residual of the weak form in (\ref{eqn:claw-weak-numflux})
corresponding to the trial space $\Vcal_{h,p}$ and test space
$\Vcal_{h',p'}$ is
\begin{equation} \label{eqn:claw-findim-dg-pg}
 \mappedTwo{r}{h',p'}^K(\mappedTwo{\stvc}{h,p},\,\map) \coloneqq
  \int_{\partial K}
   \mappedTwo{\psi}{h',p'}^+ \cdot
   \mapped\numflux(\mappedTwo{\stvc}{h,p}^+,
                   \mappedTwo{\stvc}{h,p}^-,N,\map) \, dS -
  \int_K \mapped\flux(\mappedTwo{\stvc}{h,p},\map) :
         \mapped\nabla\mappedTwo{\psi}{h',p'} \, dV.
\end{equation}
Finally, define the continuous global function space
\begin{equation} \label{eqn:gfcnsp}
  \Wcal_{h,q} \coloneqq \Vcal_{h,q} \cap C^0(\Omega),
\end{equation}
and constrain domain deformation to lie in it:
$\map \approx \map_{h,q} \in \Wcal_{h,q}$. After summing over
all elements $K \in \Ecal_{h,q}$, the final version of the
finite-dimensional Galerkin weak form is: given
$\map_{h,q} \in \Wcal_{h,q}$, find $\mappedTwo{\stvc}{h,p} \in \Vcal_{h,p}$
such that
\begin{equation} \label{eqn:claw-findim-dg}
 \sum_{K \in \Ecal_{h,q}} \mappedTwo{r}{h,p}^K(\mappedTwo{\stvc}{h,p},
                                               \map_{h,q}) = 0
\end{equation}
for all $\mappedTwo{\psi}{h,p} \in \Vcal_{h,p}$.

To establish the discrete (algebraic) form of (\ref{eqn:claw-findim-dg}),
we introduce a (nodal) basis over each element and expand the
finite-dimensional test functions ($\mappedTwo{\psi}{h,p}$),
solution ($\mappedTwo{\stvc}{h,p}$), and domain deformation
($\map_{h,q}$) in terms of these basis functions and coefficients.
Invoking arbitrariness of the test functions in
$\Vcal_{h,p}$ and assembling an algebraic system that respects the global
functions spaces in (\ref{eqn:gfcnsp}), we obtain
\begin{equation} \label{eqn:claw-disc}
 \rbm(\ubm,\,\xbm) = \zerobold,
\end{equation}
where $\ubm \in \Rbb^{N_\ubm}$ are the (assembled) coefficients of the
solution $\mappedTwo{\stvc}{h,p} \in \Vcal_{h,p}$ and $\xbm \in \Rbb^{N_\xbm}$
are the (assembled) coefficients of the domain deformation
$\map_{h,q} \in \Wcal_{h,q}$. Since we are using nodal bases,
the entries of $\xbm$ are the \textit{coordinates} of the nodes
of the mesh and $\ubm$ contains the components of the solution at
the nodes.

To close this section, we introduce an \textit{enriched} discrete
residual $\Rbm(\ubm, \xbm)$ that will be used to define the proposed
shock tracking objective function. Let $\Rbm(\ubm,\xbm)$
be the algebraic version of
\begin{equation}
 \sum_{K\in\Ecal_{h,q}} \mappedTwo{r}{h',p'}^K(\mappedTwo{\stvc}{h,p},
                                               \map_{h,q}),
\end{equation}
where $h' \leq h$ and $p' \geq p$.
That is, the enriched residual is defined by the same trial space
($\Vcal_{h,p}$) and space for the domain deformation
($\Wcal_{h,q}$) as the residual in (\ref{eqn:claw-disc}), but
uses an enriched test space $\Vcal_{h',p'}$. In this work,
we take $h' = h$ and $p' = p+1$ although other choices are possible
and will be explored in future work.

\subsection{Numerical flux function}
\label{sec:disc:numflux}


Recall the numerical flux is a quantity that replaces the flux dotted with
the outward unit normal and expect it to transform according to
\begin{equation} \label{eqn:claw-transf-numflux}
 \mapped\numflux(\mapped\stvc^+, \mapped\stvc^-, N, \map) =
 \norm{\defjac\defgrad^{-T}N} \numflux(\stvc^+, \stvc^-, n),
\end{equation}
which follows from (\ref{eqn:transf-normal}) and (\ref{eqn:claw-transf-srcflux})
as
\begin{equation} \label{eqn:claw-transf-flux-dot-normal}
 \mapped\numflux \sim
 \mapped\flux \cdot N =
 \defjac\flux\cdot\defgrad^{-T}N =
 \norm{\defjac\defgrad^{-T}N} \flux \cdot n \sim
 \norm{\defjac\defgrad^{-T}N} \numflux,
\end{equation}
where arguments have been dropped for brevity.
Therefore, the mapped numerical flux is uniquely determined from the
physical numerical flux. For the remainder of this section, we will discuss
the numerical flux in the physical domain $\numflux(\stvc^+, \stvc^-, n)$
and it will be transformed to the reference domain according to
(\ref{eqn:claw-transf-numflux}).

For DG methods to be stable, the numerical flux may be any two-point monotone
Lipschitz function that is
\begin{enumerate}[(i)]
 \item \textit{consistent} with the flux function, i.e., for any
       $\stvc \in \Rbb^M$
\begin{equation}
 \numflux(\stvc,\stvc,n) = \flux(\stvc)\cdot n
\end{equation}
 \item \textit{conservative}, i.e., for any $\stvc, \stvc' \in \Rbb^M$
\begin{equation}
 \numflux(\stvc, \stvc', n) = -\numflux(\stvc', \stvc, -n).
\end{equation}
\end{enumerate}
Conditions (i)-(ii) are satisfied by all standard numerical fluxes as
these are the minimum requirements for a stable and accurate DG method.
However, it was observed in \cite{zahr2018optimization,corrigan2019moving}
that the requirements are higher for tracking-based discretizations since
inter-element jumps do not tend to zero under refinement. As introduced in
\cite{corrigan2019moving}, this requires the additional condition:
\begin{enumerate}[(i)]
 \addtocounter{enumi}{2}
 \item \textit{preservation} of the Rankine-Hugoniot conditions, i.e.,
 given $\stvc^+, \stvc^- \in \Rbb^M$ such that
 $\flux(\stvc^+) \cdot n = \flux(\stvc^-) \cdot n$, then
\begin{equation} \label{eqn:numflux-rh-cond}
 \numflux(\stvc^+, \stvc^-, n) = \flux(\stvc^+)\cdot n = \flux(\stvc^-)\cdot n
\end{equation}
\end{enumerate}
Conditions (i)-(iii) are satisfied by exact Riemann solvers and a number of
approximate Riemann solvers including the Roe \cite{roe1981approximate} and
HLLC \cite{toro1994restoration} fluxes; however, many popular
numerical flux functions such as the local Lax-Friedrichs, Rusanov, and Roe
with entropy fix do not satisfy (iii) \cite{corrigan2019moving}.
We highlight one additional desired
property of the numerical flux function that we will show is important
for optimization-based discontinuity tracking:
\begin{enumerate}[(i)]
 \addtocounter{enumi}{3}
 \item \textit{smoothness} with respect to variations in the normal.
\end{enumerate}
Since the numerical flux must have upwind-like properties, it is
difficult to construct numerical fluxes to satisfy (i)-(iv), particularly
the smoothness property. Therefore, we choose numerical fluxes that
satisfy (i)-(iii) and replace non-smooth terms with smooth approximations
to recover smoothness (iv). Strictly speaking, this will not preserve the
Rankine-Hugoniot conditions (iii); in Section~\ref{sec:num-exp} we
demonstrate this is a desirable trade-off.

As an example, consider linear advection of a scalar field 
$\func{u}{\dom}{\Rbb}$ in a spatially varying direction
$\func{\beta}{\dom}{\Rbb^d}$ governed by conservation law
of the form (\ref{eqn:claw-phys}) (see (\ref{eqn:advec}) in
Section~\ref{sec:num-exp:advec}) with flux function
\begin{equation}
 F_\adv(u; \beta) = u \beta^T
\end{equation}
and the upwind numerical flux
\begin{equation} \label{eqn:upwind}
 \Hcal_\up(u^+, u^-, n) =
 \begin{cases}
  (\beta\cdot n) u^+ & \text{ if } \beta\cdot n \geq 0 \\
  (\beta\cdot n) u^- & \text{ if } \beta\cdot n < 0,
 \end{cases}
\end{equation}
which can equivalently be written in terms of the Heaviside function
$\func{H}{\Rbb}{\{0, 1\}}$ as
\begin{equation}
 \Hcal_\up(u^+, u^-, n) =
 (\beta\cdot n)
 \left[u^+ \cdot H(\beta\cdot n) + u^- \cdot (1-H(\beta\cdot n))\right].
\end{equation}
The upwind flux satisfies conditions (i)-(iii) for an admissible flux
for shock tracking, but fails to satisfy condition (iv). Condition (i)
follows from
\begin{equation}
 \Hcal_\up(u, u, n) = (\beta\cdot n) u = F_\adv(u)n
\end{equation}
and condition (ii) follows from
\begin{equation}
 -\Hcal(u', u, -n) = 
 (\beta\cdot n)
 \left[u' \cdot H(-\beta\cdot n) + u \cdot (1-H(-\beta\cdot n))\right] =
  \Hcal(u, u', n),
\end{equation}
for any $u, u' \in \Rbb$, where the second equality follows from the
property of the Heaviside function: $H(-s) = 1-H(s)$ for $s \in \Rbb$.
To verify condition (iii), consider $u^+, u^- \in \Rbb$ such that
$u^+ \neq u^-$ (discontinuity) and assume
$F_\adv(u^+)\cdot n = F_\adv(u^-)\cdot n$. This condition
implies $(u^+-u^-)(\beta\cdot n) = 0$, which in turn implies
$\beta\cdot n = 0$ from the assumption that $u^+\neq u^-$. Therefore
condition (iii) holds from
\begin{equation}
 \Hcal(u^+,u^-, n) = F_\adv(u^+) n = F_\adv(u^-)n = 0.
\end{equation}
Finally, it is easy to see that condition (iv) is not satisfied; the
upwind flux is continuous with respect to variations in the normal $n$,
but not smooth due to the $H(\beta \cdot n)$ terms.
A single isolated point of non-smoothness does not necessarily hinder
the optimization solver because it is unlikely to be visited during
the solution procedure. However, in this case, the kinks in the
numerical flux function lie at points where $\beta\cdot n = 0$,
which is precisely the requirement for a discontinuous solution to
satisfy the Rankine-Hugoniot conditions and will certainly be
approached as the mesh faces align with the discontinuity.

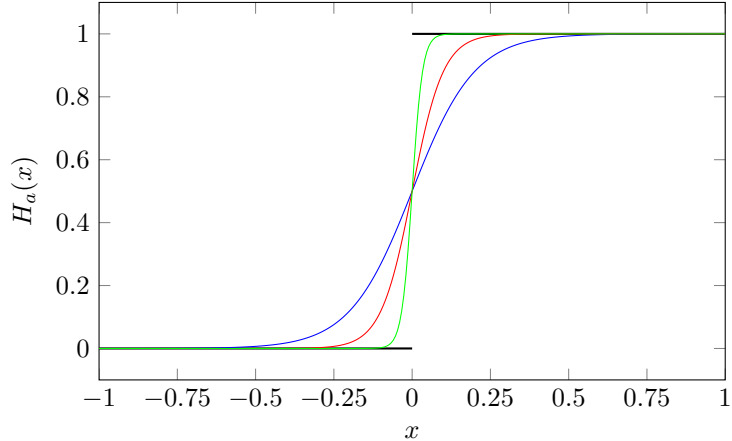
\begin{figure}
 \centering
 \input{py/smooth_heaviside.tikz}
 \caption{Smoothed Heaviside (logistic) function for
          $a =  5$ (\ref{line:smheavi5}),
          $a = 10$ (\ref{line:smheavi10}),
          $a = 30$ (\ref{line:smheavi30}),
          $a = \infty$ (\ref{line:heavi}).}
 \label{fig:smooth_heaviside}
\end{figure}

To recover condition (iv), we introduce a smoothed version of the upwind
flux where the Heaviside function is replaced with a smoothed step function
$\func{H_a}{\Rbb}{\Rbb}$
\begin{equation}
 \Hcal_\up^a(U^+, U^-, n) =
 (\beta\cdot n)
 \left[U^+ \cdot H_a(\beta\cdot n) + U^- \cdot (1-H_a(\beta\cdot n))\right].
\end{equation}
In this work, we use the logistic function as the smoothed step function
\begin{equation}
 H_a(x) \coloneqq \frac{1}{1+e^{-2 a x}},
\end{equation}
where $a \in \Rbb$ is the smoothing parameter
(Figure~\ref{fig:smooth_heaviside}).
Section~\ref{sec:num-exp} provides a detailed study of the impact
of the smoothness of the numerical flux on the convergence of the
tracking algorithm.

%



\section{Optimization formulation of $r$-adaptivity for implicit
         tracking of discontinuities}
\label{sec:optim}
In this section, we introduce the main contribution of this work:
an $r$-adaptivity framework that recasts the discrete conservation law
(\ref{eqn:claw-disc}) as an optimization problem
over the discrete solution and mesh that aims to align features
in the solution basis with features in the solution itself. In this
work these features are discontinuities since we only consider
inviscid conservation laws; however, future work will consider
steep gradients (viscous conservation laws) and interfaces.
In the present setting, this amounts to aligning element
faces with discontinuities. The method builds
upon our previous work \cite{zahr2018optimization}, where we demonstrated
that high-order methods are capable of approximating discontinuous
solutions of PDEs using extremely coarse discretizations provided
the discontinuities are tracked. In this section, we focus on the
new aspects of the optimization formulation and the next section will
address the issue of solvers for the optimization problem. In particular,
we introduce a new error-like objective function based on the DG
residual using an enriched \textit{test space} and a term to
penalize mesh distortion.

\subsection{Constrained optimization formulation}
\label{sec:optim:pdeconstr}
Following our work in \cite{zahr2018optimization}, we formulate the
problem of tracking discontinuities as a constrained optimization
problem over the PDE state and coordinates of the mesh nodes that
minimizes some objective function
$\func{f}{\Rbb^{N_\ubm}\times\Rbb^{N_\xbm}}{\Rbb}$
while enforcing the DG discretization of the conservation law
\begin{equation} \label{eqn:pde-opt}
 \optconOne{\ubm\in\Rbb^{N_\ubm},\xbm\in\Rbb^{N_\xbm}}
           {f(\ubm, \xbm)}{\rbm(\ubm, \xbm) = \zerobold.}
\end{equation}
The objective function is constructed such that the solution of
the optimization problem is a mesh that aligns
with discontinuities in the solution. The optimization-based tracking
method directly inherits the benefits of standard DG methods, i.e., high-order
accuracy and conservation, due to the constraint that exactly enforces the
DG discretization.
Finally, the optimization formulation in (\ref{eqn:pde-opt}) will
provide nonlinear stability if all discontinuities are successfully tracked,
which we will demonstrate using several examples in Section~\ref{sec:num-exp}.

The Lagrangian of the optimization problem in (\ref{eqn:pde-opt})
$\func{\Lcal}{\Rbb^{N_\ubm}\times\Rbb^{N_\xbm}\times\Rbb^{N_\ubm}}{\Rbb}$
takes the form
\begin{equation} \label{eqn:lagr}
 \Lcal(\ubm,\xbm,\lambdabold) = f(\ubm, \xbm)-\lambdabold^T\rbm(\ubm,\xbm),
\end{equation}
where $\lambdabold \in \Rbb^{N_\ubm}$ is a vector of Lagrange multipliers
associated with the DG constraint in (\ref{eqn:pde-opt}). The first-order
optimality, or Karush-Kuhn-Tucker (KKT), conditions state that the
$(\ubm^\star, \xbm^\star)$ is a first-order solution of the optimization
problem if there exists $\lambdabold^\star$ such that
\begin{equation} \label{eqn:kkt0}
 \nabla_\ubm \Lcal(\ubm^\star,\xbm^\star,\lambdabold^\star) = \zerobold, \qquad
 \nabla_\xbm \Lcal(\ubm^\star,\xbm^\star,\lambdabold^\star) = \zerobold, \qquad
 \nabla_\lambdabold\Lcal(\ubm^\star,\xbm^\star,\lambdabold^\star) = \zerobold,
\end{equation}
or equivalently,
\begin{equation} \label{eqn:kkt1}
 \pder{f}{\ubm}(\ubm^\star,\xbm^\star)^T -
 \pder{\rbm}{\ubm}(\ubm^\star,\xbm^\star)^T\lambdabold^\star = \zerobold,
 \qquad
 \pder{f}{\xbm}(\ubm^\star,\xbm^\star)^T -
 \pder{\rbm}{\xbm}(\ubm^\star,\xbm^\star)^T\lambdabold^\star = \zerobold,
 \qquad
\rbm(\ubm^\star, \xbm^\star) = \zerobold.
\end{equation}
Since the DG Jacobian with respect to the state variables $\ubm$ is assumed
to be invertible, we define the estimate of the optimal Lagrange multiplier
$\func{\hat\lambdabold}{\Rbb^{N_\ubm}\times \Rbb^{N_\xbm}}{\Rbb^{N_\ubm}}$
such that the first equation ($\nabla_\ubm \Lcal = 0$) (adjoint equation) is
always satisfied
\begin{equation} \label{eqn:lagrmult}
 \hat\lambdabold(\ubm,\xbm) =
 \pder{\rbm}{\ubm}(\ubm,\xbm)^{-T}\pder{f}{\ubm}(\ubm,\xbm)^T.
\end{equation}
Then the optimality criteria becomes
\begin{equation} \label{eqn:kkt2}
 \cbm(\ubm^\star, \xbm^\star) = \zerobold, \qquad
 \rbm(\ubm^\star, \xbm^\star) = \zerobold,
\end{equation}
where $\func{\cbm}{\Rbb^{N_\ubm}\times\Rbb^{N_\xbm}}{\Rbb^{N_\ubm}}$ is
defined as
\begin{equation}
 \cbm(\ubm,\xbm) \coloneqq
 \nabla_\xbm\Lcal(\ubm,\xbm,\hat\lambdabold(\ubm,\xbm)) = 
 \pder{f}{\xbm}(\ubm,\xbm)^T -
 \pder{\rbm}{\xbm}(\ubm,\xbm)^T\pder{\rbm}{\ubm}(\ubm,\xbm)^{-T}
 \pder{f}{\ubm}(\ubm,\xbm)^T.
\end{equation}
In Section~\ref{sec:solver:term}, $\norm{\cbm(\ubm,\xbm)}$ and
$\norm{\rbm(\ubm,\xbm)}$ will be used to define the termination criteria
for the proposed solver.

\subsection{Choice of objective function}
\label{sec:optim:obj}
We propose an objective function that consists of two terms: one term
penalizes a measure of the DG solution error
$\func{f_\err}{\Rbb^{N_\ubm}\times\Rbb^{N_\xbm}}{\Rbb}$
and the other term penalizes distortion of the mesh
$\func{f_\msh}{\Rbb^{N_\xbm}}{\Rbb}$, i.e.,
\begin{equation} \label{eqn:obj-enrich-res-msh}
 f(\ubm, \xbm) = f_\err(\ubm, \xbm) + \kappa^2 f_\msh(\xbm),
\end{equation}
where $\kappa \in \Rbb_+$ is a parameter that weights the
contribution of the two terms. Since a piecewise polynomial solution
on an aligned mesh will have much lower error than on a non-aligned mesh,
$f_\err$ promotes alignment of the mesh with discontinuities while
$f_\msh$ prevents the mesh from entangling or becoming unacceptably skewed.

For the error-like tracking term, we use the norm of the DG residual
corresponding to an enriched \textit{test space}, i.e.,
\begin{equation} \label{eqn:obj-enrich-res}
 f_\err(\ubm, \xbm) \coloneqq \frac{1}{2} \Rbm(\ubm,\xbm)^T\Rbm(\ubm,\xbm),
\end{equation}
where we enrich the test space using polynomials of one degree higher than
the trial space. This follows on a large body of work that uses residual-based
error indicators to drive $h$-, $p$-, and $r$-adaptivity
\cite{fidkowski2011output}.
This is a reasonable choice for the objective function because
the enriched test space adds additional constraints to the
solution and even though a discrete solution $\ubm$ containing significant
oscillations satisfies $\rbm(\ubm, \xbm) = \zerobold$, it will likely not
minimize $\Rbm(\,\cdot\,, \xbm)$. Furthermore, the more the test space is
enriched, the closer the DG residual comes to enforcing the true conservation
law and therefore the minimum-residual solution approaches the exact solution
of the PDE. Given that we are explicitly enforcing the constraint
$\rbm(\ubm,\xbm) = \zerobold$, which fixes $\ubm$ for a given $\xbm$
(assuming $\rbm(\cdot,\xbm) = \zerobold$ has a unique solution),
the solution of the optimization problem must deform the mesh $\xbm$
to drive the pair $(\ubm,\xbm)$ to a point where the enriched residual
is minimized, which we expect to be a mesh that tracks sharp features
in the solution. Therefore, we expect this choice of objective function
to have desirable tracking properties, which will be confirmed by our
numerical experiments (Section~\ref{sec:num-exp}). Finally, this
choice of objective function is agnostic to whether we are considering
an inviscid or viscous conservation law, which is not the case for
many popular physics-based feature indicators than rely, e.g., on the
Rankine-Hugoniot conditions.

Our work in \cite{zahr2018optimization} considered the objective function
in (\ref{eqn:obj-enrich-res}), but did not endorse it due to numerical
experiments that showed it possessed non-aligned local minima that made
it impractical. However, in that work we only considered meshes with a
fixed topology, which contributed to these local minima. In the present work,
we collapse elements with small volumes
after each major optimization iteration, which either
eliminates the local minima of (\ref{eqn:obj-enrich-res}) or suggests this
strategy less sensitive to local minima.

Other choices for the objective function are possible, but are not considered
here. The objective function proposed in \cite{zahr2018optimization} is the
elementwise deviation of the DG solution from its mean, which was shown to
have excellent tracking properties. However, it tends to move the mesh in
regions of the domain where the solution is smooth and therefore not ideal
when the entire mesh is parametrized. The work in \cite{corrigan2019moving}
uses a physics-based objective function based on the Rankine-Hugoniot
conditions at all element faces for inviscid conservation laws, which was
shown to work well when combined with a DG-like discretization in a
minimum-residual framework.

To not only prevent the degradation of the mesh quality, but actively
promote mesh smoothing, we take $\kappa \in \Rbb_+$ and define $f_\msh$
as the deviation of the distortion of the physical mesh from the
distortion of the reference mesh
\begin{equation} \label{eqn:obj-msh}
 f_\msh(\xbm) = \frac{1}{2}(\Rbm_\msh(\xbm)-\Rbm_\msh(\Xbm))^T
                           (\Rbm_\msh(\xbm)-\Rbm_\msh(\Xbm)),
\end{equation}
where $\func{\Rbm_\msh}{\Rbb^{N_\xbm}}{\Rbb^{|\Ecal_{h,q}|}}$ is the
algebraic system corresponding to the elementwise mesh distortion
used for high-order mesh generation \cite{knupp01quality,roca16distortion}
\begin{equation} \label{eqn:res-mshdist}
 r_\msh^K(\map_{h,q}) \coloneqq
 \int_{K}
 \left(
  \frac{\norm{\defgrad_{h,q}}_F^2}{(\det\defgrad_{h,q})_+^{2/d}}
 \right)^2 \, dv
\end{equation}
and $\Xbm \in \Rbb^{N_\xbm}$ are the nodal coordinates of the reference
mesh.
Equation (\ref{eqn:res-mshdist}) is similar to the distortion measure used in
our previous work \cite{zahr2018optimization}, which was shown to
maintain high-quality meshes even when tracking difficult discontinuity
surfaces. We define $\Rbm_\msh$ as the deviation from the distortion of the
reference mesh rather than the mesh distortion itself because the Hessian
approximation used for our solver performs best when the objective function
approaches zero (Section~\ref{sec:solver:hess}). Even though the
distortion term prevents the combined objective function from converging
to zero, if the mesh is high-quality, we expect $\Rbm_\msh(\xbm)$ to be
close to its minimum value (component-wise) $\Rbm_\msh(\Xbm)$.

To close this section, we define the following vector-valued function
\begin{equation}
 \Fbm(\ubm,\xbm) \coloneqq
 \begin{bmatrix} 
  \Rbm(\ubm, \xbm) \\ \kappa(\Rbm_\msh(\xbm)-\Rbm_\msh(\Xbm))
 \end{bmatrix} 
\end{equation}
and re-write the objective function as
\begin{equation} \label{eqn:obj-minres}
 f(\ubm,\xbm) = \frac{1}{2}\norm{\Fbm(\ubm,\xbm)}_2^2
              = \frac{1}{2}\Rbm(\ubm,\xbm)^T\Rbm(\ubm,\xbm) +
                \frac{\kappa^2}{2}
                (\Rbm_\msh(\xbm)-\Rbm_\msh(\Xbm))^T
                (\Rbm_\msh(\xbm)-\Rbm_\msh(\Xbm)),
\end{equation}
to emphasize the objective function is the square two-norm of a
residual function. This has implications in terms of available
solvers for the optimization problem in (\ref{eqn:pde-opt}) as will be
discussed in Section~\ref{sec:solver:hess}.


\subsection{Boundary constraint enforcement}
\label{sec:optim:bndconstr}
In order to maintain a boundary-conforming mesh, the coordinates of
all mesh nodes cannot be allowed to move freely; rather, we must add
boundary constraints to ensure nodes \textit{slide} along the domain
boundaries. To this end, we write the mesh node coordinates as the
result of a mapping $\func{\chi}{\Rbb^{N_\phibold}}{\Rbb^{N_\xbm}}$ from
the unconstrained degrees of freedom $\phibold\in\Rbb^{N_\phibold}$ that
incorporates all boundary constraints
\begin{equation}
 \xbm = \chi(\phibold),
\end{equation}
where the specific form of the mapping depends on the domain under
consideration. This constraint is incorporated into the optimization
problem (\ref{eqn:pde-opt}) as:
\begin{equation} \label{eqn:pde-opt1}
 \optconOne{\ubm\in\Rbb^{N_\ubm},\phibold\in\Rbb^{N_\phibold}}
           {f(\ubm, \chi(\phibold))}{\rbm(\ubm, \chi(\phibold)) = \zerobold.}
\end{equation}
As discussed in \cite{zahr2018optimization}, the re-parametrized
formulation can also be used to explicitly incorporate mesh
smoothing into the mesh deformation. By introducing the following
definitions
\begin{equation}
 \tilde{f}(\ubm,\phibold) \coloneqq f(\ubm,\chi(\phibold)), \qquad
 \tilde{\rbm}(\ubm,\phibold) \coloneqq \rbm(\ubm,\chi(\phibold)), \qquad
 \tilde{\Fbm}(\ubm,\phibold) \coloneqq \Fbm(\ubm,\chi(\phibold)),
\end{equation}
the optimization problem with boundary enforcement in (\ref{eqn:pde-opt1})
becomes
\begin{equation} \label{eqn:pde-opt2}
 \optconOne{\ubm\in\Rbb^{N_\ubm},\phibold\in\Rbb^{N_\phibold}}
           {\tilde{f}(\ubm, \phibold)}{\tilde\rbm(\ubm, \phibold) = \zerobold,}
\end{equation}
which has the same structure as the original optimization problem
without boundary enforcement (\ref{eqn:pde-opt}) with the following
replacements: $\tilde{f} \leftarrow f$, $\tilde{\rbm} \leftarrow \rbm$,
and $\phibold \leftarrow \xbm$.

Most of the problems in this work only require nodes to
slide along boundaries aligned with coordinate directions; the
nodes on all other boundaries are fixed. Nodes sliding along more
general boundaries has been considered in
\cite{corrigan2019convergence, zahr2020radapt}.
In this special case, we partition the mesh node coordinates into
the constrained $\xbm_c \in \Rbb^{N_\xbm-N_\phibold}$ and
unconstrained $\phibold \in \Rbb^{N_\phibold}$ coordinates
\begin{equation}
 \xbm = \begin{bmatrix} \phibold \\ \xbm_c \end{bmatrix},
\end{equation}
which implies the boundary mapping is the padded identity mapping
\begin{equation}
 \xbm =
 \chi(\phibold; \xbm_c) \coloneqq
 \begin{bmatrix} \phibold \\ \xbm_c \end{bmatrix}, \qquad
 \pder{\chi}{\phibold}(\phibold; \xbm_c) =
 \begin{bmatrix} \Ibm \\ \zerobold \end{bmatrix}.
\end{equation}



\section{Full space, minimum-residual solver for optimization-based
         discontinuity tracking}
\label{sec:solver}
With the formulation of the optimization problem given in
(\ref{eqn:pde-opt2}), this section introduces a robust, iterative
solver. For simplicity, we introduce the solver for
the optimization problem in (\ref{eqn:pde-opt}), i.e., without boundary
enforcement; however, due to the mirror structure between
(\ref{eqn:pde-opt}) and (\ref{eqn:pde-opt2}), the exact algorithm
applies to solve (\ref{eqn:pde-opt2}) with the following replacements:
$\tilde{f} \leftarrow f$, $\tilde{\rbm} \leftarrow \rbm$,
and $\phibold \leftarrow \xbm$.

For brevity, we combine the PDE solution $\ubm$ and mesh coordinates $\xbm$
into a single vector
\begin{equation}
 \zbm \coloneqq \begin{bmatrix} \ubm \\ \xbm \end{bmatrix} \in\Rbb^{N_\zbm},
\end{equation}
where $N_\zbm = N_\ubm+N_\xbm$. In the remainder, we will replace
$(\ubm,\xbm)$ with $\zbm$, and vice versa, as needed. We will also
abbreviate the partial derivatives of the objective function and
DG residual as
\begin{equation}
 \begin{aligned}
  \gbm_\zbm(\zbm) &= \pder{f}{\zbm}(\zbm)^T
  &&\gbm_\ubm(\zbm) = \pder{f}{\ubm}(\zbm)^T
  &&\gbm_\xbm(\zbm) = \pder{f}{\xbm}(\zbm)^T \\
  \Jbm_\zbm(\zbm) &= \pder{\rbm}{\zbm}(\zbm),
  &&\Jbm_\ubm(\zbm) = \pder{\rbm}{\ubm}(\zbm),
  &&\Jbm_\xbm(\zbm) = \pder{\rbm}{\xbm}(\zbm).
 \end{aligned}
\end{equation}

\subsection{Sequential quadratic programming solver}
\label{sec:solver:sqp}
The proposed solver is a sequential quadratic programming (SQP) method
\cite{boggs2000sequential} that uses a sequence of quadratic programs
to solve (\ref{eqn:pde-opt}). Let $\zbm_k$ denote the current iterate and
formulate a quadratic subproblem by linearizing the constraint about $\zbm_k$
and using a second-order Taylor series centered at $\zbm_k$ to approximate
the Lagrangian function
\begin{equation} \label{eqn:qp1}
 \optconOne{\Delta\zbm\in\Rbb^{N_\zbm}}
           {\gbm_\zbm(\zbm_k)^T\Delta\zbm +
            \frac{1}{2}\Delta\zbm^T\Bbm(\zbm_k)\Delta\zbm}
           {\rbm(\zbm_k) + \Jbm_\zbm(\zbm_k)\Delta\zbm = \zerobold,}
\end{equation}
where $\Bbm(\zbm) \in \Rbb^{N_\zbm\times N_\zbm}$ a symmetric positive definite
approximation of the Hessian of $\Lcal(\zbm)$, i.e.,
\begin{equation}
 \Bbm(\zbm) \coloneqq
 \begin{bmatrix}
  \Bbm_{\ubm\ubm}(\zbm)   & \Bbm_{\ubm\xbm}(\zbm) \\
  \Bbm_{\ubm\xbm}(\zbm)^T & \Bbm_{\xbm\xbm}(\zbm)
 \end{bmatrix}
\end{equation}
such that
\begin{equation}
 \begin{aligned}
  \pderTwo{\Lcal}{\ubm}{\ubm}(\zbm) &\approx
  \Bbm_{\ubm\ubm}(\zbm) \in \Rbb^{N_\ubm\times N_\ubm} \\
  \pderTwo{\Lcal}{\ubm}{\xbm}(\zbm) &\approx
  \Bbm_{\ubm\xbm}(\zbm) \in \Rbb^{N_\ubm\times N_\xbm} \\
  \pderTwo{\Lcal}{\xbm}{\xbm}(\zbm) &\approx
  \Bbm_{\xbm\xbm}(\zbm) \in \Rbb^{N_\xbm\times N_\xbm}.
 \end{aligned}
\end{equation}
and the dependence on the Lagrange multipliers has been dropped.
Once the quadratic subproblem is solved to obtain the current search
direction $\Delta\zbm_{k+1} \in \Rbb^{N_\zbm}$, the current iterate is updated
according to
\begin{equation}
 \zbm_{k+1} = \zbm_k + \alpha_{k+1} \Delta\zbm_{k+1}
\end{equation}
to yield a sequence of iterates $\{\zbm_k\}$ where $\alpha_{k+1} \in (0, 1]$
is a step length parameter (Section~\ref{sec:solver:lsrch}).
If $\Bbm$ is the true Hessian and $\alpha = 1$,
this method is equivalent to Newton's method applied to the first-order
optimality conditions of (\ref{eqn:pde-opt}) \cite{boggs2000sequential}
and therefore will converge quadratically provided the initial guess
is sufficiently close to a local minima.

\subsection{Linear subproblem}
\label{sec:solver:linsub}
The first-order optimality condition of the quadratic program (\ref{eqn:qp1})
leads to the following linear system of equations
\begin{equation} \label{eqn:sqp1-linsys1}
 \begin{bmatrix}
  \Bbm_{\ubm\ubm}(\zbm_k)   & \Bbm_{\ubm\xbm}(\zbm_k) & \Jbm_\ubm(\zbm_k)^T \\
  \Bbm_{\ubm\xbm}(\zbm_k)^T & \Bbm_{\xbm\xbm}(\zbm_k) & \Jbm_\xbm(\zbm_k)^T \\
  \Jbm_\ubm(\zbm_k) & \Jbm_\xbm(\zbm_k) & \zerobold
 \end{bmatrix}
 \begin{bmatrix}
  \Delta\ubm_{k+1} \\ \Delta\xbm_{k+1} \\ \etabold_{k+1}
 \end{bmatrix} =
 -
 \begin{bmatrix}
  \gbm_\ubm(\zbm_k) \\ \gbm_\xbm(\zbm_k) \\ \rbm(\zbm_k)
 \end{bmatrix},
\end{equation}
where $\etabold_{k+1} \in \Rbb^{N_\ubm}$ are the Lagrange multipliers
associated with the linearized constraint and the step
$\Delta\zbm_{k+1}$ is decomposed as a step in the PDE state $\Delta\ubm_{k+1}$
and the nodal mesh coordinates $\Delta\xbm_{k+1}$
\begin{equation}
 \Delta\zbm_{k+1} =
 \begin{bmatrix} \Delta\ubm_{k+1} \\ \Delta\xbm_{k+1} \end{bmatrix}.
\end{equation}
The state update $\Delta\ubm$ and Lagrange multipliers $\etabold$ can
be eliminated from the linear system to obtain an explicit expression
for the state and mesh updates
\begin{equation} \label{eqn:sqp1-linsys2}
 \Delta\xbm_{k+1} = -\Abm(\zbm_k)^{-1}\bbm(\zbm_k), \qquad
 \Delta\ubm_{k+1} = \Delta\ubm_0(\zbm_k) + \Cbm(\zbm_k)\Delta\xbm,
\end{equation}
where
\begin{equation}
 \begin{aligned}
  \Delta\ubm_0(\zbm) &= -\Jbm_\ubm(\zbm)^{-1}\rbm(\zbm) \\
  \Cbm(\zbm) &= -\Jbm_\ubm(\zbm)^{-1}\Jbm_\xbm(\zbm) \\
  \Abm(\zbm) &= \Bbm_{\xbm\xbm}(\zbm)+2\Cbm(\zbm)^T\Bbm_{\ubm\xbm}(\zbm)+
                \Cbm(\zbm)^T\Bbm_{\ubm\ubm}(\zbm)\Cbm(\zbm) \\
  \bbm(\zbm) &= \gbm_\xbm(\zbm)+\Cbm(\zbm)^T\gbm_\ubm(\zbm)+
          \left(\Bbm_{\ubm\xbm}(\zbm)^T+\Cbm(\zbm)^T\Bbm_{\ubm\ubm}(\zbm)\right)
          \Delta\ubm_0.
 \end{aligned}
\end{equation}
Notice $\Delta\ubm_0(\zbm)$ is the Newton step for the nonlinear system
$\rbm(\ubm,\xbm) = \zerobold$ if $\xbm$ is fixed and $\Cbm(\zbm)$ is
the sensitivity of the solution $\ubm$ with respect to the mesh nodes
$\xbm$.

The linear systems in (\ref{eqn:sqp1-linsys1}) and (\ref{eqn:sqp1-linsys2})
are two mathematically equivalent options to compute the solution to the SQP
subproblem that place different requirements on the linear solver. The
system in (\ref{eqn:sqp1-linsys1}) is large (size: $2N_\ubm+N_\xbm$)
since it simultaneously computes the steps $\Delta\ubm$, $\Delta\xbm$
and Lagrange multipliers $\etabold$, but each (block) entry is relatively
easy to form. On the other hand, the system in
(\ref{eqn:sqp1-linsys2}) is smaller (size: $N_\xbm$), but formation of
$\Abm(\zbm_k)$ requires computation of $\Cbm(\zbm_k)$, which is impractical
if a direct solver is not available. In this work, we use a direct solver
to solve (\ref{eqn:sqp1-linsys1}); however, this approach is not
practical for large-scale problems. In future work, we will develop iterative
solvers and preconditioners for the full system in (\ref{eqn:sqp1-linsys1}).

\subsection{Line search globalization}
\label{sec:solver:lsrch}
To ensure the sequence $\{\zbm_k\}$ converges to a first-order critical
point of (\ref{eqn:pde-opt}) from an arbitrary initial guess, the
SQP algorithm must be globalized with e.g., a trust region strategy
or line search. We choose to globalize the SQP method with a line search
that computes the step length such that it minimizes a merit
function that combines the objective function and a measure of constraint
violation into a scalar function. In this work, we use the $\ell_1$
penalty function $\func{\varphi_k}{\Rbb}{\Rbb}$
\begin{equation} \label{eqn:l1merit}
 \varphi_k(\alpha) \coloneqq
 f(\zbm_k+\alpha\Delta\zbm_k)+
 \mu \norm{\rbm(\zbm_k+\alpha\Delta\zbm_{k+1})}_1
\end{equation}
where $\mu > 0$ is the penalty parameter. This is an exact merit function
in the sense that there exists a positive scalar $\hat\mu$ such that for
any $\mu > \hat\mu$, any local solution of (\ref{eqn:pde-opt}) is a
local minimizer of (\ref{eqn:l1merit}) \cite{nocedal2006numerical}.
Furthermore, 
$\hat\mu = \norm{\lambdabold^\star}_\infty$, where $\lambdabold^\star$
are the Lagrange multipliers associated with the optimal solution
$(\ubm^\star,\xbm^\star)$ \cite{nocedal2006numerical}. Therefore, in the present
setting we take $\mu = 2\norm{\hat\lambdabold(\zbm_k)}_\infty$, where
$\hat\lambdabold$ is defined in (\ref{eqn:lagrmult}).

It is well-known that it is not necessary to find the exact minimizer
of the merit function to obtain a convergent algorithm. Instead, we
search for $\alpha_{k+1} \in (0, 1]$ that satisfies \textit{sufficient decrease}
\begin{equation} \label{eqn:suffdec}
 \varphi_k(\alpha_{k+1}) \leq \varphi_k(0) + c\alpha_{k+1}\varphi_k'(0),
\end{equation}
where $c \in (0, 1)$. We use a backtracking strategy
\cite{nocedal2006numerical} to determine $\alpha_{k+1}$: define
$\alpha_{k+1} = \tau^{n-1}$ for $\tau \in (0, 1)$ and let $n \in \Nbb$
be the smallest number such that sufficient decrease (\ref{eqn:suffdec}) holds.
In this work we make standard choices for these parameters
\cite{nocedal2006numerical}: $c = 10^{-4}$ and $\tau = 0.5$.

\subsection{Levenberg-Marquardt Hessian approximation}
\label{sec:solver:hess}
The true Hessian of the Lagrangian in (\ref{eqn:lagr}) is
\begin{equation}
 \Hbm(\zbm) = \pder{\Fbm}{\zbm}(\zbm)^T\pder{\Fbm}{\zbm}(\zbm) +
              \Fbm_i(\zbm) \pderTwo{\Fbm_i}{\zbm}{\zbm}(\zbm) -
              \lambdabold_i\pderTwo{\rbm_i}{\zbm}{\zbm}(\zbm)
\end{equation}
where summation is implied over the repeated index. The second
and third term involves the Hessian of the DG residual, which
is a complicated third-order tensor that is rarely available in
computational mechanics codes. Therefore, we use the Gauss-Newton
assumption and approximate the Hessian as
\begin{equation} \label{eqn:gnhess}
 \Hbm(\zbm)\approx\pder{\Fbm}{\zbm}(\zbm)^T\pder{\Fbm}{\zbm}(\zbm),
\end{equation}
which is justified if the combined enriched residual and mesh distortion
$\Fbm(\zbm)$ is small.
While this approximation is convenient, it could also lead to singular
or ill-conditioned Hessian approximations, which would in turn lead to
poor search directions. Therefore, we use the Levenberg-Marquardt approach
that adds a scaled multiple of a symmetric positive definite matrix to
regularize the system.
Based on the observations in \cite{corrigan2019moving} that
Levenberg-Marquardt regularization improves the mesh motion and
is not needed for the state, we only regularize the mesh components
of the Hessian, i.e.,
\begin{equation} \label{eqn:hess-approx}
 \Bbm_{\ubm\ubm}(\zbm) =
 \pder{\Fbm}{\ubm}(\zbm)^T\pder{\Fbm}{\ubm}(\zbm), \quad
 \Bbm_{\ubm\xbm}(\zbm) =
 \pder{\Fbm}{\ubm}(\zbm)^T\pder{\Fbm}{\xbm}(\zbm), \quad
 \Bbm_{\xbm\xbm}(\zbm) =
 \pder{\Fbm}{\xbm}(\zbm)^T\pder{\Fbm}{\xbm}(\zbm) +
 \gamma \Dbm,
\end{equation}
where $\gamma \in \Rbb_+$ and $\Dbm \in \Rbb^{N_\xbm\times N_\xbm}$
is a symmetric positive definite (SPD) matrix
(Section~\ref{sec:solver:regmat}).

\subsection{Choice of regularization matrix}
\label{sec:solver:regmat}
The classic Levenberg-Marquardt algorithm uses the identity matrix as
the regularization matrix $\Dbm = \Ibm_{N_\xbm}$ to guard against poor
search directions that could result if the Jacobian of $\Fbm$
with respect to $\xbm$ is rank deficient or ill-conditioned. However,
the same result can be achieved for any SPD matrix $\Dbm$. We consider a
regularization matrix that is known to possess smoothing properties: the
stiffness matrix of a linear elliptic PDE. To this end,
define a vector-valued function $\func{v}{\refdom}{\Rbb^d}$ in which each
component satisfies the elliptic PDE with homogeneous Neumann boundary
conditions
\begin{equation} \label{eqn:poi1}
 \nabla_X \cdot (k \nabla_X v_i) = 0 \quad \text{in}~~\refdom, \qquad
 \nabla_X v_i \cdot N = 0 \quad \text{on}~~\partial\refdom
\end{equation}
for $i = 1, \dots, d$ and the coefficient $\func{k}{\refdom}{\Rbb_+}$
is piecewise constant over each element in $\Ecal_{h,q}$.
Numerical experimentation on problems where elements in the reference
mesh significantly vary in size indicate that the natural scaling of
the stiffness matrix ($k = 1$) is not sufficient to provide
search directions that are scaled according to the mesh resolution. Instead,
we take the piecewise constant coefficient  as
\begin{equation} \label{eqn:poi2}
 k(x) = \omega_K, \quad x \in K
\end{equation}
for each $K \in \Ecal_{h,q}$, where $\omega_K \in \Rbb_+$ is inversely
proportional to the size of the element
\begin{equation}
 \omega_K \coloneqq \frac{\underset{K'\in\Ecal_{h,q}}{\min} |K'|}{|K|}.
\end{equation}
The regularization matrix is the assembled stiffness matrix associated with
the elliptic system in (\ref{eqn:poi1})-(\ref{eqn:poi2}) over the global finite
element space $\Wcal_{h,q}$. This is similar to the regularization
matrix used in \cite{corrigan2019convergence} without volume-based scaling
of the elliptic coefficient, which is important for problems with
reference domain elements of varying size
(Section~\ref{sec:num-exp:euler:naca0}).
Boundary constraints are incorporated into the elliptic regularization
matrix as homogeneous Dirichlet boundary conditions.

\subsection{Adaptive regularization parameter}
\label{sec:solver:regparam}
The ideal value for the regularization parameter $\gamma$ in
(\ref{eqn:hess-approx}) is difficult to know \textit{a priori}. Large values
of $\gamma$ produce highly regularized search directions $\Delta\xbm$,
e.g., the search direction becomes the solution of the elliptic PDE
with the steepest descent direction as the right-hand side; however, in
this case, the Hessian approximation is poor. Alternatively, small values
of $\gamma$, particularly in early iterations, can leads to bad search
directions due to ill-conditioning of the Hessian approximation.
Following the work in \cite{corrigan2019moving}, we introduce
an adaptive algorithm to minimize the importance of choosing the appropriate
value for $\gamma$ \textit{a priori}. Our algorithm is based on
the heuristic that $\gamma$ should control the size of the mesh deformation
at a given iteration. Therefore, if a step is too large,
we increase the value of $\gamma$ for the next iteration and vice
versa. Let $\gamma_k > 0$ be the value of the regularization parameter
at iteration $k$ and define constants $\kappa_1, \kappa_2 > 0$,
$\sigma \in (0, 1)$, $\gamma_{min} > 0$ then
\begin{equation}
 \gamma_{k+1} = \max\{\bar\gamma_{k+1}, \gamma_{min}\}, \qquad
 \bar\gamma_{k+1} =
 \begin{cases}
  \sigma^{-1} \gamma_k & \text{if } \norm{\Delta \xbm_k} < \kappa_1 \\
  \sigma \gamma_k & \text{if } \norm{\Delta \xbm_k} > \kappa_2 \\
  \gamma_k & \text{otherwise}.
 \end{cases}
\end{equation}
For domains with dimensions $\Ocal(1)$, we take $\kappa_1 = 10^{-2}$ and
$\kappa_2 = 10^{-1}$. In this work, we choose a conservative value
$\sigma = 0.5$ to prevent $\gamma$ from changing significantly between
iterations. The appropriate value of $\gamma_{min}$ is problem-dependent
and will be addressed in Section~\ref{sec:num-exp}.


\subsection{Termination criteria}
\label{sec:solver:term}
The termination criteria for the SQP method comes from the first-order
optimality criteria discussed in Section~\ref{sec:optim:pdeconstr}
(\ref{eqn:kkt2}).
That is, given tolerances $\epsilon_1, \epsilon_2 > 0$, $\zbm_k$
is a considered a numerical solution of (\ref{eqn:pde-opt}) if
\begin{equation} \label{eqn:conv}
 \norm{\cbm(\zbm_k)} < \epsilon_1, \qquad
 \norm{\rbm(\zbm_k)} < \epsilon_2.
\end{equation}
We have empirically observed that feasibility can be driven to near
machine tolerance (take $\epsilon_2 = 10^{-10}$); however, the
optimality condition is a more difficult condition
(take $\epsilon_1 = 10^{-5}$). For difficult problems, it may require
a large number of iterations, so we safeguard the algorithm
with a maximum number of iterations ($N_{max}$), i.e., iterations
terminate when either the convergence criteria (\ref{eqn:conv}) are met
or a $N_{max}$ iterations have been completed. If termination is based
on $N_{max}$, it is likely the DG constraint is not satisfied, i.e.,
\begin{equation}
 \rbm(\zbm_{N_{max}}) > \epsilon_2,
\end{equation}
which could lead to a non-conservative scheme. In this case, we
solve the DG equations on the fixed mesh with nodal coordinates
$\xbm_{N_{max}}$ using Newton's method.

\section{Practical considerations}
\label{sec:practical}

\subsection{Solution and mesh initialization}
\label{sec:practical:init}
The discontinuity-tracking optimization problem in (\ref{eqn:pde-opt})
is non-convex and therefore the initial guess for the SQP solver
is critical to obtain a good solution. In the present context, this
means we must provide a reasonable initial guess for the mesh
coordinates $\xbm_0$ and DG solution $\ubm_0$. 

\subsubsection{Special case: straight-sided mesh ($q = 1$)}
\label{sec:practical:init1}
First consider the case of a straight-sided mesh ($q = 1$) combined
with any finite element space for the solution ($p \geq 0$). The
mesh is always initialized from the reference mesh
which comes from mesh generation agnostic to the discontinuity because
we usually do not have an estimate of the discontinuity surface. In special
cases where an estimate of the discontinuity surface is available, this could
be used to drive generation of the reference mesh. The DG solution is
initialized from the DG ($p = 0$) solution on the reference mesh; $p = 0$
is used because nonlinear instabilities resulting from oscillations about
discontinuities cannot arise with a piecewise constant solution field. For
nonlinear problems the $p = 0$ solution is obtained using pseudo-transient
continuation \cite{kelley1998convergence} with adaptive time steps,
initialized from uniform flow.

\subsubsection{General case: high-order meshes ($q > 1$)}
\label{sec:practical:init2}
To avoid local minima in the optimization problem (\ref{eqn:pde-opt})
that arise when high-order meshes are used, we usually initialize the
tracking problem for $p \geq 0$, $q > 1$ from the solution of the
tracking problem for $p' \leq p$, $q' = 1$. That is, we solve the
the tracking problem using a solution space with polynomial degree
$p' \leq p$ and mesh deformation with polynomial degree $q' = 1$ and
initialize the desired tracking problem (solution space of degree $p$ and
mesh of degree $q > 1$) from the resulting DG solution and mesh. This
strategy comes from our observations that $q = 1$ tracking is quite robust
and convergence of the $q > 1$ solution from a straight-sided tracking mesh
is rapid. For difficult problems, it may be helpful to use continuation
on the polynomial degree for the solution as well, i.e., take $p' < p$.

\subsection{Edge collapses and solution transfer}
\label{sec:practical:meshop}

As we will show in our numerical experiments in Section~\ref{sec:num-exp}, 
we typically start with an initial mesh that is far from alignment with the shock 
and iterate using an SQP solver which attempts to move nodes onto the shock. 
Since the deformation of the initial mesh can be quite large, this can result in 
severely ill-conditioned elements that drastically degrade the quality of our solution. 
To address this, we follow the approach in \cite{corrigan2019moving} 
and collapse elements once they become problematic. 
In particular, after each SQP iteration, we compute the
volume of each element $K$ in the physical domain, $v^K$, and compare
it to the volume of the corresponding element in the reference domain,
$v_0^K$, where
\begin{equation}
 v_0^K = \int_{K} \, dV, \qquad
 v^K = \int_{\Gcal_{h,q}(K)} \, dV.
\end{equation}
If the volume of an element has decreased more than a certain factor,
i.e., $v^K < \epsilon v_0^K$ (where $\epsilon = 0.2$ in this work),
the element is tagged for removal.
In principle, $\kappa$ in the objective function (\ref{eqn:obj-enrich-res-msh}) could be 
chosen large enough to prevent these ill-conditioned elements from ever appearing. 
However, such a $\kappa$ would weight $f_\msh$ too heavily over $f_\err$ 
and interfere with the tracking capabilities of the method. Therefore, this weighting term 
$\kappa$ is only chosen large enough to prevent unacceptably bad elements 
(such as tangled ones) which would cause the mesh distortion term 
$f_\msh$ to blow up dramatically. As long as we avoid this situation, 
we can handle the remainder of the ill-conditioned elements through collapses. 

\begin{figure}[H]
 \centering
 \includegraphics[width=0.45\textwidth]{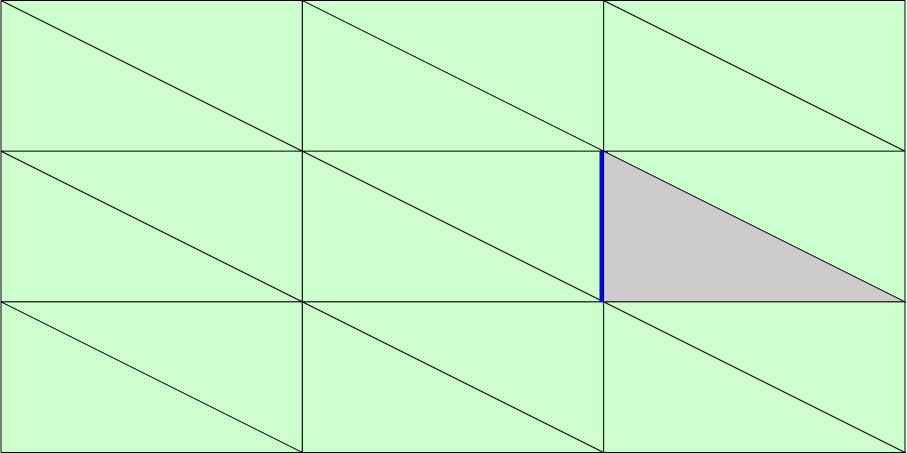} \quad
 \includegraphics[width=0.45\textwidth]{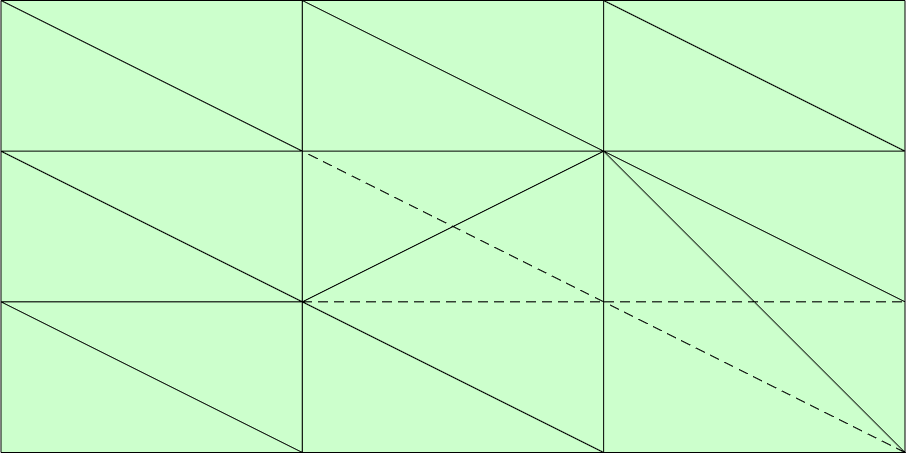}
 \caption{Demonstration of edge collapse algorithm: the element identified
          in the original mesh (\textit{left}) is collapsed along the
          highlighted edge to produce the new mesh (\textit{right})
          with the original elements shown in dashed lines for reference.}
 \label{fig:edge_collapse_demo}
\end{figure}

Each element in the physical domain $\xbm$ tagged for removal 
is eliminated by an edge collapse \cite{lohner2008applied}. 
For a given tagged element, we choose to collapse the shortest edge 
into the longest edge (Figure~\ref{fig:edge_collapse_demo}). 
This is based off the principle that the shortest edge is likely 
to be transverse to the shock and the longest edge is likely 
to be aligned with the shock. 
A special case to note is when tagged elements have nodes or edges in common. 
Additional care must be taken to choose an edge to collapse that is consistent 
with all affected elements. For elements on the boundary, 
this logic is slightly modified to ensure that an edge collapse does not move nodes
 off the boundary. Note that an edge collapse must also be applied 
 to the corresponding element in the reference domain $\Xbm$ 
 to ensure the physical and reference domains always have the same topology. 

The solution $\ubm$ is transferred to the new mesh by removing the entries 
corresponding to the degrees of freedom in the collapsed elements. 
To update the data structures for both the mesh and the solution, 
we simply need to delete the entries in the element connectivity matrix 
and solution vector corresponding to the collapsed elements and renumber 
based off the new node numbering. Because of our choice of a DG discretization, 
this update is particularly easy to do, and furthermore, removes the need to modify 
the degrees of freedom in the neighboring elements. In a CG framework, 
the degrees of freedom in the neighboring elements would need to be modified 
in order to guarantee continuity on the new mesh. 

While these edge collapses do modify the objective function and formally result 
in a new optimization problem, in practice this does not lead to issues in the 
convergence of the solver. We observe that these collapses mostly occur in the initial 
iterations of the SQP solver when the mesh is far from convergence and each iteration 
results in a large update to the mesh. At some point, the line search described in 
Section~\ref{sec:solver:lsrch} will reject steps that would cause large 
deformations in the mesh, hence precluding the need 
for further collapses and ensuring that they do not occur indefinitely. 

\section{Numerical experiments}
\label{sec:num-exp}
In this section we introduce three inviscid conservation laws
and demonstrate the tracking framework on six problems with
discontinuous solutions of varying difficulty. We also provide
a detailed study of the various algorithmic parameters introduced
in Section~\ref{sec:disc}-\ref{sec:optim}, in particular, the choice
of numerical flux ($\Hcal$) and mesh distortion parameter ($\kappa$).

\subsection{Linear advection}
\label{sec:num-exp:advec}
The first conservation law we consider is steady linear advection of a scalar
quantity $\func{U}{\dom}{\Rbb}$ through a domain $\dom \subset \Rbb^d$
\begin{equation} \label{eqn:advec}
 \nabla\cdot(\beta U) = 0 \quad \text{in } \dom, \qquad
 U = U_\infty \quad \text{on } \Gamma_i,
\end{equation}
where $\func{\beta}{\dom}{\Rbb^d}$ is the local flow direction,
$\Gamma_i \coloneqq \{ x \in \partial\dom \mid \beta\cdot n(x) < 0\}$
is the inflow boundary, $\func{\partial\dom}{\Rbb^d}$ is the unit outward
normal to the boundary, and $\func{U_\infty}{\Gamma_i}{\Rbb}$ is the
inflow boundary condition. We consider the pure upwind numerical flux
$\Hcal_\up$ introduced in (\ref{eqn:upwind}) because it satisfies conditions
(i)-(ii) (Section~\ref{sec:disc:numflux}) for a stable DG discretization
and condition (iii) to be suitable for tracking, as well as its smoothed
version $\Hcal_\up^a$ to recover smoothness with respect to variations
of the mesh nodes (iv). In this section we take the smoothness parameter
to be $a = 10$ to trade-off between smoothness and appropriate upwinding.
The boundary condition is enforced via the pure upwind numerical flux
evaluated at the interior state $U^+$ and the boundary state
$U^\partial = U_\infty$.

\subsubsection{Straight shock, piecewise constant solution}
\label{sec:num-exp:advec:straight}
First we consider a two-dimensional ($d = 2$) domain
$\dom \coloneqq (-1, 1) \times (0, 1)$ with constant
advection field and piecewise constant boundary condition
\begin{equation} \label{eqn:advec-bndcnd0}
 \beta(x) = \begin{bmatrix} -1.25 \\ 1 \end{bmatrix}, \quad U_\infty(x) = H(x).
\end{equation}
This leads to a linear discontinuity surface $\Gamma_s \subset \dom$
\begin{equation}
 \Gamma_s \coloneqq \{(-1.25s, s) \mid s \in (0, 0.8)\}
\end{equation}
and piecewise constant solution field
\begin{equation}
 U(x) = H(x_1+1.25 x_2).
\end{equation}
To ensure the Heaviside function on the bottom (inflow) boundary is
accurately represented and integrated in the weak form, we require
our computational mesh to have an element face that intersects $(0, 0)$
and do not allow the corresponding node to move throughout iterations
using the boundary mapping described in Section~\ref{sec:optim:bndconstr}.

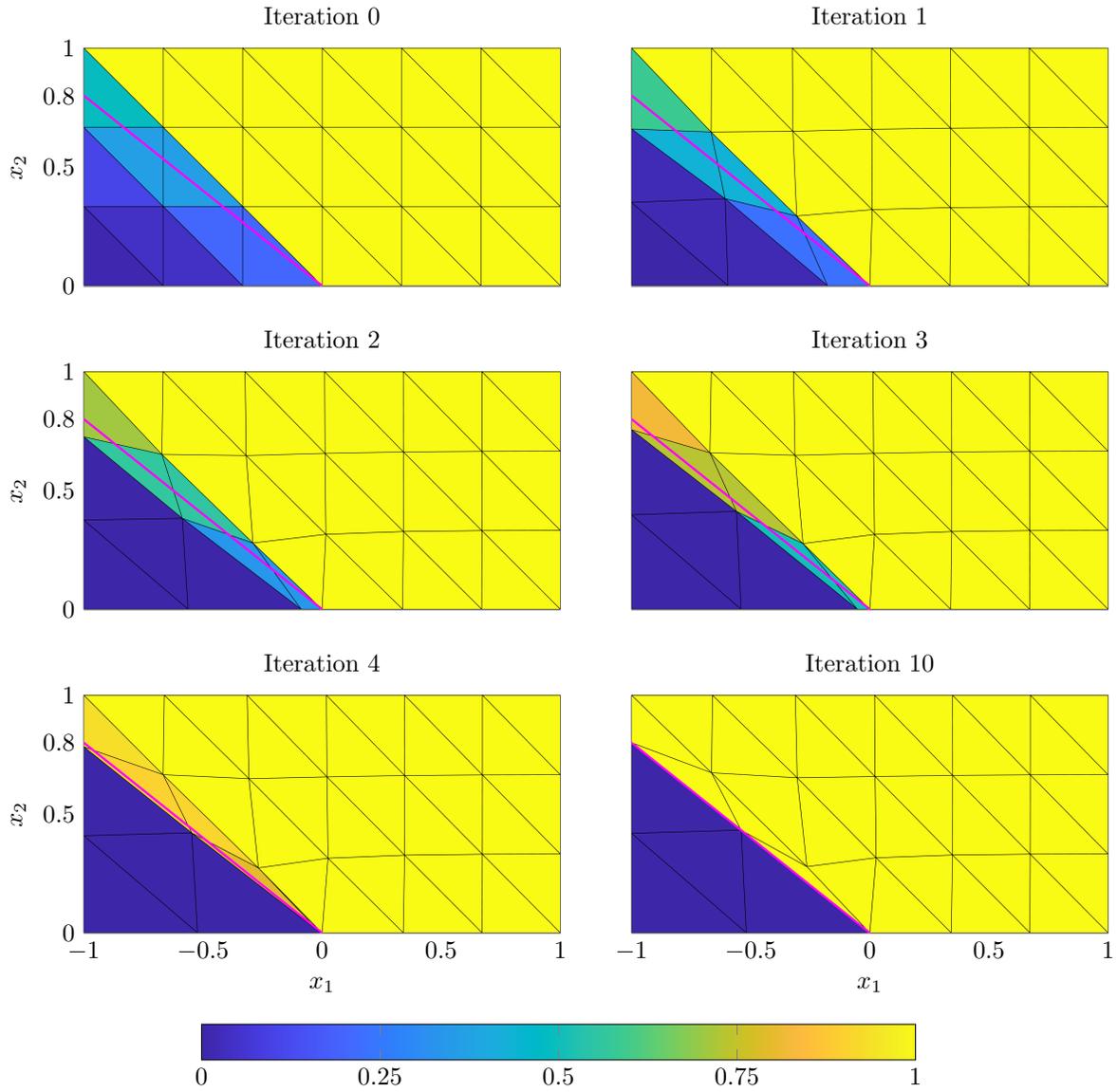
\begin{figure}
 \centering
 \input{tikz/advec2d_beta0_c1s1_nel6x3_porder0x1.tikz}
 \input{tikz/cbar_parula_0x1.tikz}
 \caption{Solution of advection equation with a straight shock (constant
          advection field) using the
          tracking method at various iterations throughout the solution
          procedure using $p = 0$ basis for the solution and $q = 1$
          basis for the mesh. The magenta line highlights the discontinuity
          surface of the exact solution. The method converges to the exact
          solution, which lies in the chosen finite element subspace, to
          near machine precision in only $10$ iterations.}
 \label{fig:advec2d_beta0_c1s1}
\end{figure}

\begin{figure}
 \centering
 \input{py/advec2d_beta0_c1s1_nel6x3_porder0x1_conv.tikz}
 \caption{Convergence of the DG residual $\norm{\rbm(\ubm,\xbm)}$
          (\ref{line:advec2d_beta0_c1s1:R0}), enriched DG residual
          $\norm{\Rbm(\ubm,\xbm)}$ (\ref{line:advec2d_beta0_c1s1:R1}),
          optimality condition $\norm{\cbm(\ubm,\xbm)}$
          (\ref{line:advec2d_beta0_c1s1:dLdY}), and control of the
          regularization parameter (\ref{line:advec2d_beta0_c1s1:lam})
          for the tracking method applied to the advection equation
          with a straight shock (constant advection field). For this
          simple problem, Newton-like convergence is achieved.}
 \label{fig:advec2d_beta0_c1s1_conv}
\end{figure}
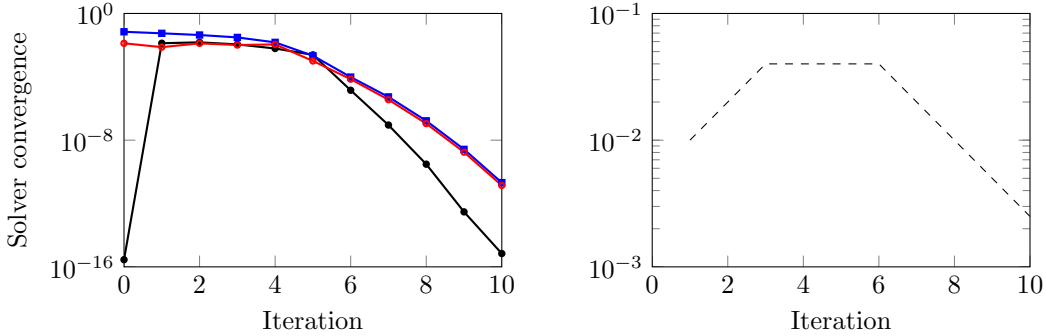

Given the piecewise constant solution and linear shock surface, the
solution lies in a $p = 0$ polynomial basis on a piecewise linear
mesh $q = 1$, provided the mesh tracks the discontinuity. To
demonstrate the performance of the tracking method, we use this
minimal basis with $p = 0$, $q = 1$. The reference mesh is taken
to be a uniform triangular mesh of the domain with $36$ elements.
The mesh and solution are initialized according to
Section~\ref{sec:practical:init}.  The SQP solver is used with
$\lambda$ chosen adaptively (Section~\ref{sec:solver:regparam}).
The various DG/tracking parameters are set as follows:
$\kappa = 0$ (no mesh smoothing), pure upwind numerical flux ($a=\infty$),
$\gamma_0 = 10^{-2}$ and $\gamma_{min} = 10^{-8}$ (regularization parameter
adaptivity), and $\epsilon_1 = 10^{-10}$, and $\epsilon_2 = 10^{-12}$
(termination criteria). After only $10$ iterations, the mesh perfectly
tracks the discontinuity (Figure~\ref{fig:advec2d_beta0_c1s1}) and the DG
solution closely matches the exact solution; the $L^1$ error of the
solution is $3.84\times 10^{-11}$.
For this problem where the finite element subspace contains the exact
solution of the problem, the tracking method exhibits Newton-like
convergence (Figure~\ref{fig:advec2d_beta0_c1s1_conv}).
Furthermore, the tracking solver is robust with regard to the
various algorithmic parameters introduced, i.e., fast convergence to
the exact solution was obtained without mesh smoothing, with the upwind
numerical flux (not smooth with respect to domain deformation), and
nearly independent of the choice of regularization adaptivity
($\gamma_0$, $\gamma_{min}$). The choice of these parameters becomes
significant for curved discontinuities and nontrivial flows as we
will demonstrate in subsequent sections.

\subsubsection{Curved shock, piecewise constant solution}
\label{sec:num-exp:advec:curved}
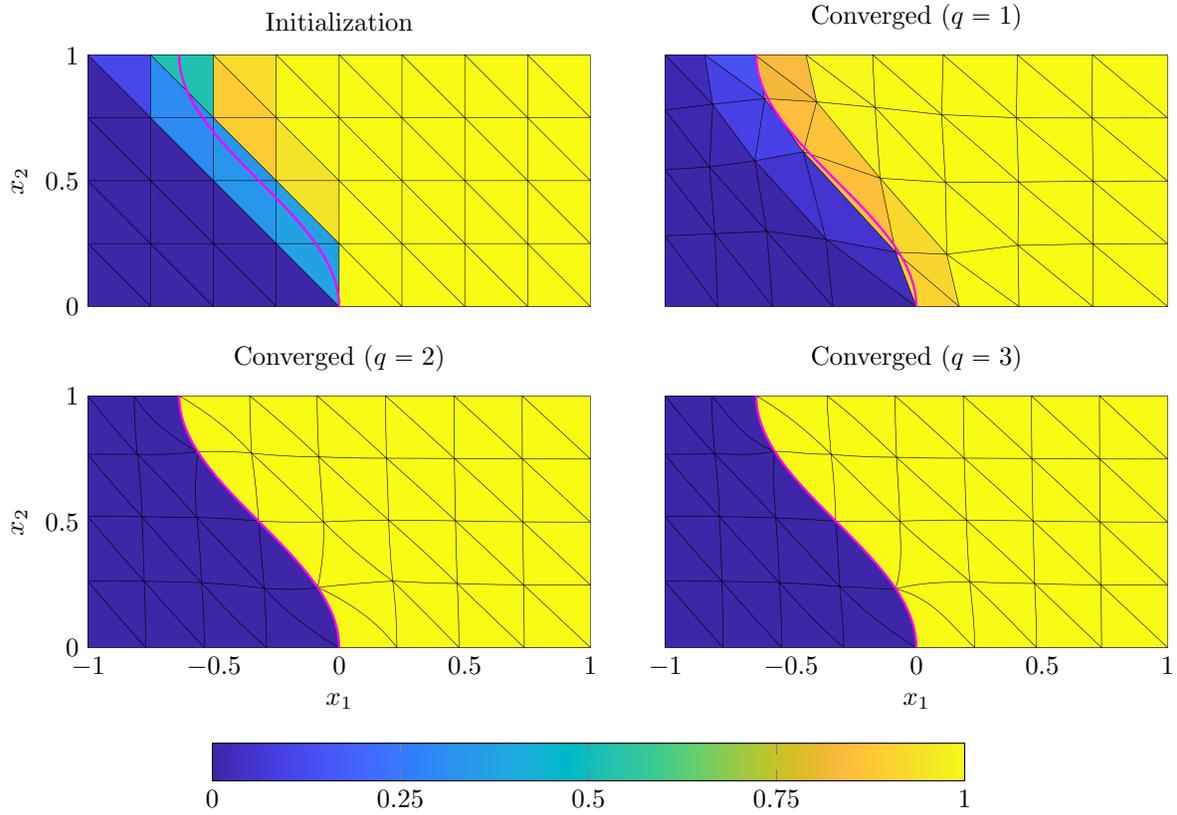
\begin{figure}
 \centering
 \input{tikz/advec2d_beta0_c2s1_nel8x4.tikz}
 \input{tikz/cbar_parula_0x1.tikz}
 \caption{Solution of advection equation with the trigonometric shock
          using the using a $p = 0$ basis for the solution and $q = 1$
          (\textit{top right}), $q = 2$ (\textit{bottom left}), and
          $q = 3$ (\textit{bottom right}) basis for the mesh. The
          DG solution on a uniform triangular mesh with $64$ elements is used
          to initialize the tracking method (\textit{top left}).
          The magenta line highlights
          the discontinuity surface of the exact solution. The method
          converges to a mesh and solution that approximates the true solution
          as good as can be expected given the resolution of the finite element
          space; however, only the high-order elements provide a reasonable
          approximation on this coarse mesh.}
 \label{fig:advec2d_beta0_c2s1}
\end{figure}

Next we consider linear advection through a two-dimensional ($d = 2$) domain
$\dom \coloneqq (-1, 1) \times (0, 1)$ with a spatially varying advection
field and piecewise constant boundary condition
\begin{equation} \label{eqn:advec-bndcnd1}
 \beta(x) = \begin{bmatrix} -\sin(\pi x_2) \\ 1 \end{bmatrix}, \quad
 U_\infty(x) = H(x).
\end{equation}
This leads to a trigonometric discontinuity surface $\Gamma_s \subset \Omega$
\begin{equation} \label{eqn:advec-discont-bndcnd1}
 \Gamma_s \coloneqq
 \left\{\left(\frac{\cos(\pi s)-1}{\pi}, s\right) \mid s \in (0, 1)\right\}
\end{equation}
and piecewise constant solution field
\begin{equation}
 U(x) = H(\pi x_1 - \cos(\pi x_2)+1).
\end{equation}
Unlike the previous problem, this solution cannot be represented
exactly using polynomial basis functions since the discontinuity
surface is non-polynomial (trigonometric). We use a piecewise
constant solution basis ($p = 0$) and piecewise linear, quadratic,
and cubic basis for the mesh ($q = 1,2,3$).
The reference mesh is taken to be a uniform triangular
mesh of the domain with $64$ elements. The mesh and solution are initialized
according to Section~\ref{sec:practical:init}; however, for this problem we
do not use continuation in the polynomial degree.  The SQP solver is used with
$\lambda$ chosen adaptively (Section~\ref{sec:solver:regparam}).
The various DG/tracking parameters are set as follows:
$\kappa = 0$ (mesh smoothing), smoothed upwind numerical flux
($a=10$), $\gamma_0 = 10^{-2}$ and $\gamma_{min} = 10^{-3}$ (regularization
parameter adaptivity), and $N_{max} = 80$, $\epsilon_1 = 10^{-10}$,
$\epsilon_2 = 10^{-12}$ (termination criteria).
For all polynomial degrees considered, the mesh
tracks the discontinuity as accurately as possible given the resolution
in the finite element space (Figure~\ref{fig:advec2d_beta0_c2s1}).
The $q = 1$ solution is under-resolved since the combination of the coarse mesh
and straight-sided elements is not sufficient to resolve the discontinuity
structure and as a result the solution exhibits over- and under-shoot; the
$L^1$ error associated with this solution is $5.79\times 10^{-2}$.
However, the $q = 2$ and $q = 3$ solutions are extremely accurate even on
the coarse mesh due to the high-order elements that curve to conform to
the discontinuity structure; the $L^1$ error associated with each solution
is $1.15\times 10^{-3}$ and $5.50\times 10^{-4}$, respectively.

\begin{figure}
 \centering
 \input{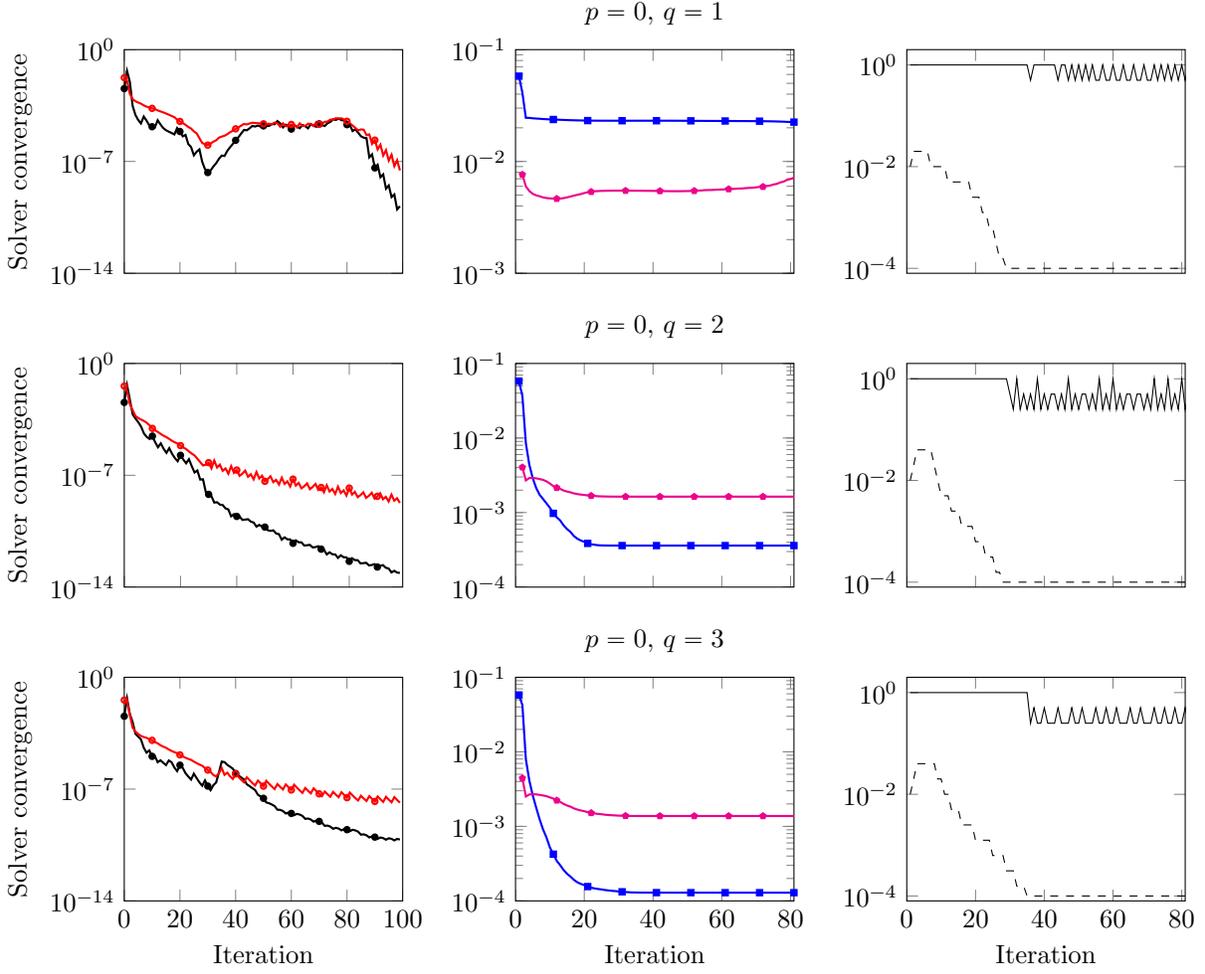}
 \caption{Convergence of the DG residual $\norm{\rbm(\ubm,\xbm)}$
          (\ref{line:advec2d_beta0_c2s1:R0}), enriched DG residual
          $\norm{\Rbm(\ubm,\xbm)}$ (\ref{line:advec2d_beta0_c2s1:R1a}),
          optimality condition $\norm{\cbm(\ubm,\xbm)}$
          (\ref{line:advec2d_beta0_c2s1:dLdY}), mesh distortion
          $\norm{\kappa\Rbm_\msh(\xbm)}$ (\ref{line:advec2d_beta0_c2s1:R1b}),
          and control of the regularization parameter
          (\ref{line:advec2d_beta0_c2s1:lam}) and step size
          (\ref{line:advec2d_beta0_c2s1:alpha}) for the
          tracking method applied to the advection equation with the
          trigonometric shock. Even though Newton convergence is
          not achieved, the proposed solver converges the KKT system 
          to tight tolerances in a reasonable number of iterations.}
 \label{fig:advec2d_beta0_c2s1_conv}
\end{figure}

For all polynomial degrees considered, the solver is able to converge the
KKT system to relatively tight tolerances ($\norm{\rbm(\ubm,\xbm)} < 10^{-10}$
and $\norm{\cbm(\ubm,\xbm)}<10^{-7}$). As the polynomial degree increases,
the enriched residual and mesh distortion converge to increasingly small
values and the overall convergence of the solver becomes cleaner
(Figure~\ref{fig:advec2d_beta0_c2s1_conv}). This comes from the
improved Hessian approximation, which comes from the mesh distortion
and enriched residual converging to smaller values and justifies
dropping of the second term in (\ref{eqn:gnhess}). The regularization
parameter is adapted to control the size of $\Delta\xbm$ produced
from the linear solve and the line search ensures sufficient progress
is made with respect to the $\ell_1$ merit function
(Figure~\ref{fig:advec2d_beta0_c2s1_conv}). For this problem, a
non-unity step size is only required once the regularization parameter
is small, which is required near convergence to have a decent approximation
to the Hessian. 

\begin{figure}
 \centering
 \input{py/advec2d_beta0_c2s2_nel8x4_porder0x2_conv.tikz}
 \caption{Convergence of the constraint (\textit{left}),
          optimality condition (\text{middle}), and objective
          function (\textit{right}) of the shock tracking optimization
          problem for linear advection with the trigonometric shock
          ($p = 0$, $q = 2$) when the non-smooth upwind flux $\Hcal_\up$
          (\ref{line:advec2d_beta0_c2s2_nel8x4_porder0x2:upwind})
          and smoothed ($a = 10$) upwind flux $\Hcal_\up^a$
          (\ref{line:advec2d_beta0_c2s2_nel8x4_porder0x2:upwind_smooth})
          are used as the numerical flux function.}
 \label{fig:advec2d_beta0_c2s2_nel8x4_porder0x2_conv}
\end{figure}
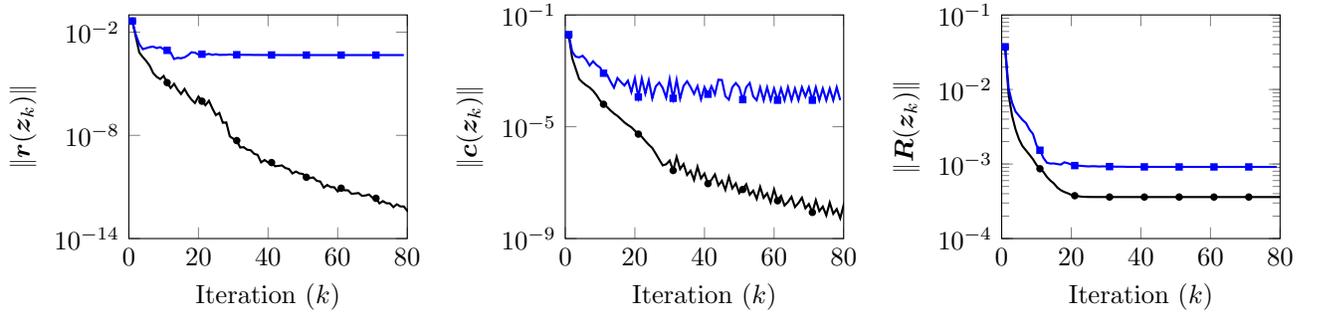

For this problem, the mesh smoothing parameter ($\kappa$) and minimum
value for the regularization parameter ($\gamma_{min}$) were important
for the $q = 1$ tracking problem; for $\kappa < 10^{-2}$ or
$\gamma_{min} < 10^{-3}$, the tracking algorithm would continually
collapse elements near the discontinuity as it pushed nodes to this
region to improve the faceted approximation of the trigonometric
discontinuity surface.
For $q>1$, the smoothing parameter played no role as the mesh is
well-conditioned even with $\kappa = 0$ and $\gamma_{min}$ can be
taken much smaller, e.g., $\gamma_{min} = 10^{-8}$, without adversely
affecting the performance of the solver or the final solution. The
numerical flux plays a more significant role with regard to the
behavior of the SQP solver, even for higher order elements. For
the non-smooth (with respect to domain deformation) numerical flux,
the solver does not converge because the first-order information is
meaningless in regions near the kink in the numerical flux, whereas
the convergence is much faster and cleaner for the smoothed upwind flux
(Figure~\ref{fig:advec2d_beta0_c2s2_nel8x4_porder0x2_conv}).

\subsection{Time-dependent, inviscid Burgers' equation}
\label{sec:num-exp:burg}
Next, we consider the time-dependent, inviscid Burgers' equation that
governs nonlinear advection of a scalar quantity
$\func{U}{\dom}{\Rbb}$ through the space-time domain
$\dom \coloneqq \Tcal \times \bar\dom$,
\begin{equation} \label{eqn:burg-sptm}
 \pder{U}{t} + U\pder{U}{x} = 0 \quad \text{in } \Omega, \qquad
 U = U_\infty \quad \text{on } \Gamma_i(U)
\end{equation}
where $\Tcal \coloneqq (0, 1)$ is the time domain,
$\bar\dom \coloneqq (-1,1)$ is the one-dimensional spatial domain,
$\Gamma_i(U)\coloneqq\{(t,x)\in\partial\Omega \mid \beta(U) \cdot n(x) < 0\}$
is the inflow boundary, $\beta(U)$ is the space-time flow direction defined as
\begin{equation} \label{eqn:burg-advec}
 \beta(U) = \begin{bmatrix} 1 \\ U \end{bmatrix},
\end{equation}
$\func{n}{\partial\dom}{\Rbb^2}$ is the space-time unit outward normal,
and $\func{U_\infty}{\partial\dom}{\Rbb}$ is the inflow boundary condition.
This fits the form of a general, steady conservation law
 over the space-time domain $\Omega$ with flux function
\begin{equation}
 F(U) = \begin{bmatrix} U & \ds{\frac{U^2}{2}} \end{bmatrix}.
\end{equation}
In the space-time setting, Burgers' equation has the same form as linear
advection with the solution-dependent advection field (\ref{eqn:burg-advec}).
Similar to the advection equation, we consider a pure space-time upwind
numerical flux since it satisfies conditions (i)-(iii)
(Section~\ref{sec:disc:numflux})
\begin{equation}
 \Hcal(U^+, U^-, n) = \Hcal_\up\left(U^+, U^-, n;
                                      \beta\left(\frac{U^++U^-}{2}\right)\right)
\end{equation}
and yields a stable DG discretization suitable for tracking. We also consider
the smoothed version to satisfy condition (iv)
\begin{equation}
 \Hcal(U^+, U^-, n) = \Hcal_\up^a\left(U^+, U^-, n;
                                     \beta\left(\frac{U^++U^-}{2}\right)\right).
\end{equation}
We consider the following piecewise quadratic boundary condition
\begin{equation} \label{eqn:burg-sptm-bndcnd}
 U_\infty(t, x) = 2(x+1)^2 (1-H(x)),
\end{equation}
which is enforced via the pure upwind numerical flux
evaluated at the interior state $U^+$ and boundary state
$U^\partial = U_\infty$.

\begin{figure}
 \centering
 \input{tikz/burg1d_sptm_c1s1_nel8x4_porder1x1.tikz}
 \caption{Space-time solution of one-dimensional, inviscid Burgers' equation
          using the tracking method at various iterations throughout the
          solution procedure using a $p = q = 1$ basis for the solution and
          mesh. The method collapses $7$ elements throughout the solution
          procedure and tracks a faceted approximation of the shock trajectory
          using the $q=1$ mesh in only $40$ iterations.}
 \label{fig:burg1d_sptm_c1s1_porder1x1}
\end{figure}
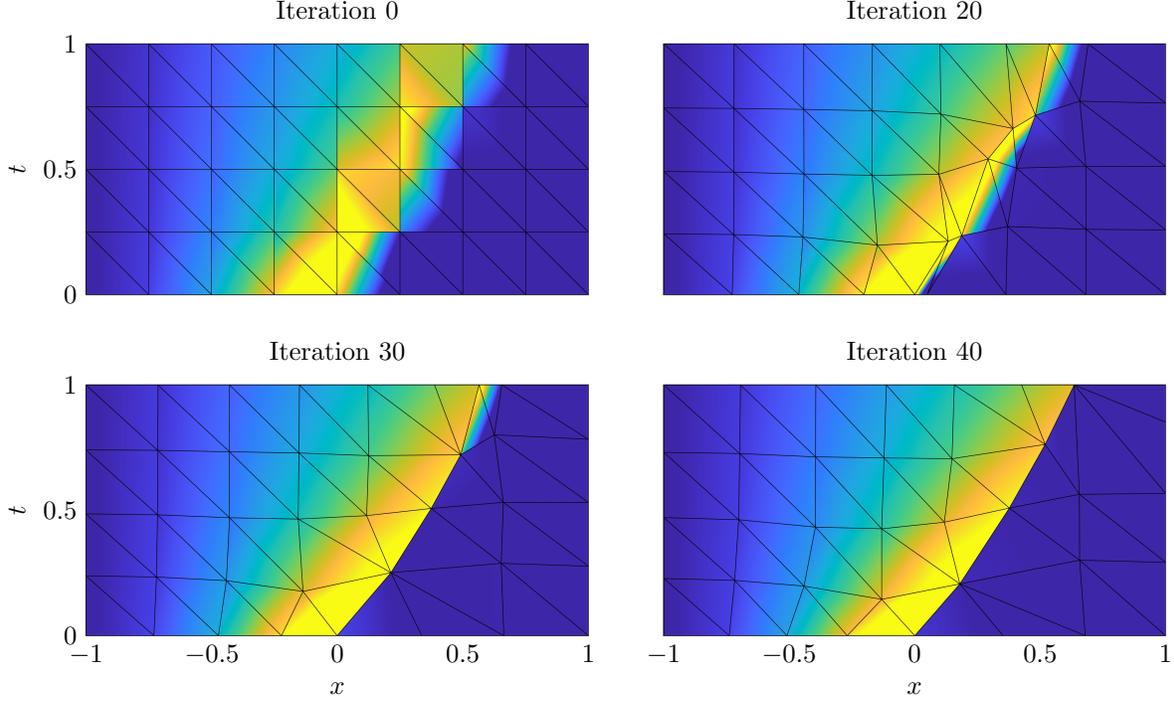

For this problem, we consider a solution and mesh basis of equal
polynomial degree $p = q$ up to $p = q = 4$. The reference mesh is taken
to be a uniform triangular mesh of the domain with $64$ elements. The mesh
and solution are initialized according to Section~\ref{sec:practical:init},
including continuation in the polynomial degree. The SQP solver is used with
$\lambda$ chosen adaptively (Section~\ref{sec:solver:regparam}).
The various DG/tracking parameters are set as follows:
$\kappa = 10^{-4}$ (mesh smoothing), smoothed upwind numerical flux
($a=10$), $\gamma_0 = 10^{-1}$ and $\gamma_{min} = 10^{-4}$ (regularization
parameter adaptivity), and $N_{max} = 100$, $\epsilon_1 = 10^{-6}$,
$\epsilon_2 = 10^{-10}$ (termination criteria). We initialize the
$p = q = 1$ simulation from the $p = 1$ DG solution on the reference mesh.
Even though the initial mesh is far from tracking the discontinuity
(some faces are nearly orthogonal to the discontinuity, rather than
parallel to it), our method tracks a faceted approximation to the
discontinuity in only $40$ iterations, requiring $7$ element collapses
(Figure~\ref{fig:burg1d_sptm_c1s1_porder1x1}).

\begin{figure}
 \centering
 \input{tikz/burg1d_sptm_c1s1_nel8x4.tikz}
 \caption{Space-time solution of one-dimensional, inviscid Burgers' equation
          using the proposed tracking method with a
          $p = q = 1$ (\textit{top}),
          $p = q = 2$ (\textit{middle}), and
          $p = q = 3$ (\textit{bottom}) basis for the solution and mesh
          with (\textit{left}) and without (\textit{right}) element boundaries.
          In all cases, the finite element solution provides a high-quality
          approximation to the true solution by tracking the discontinuity
          with a well-conditioned mesh. This is particularly true for the
          high-order elements that curve to the space-time trajectory of the
          discontinuity.}
 \label{fig:burg1d_sptm_c1s1}
\end{figure}
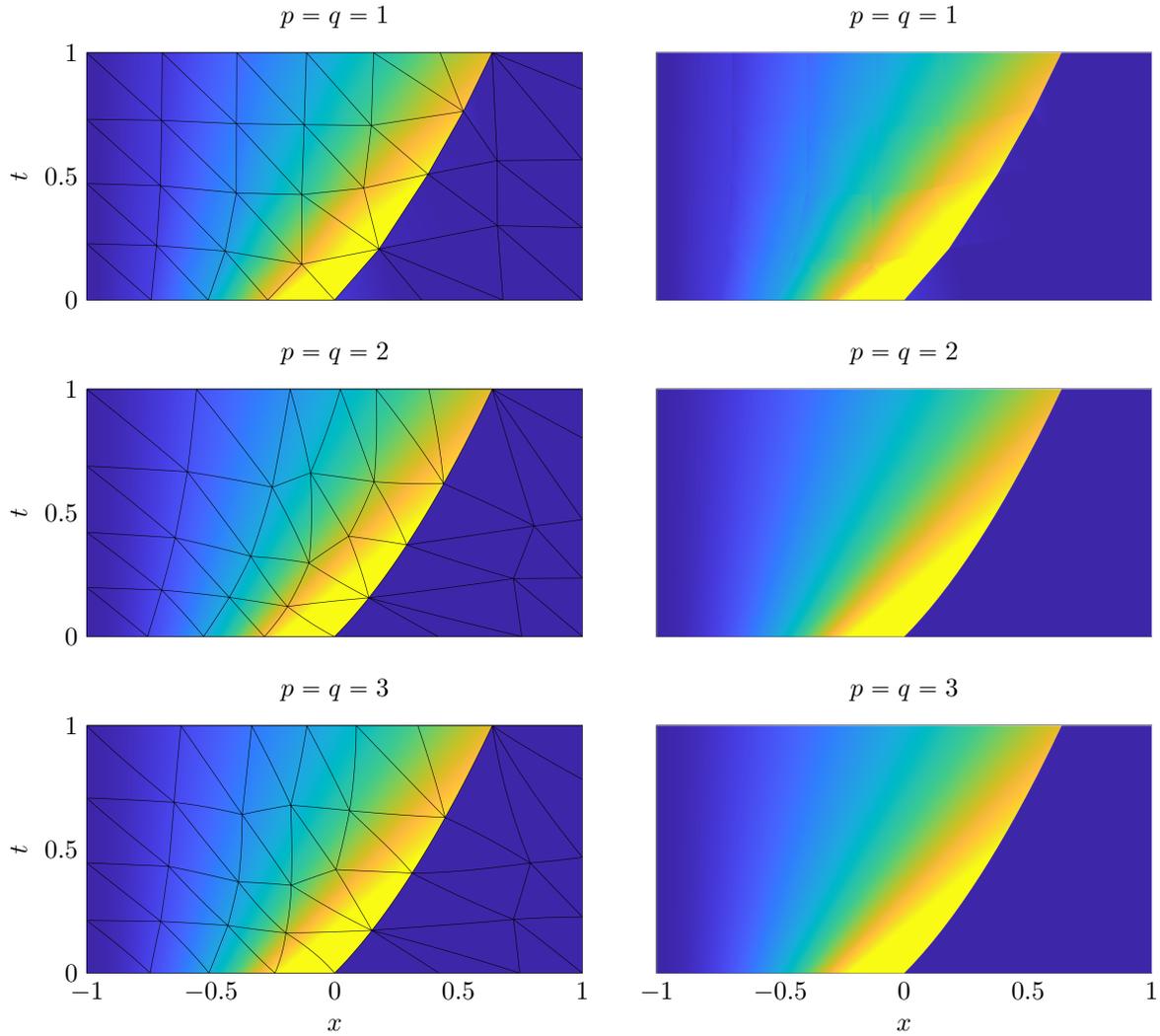

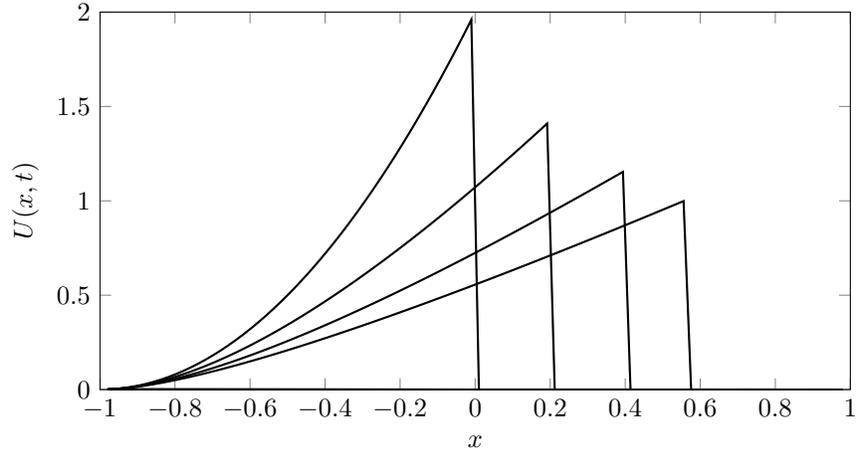
\begin{figure}
 \centering
 \input{py/burg1d_sptm_c1s1_nel8x4_porder4x4_slice.tikz}
 \caption{Temporal slices of the $p = q = 4$ tracking solution to the
          inviscid Burgers' equation at times
          $t = 0.05, 0.35, 0.65, 0.95$.}
 \label{fig:burg1d_sptm_c1s1_slice}
\end{figure}

The mesh and solution for the high-order elements ($p = q > 1$) are initialized
from the $p = q = 1$ tracking mesh and solution. These high-order
approximations provide high-quality approximations of the discontinuous
space-time solution on the coarse mesh
(Figure~\ref{fig:burg1d_sptm_c1s1}).
This can further be verified
from the temporal slices (Figure~\ref{fig:burg1d_sptm_c1s1_slice}),
which show the moving discontinuity is perfectly tracked and the
solution is smooth and non-oscillatory away from the discontinuity.

\begin{figure}
 \centering
 \input{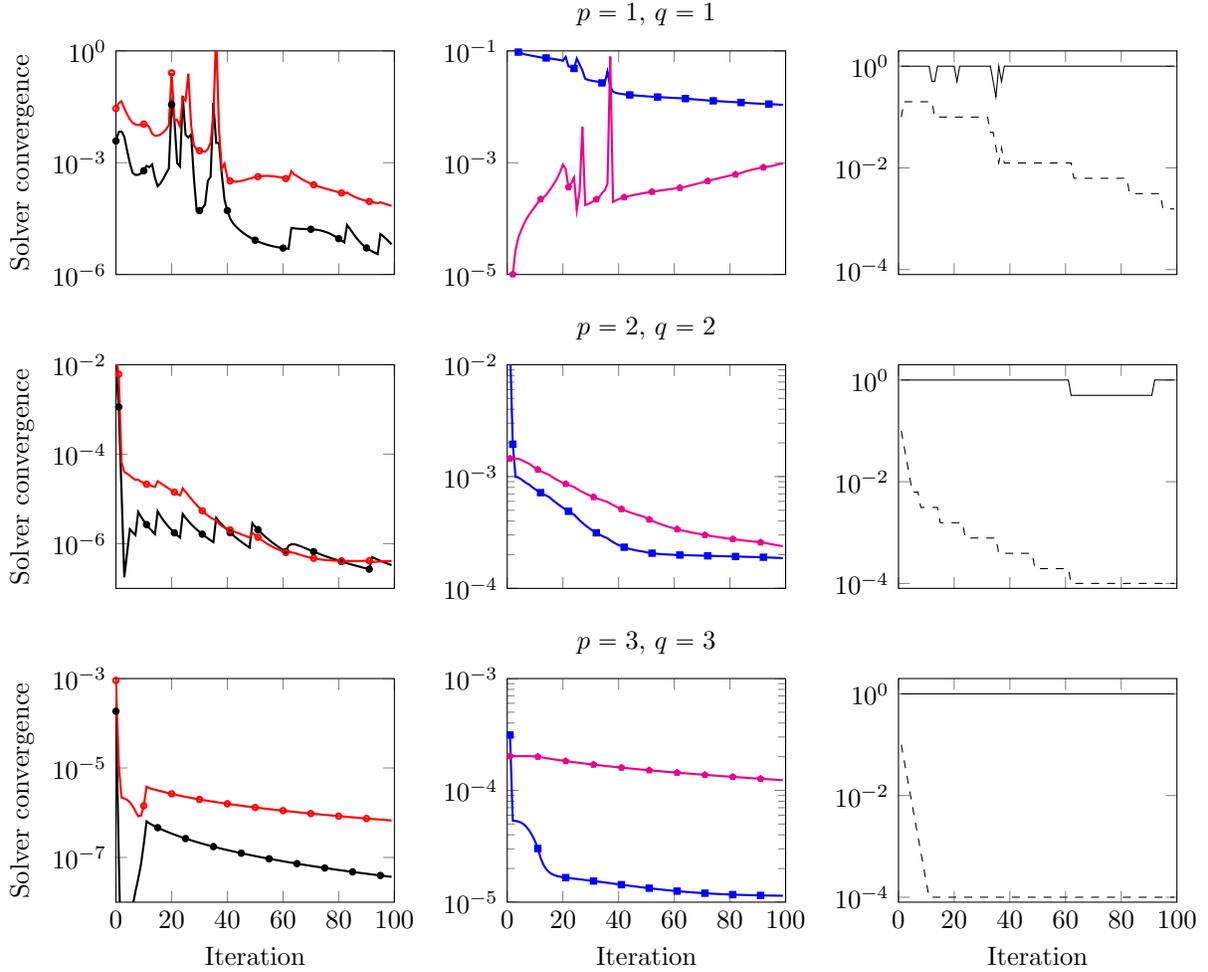}
 \caption{Convergence of the DG residual $\norm{\rbm(\ubm,\xbm)}$
          (\ref{line:burg1d_sptm_c1s1:R0}), enriched DG residual
          $\norm{\Rbm(\ubm,\xbm)}$ (\ref{line:burg1d_sptm_c1s1:R1a}),
          optimality condition $\norm{\cbm(\ubm,\xbm)}$
          (\ref{line:burg1d_sptm_c1s1:dLdY}), mesh distortion
          $\norm{\kappa\Rbm_\msh(\xbm)}$ (\ref{line:burg1d_sptm_c1s1:R1b})
          and control of the regularization parameter
          (\ref{line:burg1d_sptm_c1s1:lam}) and step size
          (\ref{line:burg1d_sptm_c1s1:alpha}) for the
          tracking method applied to the inviscid Burgers' equation.}
 \label{fig:burg1d_sptm_c1s1_conv}
\end{figure}

The convergence of the SQP solver in general is similar to convergence behavior 
observed with the linear advection equation with a trigonometric shock in that the
higher the polynomial degree, the cleaner the convergence, and the
enriched residual and mesh distortion converge to smaller absolute
values, indicating a solution that provides a better approximation
to the continuous weak form on a higher quality mesh
(Figure~\ref{fig:burg1d_sptm_c1s1_conv}). The line search
is active at some intermediate iterations, but plays less of a role
near convergence. Furthermore, the same observations regarding the
tracking parameters made for the advection equation with the trigonometric
shock hold for this problem: the smoothed upwind flux
significantly improves convergence of the solver, and if $\kappa$ or
$\gamma_{min}$ are any smaller, too many elements will collapse onto
the discontinuity surface for $p = q = 1$. 

\subsection{Steady, compressible Euler equations}
\label{sec:num-exp:euler}
The Euler equations govern the steady flow of an inviscid, compressible fluid
through a domain $\Omega \subset \Rbb^d$
\begin{equation} \label{eqn:euler}
  (\rho v_j)_{,j} = 0, \quad (\rho v_i v_j + p \delta_{ij})_{,j} = 0, \quad
  (\rho H v_i)_{,j} = 0 \quad \text{ in } \Omega
\end{equation}
where $\func{\rho}{\Omega\times (0, T)}{\Rbb_+}$ is the density of the
fluid, $\func{v_i}{\Omega\times (0, T)}{\Rbb}$ for $i = 1, \dots, d$
is the velocity of the fluid in the $i$th coordinate direction, and
$\func{E}{\Omega\times (0, T)}{\Rbb_+}$ is the total energy of the
fluid. The square velocity $\func{q^2}{\Omega\times (0, T)}{\Rbb_+}$ and
kinetic energy $\func{k_e}{\Omega\times (0, T)}{\Rbb_+}$ of the fluid
are
\begin{equation}
 q^2 = v_i v_i, \quad
 k_e = \frac{1}{2}\rho q^2.
\end{equation}
The enthalpy of the fluid
$\func{H}{\Omega\times (0, T)}{\Rbb_+}$ is defined as
\begin{equation}
 \rho H = \rho E + P,
\end{equation}
where $\func{P}{\Omega\times (0, T)}{\Rbb_+}$ is the pressure.
For a calorically ideal fluid, the pressure and energy are related via
the ideal gas law
\begin{equation}
 P = (\gamma-1)(\rho E - k_e)
\end{equation}
and the speed of sound is
\begin{equation}
 c = \sqrt{\gamma P/\rho}.
\end{equation}
The density, velocity, and energy are combined into a vector of conservative
variables $U$ and the Euler equations take the form of an inviscid conservation
law (\ref{eqn:claw-phys}) with flux function $F(U)$
\begin{equation}
 U = \begin{bmatrix} \rho \\ \rho v \\ \rho E \end{bmatrix}, \qquad
 F(U) = \begin{bmatrix} \rho v^T \\ \rho vv^T + P I_2 \\ \rho H v^T \end{bmatrix},
\end{equation}
where $I_2$ is the $2\times 2$ identity matrix.

We use Roe's flux \cite{roe1981approximate} as the numerical flux function
to ensure a stable discretization suitable for tracking since it satisfies
(i)-(iii) (Section~\ref{sec:disc:numflux}); however, it is not smooth with
respect to variations in the domain deformation. To eliminate entropy
violating rarefaction shocks and improve the smoothness with respect to
the domain deformations, we modify the interior numerical fluxes with the
Harten-Hyman entropy fix \cite{harten1983self}. This leads
to a numerical flux function that violates (iii), as observed in
\cite{corrigan2019moving}, but recovers smoothness with respect
to mesh deformation, which proves to be a reasonable trade-off.

We consider three types of boundaries: slip wall ($\Gamma_w$), supersonic
inflow ($\Gamma_D$), and supersonic outflow ($\Gamma_N$). The supersonic
inflow is also known as a farfield or Dirichlet condition and the supersonic
outflow is a Neumann condition. For a slip wall ($v\cdot n = 0$), the boundary
state is defined as
\begin{equation}
 U_w^\partial(U,n) \coloneqq
 \begin{bmatrix}
  \rho \\ \rho v^- \\ \rho E,
 \end{bmatrix}
\end{equation}
where $v^- = (v - 2 v\cdot n)n$ is the velocity reflected about the normal.
For a supersonic inflow, all characteristics are coming into the domain
and the boundary state depends solely on the prescribed density $\rho_\infty$,
velocity $v_\infty$, and pressure $p_\infty$
\begin{equation}
 U_D^\partial(U) \coloneqq
 \begin{bmatrix}
  \rho_\infty \\ \rho_\infty v_\infty \\
  \frac{p_\infty}{\gamma-1}+\frac{\rho_\infty}{2} v_\infty \cdot v_\infty
 \end{bmatrix}.
\end{equation}
Finally, at a supersonic outflow, all characteristics are leaving the domain
and the boundary state is taken from the interior
\begin{equation}
 U_N^\partial(U) \coloneqq U.
\end{equation}
All boundary conditions are enforced via the pure Roe flux (without entropy
fix) evaluated at the interior state $U^+$ and appropriate boundary
state
\begin{equation}
 U^\partial(U, n) = 
 \begin{cases}
  U_w^\partial(U,n) & \text{on } \Gamma_w \\
  U_D^\partial(U) & \text{on } \Gamma_D \\
  U_N^\partial(U) & \text{on } \Gamma_N.
 \end{cases}
\end{equation}

A useful property of the inviscid flows is the enthalpy is constant
throughout the domain. Therefore, to quantify the error in the numerical
solution obtained using the proposed method, we will use the deviation
of the flow enthalpy from the inflow enthalpy,
$H_\infty \coloneqq \frac{\gamma}{\gamma-1}\frac{p_\infty}{\rho_\infty} +
                    \frac{1}{2} v_\infty \cdot v_\infty$
\begin{equation}
 e_H(U) \coloneqq
 \sqrt{\frac{\int_\Omega (H(U) - H_\infty)^2}{\int_\Omega \, dV}}.
\end{equation}

\subsubsection{Supersonic flow over wedge}
\label{sec:num-exp:euler:wedge}
\begin{figure}[H]
 \centering
 \input{py/wedge0_geom.tikz}
 \caption{Geometry and boundary conditions of the wedge problem. Boundary
          conditions: slip wall (\ref{line:wedge0:wall}), supersonic inflow
          with $\rho_\infty = 1.4$, $v_\infty = (2, 0)$, $p_\infty = 1$
          ($M_\infty = 2$) (\ref{line:wedge0:supin}), and supersonic outflow
          (\ref{line:wedge0:supout}).}
 \label{fig:wedge0_geom}
\end{figure}
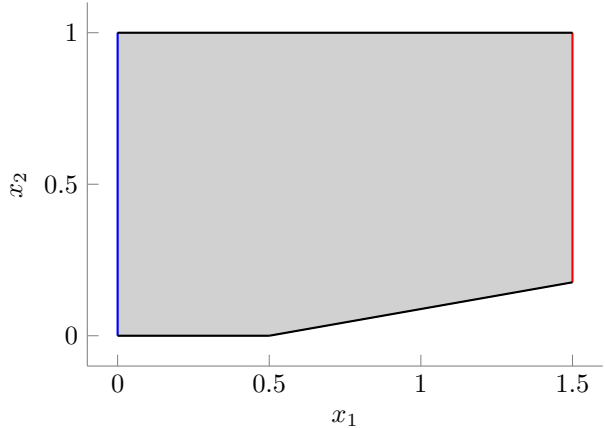

\begin{figure}
 \centering
 \input{tikz/euler2d0_wedge0_c1s1_nel6x4_porder0x1.tikz} \\
 \input{tikz/cbar_parula_1p65x2.tikz}
 \caption{Solution (Mach) of Euler equations over a wedge (supersonic regime)
          using the tracking method at various iterations throughout the
          solution procedure using $p = 0$ basis for the solution and $q = 1$
          basis for the mesh. The method converges to nearly the exact
          solution (enthalpy error $e_H = 7.94\times 10^{-10}$), in only
          $20$ iterations.}
 \label{fig:euler2d0_wedge0_c1s1}
\end{figure}
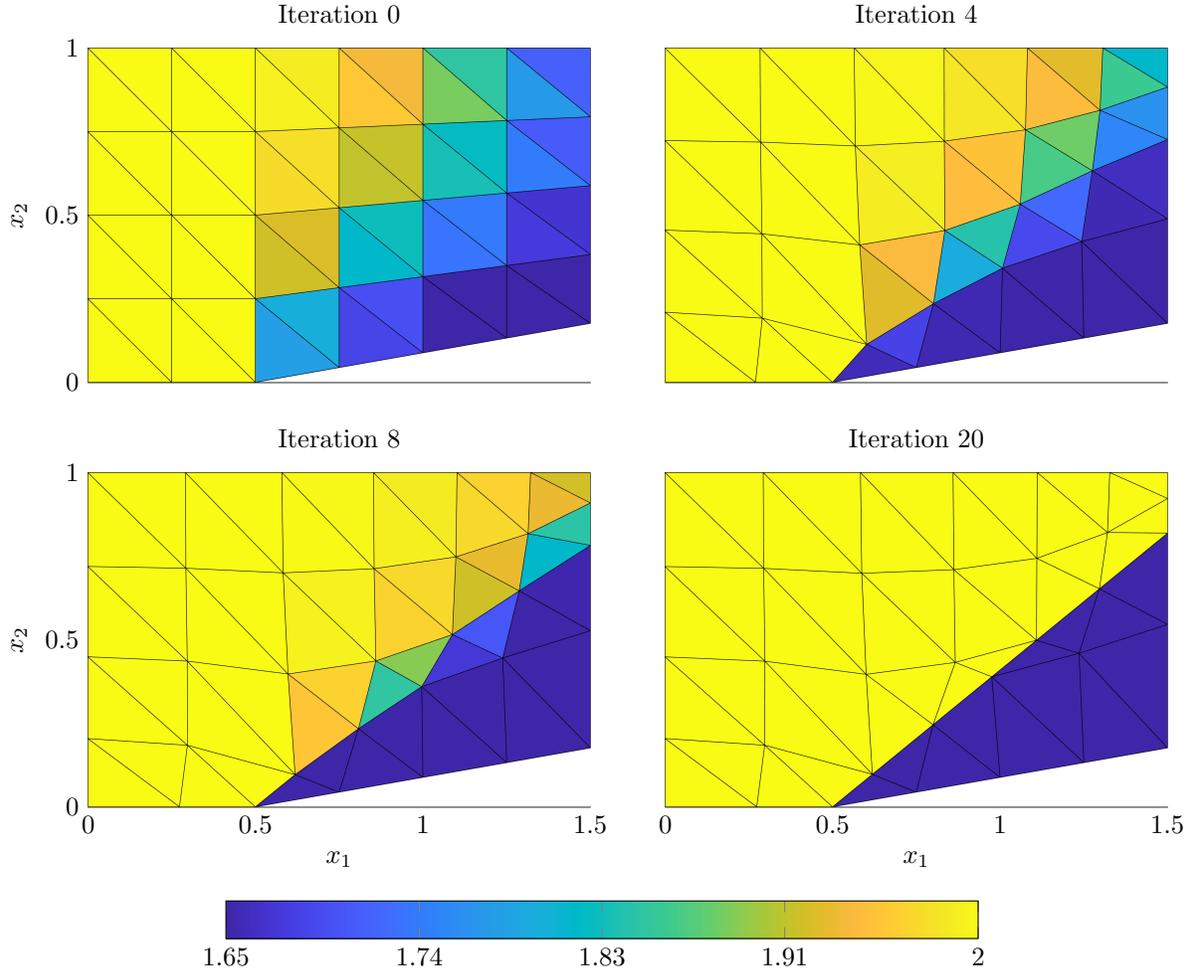

First we consider supersonic flow
($M_\infty = 2$) over a $\theta=10^\circ$ inclined plane
(Figure~\ref{fig:wedge0_geom}). Since all wall boundaries are
straight-sided and the incoming flow is uniform, the flow is
piecewise constant and the shocks are straight.
Therefore, we apply the tracking method with a $p = 0$ polynomial basis
for the solution and $q = 1$ basis for the mesh. The reference mesh is taken
to be a uniform triangular mesh of the domain with $48$ elements.
The solver is initialized with the $p=0$ DG solution on the reference
mesh (Section~\ref{sec:practical:init}).  The SQP solver is used with
$\lambda$ chosen adaptively (Section~\ref{sec:solver:regparam}).
The various DG/tracking parameters are set as follows:
$\kappa = 0$ (no mesh smoothing), Roe flux (without entropy fix),
$\gamma_0 = 1$ and $\gamma_{min} = 10^{-8}$ (regularization parameter
adaptivity), and $\epsilon_1 = 10^{-8}$, and $\epsilon_2 = 10^{-12}$
(termination criteria). After only $20$ iterations, the mesh perfectly
tracks the shock (Figure~\ref{fig:euler2d0_wedge0_c1s1}) and the DG
solution closely matches the exact solution; the enthalpy error of the
solution is $e_H = 7.94\times 10^{-10}$. The solver performs similarly
to the advection equation with a straight discontinuity and therefore
a full discussion is omitted for brevity.

\subsubsection{Supersonic flow over airfoil}
\label{sec:num-exp:euler:naca0}
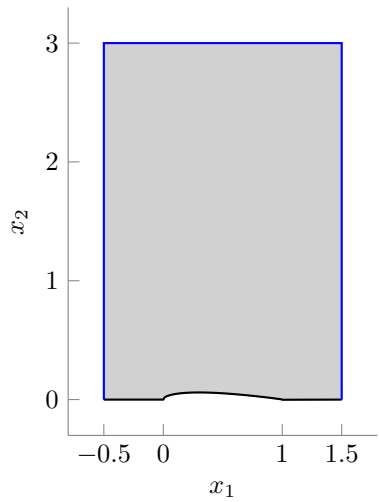
\begin{figure}
 \centering
 \input{py/naca0012cfg1_geom.tikz}
 \caption{Geometry and boundary conditions of the airfoil problem. Boundary
          conditions: slip wall (\ref{line:naca0012cfg1:wall}) and farfield
          (characteristic) conditions with $\rho_\infty = 1.4$,
          $v_\infty = (M_\infty, 0)$, $p_\infty = 1$ ($M_\infty = 0.85$
          for the transonic case and $M_\infty = 1.5$ for the supersonic case).
          (\ref{line:naca0012cfg1:farf}).}
 \label{fig:naca0012cfg1_geom}
\end{figure}

Next, we apply the proposed tracking method to solve for supersonic
flow over a NACA0012 airfoil (Figure~\ref{fig:naca0012cfg1_geom}).
This is a difficult problem because there are two distinct shocks that
must be resolved: a bow shock ahead of the leading edge and an oblique
shock off the tail. The reference mesh is taken as an unstructured
triangular mesh of the domain with $160$ elements generated by
DistMesh \cite{persson2004simple}. The mesh and solution are initialized
according to Section~\ref{sec:practical:init}, including continuation in the
polynomial degree. The SQP solver is used with $\lambda$ chosen
adaptively (Section~\ref{sec:solver:regparam}).
The various DG/tracking parameters are set as follows:
$\kappa = 10^{-4}$ (mesh smoothing), Roe flux with entropy fix,
$\gamma_0 = 1$ and $\gamma_{min} = 1$ (regularization
parameter adaptivity), and $N_{max} = 100$, $\epsilon_1 = 10^{-6}$,
$\epsilon_2 = 10^{-10}$ (termination criteria).

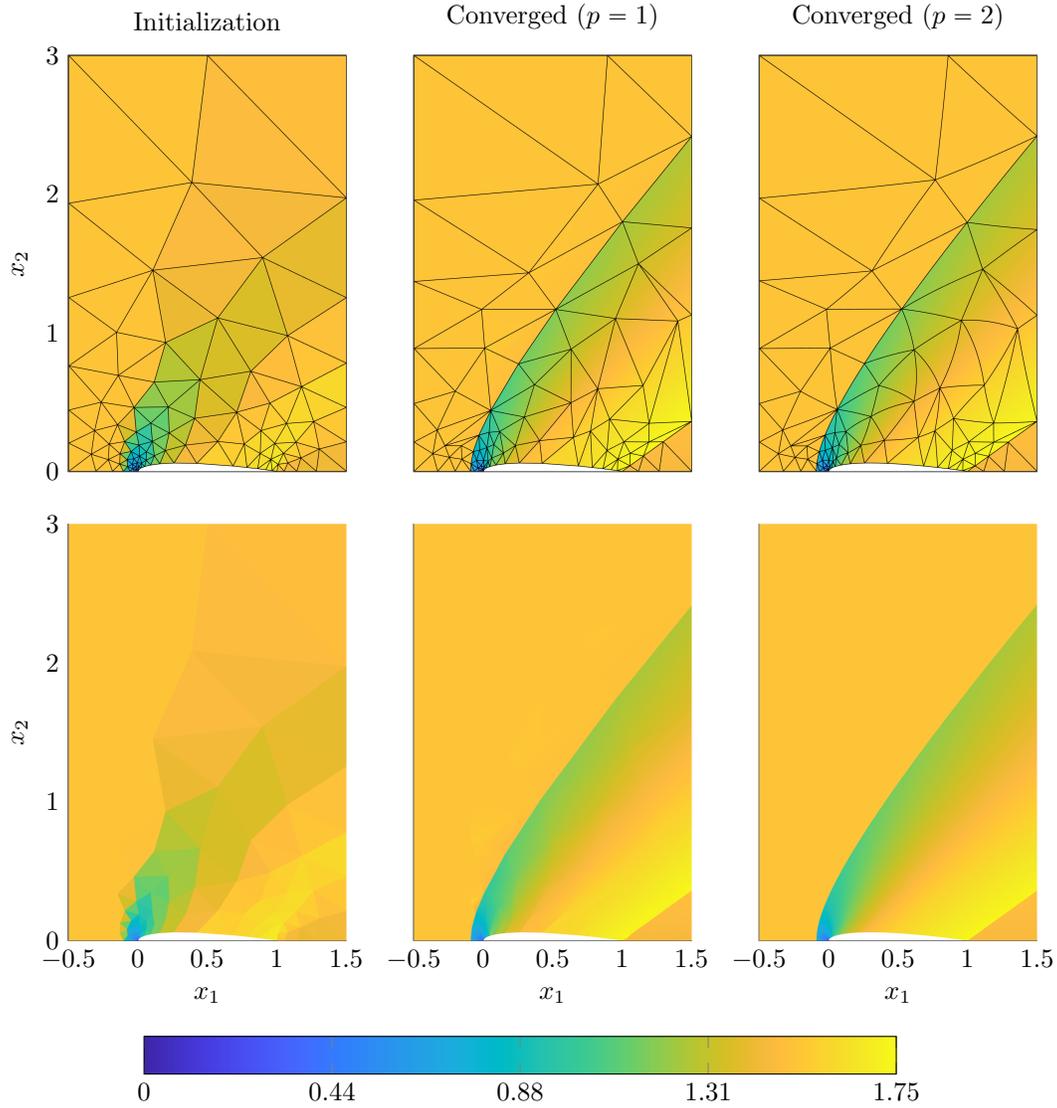
\begin{figure}
 \centering
 \input{tikz/euler2d0_naca0012cfg1_nref1_c1s1.tikz} \\
 \input{tikz/cbar_parula_0x1p75.tikz}
 \caption{Solution (Mach) of Euler equations over the NACA0012 airfoil
          ($M_\infty = 1.5$) using the proposed tracking method with a
          $p = q = 1$ (\textit{center}) and
          $p = q = 2$ (\textit{right}) basis for the solution and mesh
          with (\textit{top}) and without (\textit{bottom}) element boundaries.
          In both cases, the tracking procedure successfully tracks the
          shocks given the resolution in the finite element space, despite
          the initial mesh and solution (\textit{left}) being far
          from aligned with the shock. The high-order ($p = q = 2$)
          basis yields an accurate approximation to the flow on
          the coarse mesh while the low-order ($p = q = 1$) solution is
          under-resolved.}
 \label{fig:euler2d0_naca0012cfg1_nref1_c1s1}
\end{figure}

\begin{figure}
 \centering
 \input{py/euler2d0_naca0012cfg1_nref1p3_c1s1_slice.tikz}
 \caption{Slices of density (\textit{left}), Mach number (\textit{center}),
          and pressure (\textit{right}) of the $p = q = 3$ tracking solution
          along the curve
          $\Gamma \coloneqq \{(s, 0.14) \mid s \in (-0.5,1.5)\}$ for the
          NACA problem. The discontinuities are captured perfectly between
          DG elements and the solution is smooth and non-oscillatory away
          from the discontinuities, indicating that the solution is
          well-resolved and the discontinuities are successfully tracked.}
 \label{fig:euler2d0_naca0012cfg1_nref1p3_c1s1_slice}
\end{figure}
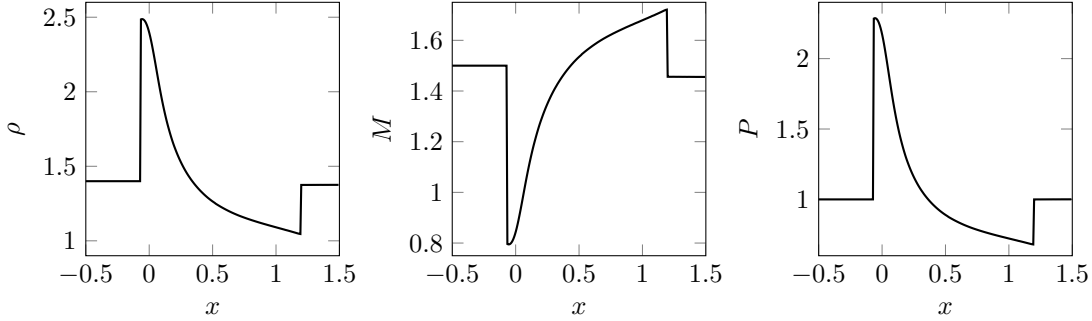

For the $p = q = 1$ simulation, the proposed method tracks a faceted
approximation to the discontinuity in only $100$ iterations, requiring
$11$ element collapses and maintain high-quality elements, despite the
coarse elements in the initial mesh that do not conform to either shock
(Figure~\ref{fig:euler2d0_naca0012cfg1_nref1_c1s1}).
Even though the discontinuities are successfully tracked, the solution
is under-resolved, particularly near the body, since the large $p = 1$
elements are not sufficient to resolve the solution in this region.
As a result, the enthalpy error is large $e_H = 1.30\times 10^{-3}$.
The $p = q = 2$ and $p = q = 3$ tracking solutions provide highly accurate
approximations to the true flow even on this coarse mesh due to the high-order
resolution of the discontinuity surface with curved elements and high-order
approximation of the flow (Figure~\ref{fig:euler2d0_naca0012cfg1_nref1_c1s1});
the enthalpy error for the $p = q = 2$ tracking solution is
$e_H = 6.73\times 10^{-5}$ and for $p = q = 3$ is $e_H = 1.02\times 10^{-5}$.
The perfect (zero-thickness) capturing of the shocks and high-order
approximation of the solution can further be seen from the solution slices
in Figure~\ref{fig:euler2d0_naca0012cfg1_nref1p3_c1s1_slice}.

For this problem, the convergence of the solver is not as clean as the
other problems (Figure~\ref{fig:euler2d0_naca0012cfg1_nref1_c1s1_conv})
and heavily relies on the line search throughout the solution procedure
to ensure sufficient decrease in the $\ell_1$ merit function. However,
the solver still drives the KKT conditions to reasonable tolerances within
$100$ iterations for the $p = q = 3$ simulation. The spikes in the KKT
residuals around iterations $20$ and $60$ in the $p = q = 1$ simulation
are due to element collapses. Our choice of algorithmic parameters were
critical for this case. If $\kappa$ is significantly larger, the element
collapses would not occur and the mesh would not be able to track the
shocks, particularly the weaker shock at the trailing edge. Furthermore,
it was important to use the relatively large value of $\gamma_{min} = 1$,
otherwise the line search would be taxed more heavily in the later iterations
and completely stalled convergence.
\begin{figure}
 \centering
 \input{py/euler2d0_naca0012cfg1_nref1_c1s1_conv.tikz}
 \caption{Convergence of the DG residual $\norm{\rbm(\ubm,\xbm)}$
          (\ref{line:euler2d0_naca0012cfg1_c1s1:R0}), enriched DG residual
          $\norm{\Rbm(\ubm,\xbm)}$ (\ref{line:euler2d0_naca0012cfg1_c1s1:R1a}),
          optimality condition $\norm{\cbm(\ubm,\xbm)}$
          (\ref{line:euler2d0_naca0012cfg1_c1s1:dLdY}), mesh distortion
          $\norm{\kappa\Rbm_\msh(\xbm)}$
          (\ref{line:euler2d0_naca0012cfg1_c1s1:R1b})
          and control of the regularization parameter
          (\ref{line:euler2d0_naca0012cfg1_c1s1:lam}) and step size
          (\ref{line:euler2d0_naca0012cfg1_c1s1:alpha}) for the
          tracking method applied to solve the supersonic flow around
          the NACA airfoil. For this difficult problem, the convergence
          of the KKT system is not as clean as the other problems and
          relies heavily on the line search; however, the solver still
          tracks the shocks and returns an accurate flow as seen from
          the solution plots
          (Figures~\ref{fig:euler2d0_naca0012cfg1_nref1_c1s1}-%
                   \ref{fig:euler2d0_naca0012cfg1_nref1p3_c1s1_slice}).
          The spikes in the KKT conditions for $p = q = 1$ are due to
          element collapses.}
 \label{fig:euler2d0_naca0012cfg1_nref1_c1s1_conv}
\end{figure}

Finally, we mention that this is the first problem where our choice of
the elliptic regularization matrix with the coefficients chosen inversely
proportional to the size of the element in the reference domain
(Section~\ref{sec:solver:regmat}) is significant. The other choices
we explored included the identity matrix \cite{corrigan2019moving},
elliptic PDE stiffness matrix without volume-based weighting
\cite{corrigan2019convergence}, and other choices involving the
finite element mass matrix. All of these options performed similarly
in creating unacceptably large search direction in regions near the
leading and trailing edge (small elements) relative to other regions
in the mesh (large elements). Without a line search, these steps would
cause the mesh to entangle and the simulation to crash. With the line
search, the step size in regions with larger elements would be driven
nearly to zero causing the solver to fail to track the discontinuity
in these regions. From this, we conclude the combination of the elliptic
PDE stiffness matrix and the weighting of the coefficients inversely with
respect to the element size in the reference domain is important for
meshes with elements of significantly varying size.

\subsubsection{Transonic flow over airfoil}
\label{sec:num-exp:euler:naca1}
Next, we consider transonic flow ($M_\infty = 0.85$) over the same NACA0012
airfoil and domain from the previous section
(Figure~\ref{fig:naca0012cfg1_geom}). This problem has a shock attached
to the curved airfoil profile, which requires a nontrivial boundary mapping
(Section~\ref{sec:optim:bndconstr}) to ensure nodes slide along the airfoil
surface. We construct the mapping $\chi$ using the procedure in
\cite{zahr2020radapt}, i.e.,
the $x$-coordinates of the nodes on the surface are taken as optimization
parameters and the $y$-coordinates are determined from the expression
for the airfoil profile
\begin{equation}
 y(x) = 0.6\left(0.2960\sqrt{x}-0.126x-0.3516x^2+0.2843x^3-0.1036x^4\right).
\end{equation}
The reference mesh is taken as an unstructured triangular mesh of the domain
with $127$ elements generated by DistMesh \cite{persson2004simple}. The mesh
and solution are initialized according to Section~\ref{sec:practical:init},
including continuation in the polynomial degree. The SQP solver is used with
$\lambda$ chosen adaptively (Section~\ref{sec:solver:regparam}). The
various DG/tracking parameters are set as follows:
$\kappa = 10^{-3}$ (mesh smoothing), Roe flux with entropy fix,
$\gamma_0 = 10^{-2}$ and $\gamma_{min} = 1$ (regularization
parameter adaptivity), and $N_{max} = 30$, $\epsilon_1 = 10^{-6}$,
$\epsilon_2 = 10^{-10}$ (termination criteria).

For polynomial degrees $p = q = 1, 2, 3$, the tracking method accurately
tracks the attached discontinuity without requiring element collapses and
produces an accurate flow solution given the resolution of the approximation
space (Figure~\ref{fig:euler2d0_naca0012cfg1_nref0_c2s1}). The high resolution
of the $p = q = 3$ simulation can be seen through the slices of the solution
(Figure~\ref{fig:euler2d0_naca0012cfg1_nref0p3_c2s1_slice}): the discontinuity
is captured perfectly and the solution away from the discontinuity is smooth
and non-oscillatory. The convergence of the method has been thoroughly studied
in the previous sections and omitted for brevity.
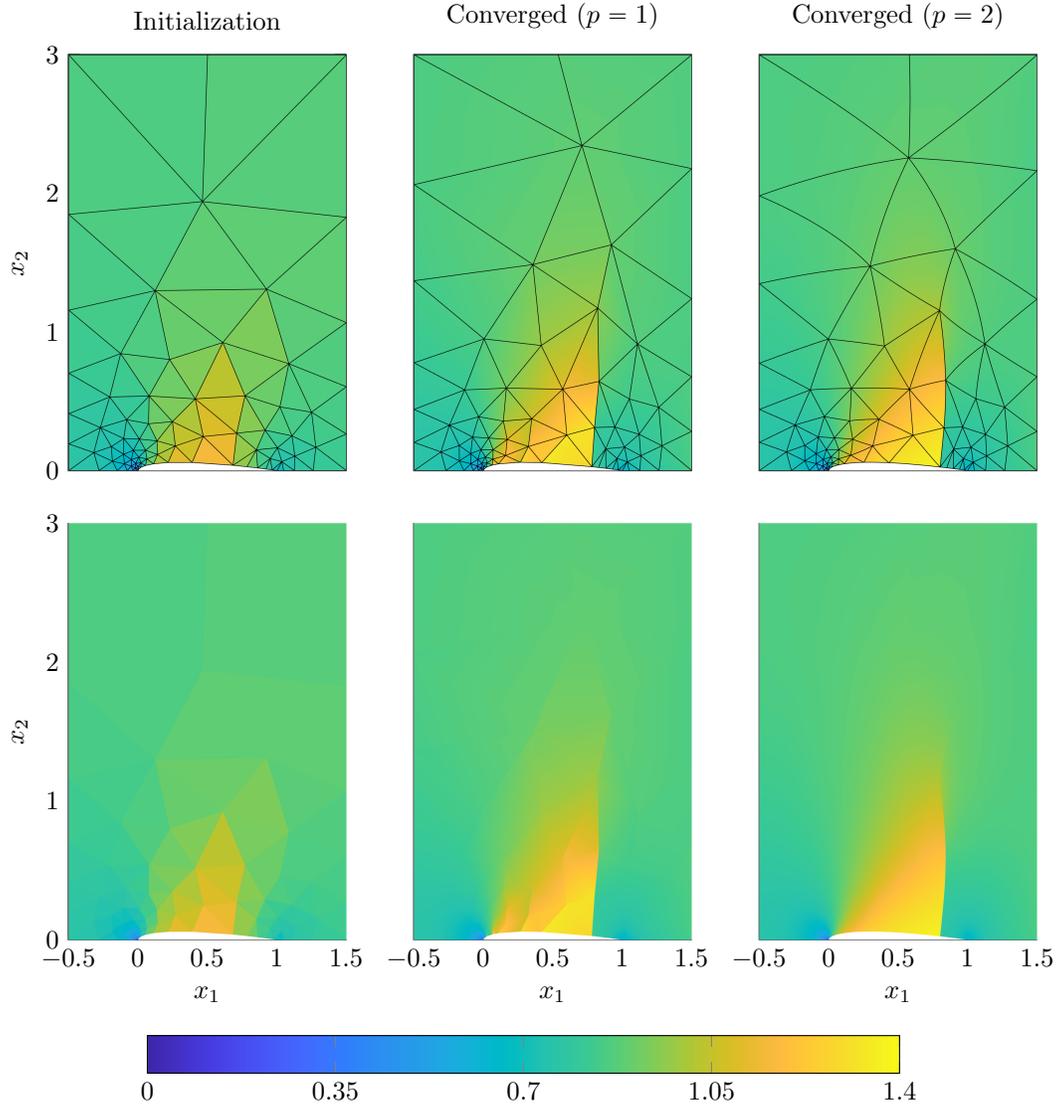
\begin{figure}
 \centering
 \input{tikz/euler2d0_naca0012cfg1_nref0_c2s1.tikz} \\
 \input{tikz/cbar_parula_0x1p4.tikz}
 \caption{Solution (Mach) of Euler equations over the NACA0012 airfoil
          ($M_\infty = 0.85$) using the proposed tracking method with a
          $p = q = 1$ (\textit{center}) and
          $p = q = 2$ (\textit{right}) basis for the solution and mesh
          with (\textit{top}) and without (\textit{bottom}) element boundaries.
          In both cases, the tracking procedure successfully tracks the
          shocks given the resolution in the finite element space, despite
          the initial mesh and solution (\textit{left}) being far
          from aligned with the shock.}
 \label{fig:euler2d0_naca0012cfg1_nref0_c2s1}
\end{figure}
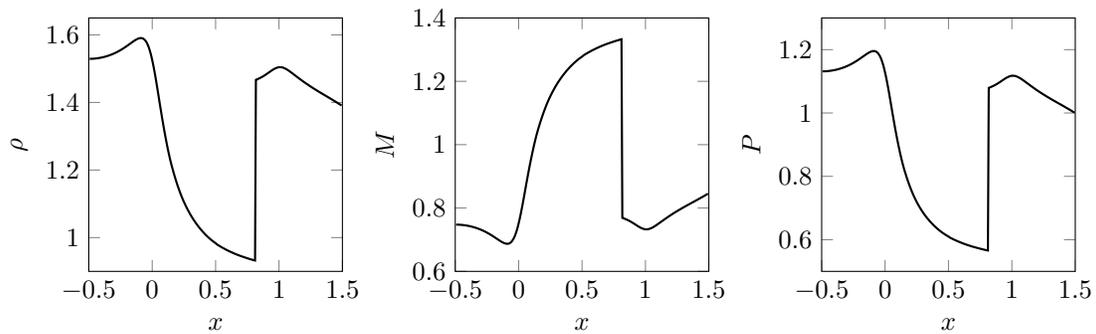
\begin{figure}
 \centering
 \input{py/euler2d0_naca0012cfg1_nref0p3_c2s1_slice.tikz}
 \caption{Slices of density (\textit{left}), Mach number (\textit{center}),
          and pressure (\textit{right}) of the $p = q = 3$ tracking solution
          along the curve
          $\Gamma \coloneqq \{(s, 0.14) \mid s \in (-0.5,1.5)\}$ for the
          NACA problem ($M_\infty=0.85$). The discontinuity is captured
          perfectly between DG elements and the solution is smooth and
          non-oscillatory away from the discontinuity, indicating that
          the solution is well-resolved.}
 \label{fig:euler2d0_naca0012cfg1_nref0p3_c2s1_slice}
\end{figure}

%


\section{Conclusions}
\label{sec:conclude}
We introduced an improved formulation of the optimization-based implicit
shock tracking method proposed in \cite{zahr2018optimization} and an
associated solver
that leverages the structure of the problem. The proposed optimization
problem minimizes the DG residual in an \textit{enriched} test and the
distortion of the mesh, constrained by the standard DG residual
(equal trial and test spaces).
The enriched residual is a practical
surrogate for the violation of the weak formulation of the conservation
law; its magnitude serves an error indicator for a DG solution.
Therefore, penalizing the enriched DG residual
promotes alignment of the element faces with discontinuities; otherwise,
the DG solution would oscillate about the discontinuities and provide a
poor approximation to the conservation law. The proposed solver for the
constrained optimization problem over the DG solution and mesh coordinates
is an SQP method that uses a Levenberg-Marquardt Hessian approximation and
is globalized via a line search on the $\ell_1$ merit function. The Hessian
approximation is regularized with the stiffness matrix corresponding to a
linear elliptic partial equation with coefficients inversely proportional
to the local element size and adaptively chosen regularization parameter,
which turns out to be significant for meshes containing elements of
significantly different size.

For problems where the finite element space contains the exact solution,
the SQP method exhibits Newton-like convergence to the exact solution,
which was demonstrated via linear advection of a scalar through a constant
advection field and supersonic, inviscid flow over a wedge. The framework
was also shown to be effective in accurately resolving more complex flows
with intricate shock structures using coarse, high-order meshes as demonstrated
using the inviscid Burgers' equation (time-dependent), and transonic and
supersonic flow over a NACA0012 airfoil.
For these problems, the convergence of the solver is slower, but still
drives the first-order optimality conditions to tight tolerances in
a reasonable number of iterations. 
 
Future work will develop iterative solvers and preconditioners for the
SQP linear system in (\ref{eqn:sqp1-linsys1}) to make the approach practical
for large-scale problems. We also intend to further improve the robustness of
the solver by incorporating pseudo-transient continuation to avoid the need to
initialize from a $p = 0$ DG solution and use continuation on the polynomial
degrees $p$ (DG solution) and $q$ (domain deformation). We will also develop
a method of lines and adaptive space-time approaches to handle more complex
time-dependent cases where a single space-time discretization may not be
feasible. Finally, we will demonstrate high-order convergence
of the method for inviscid flows in two- and three-dimensions and consider
more complex flows, including viscous and relevant 3D problems.

\section*{Acknowledgments}
This work was supported in part by the Director, Office of Science, Office of
Advanced Scientific Computing Research, U.S. Department of Energy under
Contract No. DE-AC02-05CH11231. The content of this publication does not
necessarily reflect the position or policy of any of these supporters, and
no official endorsement should be inferred.

\bibliographystyle{plain}
\bibliography{biblio,biblio_shock,biblio_track}

\end{document}

%% file: preamble.pap.tex
\usepackage[margin=1in]{geometry}

\usepackage{graphicx}
\usepackage{pstool}
\usepackage{wrapfig}
\usepackage{caption}
\usepackage{subcaption}
\usepackage[section]{placeins} 

\usepackage{fancyhdr} 
\usepackage{lastpage} 
\usepackage{extramarks} 

\usepackage{xcolor} 
\usepackage{enumerate} 
\usepackage{paralist} 
\usepackage{amsmath, amsthm, amssymb, mathtools}
\usepackage{mathabx, pifont, stmaryrd} 
\usepackage[explicit]{titlesec} 
\usepackage{etoolbox} 
\usepackage{relsize}
\usepackage{bibentry} 
\makeatletter\let\saved@bibitem\@bibitem\makeatother 
\usepackage[colorlinks, bookmarksopen, bookmarksnumbered,
            citecolor=red,urlcolor=red]{hyperref} 
\makeatletter\let\@bibitem\saved@bibitem\makeatother 

\usepackage{algorithm}
\usepackage{algorithmic}


%% file: preamble.macros.tex
\usepackage{bm} 
\usepackage{amsfonts} 

\newcommand{\ds}[1]{\ensuremath{\displaystyle{#1}}}
\newcommand{\func}[3]{\ensuremath{#1 : #2 \rightarrow #3}}

\newcommand{\norm}[1]{\ensuremath{\left\| #1 \right\|}}

\newcommand{\optconOne}[3]{
\begin{aligned}
& \underset{#1}{\text{minimize}}
& & #2 \\
& \text{subject to} & & #3
\end{aligned}}

\newcommand{\pder}[2]{\ensuremath{\frac{\partial #1}{\partial #2}}}

\newcommand{\pderTwo}[3]{\ensuremath{\frac{\partial^2 #1}{\partial #2 \partial #3}}}

\newcommand{\Ecal}{\ensuremath{\mathcal{E}}}

\newcommand{\Gcal}{\ensuremath{\mathcal{G}}}
\newcommand{\Hcal}{\ensuremath{\mathcal{H}}}

\newcommand{\Lcal}{\ensuremath{\mathcal{L}}}

\newcommand{\Ocal}{\ensuremath{\mathcal{O}}}
\newcommand{\Pcal}{\ensuremath{\mathcal{P}}}

\newcommand{\Tcal}{\ensuremath{\mathcal{T}}}

\newcommand{\Vcal}{\ensuremath{\mathcal{V}}}
\newcommand{\Wcal}{\ensuremath{\mathcal{W}}}

\newcommand{\Nbb}{\ensuremath{\mathbb{N} }}

\newcommand{\Rbb}{\ensuremath{\mathbb{R} }}

\newcommand\Abm{{\ensuremath{\bm{A}}}}
\newcommand\Bbm{{\ensuremath{\bm{B}}}}
\newcommand\Cbm{{\ensuremath{\bm{C}}}}
\newcommand\Dbm{{\ensuremath{\bm{D}}}}

\newcommand\Fbm{{\ensuremath{\bm{F}}}}

\newcommand\Hbm{{\ensuremath{\bm{H}}}}
\newcommand\Ibm{{\ensuremath{\bm{I}}}}
\newcommand\Jbm{{\ensuremath{\bm{J}}}}

\newcommand\Rbm{{\ensuremath{\bm{R}}}}

\newcommand\Xbm{{\ensuremath{\bm{X}}}}

\newcommand\bbm{{\ensuremath{\bm{b}}}}
\newcommand\cbm{{\ensuremath{\bm{c}}}}

\newcommand\gbm{{\ensuremath{\bm{g}}}}

\newcommand\rbm{{\ensuremath{\bm{r}}}}

\newcommand\ubm{{\ensuremath{\bm{u}}}}

\newcommand\xbm{{\ensuremath{\bm{x}}}}

\newcommand\zbm{{\ensuremath{\bm{z}}}}

\newcommand\lambdabold{{\ensuremath{\boldsymbol{\lambda}}}}

\newcommand\etabold{{\ensuremath{\boldsymbol{\eta}}}}

\newcommand\phibold{{\ensuremath{\boldsymbol{\phi}}}}

\newcommand\zerobold{\ensuremath{\mathbf{0}}}

%% file: preamble.tikz.tex
\usepackage{tikz}
\usepackage{pgfplots}
\usepackage{pgfplotstable, filecontents, booktabs}
\pgfplotsset{compat=1.9}

\usetikzlibrary{pgfplots.groupplots}
\usepgfplotslibrary{fillbetween}
\usetikzlibrary{calc,fit,matrix,arrows,automata,positioning,shapes}
\usetikzlibrary{arrows.meta}

\pgfplotsset{select coords between index/.style 2 args={
    x filter/.code={
        \ifnum\coordindex<#1\fi
        \ifnum\coordindex>#2\fi
    }
}}

\tikzset{
 invisible/.style={opacity=0},
 visible on/.style={alt={#1{}{invisible}}},
 alt/.code args={<#1>#2#3}{%
   \alt<#1>{\pgfkeysalso{#2}}{\pgfkeysalso{#3}}
 },
}


%% file: py/smooth_heaviside.tikz
\begin{tikzpicture}
\begin{axis}[
width=0.6\textwidth,
xtick={-1, -0.75, -0.5, -0.25, 0, 0.25, 0.5, 0.75, 1},
xlabel={$x$},
ymax=1.1,
xmax=1,
ylabel={$H_a(x)$},
xmin=-1,
ymin=-0.1,
height=0.4\textwidth]
\addplot [solid, thick, color=black]
coordinates {
(-1.00000000e+00,  0.00000000e+00)
( 0.00000000e+00,  0.00000000e+00)};\label{line:heavi}

\addplot [solid, thick, color=black, forget plot]
coordinates {
( 0.00000000e+00,  1.00000000e+00)
( 1.00000000e+00,  1.00000000e+00)};

\addplot [solid, color=blue]
coordinates {
(-1.00000000e+00,  4.53978687e-05)
(-9.93333333e-01,  4.85274053e-05)
(-9.86666667e-01,  5.18726677e-05)
(-9.80000000e-01,  5.54485247e-05)
(-9.73333333e-01,  5.92708699e-05)
(-9.66666667e-01,  6.33566917e-05)
(-9.60000000e-01,  6.77241496e-05)
(-9.53333333e-01,  7.23926539e-05)
(-9.46666667e-01,  7.73829526e-05)
(-9.40000000e-01,  8.27172229e-05)
(-9.33333333e-01,  8.84191700e-05)
(-9.26666667e-01,  9.45141324e-05)
(-9.20000000e-01,  1.01029194e-04)
(-9.13333333e-01,  1.07993304e-04)
(-9.06666667e-01,  1.15437406e-04)
(-9.00000000e-01,  1.23394576e-04)
(-8.93333333e-01,  1.31900166e-04)
(-8.86666667e-01,  1.40991963e-04)
(-8.80000000e-01,  1.50710358e-04)
(-8.73333333e-01,  1.61098522e-04)
(-8.66666667e-01,  1.72202597e-04)
(-8.60000000e-01,  1.84071905e-04)
(-8.53333333e-01,  1.96759161e-04)
(-8.46666667e-01,  2.10320708e-04)
(-8.40000000e-01,  2.24816770e-04)
(-8.33333333e-01,  2.40311713e-04)
(-8.26666667e-01,  2.56874332e-04)
(-8.20000000e-01,  2.74578156e-04)
(-8.13333333e-01,  2.93501773e-04)
(-8.06666667e-01,  3.13729175e-04)
(-8.00000000e-01,  3.35350130e-04)
(-7.93333333e-01,  3.58460581e-04)
(-7.86666667e-01,  3.83163065e-04)
(-7.80000000e-01,  4.09567165e-04)
(-7.73333333e-01,  4.37789998e-04)
(-7.66666667e-01,  4.67956726e-04)
(-7.60000000e-01,  5.00201107e-04)
(-7.53333333e-01,  5.34666087e-04)
(-7.46666667e-01,  5.71504424e-04)
(-7.40000000e-01,  6.10879359e-04)
(-7.33333333e-01,  6.52965338e-04)
(-7.26666667e-01,  6.97948766e-04)
(-7.20000000e-01,  7.46028834e-04)
(-7.13333333e-01,  7.97418383e-04)
(-7.06666667e-01,  8.52344837e-04)
(-7.00000000e-01,  9.11051194e-04)
(-6.93333333e-01,  9.73797088e-04)
(-6.86666667e-01,  1.04085991e-03)
(-6.80000000e-01,  1.11253603e-03)
(-6.73333333e-01,  1.18914207e-03)
(-6.66666667e-01,  1.27101626e-03)
(-6.60000000e-01,  1.35851995e-03)
(-6.53333333e-01,  1.45203911e-03)
(-6.46666667e-01,  1.55198603e-03)
(-6.40000000e-01,  1.65880108e-03)
(-6.33333333e-01,  1.77295459e-03)
(-6.26666667e-01,  1.89494889e-03)
(-6.20000000e-01,  2.02532039e-03)
(-6.13333333e-01,  2.16464193e-03)
(-6.06666667e-01,  2.31352517e-03)
(-6.00000000e-01,  2.47262316e-03)
(-5.93333333e-01,  2.64263312e-03)
(-5.86666667e-01,  2.82429935e-03)
(-5.80000000e-01,  3.01841632e-03)
(-5.73333333e-01,  3.22583199e-03)
(-5.66666667e-01,  3.44745130e-03)
(-5.60000000e-01,  3.68423990e-03)
(-5.53333333e-01,  3.93722809e-03)
(-5.46666667e-01,  4.20751504e-03)
(-5.40000000e-01,  4.49627316e-03)
(-5.33333333e-01,  4.80475289e-03)
(-5.26666667e-01,  5.13428760e-03)
(-5.20000000e-01,  5.48629890e-03)
(-5.13333333e-01,  5.86230220e-03)
(-5.06666667e-01,  6.26391256e-03)
(-5.00000000e-01,  6.69285092e-03)
(-4.93333333e-01,  7.15095064e-03)
(-4.86666667e-01,  7.64016433e-03)
(-4.80000000e-01,  8.16257115e-03)
(-4.73333333e-01,  8.72038437e-03)
(-4.66666667e-01,  9.31595935e-03)
(-4.60000000e-01,  9.95180187e-03)
(-4.53333333e-01,  1.06305769e-02)
(-4.46666667e-01,  1.13551176e-02)
(-4.40000000e-01,  1.21284350e-02)
(-4.33333333e-01,  1.29537275e-02)
(-4.26666667e-01,  1.38343916e-02)
(-4.20000000e-01,  1.47740317e-02)
(-4.13333333e-01,  1.57764718e-02)
(-4.06666667e-01,  1.68457659e-02)
(-4.00000000e-01,  1.79862100e-02)
(-3.93333333e-01,  1.92023532e-02)
(-3.86666667e-01,  2.04990098e-02)
(-3.80000000e-01,  2.18812709e-02)
(-3.73333333e-01,  2.33545165e-02)
(-3.66666667e-01,  2.49244266e-02)
(-3.60000000e-01,  2.65969936e-02)
(-3.53333333e-01,  2.83785324e-02)
(-3.46666667e-01,  3.02756920e-02)
(-3.40000000e-01,  3.22954647e-02)
(-3.33333333e-01,  3.44451957e-02)
(-3.26666667e-01,  3.67325907e-02)
(-3.20000000e-01,  3.91657228e-02)
(-3.13333333e-01,  4.17530374e-02)
(-3.06666667e-01,  4.45033553e-02)
(-3.00000000e-01,  4.74258732e-02)
(-2.93333333e-01,  5.05301622e-02)
(-2.86666667e-01,  5.38261628e-02)
(-2.80000000e-01,  5.73241759e-02)
(-2.73333333e-01,  6.10348511e-02)
(-2.66666667e-01,  6.49691691e-02)
(-2.60000000e-01,  6.91384203e-02)
(-2.53333333e-01,  7.35541769e-02)
(-2.46666667e-01,  7.82282595e-02)
(-2.40000000e-01,  8.31726965e-02)
(-2.33333333e-01,  8.83996772e-02)
(-2.26666667e-01,  9.39214960e-02)
(-2.20000000e-01,  9.97504891e-02)
(-2.13333333e-01,  1.05898962e-01)
(-2.06666667e-01,  1.12379109e-01)
(-2.00000000e-01,  1.19202922e-01)
(-1.93333333e-01,  1.26382089e-01)
(-1.86666667e-01,  1.33927888e-01)
(-1.80000000e-01,  1.41851065e-01)
(-1.73333333e-01,  1.50161706e-01)
(-1.66666667e-01,  1.58869105e-01)
(-1.60000000e-01,  1.67981615e-01)
(-1.53333333e-01,  1.77506503e-01)
(-1.46666667e-01,  1.87449793e-01)
(-1.40000000e-01,  1.97816111e-01)
(-1.33333333e-01,  2.08608527e-01)
(-1.26666667e-01,  2.19828398e-01)
(-1.20000000e-01,  2.31475217e-01)
(-1.13333333e-01,  2.43546471e-01)
(-1.06666667e-01,  2.56037509e-01)
(-1.00000000e-01,  2.68941421e-01)
(-9.33333333e-02,  2.82248943e-01)
(-8.66666667e-02,  2.95948372e-01)
(-8.00000000e-02,  3.10025519e-01)
(-7.33333333e-02,  3.24463677e-01)
(-6.66666667e-02,  3.39243631e-01)
(-6.00000000e-02,  3.54343694e-01)
(-5.33333333e-02,  3.69739777e-01)
(-4.66666667e-02,  3.85405502e-01)
(-4.00000000e-02,  4.01312340e-01)
(-3.33333333e-02,  4.17429794e-01)
(-2.66666667e-02,  4.33725606e-01)
(-2.00000000e-02,  4.50166003e-01)
(-1.33333333e-02,  4.66715962e-01)
(-6.66666667e-03,  4.83339503e-01)
( 0.00000000e+00,  5.00000000e-01)
( 6.66666667e-03,  5.16660497e-01)
( 1.33333333e-02,  5.33284038e-01)
( 2.00000000e-02,  5.49833997e-01)
( 2.66666667e-02,  5.66274394e-01)
( 3.33333333e-02,  5.82570206e-01)
( 4.00000000e-02,  5.98687660e-01)
( 4.66666667e-02,  6.14594498e-01)
( 5.33333333e-02,  6.30260223e-01)
( 6.00000000e-02,  6.45656306e-01)
( 6.66666667e-02,  6.60756369e-01)
( 7.33333333e-02,  6.75536323e-01)
( 8.00000000e-02,  6.89974481e-01)
( 8.66666667e-02,  7.04051628e-01)
( 9.33333333e-02,  7.17751057e-01)
( 1.00000000e-01,  7.31058579e-01)
( 1.06666667e-01,  7.43962491e-01)
( 1.13333333e-01,  7.56453529e-01)
( 1.20000000e-01,  7.68524783e-01)
( 1.26666667e-01,  7.80171602e-01)
( 1.33333333e-01,  7.91391473e-01)
( 1.40000000e-01,  8.02183889e-01)
( 1.46666667e-01,  8.12550207e-01)
( 1.53333333e-01,  8.22493497e-01)
( 1.60000000e-01,  8.32018385e-01)
( 1.66666667e-01,  8.41130895e-01)
( 1.73333333e-01,  8.49838294e-01)
( 1.80000000e-01,  8.58148935e-01)
( 1.86666667e-01,  8.66072112e-01)
( 1.93333333e-01,  8.73617911e-01)
( 2.00000000e-01,  8.80797078e-01)
( 2.06666667e-01,  8.87620891e-01)
( 2.13333333e-01,  8.94101038e-01)
( 2.20000000e-01,  9.00249511e-01)
( 2.26666667e-01,  9.06078504e-01)
( 2.33333333e-01,  9.11600323e-01)
( 2.40000000e-01,  9.16827304e-01)
( 2.46666667e-01,  9.21771741e-01)
( 2.53333333e-01,  9.26445823e-01)
( 2.60000000e-01,  9.30861580e-01)
( 2.66666667e-01,  9.35030831e-01)
( 2.73333333e-01,  9.38965149e-01)
( 2.80000000e-01,  9.42675824e-01)
( 2.86666667e-01,  9.46173837e-01)
( 2.93333333e-01,  9.49469838e-01)
( 3.00000000e-01,  9.52574127e-01)
( 3.06666667e-01,  9.55496645e-01)
( 3.13333333e-01,  9.58246963e-01)
( 3.20000000e-01,  9.60834277e-01)
( 3.26666667e-01,  9.63267409e-01)
( 3.33333333e-01,  9.65554804e-01)
( 3.40000000e-01,  9.67704535e-01)
( 3.46666667e-01,  9.69724308e-01)
( 3.53333333e-01,  9.71621468e-01)
( 3.60000000e-01,  9.73403006e-01)
( 3.66666667e-01,  9.75075573e-01)
( 3.73333333e-01,  9.76645484e-01)
( 3.80000000e-01,  9.78118729e-01)
( 3.86666667e-01,  9.79500990e-01)
( 3.93333333e-01,  9.80797647e-01)
( 4.00000000e-01,  9.82013790e-01)
( 4.06666667e-01,  9.83154234e-01)
( 4.13333333e-01,  9.84223528e-01)
( 4.20000000e-01,  9.85225968e-01)
( 4.26666667e-01,  9.86165608e-01)
( 4.33333333e-01,  9.87046272e-01)
( 4.40000000e-01,  9.87871565e-01)
( 4.46666667e-01,  9.88644882e-01)
( 4.53333333e-01,  9.89369423e-01)
( 4.60000000e-01,  9.90048198e-01)
( 4.66666667e-01,  9.90684041e-01)
( 4.73333333e-01,  9.91279616e-01)
( 4.80000000e-01,  9.91837429e-01)
( 4.86666667e-01,  9.92359836e-01)
( 4.93333333e-01,  9.92849049e-01)
( 5.00000000e-01,  9.93307149e-01)
( 5.06666667e-01,  9.93736087e-01)
( 5.13333333e-01,  9.94137698e-01)
( 5.20000000e-01,  9.94513701e-01)
( 5.26666667e-01,  9.94865712e-01)
( 5.33333333e-01,  9.95195247e-01)
( 5.40000000e-01,  9.95503727e-01)
( 5.46666667e-01,  9.95792485e-01)
( 5.53333333e-01,  9.96062772e-01)
( 5.60000000e-01,  9.96315760e-01)
( 5.66666667e-01,  9.96552549e-01)
( 5.73333333e-01,  9.96774168e-01)
( 5.80000000e-01,  9.96981584e-01)
( 5.86666667e-01,  9.97175701e-01)
( 5.93333333e-01,  9.97357367e-01)
( 6.00000000e-01,  9.97527377e-01)
( 6.06666667e-01,  9.97686475e-01)
( 6.13333333e-01,  9.97835358e-01)
( 6.20000000e-01,  9.97974680e-01)
( 6.26666667e-01,  9.98105051e-01)
( 6.33333333e-01,  9.98227045e-01)
( 6.40000000e-01,  9.98341199e-01)
( 6.46666667e-01,  9.98448014e-01)
( 6.53333333e-01,  9.98547961e-01)
( 6.60000000e-01,  9.98641480e-01)
( 6.66666667e-01,  9.98728984e-01)
( 6.73333333e-01,  9.98810858e-01)
( 6.80000000e-01,  9.98887464e-01)
( 6.86666667e-01,  9.98959140e-01)
( 6.93333333e-01,  9.99026203e-01)
( 7.00000000e-01,  9.99088949e-01)
( 7.06666667e-01,  9.99147655e-01)
( 7.13333333e-01,  9.99202582e-01)
( 7.20000000e-01,  9.99253971e-01)
( 7.26666667e-01,  9.99302051e-01)
( 7.33333333e-01,  9.99347035e-01)
( 7.40000000e-01,  9.99389121e-01)
( 7.46666667e-01,  9.99428496e-01)
( 7.53333333e-01,  9.99465334e-01)
( 7.60000000e-01,  9.99499799e-01)
( 7.66666667e-01,  9.99532043e-01)
( 7.73333333e-01,  9.99562210e-01)
( 7.80000000e-01,  9.99590433e-01)
( 7.86666667e-01,  9.99616837e-01)
( 7.93333333e-01,  9.99641539e-01)
( 8.00000000e-01,  9.99664650e-01)
( 8.06666667e-01,  9.99686271e-01)
( 8.13333333e-01,  9.99706498e-01)
( 8.20000000e-01,  9.99725422e-01)
( 8.26666667e-01,  9.99743126e-01)
( 8.33333333e-01,  9.99759688e-01)
( 8.40000000e-01,  9.99775183e-01)
( 8.46666667e-01,  9.99789679e-01)
( 8.53333333e-01,  9.99803241e-01)
( 8.60000000e-01,  9.99815928e-01)
( 8.66666667e-01,  9.99827797e-01)
( 8.73333333e-01,  9.99838901e-01)
( 8.80000000e-01,  9.99849290e-01)
( 8.86666667e-01,  9.99859008e-01)
( 8.93333333e-01,  9.99868100e-01)
( 9.00000000e-01,  9.99876605e-01)
( 9.06666667e-01,  9.99884563e-01)
( 9.13333333e-01,  9.99892007e-01)
( 9.20000000e-01,  9.99898971e-01)
( 9.26666667e-01,  9.99905486e-01)
( 9.33333333e-01,  9.99911581e-01)
( 9.40000000e-01,  9.99917283e-01)
( 9.46666667e-01,  9.99922617e-01)
( 9.53333333e-01,  9.99927607e-01)
( 9.60000000e-01,  9.99932276e-01)
( 9.66666667e-01,  9.99936643e-01)
( 9.73333333e-01,  9.99940729e-01)
( 9.80000000e-01,  9.99944551e-01)
( 9.86666667e-01,  9.99948127e-01)
( 9.93333333e-01,  9.99951473e-01)
( 1.00000000e+00,  9.99954602e-01)};\label{line:smheavi5}

\addplot [solid, color=red]
coordinates {
(-1.00000000e+00,  2.06115362e-09)
(-9.93333333e-01,  2.35513763e-09)
(-9.86666667e-01,  2.69105282e-09)
(-9.80000000e-01,  3.07487987e-09)
(-9.73333333e-01,  3.51345248e-09)
(-9.66666667e-01,  4.01457906e-09)
(-9.60000000e-01,  4.58718173e-09)
(-9.53333333e-01,  5.24145518e-09)
(-9.46666667e-01,  5.98904818e-09)
(-9.40000000e-01,  6.84327098e-09)
(-9.33333333e-01,  7.81933226e-09)
(-9.26666667e-01,  8.93460996e-09)
(-9.20000000e-01,  1.02089606e-08)
(-9.13333333e-01,  1.16650729e-08)
(-9.06666667e-01,  1.33288717e-08)
(-9.00000000e-01,  1.52299795e-08)
(-8.93333333e-01,  1.74022438e-08)
(-8.86666667e-01,  1.98843399e-08)
(-8.80000000e-01,  2.27204594e-08)
(-8.73333333e-01,  2.59610969e-08)
(-8.66666667e-01,  2.96639491e-08)
(-8.60000000e-01,  3.38949421e-08)
(-8.53333333e-01,  3.87294050e-08)
(-8.46666667e-01,  4.42534113e-08)
(-8.40000000e-01,  5.05653109e-08)
(-8.33333333e-01,  5.77774819e-08)
(-8.26666667e-01,  6.60183305e-08)
(-8.20000000e-01,  7.54345778e-08)
(-8.13333333e-01,  8.61938720e-08)
(-8.06666667e-01,  9.84877727e-08)
(-8.00000000e-01,  1.12535162e-07)
(-7.93333333e-01,  1.28586142e-07)
(-7.86666667e-01,  1.46926485e-07)
(-7.80000000e-01,  1.67882725e-07)
(-7.73333333e-01,  1.91827970e-07)
(-7.66666667e-01,  2.19188543e-07)
(-7.60000000e-01,  2.50451575e-07)
(-7.53333333e-01,  2.86173676e-07)
(-7.46666667e-01,  3.26990846e-07)
(-7.40000000e-01,  3.73629798e-07)
(-7.33333333e-01,  4.26920897e-07)
(-7.26666667e-01,  4.87812942e-07)
(-7.20000000e-01,  5.57390059e-07)
(-7.13333333e-01,  6.36891004e-07)
(-7.06666667e-01,  7.27731219e-07)
(-7.00000000e-01,  8.31528028e-07)
(-6.93333333e-01,  9.50129433e-07)
(-6.86666667e-01,  1.08564702e-06)
(-6.80000000e-01,  1.24049354e-06)
(-6.73333333e-01,  1.41742589e-06)
(-6.66666667e-01,  1.61959417e-06)
(-6.60000000e-01,  1.85059777e-06)
(-6.53333333e-01,  2.11454948e-06)
(-6.46666667e-01,  2.41614866e-06)
(-6.40000000e-01,  2.76076495e-06)
(-6.33333333e-01,  3.15453385e-06)
(-6.26666667e-01,  3.60446596e-06)
(-6.20000000e-01,  4.11857174e-06)
(-6.13333333e-01,  4.70600421e-06)
(-6.06666667e-01,  5.37722180e-06)
(-6.00000000e-01,  6.14417460e-06)
(-5.93333333e-01,  7.02051706e-06)
(-5.86666667e-01,  8.02185108e-06)
(-5.80000000e-01,  9.16600372e-06)
(-5.73333333e-01,  1.04733446e-05)
(-5.66666667e-01,  1.19671483e-05)
(-5.60000000e-01,  1.36740091e-05)
(-5.53333333e-01,  1.56243136e-05)
(-5.46666667e-01,  1.78527824e-05)
(-5.40000000e-01,  2.03990873e-05)
(-5.33333333e-01,  2.33085578e-05)
(-5.26666667e-01,  2.66329878e-05)
(-5.20000000e-01,  3.04315569e-05)
(-5.13333333e-01,  3.47718836e-05)
(-5.06666667e-01,  3.97312286e-05)
(-5.00000000e-01,  4.53978687e-05)
(-4.93333333e-01,  5.18726677e-05)
(-4.86666667e-01,  5.92708699e-05)
(-4.80000000e-01,  6.77241496e-05)
(-4.73333333e-01,  7.73829526e-05)
(-4.66666667e-01,  8.84191700e-05)
(-4.60000000e-01,  1.01029194e-04)
(-4.53333333e-01,  1.15437406e-04)
(-4.46666667e-01,  1.31900166e-04)
(-4.40000000e-01,  1.50710358e-04)
(-4.33333333e-01,  1.72202597e-04)
(-4.26666667e-01,  1.96759161e-04)
(-4.20000000e-01,  2.24816770e-04)
(-4.13333333e-01,  2.56874332e-04)
(-4.06666667e-01,  2.93501773e-04)
(-4.00000000e-01,  3.35350130e-04)
(-3.93333333e-01,  3.83163065e-04)
(-3.86666667e-01,  4.37789998e-04)
(-3.80000000e-01,  5.00201107e-04)
(-3.73333333e-01,  5.71504424e-04)
(-3.66666667e-01,  6.52965338e-04)
(-3.60000000e-01,  7.46028834e-04)
(-3.53333333e-01,  8.52344837e-04)
(-3.46666667e-01,  9.73797088e-04)
(-3.40000000e-01,  1.11253603e-03)
(-3.33333333e-01,  1.27101626e-03)
(-3.26666667e-01,  1.45203911e-03)
(-3.20000000e-01,  1.65880108e-03)
(-3.13333333e-01,  1.89494889e-03)
(-3.06666667e-01,  2.16464193e-03)
(-3.00000000e-01,  2.47262316e-03)
(-2.93333333e-01,  2.82429935e-03)
(-2.86666667e-01,  3.22583199e-03)
(-2.80000000e-01,  3.68423990e-03)
(-2.73333333e-01,  4.20751504e-03)
(-2.66666667e-01,  4.80475289e-03)
(-2.60000000e-01,  5.48629890e-03)
(-2.53333333e-01,  6.26391256e-03)
(-2.46666667e-01,  7.15095064e-03)
(-2.40000000e-01,  8.16257115e-03)
(-2.33333333e-01,  9.31595935e-03)
(-2.26666667e-01,  1.06305769e-02)
(-2.20000000e-01,  1.21284350e-02)
(-2.13333333e-01,  1.38343916e-02)
(-2.06666667e-01,  1.57764718e-02)
(-2.00000000e-01,  1.79862100e-02)
(-1.93333333e-01,  2.04990098e-02)
(-1.86666667e-01,  2.33545165e-02)
(-1.80000000e-01,  2.65969936e-02)
(-1.73333333e-01,  3.02756920e-02)
(-1.66666667e-01,  3.44451957e-02)
(-1.60000000e-01,  3.91657228e-02)
(-1.53333333e-01,  4.45033553e-02)
(-1.46666667e-01,  5.05301622e-02)
(-1.40000000e-01,  5.73241759e-02)
(-1.33333333e-01,  6.49691691e-02)
(-1.26666667e-01,  7.35541769e-02)
(-1.20000000e-01,  8.31726965e-02)
(-1.13333333e-01,  9.39214960e-02)
(-1.06666667e-01,  1.05898962e-01)
(-1.00000000e-01,  1.19202922e-01)
(-9.33333333e-02,  1.33927888e-01)
(-8.66666667e-02,  1.50161706e-01)
(-8.00000000e-02,  1.67981615e-01)
(-7.33333333e-02,  1.87449793e-01)
(-6.66666667e-02,  2.08608527e-01)
(-6.00000000e-02,  2.31475217e-01)
(-5.33333333e-02,  2.56037509e-01)
(-4.66666667e-02,  2.82248943e-01)
(-4.00000000e-02,  3.10025519e-01)
(-3.33333333e-02,  3.39243631e-01)
(-2.66666667e-02,  3.69739777e-01)
(-2.00000000e-02,  4.01312340e-01)
(-1.33333333e-02,  4.33725606e-01)
(-6.66666667e-03,  4.66715962e-01)
( 0.00000000e+00,  5.00000000e-01)
( 6.66666667e-03,  5.33284038e-01)
( 1.33333333e-02,  5.66274394e-01)
( 2.00000000e-02,  5.98687660e-01)
( 2.66666667e-02,  6.30260223e-01)
( 3.33333333e-02,  6.60756369e-01)
( 4.00000000e-02,  6.89974481e-01)
( 4.66666667e-02,  7.17751057e-01)
( 5.33333333e-02,  7.43962491e-01)
( 6.00000000e-02,  7.68524783e-01)
( 6.66666667e-02,  7.91391473e-01)
( 7.33333333e-02,  8.12550207e-01)
( 8.00000000e-02,  8.32018385e-01)
( 8.66666667e-02,  8.49838294e-01)
( 9.33333333e-02,  8.66072112e-01)
( 1.00000000e-01,  8.80797078e-01)
( 1.06666667e-01,  8.94101038e-01)
( 1.13333333e-01,  9.06078504e-01)
( 1.20000000e-01,  9.16827304e-01)
( 1.26666667e-01,  9.26445823e-01)
( 1.33333333e-01,  9.35030831e-01)
( 1.40000000e-01,  9.42675824e-01)
( 1.46666667e-01,  9.49469838e-01)
( 1.53333333e-01,  9.55496645e-01)
( 1.60000000e-01,  9.60834277e-01)
( 1.66666667e-01,  9.65554804e-01)
( 1.73333333e-01,  9.69724308e-01)
( 1.80000000e-01,  9.73403006e-01)
( 1.86666667e-01,  9.76645484e-01)
( 1.93333333e-01,  9.79500990e-01)
( 2.00000000e-01,  9.82013790e-01)
( 2.06666667e-01,  9.84223528e-01)
( 2.13333333e-01,  9.86165608e-01)
( 2.20000000e-01,  9.87871565e-01)
( 2.26666667e-01,  9.89369423e-01)
( 2.33333333e-01,  9.90684041e-01)
( 2.40000000e-01,  9.91837429e-01)
( 2.46666667e-01,  9.92849049e-01)
( 2.53333333e-01,  9.93736087e-01)
( 2.60000000e-01,  9.94513701e-01)
( 2.66666667e-01,  9.95195247e-01)
( 2.73333333e-01,  9.95792485e-01)
( 2.80000000e-01,  9.96315760e-01)
( 2.86666667e-01,  9.96774168e-01)
( 2.93333333e-01,  9.97175701e-01)
( 3.00000000e-01,  9.97527377e-01)
( 3.06666667e-01,  9.97835358e-01)
( 3.13333333e-01,  9.98105051e-01)
( 3.20000000e-01,  9.98341199e-01)
( 3.26666667e-01,  9.98547961e-01)
( 3.33333333e-01,  9.98728984e-01)
( 3.40000000e-01,  9.98887464e-01)
( 3.46666667e-01,  9.99026203e-01)
( 3.53333333e-01,  9.99147655e-01)
( 3.60000000e-01,  9.99253971e-01)
( 3.66666667e-01,  9.99347035e-01)
( 3.73333333e-01,  9.99428496e-01)
( 3.80000000e-01,  9.99499799e-01)
( 3.86666667e-01,  9.99562210e-01)
( 3.93333333e-01,  9.99616837e-01)
( 4.00000000e-01,  9.99664650e-01)
( 4.06666667e-01,  9.99706498e-01)
( 4.13333333e-01,  9.99743126e-01)
( 4.20000000e-01,  9.99775183e-01)
( 4.26666667e-01,  9.99803241e-01)
( 4.33333333e-01,  9.99827797e-01)
( 4.40000000e-01,  9.99849290e-01)
( 4.46666667e-01,  9.99868100e-01)
( 4.53333333e-01,  9.99884563e-01)
( 4.60000000e-01,  9.99898971e-01)
( 4.66666667e-01,  9.99911581e-01)
( 4.73333333e-01,  9.99922617e-01)
( 4.80000000e-01,  9.99932276e-01)
( 4.86666667e-01,  9.99940729e-01)
( 4.93333333e-01,  9.99948127e-01)
( 5.00000000e-01,  9.99954602e-01)
( 5.06666667e-01,  9.99960269e-01)
( 5.13333333e-01,  9.99965228e-01)
( 5.20000000e-01,  9.99969568e-01)
( 5.26666667e-01,  9.99973367e-01)
( 5.33333333e-01,  9.99976691e-01)
( 5.40000000e-01,  9.99979601e-01)
( 5.46666667e-01,  9.99982147e-01)
( 5.53333333e-01,  9.99984376e-01)
( 5.60000000e-01,  9.99986326e-01)
( 5.66666667e-01,  9.99988033e-01)
( 5.73333333e-01,  9.99989527e-01)
( 5.80000000e-01,  9.99990834e-01)
( 5.86666667e-01,  9.99991978e-01)
( 5.93333333e-01,  9.99992979e-01)
( 6.00000000e-01,  9.99993856e-01)
( 6.06666667e-01,  9.99994623e-01)
( 6.13333333e-01,  9.99995294e-01)
( 6.20000000e-01,  9.99995881e-01)
( 6.26666667e-01,  9.99996396e-01)
( 6.33333333e-01,  9.99996845e-01)
( 6.40000000e-01,  9.99997239e-01)
( 6.46666667e-01,  9.99997584e-01)
( 6.53333333e-01,  9.99997885e-01)
( 6.60000000e-01,  9.99998149e-01)
( 6.66666667e-01,  9.99998380e-01)
( 6.73333333e-01,  9.99998583e-01)
( 6.80000000e-01,  9.99998760e-01)
( 6.86666667e-01,  9.99998914e-01)
( 6.93333333e-01,  9.99999050e-01)
( 7.00000000e-01,  9.99999168e-01)
( 7.06666667e-01,  9.99999272e-01)
( 7.13333333e-01,  9.99999363e-01)
( 7.20000000e-01,  9.99999443e-01)
( 7.26666667e-01,  9.99999512e-01)
( 7.33333333e-01,  9.99999573e-01)
( 7.40000000e-01,  9.99999626e-01)
( 7.46666667e-01,  9.99999673e-01)
( 7.53333333e-01,  9.99999714e-01)
( 7.60000000e-01,  9.99999750e-01)
( 7.66666667e-01,  9.99999781e-01)
( 7.73333333e-01,  9.99999808e-01)
( 7.80000000e-01,  9.99999832e-01)
( 7.86666667e-01,  9.99999853e-01)
( 7.93333333e-01,  9.99999871e-01)
( 8.00000000e-01,  9.99999887e-01)
( 8.06666667e-01,  9.99999902e-01)
( 8.13333333e-01,  9.99999914e-01)
( 8.20000000e-01,  9.99999925e-01)
( 8.26666667e-01,  9.99999934e-01)
( 8.33333333e-01,  9.99999942e-01)
( 8.40000000e-01,  9.99999949e-01)
( 8.46666667e-01,  9.99999956e-01)
( 8.53333333e-01,  9.99999961e-01)
( 8.60000000e-01,  9.99999966e-01)
( 8.66666667e-01,  9.99999970e-01)
( 8.73333333e-01,  9.99999974e-01)
( 8.80000000e-01,  9.99999977e-01)
( 8.86666667e-01,  9.99999980e-01)
( 8.93333333e-01,  9.99999983e-01)
( 9.00000000e-01,  9.99999985e-01)
( 9.06666667e-01,  9.99999987e-01)
( 9.13333333e-01,  9.99999988e-01)
( 9.20000000e-01,  9.99999990e-01)
( 9.26666667e-01,  9.99999991e-01)
( 9.33333333e-01,  9.99999992e-01)
( 9.40000000e-01,  9.99999993e-01)
( 9.46666667e-01,  9.99999994e-01)
( 9.53333333e-01,  9.99999995e-01)
( 9.60000000e-01,  9.99999995e-01)
( 9.66666667e-01,  9.99999996e-01)
( 9.73333333e-01,  9.99999996e-01)
( 9.80000000e-01,  9.99999997e-01)
( 9.86666667e-01,  9.99999997e-01)
( 9.93333333e-01,  9.99999998e-01)
( 1.00000000e+00,  9.99999998e-01)};\label{line:smheavi10}

\addplot [solid, color=green]
coordinates {
(-1.00000000e+00,  8.75651076e-27)
(-9.93333333e-01,  1.30631790e-26)
(-9.86666667e-01,  1.94879731e-26)
(-9.80000000e-01,  2.90726396e-26)
(-9.73333333e-01,  4.33712817e-26)
(-9.66666667e-01,  6.47023493e-26)
(-9.60000000e-01,  9.65245626e-26)
(-9.53333333e-01,  1.43997726e-25)
(-9.46666667e-01,  2.14819365e-25)
(-9.40000000e-01,  3.20472834e-25)
(-9.33333333e-01,  4.78089288e-25)
(-9.26666667e-01,  7.13225408e-25)
(-9.20000000e-01,  1.06400728e-24)
(-9.13333333e-01,  1.58731234e-24)
(-9.06666667e-01,  2.36799175e-24)
(-9.00000000e-01,  3.53262857e-24)
(-8.93333333e-01,  5.27006255e-24)
(-8.86666667e-01,  7.86200947e-24)
(-8.80000000e-01,  1.17287399e-23)
(-8.73333333e-01,  1.74972239e-23)
(-8.66666667e-01,  2.61027907e-23)
(-8.60000000e-01,  3.89407878e-23)
(-8.53333333e-01,  5.80928290e-23)
(-8.46666667e-01,  8.66643171e-23)
(-8.40000000e-01,  1.29287969e-22)
(-8.33333333e-01,  1.92874985e-22)
(-8.26666667e-01,  2.87735666e-22)
(-8.20000000e-01,  4.29251173e-22)
(-8.13333333e-01,  6.40367501e-22)
(-8.06666667e-01,  9.55316054e-22)
(-8.00000000e-01,  1.42516408e-21)
(-7.93333333e-01,  2.12609498e-21)
(-7.86666667e-01,  3.17176100e-21)
(-7.80000000e-01,  4.73171139e-21)
(-7.73333333e-01,  7.05888391e-21)
(-7.66666667e-01,  1.05306174e-20)
(-7.60000000e-01,  1.57098351e-20)
(-7.53333333e-01,  2.34363199e-20)
(-7.46666667e-01,  3.49628809e-20)
(-7.40000000e-01,  5.21584892e-20)
(-7.33333333e-01,  7.78113224e-20)
(-7.26666667e-01,  1.16080853e-19)
(-7.20000000e-01,  1.73172283e-19)
(-7.13333333e-01,  2.58342688e-19)
(-7.06666667e-01,  3.85402003e-19)
(-7.00000000e-01,  5.74952226e-19)
(-6.93333333e-01,  8.57727931e-19)
(-6.86666667e-01,  1.27957971e-18)
(-6.80000000e-01,  1.90890862e-18)
(-6.73333333e-01,  2.84775702e-18)
(-6.66666667e-01,  4.24835426e-18)
(-6.60000000e-01,  6.33779980e-18)
(-6.53333333e-01,  9.45488627e-18)
(-6.46666667e-01,  1.41050329e-17)
(-6.40000000e-01,  2.10422364e-17)
(-6.33333333e-01,  3.13913279e-17)
(-6.26666667e-01,  4.68303583e-17)
(-6.20000000e-01,  6.98626851e-17)
(-6.13333333e-01,  1.04222879e-16)
(-6.06666667e-01,  1.55482265e-16)
(-6.00000000e-01,  2.31952283e-16)
(-5.93333333e-01,  3.46032144e-16)
(-5.86666667e-01,  5.16219299e-16)
(-5.80000000e-01,  7.70108700e-16)
(-5.73333333e-01,  1.14886718e-15)
(-5.66666667e-01,  1.71390843e-15)
(-5.60000000e-01,  2.55685093e-15)
(-5.53333333e-01,  3.81437336e-15)
(-5.46666667e-01,  5.69037639e-15)
(-5.40000000e-01,  8.48904403e-15)
(-5.33333333e-01,  1.26641655e-14)
(-5.26666667e-01,  1.88927149e-14)
(-5.20000000e-01,  2.81846188e-14)
(-5.13333333e-01,  4.20465104e-14)
(-5.06666667e-01,  6.27260226e-14)
(-5.00000000e-01,  9.35762297e-14)
(-4.93333333e-01,  1.39599331e-13)
(-4.86666667e-01,  2.08257729e-13)
(-4.80000000e-01,  3.10684024e-13)
(-4.73333333e-01,  4.63486100e-13)
(-4.66666667e-01,  6.91440011e-13)
(-4.60000000e-01,  1.03150728e-12)
(-4.53333333e-01,  1.53882804e-12)
(-4.46666667e-01,  2.29566168e-12)
(-4.40000000e-01,  3.42472479e-12)
(-4.33333333e-01,  5.10908903e-12)
(-4.26666667e-01,  7.62186519e-12)
(-4.20000000e-01,  1.13704867e-11)
(-4.13333333e-01,  1.69627729e-11)
(-4.06666667e-01,  2.53054836e-11)
(-4.00000000e-01,  3.77513454e-11)
(-3.93333333e-01,  5.63183895e-11)
(-3.86666667e-01,  8.40171644e-11)
(-3.80000000e-01,  1.25338881e-10)
(-3.73333333e-01,  1.86983638e-10)
(-3.66666667e-01,  2.78946809e-10)
(-3.60000000e-01,  4.16139739e-10)
(-3.53333333e-01,  6.20807541e-10)
(-3.46666667e-01,  9.26136021e-10)
(-3.40000000e-01,  1.38163259e-09)
(-3.33333333e-01,  2.06115362e-09)
(-3.26666667e-01,  3.07487987e-09)
(-3.20000000e-01,  4.58718173e-09)
(-3.13333333e-01,  6.84327098e-09)
(-3.06666667e-01,  1.02089606e-08)
(-3.00000000e-01,  1.52299795e-08)
(-2.93333333e-01,  2.27204594e-08)
(-2.86666667e-01,  3.38949421e-08)
(-2.80000000e-01,  5.05653109e-08)
(-2.73333333e-01,  7.54345778e-08)
(-2.66666667e-01,  1.12535162e-07)
(-2.60000000e-01,  1.67882725e-07)
(-2.53333333e-01,  2.50451575e-07)
(-2.46666667e-01,  3.73629798e-07)
(-2.40000000e-01,  5.57390059e-07)
(-2.33333333e-01,  8.31528028e-07)
(-2.26666667e-01,  1.24049354e-06)
(-2.20000000e-01,  1.85059777e-06)
(-2.13333333e-01,  2.76076495e-06)
(-2.06666667e-01,  4.11857174e-06)
(-2.00000000e-01,  6.14417460e-06)
(-1.93333333e-01,  9.16600372e-06)
(-1.86666667e-01,  1.36740091e-05)
(-1.80000000e-01,  2.03990873e-05)
(-1.73333333e-01,  3.04315569e-05)
(-1.66666667e-01,  4.53978687e-05)
(-1.60000000e-01,  6.77241496e-05)
(-1.53333333e-01,  1.01029194e-04)
(-1.46666667e-01,  1.50710358e-04)
(-1.40000000e-01,  2.24816770e-04)
(-1.33333333e-01,  3.35350130e-04)
(-1.26666667e-01,  5.00201107e-04)
(-1.20000000e-01,  7.46028834e-04)
(-1.13333333e-01,  1.11253603e-03)
(-1.06666667e-01,  1.65880108e-03)
(-1.00000000e-01,  2.47262316e-03)
(-9.33333333e-02,  3.68423990e-03)
(-8.66666667e-02,  5.48629890e-03)
(-8.00000000e-02,  8.16257115e-03)
(-7.33333333e-02,  1.21284350e-02)
(-6.66666667e-02,  1.79862100e-02)
(-6.00000000e-02,  2.65969936e-02)
(-5.33333333e-02,  3.91657228e-02)
(-4.66666667e-02,  5.73241759e-02)
(-4.00000000e-02,  8.31726965e-02)
(-3.33333333e-02,  1.19202922e-01)
(-2.66666667e-02,  1.67981615e-01)
(-2.00000000e-02,  2.31475217e-01)
(-1.33333333e-02,  3.10025519e-01)
(-6.66666667e-03,  4.01312340e-01)
( 0.00000000e+00,  5.00000000e-01)
( 6.66666667e-03,  5.98687660e-01)
( 1.33333333e-02,  6.89974481e-01)
( 2.00000000e-02,  7.68524783e-01)
( 2.66666667e-02,  8.32018385e-01)
( 3.33333333e-02,  8.80797078e-01)
( 4.00000000e-02,  9.16827304e-01)
( 4.66666667e-02,  9.42675824e-01)
( 5.33333333e-02,  9.60834277e-01)
( 6.00000000e-02,  9.73403006e-01)
( 6.66666667e-02,  9.82013790e-01)
( 7.33333333e-02,  9.87871565e-01)
( 8.00000000e-02,  9.91837429e-01)
( 8.66666667e-02,  9.94513701e-01)
( 9.33333333e-02,  9.96315760e-01)
( 1.00000000e-01,  9.97527377e-01)
( 1.06666667e-01,  9.98341199e-01)
( 1.13333333e-01,  9.98887464e-01)
( 1.20000000e-01,  9.99253971e-01)
( 1.26666667e-01,  9.99499799e-01)
( 1.33333333e-01,  9.99664650e-01)
( 1.40000000e-01,  9.99775183e-01)
( 1.46666667e-01,  9.99849290e-01)
( 1.53333333e-01,  9.99898971e-01)
( 1.60000000e-01,  9.99932276e-01)
( 1.66666667e-01,  9.99954602e-01)
( 1.73333333e-01,  9.99969568e-01)
( 1.80000000e-01,  9.99979601e-01)
( 1.86666667e-01,  9.99986326e-01)
( 1.93333333e-01,  9.99990834e-01)
( 2.00000000e-01,  9.99993856e-01)
( 2.06666667e-01,  9.99995881e-01)
( 2.13333333e-01,  9.99997239e-01)
( 2.20000000e-01,  9.99998149e-01)
( 2.26666667e-01,  9.99998760e-01)
( 2.33333333e-01,  9.99999168e-01)
( 2.40000000e-01,  9.99999443e-01)
( 2.46666667e-01,  9.99999626e-01)
( 2.53333333e-01,  9.99999750e-01)
( 2.60000000e-01,  9.99999832e-01)
( 2.66666667e-01,  9.99999887e-01)
( 2.73333333e-01,  9.99999925e-01)
( 2.80000000e-01,  9.99999949e-01)
( 2.86666667e-01,  9.99999966e-01)
( 2.93333333e-01,  9.99999977e-01)
( 3.00000000e-01,  9.99999985e-01)
( 3.06666667e-01,  9.99999990e-01)
( 3.13333333e-01,  9.99999993e-01)
( 3.20000000e-01,  9.99999995e-01)
( 3.26666667e-01,  9.99999997e-01)
( 3.33333333e-01,  9.99999998e-01)
( 3.40000000e-01,  9.99999999e-01)
( 3.46666667e-01,  9.99999999e-01)
( 3.53333333e-01,  9.99999999e-01)
( 3.60000000e-01,  1.00000000e+00)
( 3.66666667e-01,  1.00000000e+00)
( 3.73333333e-01,  1.00000000e+00)
( 3.80000000e-01,  1.00000000e+00)
( 3.86666667e-01,  1.00000000e+00)
( 3.93333333e-01,  1.00000000e+00)
( 4.00000000e-01,  1.00000000e+00)
( 4.06666667e-01,  1.00000000e+00)
( 4.13333333e-01,  1.00000000e+00)
( 4.20000000e-01,  1.00000000e+00)
( 4.26666667e-01,  1.00000000e+00)
( 4.33333333e-01,  1.00000000e+00)
( 4.40000000e-01,  1.00000000e+00)
( 4.46666667e-01,  1.00000000e+00)
( 4.53333333e-01,  1.00000000e+00)
( 4.60000000e-01,  1.00000000e+00)
( 4.66666667e-01,  1.00000000e+00)
( 4.73333333e-01,  1.00000000e+00)
( 4.80000000e-01,  1.00000000e+00)
( 4.86666667e-01,  1.00000000e+00)
( 4.93333333e-01,  1.00000000e+00)
( 5.00000000e-01,  1.00000000e+00)
( 5.06666667e-01,  1.00000000e+00)
( 5.13333333e-01,  1.00000000e+00)
( 5.20000000e-01,  1.00000000e+00)
( 5.26666667e-01,  1.00000000e+00)
( 5.33333333e-01,  1.00000000e+00)
( 5.40000000e-01,  1.00000000e+00)
( 5.46666667e-01,  1.00000000e+00)
( 5.53333333e-01,  1.00000000e+00)
( 5.60000000e-01,  1.00000000e+00)
( 5.66666667e-01,  1.00000000e+00)
( 5.73333333e-01,  1.00000000e+00)
( 5.80000000e-01,  1.00000000e+00)
( 5.86666667e-01,  1.00000000e+00)
( 5.93333333e-01,  1.00000000e+00)
( 6.00000000e-01,  1.00000000e+00)
( 6.06666667e-01,  1.00000000e+00)
( 6.13333333e-01,  1.00000000e+00)
( 6.20000000e-01,  1.00000000e+00)
( 6.26666667e-01,  1.00000000e+00)
( 6.33333333e-01,  1.00000000e+00)
( 6.40000000e-01,  1.00000000e+00)
( 6.46666667e-01,  1.00000000e+00)
( 6.53333333e-01,  1.00000000e+00)
( 6.60000000e-01,  1.00000000e+00)
( 6.66666667e-01,  1.00000000e+00)
( 6.73333333e-01,  1.00000000e+00)
( 6.80000000e-01,  1.00000000e+00)
( 6.86666667e-01,  1.00000000e+00)
( 6.93333333e-01,  1.00000000e+00)
( 7.00000000e-01,  1.00000000e+00)
( 7.06666667e-01,  1.00000000e+00)
( 7.13333333e-01,  1.00000000e+00)
( 7.20000000e-01,  1.00000000e+00)
( 7.26666667e-01,  1.00000000e+00)
( 7.33333333e-01,  1.00000000e+00)
( 7.40000000e-01,  1.00000000e+00)
( 7.46666667e-01,  1.00000000e+00)
( 7.53333333e-01,  1.00000000e+00)
( 7.60000000e-01,  1.00000000e+00)
( 7.66666667e-01,  1.00000000e+00)
( 7.73333333e-01,  1.00000000e+00)
( 7.80000000e-01,  1.00000000e+00)
( 7.86666667e-01,  1.00000000e+00)
( 7.93333333e-01,  1.00000000e+00)
( 8.00000000e-01,  1.00000000e+00)
( 8.06666667e-01,  1.00000000e+00)
( 8.13333333e-01,  1.00000000e+00)
( 8.20000000e-01,  1.00000000e+00)
( 8.26666667e-01,  1.00000000e+00)
( 8.33333333e-01,  1.00000000e+00)
( 8.40000000e-01,  1.00000000e+00)
( 8.46666667e-01,  1.00000000e+00)
( 8.53333333e-01,  1.00000000e+00)
( 8.60000000e-01,  1.00000000e+00)
( 8.66666667e-01,  1.00000000e+00)
( 8.73333333e-01,  1.00000000e+00)
( 8.80000000e-01,  1.00000000e+00)
( 8.86666667e-01,  1.00000000e+00)
( 8.93333333e-01,  1.00000000e+00)
( 9.00000000e-01,  1.00000000e+00)
( 9.06666667e-01,  1.00000000e+00)
( 9.13333333e-01,  1.00000000e+00)
( 9.20000000e-01,  1.00000000e+00)
( 9.26666667e-01,  1.00000000e+00)
( 9.33333333e-01,  1.00000000e+00)
( 9.40000000e-01,  1.00000000e+00)
( 9.46666667e-01,  1.00000000e+00)
( 9.53333333e-01,  1.00000000e+00)
( 9.60000000e-01,  1.00000000e+00)
( 9.66666667e-01,  1.00000000e+00)
( 9.73333333e-01,  1.00000000e+00)
( 9.80000000e-01,  1.00000000e+00)
( 9.86666667e-01,  1.00000000e+00)
( 9.93333333e-01,  1.00000000e+00)
( 1.00000000e+00,  1.00000000e+00)};\label{line:smheavi30}

\end{axis}
\end{tikzpicture}

%% file: tikz/advec2d_beta0_c1s1_nel6x3_porder0x1.tikz
\begin{tikzpicture}
\begin{groupplot}[
  group style={
      group size=2 by 3,
      horizontal sep=1cm,
      vertical sep=1.2cm
  },
  width=0.5\textwidth,
  axis equal image,
  xlabel={$x_1$},
  ylabel={$x_2$},
  xtick = {-1, -0.5, 0.0, 0.5, 1.0},
  ytick = {0.0, 0.5, 0.8, 1.0},
  xmin=-1, xmax=1,
  ymin=0, ymax=1
]

\nextgroupplot[title={Iteration 0}, xlabel={}, xtick=\empty]
\addplot graphics [xmin=-1, xmax=1, ymin=0, ymax=1] {{img/advec2d_beta0_c1s1_nel6x3_porder0x1_iter000}.png};

\nextgroupplot[title={Iteration 1}, xlabel={}, ylabel={}, xtick=\empty, ytick=\empty]
\addplot graphics [xmin=-1, xmax=1, ymin=0, ymax=1] {{img/advec2d_beta0_c1s1_nel6x3_porder0x1_iter001}.png};

\nextgroupplot[title={Iteration 2}, xlabel={}, xtick=\empty]
\addplot graphics [xmin=-1, xmax=1, ymin=0, ymax=1] {{img/advec2d_beta0_c1s1_nel6x3_porder0x1_iter002}.png};

\nextgroupplot[title={Iteration 3}, xlabel={}, ylabel={}, xtick=\empty, ytick=\empty]
\addplot graphics [xmin=-1, xmax=1, ymin=0, ymax=1] {{img/advec2d_beta0_c1s1_nel6x3_porder0x1_iter003}.png};

\nextgroupplot[title={Iteration 4}]
\addplot graphics [xmin=-1, xmax=1, ymin=0, ymax=1] {{img/advec2d_beta0_c1s1_nel6x3_porder0x1_iter004}.png};

\nextgroupplot[title={Iteration 10}, ylabel={}, ytick=\empty]
\addplot graphics [xmin=-1, xmax=1, ymin=0, ymax=1] {{img/advec2d_beta0_c1s1_nel6x3_porder0x1_iter010}.png};

\end{groupplot}
\end{tikzpicture}

%% file: tikz/cbar_parula_0x1.tikz
\begin{tikzpicture}
\begin{axis}[
   hide axis, scale only axis,
   height=0pt, width=0pt,
   colormap={parula}{rgb255=(62,38,168) rgb255=(62,39,172) rgb255=(63,40,175) rgb255=(63,41,178) rgb255=(64,42,180) rgb255=(64,43,183) rgb255=(65,44,186) rgb255=(65,45,189) rgb255=(66,46,191) rgb255=(66,47,194) rgb255=(67,48,197) rgb255=(67,49,200) rgb255=(67,50,202) rgb255=(68,51,205) rgb255=(68,52,208) rgb255=(69,53,210) rgb255=(69,55,213) rgb255=(69,56,215) rgb255=(70,57,217) rgb255=(70,58,220) rgb255=(70,59,222) rgb255=(70,61,224) rgb255=(71,62,225) rgb255=(71,63,227) rgb255=(71,65,229) rgb255=(71,66,230) rgb255=(71,68,232) rgb255=(71,69,233) rgb255=(71,70,235) rgb255=(72,72,236) rgb255=(72,73,237) rgb255=(72,75,238) rgb255=(72,76,240) rgb255=(72,78,241) rgb255=(72,79,242) rgb255=(72,80,243) rgb255=(72,82,244) rgb255=(72,83,245) rgb255=(72,84,246) rgb255=(71,86,247) rgb255=(71,87,247) rgb255=(71,89,248) rgb255=(71,90,249) rgb255=(71,91,250) rgb255=(71,93,250) rgb255=(70,94,251) rgb255=(70,96,251) rgb255=(70,97,252) rgb255=(69,98,252) rgb255=(69,100,253) rgb255=(68,101,253) rgb255=(67,103,253) rgb255=(67,104,254) rgb255=(66,106,254) rgb255=(65,107,254) rgb255=(64,109,254) rgb255=(63,110,255) rgb255=(62,112,255) rgb255=(60,113,255) rgb255=(59,115,255) rgb255=(57,116,255) rgb255=(56,118,254) rgb255=(54,119,254) rgb255=(53,121,253) rgb255=(51,122,253) rgb255=(50,124,252) rgb255=(49,125,252) rgb255=(48,127,251) rgb255=(47,128,250) rgb255=(47,130,250) rgb255=(46,131,249) rgb255=(46,132,248) rgb255=(46,134,248) rgb255=(46,135,247) rgb255=(45,136,246) rgb255=(45,138,245) rgb255=(45,139,244) rgb255=(45,140,243) rgb255=(45,142,242) rgb255=(44,143,241) rgb255=(44,144,240) rgb255=(43,145,239) rgb255=(42,147,238) rgb255=(41,148,237) rgb255=(40,149,236) rgb255=(39,151,235) rgb255=(39,152,234) rgb255=(38,153,233) rgb255=(38,154,232) rgb255=(37,155,232) rgb255=(37,156,231) rgb255=(36,158,230) rgb255=(36,159,229) rgb255=(35,160,229) rgb255=(35,161,228) rgb255=(34,162,228) rgb255=(33,163,227) rgb255=(32,165,227) rgb255=(31,166,226) rgb255=(30,167,225) rgb255=(29,168,225) rgb255=(29,169,224) rgb255=(28,170,223) rgb255=(27,171,222) rgb255=(26,172,221) rgb255=(25,173,220) rgb255=(23,174,218) rgb255=(22,175,217) rgb255=(20,176,216) rgb255=(18,177,214) rgb255=(16,178,213) rgb255=(14,179,212) rgb255=(11,179,210) rgb255=(8,180,209) rgb255=(6,181,207) rgb255=(4,182,206) rgb255=(2,183,204) rgb255=(1,183,202) rgb255=(0,184,201) rgb255=(0,185,199) rgb255=(0,186,198) rgb255=(1,186,196) rgb255=(2,187,194) rgb255=(4,187,193) rgb255=(6,188,191) rgb255=(9,189,189) rgb255=(13,189,188) rgb255=(16,190,186) rgb255=(20,190,184) rgb255=(23,191,182) rgb255=(26,192,181) rgb255=(29,192,179) rgb255=(32,193,177) rgb255=(35,193,175) rgb255=(37,194,174) rgb255=(39,194,172) rgb255=(41,195,170) rgb255=(43,195,168) rgb255=(44,196,166) rgb255=(46,196,165) rgb255=(47,197,163) rgb255=(49,197,161) rgb255=(50,198,159) rgb255=(51,199,157) rgb255=(53,199,155) rgb255=(54,200,153) rgb255=(56,200,150) rgb255=(57,201,148) rgb255=(59,201,146) rgb255=(61,202,144) rgb255=(64,202,141) rgb255=(66,202,139) rgb255=(69,203,137) rgb255=(72,203,134) rgb255=(75,203,132) rgb255=(78,204,129) rgb255=(81,204,127) rgb255=(84,204,124) rgb255=(87,204,122) rgb255=(90,204,119) rgb255=(94,205,116) rgb255=(97,205,114) rgb255=(100,205,111) rgb255=(103,205,108) rgb255=(107,205,105) rgb255=(110,205,102) rgb255=(114,205,100) rgb255=(118,204,97) rgb255=(121,204,94) rgb255=(125,204,91) rgb255=(129,204,89) rgb255=(132,204,86) rgb255=(136,203,83) rgb255=(139,203,81) rgb255=(143,203,78) rgb255=(147,202,75) rgb255=(150,202,72) rgb255=(154,201,70) rgb255=(157,201,67) rgb255=(161,200,64) rgb255=(164,200,62) rgb255=(167,199,59) rgb255=(171,199,57) rgb255=(174,198,55) rgb255=(178,198,53) rgb255=(181,197,51) rgb255=(184,196,49) rgb255=(187,196,47) rgb255=(190,195,45) rgb255=(194,195,44) rgb255=(197,194,42) rgb255=(200,193,41) rgb255=(203,193,40) rgb255=(206,192,39) rgb255=(208,191,39) rgb255=(211,191,39) rgb255=(214,190,39) rgb255=(217,190,40) rgb255=(219,189,40) rgb255=(222,188,41) rgb255=(225,188,42) rgb255=(227,188,43) rgb255=(230,187,45) rgb255=(232,187,46) rgb255=(234,186,48) rgb255=(236,186,50) rgb255=(239,186,53) rgb255=(241,186,55) rgb255=(243,186,57) rgb255=(245,186,59) rgb255=(247,186,61) rgb255=(249,186,62) rgb255=(251,187,62) rgb255=(252,188,62) rgb255=(254,189,61) rgb255=(254,190,60) rgb255=(254,192,59) rgb255=(254,193,58) rgb255=(254,194,57) rgb255=(254,196,56) rgb255=(254,197,55) rgb255=(254,199,53) rgb255=(254,200,52) rgb255=(254,202,51) rgb255=(253,203,50) rgb255=(253,205,49) rgb255=(253,206,49) rgb255=(252,208,48) rgb255=(251,210,47) rgb255=(251,211,46) rgb255=(250,213,46) rgb255=(249,214,45) rgb255=(249,216,44) rgb255=(248,217,43) rgb255=(247,219,42) rgb255=(247,221,42) rgb255=(246,222,41) rgb255=(246,224,40) rgb255=(245,225,40) rgb255=(245,227,39) rgb255=(245,229,38) rgb255=(245,230,38) rgb255=(245,232,37) rgb255=(245,233,36) rgb255=(245,235,35) rgb255=(245,236,34) rgb255=(245,238,33) rgb255=(246,239,32) rgb255=(246,241,31) rgb255=(246,242,30) rgb255=(247,244,28) rgb255=(247,245,27) rgb255=(248,247,26) rgb255=(248,248,24) rgb255=(249,249,22) rgb255=(249,251,21) },
   colorbar horizontal,
   point meta min=0.000000e+00, point meta max=1.000000e+00,
   colorbar style={width=10cm, xtick={0.000000e+00,2.500000e-01,5.000000e-01,7.500000e-01,1.000000e+00}}
]
\addplot [draw=none] coordinates {(0,0)};
\end{axis}
\end{tikzpicture}

%% file: py/advec2d_beta0_c1s1_nel6x3_porder0x1_conv.tikz
\begin{tikzpicture}
\begin{groupplot} [
group style={group size = 3 by 1, horizontal sep = 2cm, vertical sep = 1.5cm}]
\nextgroupplot[width=0.4\textwidth, xlabel={Iteration}, ymax=1, xmax=10, ylabel={Solver convergence}, xmin=0, ymode=log, ymin=1e-16, height=0.3\textwidth]
\addplot [solid, thick, color=black, mark=mark options={solid, thin}, mark=*, mark size=1, color=black]
coordinates {
( 0.00000000e+00,  2.84471764e-16)
( 1.00000000e+00,  1.32056958e-02)
( 2.00000000e+00,  1.50172327e-02)
( 3.00000000e+00,  1.11275183e-02)
( 4.00000000e+00,  6.14424825e-03)
( 5.00000000e+00,  2.43314119e-03)
( 6.00000000e+00,  1.41769544e-05)
( 7.00000000e+00,  9.01144294e-08)
( 8.00000000e+00,  3.04319458e-10)
( 9.00000000e+00,  2.97416704e-13)
( 1.00000000e+01,  6.91560954e-16)};\label{line:advec2d_beta0_c1s1:R0}

\addplot [solid, thick, color=blue, mark=mark options={solid, thin}, mark=square*, mark size=1, color=blue]
coordinates {
( 0.00000000e+00,  6.98188131e-02)
( 1.00000000e+00,  5.57289981e-02)
( 2.00000000e+00,  4.35751105e-02)
( 3.00000000e+00,  3.07106857e-02)
( 4.00000000e+00,  1.49780895e-02)
( 5.00000000e+00,  2.15432189e-03)
( 6.00000000e+00,  9.62280164e-05)
( 7.00000000e+00,  5.40771340e-06)
( 8.00000000e+00,  1.66733401e-07)
( 9.00000000e+00,  2.62803750e-09)
( 1.00000000e+01,  2.08911178e-11)};\label{line:advec2d_beta0_c1s1:R1}

\addplot [solid, thick, color=red, mark=mark options={solid, thin}, mark=o, mark size=1, color=red]
coordinates {
( 0.00000000e+00,  1.30499751e-02)
( 1.00000000e+00,  7.48835482e-03)
( 2.00000000e+00,  1.29386264e-02)
( 3.00000000e+00,  1.01668305e-02)
( 4.00000000e+00,  1.11066407e-02)
( 5.00000000e+00,  1.00280900e-03)
( 6.00000000e+00,  7.12058363e-05)
( 7.00000000e+00,  3.60691833e-06)
( 8.00000000e+00,  1.10998319e-07)
( 9.00000000e+00,  1.75124249e-09)
( 1.00000000e+01,  1.39259689e-11)};\label{line:advec2d_beta0_c1s1:dLdY}

\nextgroupplot[width=0.4\textwidth, xlabel={Iteration}, ymax=0.1, xmax=10, xmin=0, ymode=log, ymin=0.001, height=0.3\textwidth]
\addplot [dashed, color=black]
coordinates {
( 1.00000000e+00,  1.00000000e-02)
( 2.00000000e+00,  2.00000000e-02)
( 3.00000000e+00,  4.00000000e-02)
( 4.00000000e+00,  4.00000000e-02)
( 5.00000000e+00,  4.00000000e-02)
( 6.00000000e+00,  4.00000000e-02)
( 7.00000000e+00,  2.00000000e-02)
( 8.00000000e+00,  1.00000000e-02)
( 9.00000000e+00,  5.00000000e-03)
( 1.00000000e+01,  2.50000000e-03)};\label{line:advec2d_beta0_c1s1:lam}

\end{groupplot}\end{tikzpicture}

%% file: tikz/advec2d_beta0_c2s1_nel8x4.tikz
\begin{tikzpicture}
\begin{groupplot}[
  group style={
      group size=2 by 2,
      horizontal sep=1.0cm,
      vertical sep=1.2cm
  },
  width=0.5\textwidth,
  axis equal image,
  axis equal image,
  xlabel={$x_1$},
  ylabel={$x_2$},
  xtick = {-1, -0.5, 0.0, 0.5, 1.0},
  ytick = {0.0, 0.5, 1.0},
  xmin=-1, xmax=1,
  ymin=0, ymax=1
]

\nextgroupplot[title={Initialization}, xlabel={}, xtick=\empty]
\addplot graphics [xmin=-1, xmax=1, ymin=0, ymax=1] {{img/advec2d_beta0_c2s1_nel8x4_porder0x1_iter000}.png};

\nextgroupplot[title={Converged ($q = 1$)}, xlabel={}, ylabel={}, xtick=\empty, ytick=\empty]
\addplot graphics [xmin=-1, xmax=1, ymin=0, ymax=1] {{img/advec2d_beta0_c2s1_nel8x4_porder0x1_iter100}.png};

\nextgroupplot[title={Converged ($q = 2$)}]
\addplot graphics [xmin=-1, xmax=1, ymin=0, ymax=1] {{img/advec2d_beta0_c2s1_nel8x4_porder0x2_iter045}.png};

\nextgroupplot[title={Converged ($q = 3$)}, ylabel={}, ytick=\empty]
\addplot graphics [xmin=-1, xmax=1, ymin=0, ymax=1] {{img/advec2d_beta0_c2s1_nel8x4_porder0x3_iter090}.png};

\end{groupplot}
\end{tikzpicture}

%% file: py/advec2d_beta0_c2s2_nel8x4_porder0x2_conv.tikz
\begin{tikzpicture}
\begin{groupplot} [
group style={group size = 3 by 1, horizontal sep = 2.1cm}]
\nextgroupplot[width=0.32\textwidth, ymin=1e-14, xlabel={Iteration ($k$)}, ymax=0.1, xmax=80, mark repeat={10}, xmin=0, ymode=log, ylabel={$\norm{\rbm(\zbm_k)}$}]
\addplot [solid, thick, color=black, mark=mark options={solid, thin}, mark=*, mark size=1, color=black]
coordinates {
( 1.00000000e+00,  4.31364648e-02)
( 2.00000000e+00,  5.33675085e-03)
( 3.00000000e+00,  6.51592882e-04)
( 4.00000000e+00,  3.55296366e-04)
( 5.00000000e+00,  2.09714444e-04)
( 6.00000000e+00,  1.09898243e-04)
( 7.00000000e+00,  4.16157400e-05)
( 8.00000000e+00,  2.51279208e-05)
( 9.00000000e+00,  1.87045746e-05)
( 1.00000000e+01,  2.83348692e-05)
( 1.10000000e+01,  1.16430605e-05)
( 1.20000000e+01,  5.67777319e-06)
( 1.30000000e+01,  9.93357446e-06)
( 1.40000000e+01,  4.53342274e-06)
( 1.50000000e+01,  2.52191815e-06)
( 1.60000000e+01,  4.99417351e-06)
( 1.70000000e+01,  2.20547251e-06)
( 1.80000000e+01,  1.19326879e-06)
( 1.90000000e+01,  7.55975550e-07)
( 2.00000000e+01,  1.78726550e-06)
( 2.10000000e+01,  9.64829376e-07)
( 2.20000000e+01,  5.60661164e-07)
( 2.30000000e+01,  1.02061825e-06)
( 2.40000000e+01,  3.69929519e-07)
( 2.50000000e+01,  1.49287084e-07)
( 2.60000000e+01,  1.99873812e-07)
( 2.70000000e+01,  5.37434314e-08)
( 2.80000000e+01,  5.63117062e-08)
( 2.90000000e+01,  8.74949106e-09)
( 3.00000000e+01,  6.36645316e-09)
( 3.10000000e+01,  4.94207163e-09)
( 3.20000000e+01,  1.95049172e-09)
( 3.30000000e+01,  2.05438123e-09)
( 3.40000000e+01,  1.45503273e-09)
( 3.50000000e+01,  1.36948622e-09)
( 3.60000000e+01,  8.42802794e-10)
( 3.70000000e+01,  7.32183399e-10)
( 3.80000000e+01,  2.69062125e-10)
( 3.90000000e+01,  3.69064353e-10)
( 4.00000000e+01,  2.66266419e-10)
( 4.10000000e+01,  2.61829131e-10)
( 4.20000000e+01,  1.60415099e-10)
( 4.30000000e+01,  1.67820708e-10)
( 4.40000000e+01,  1.78939951e-10)
( 4.50000000e+01,  1.24897550e-10)
( 4.60000000e+01,  1.10430140e-10)
( 4.70000000e+01,  6.07684189e-11)
( 4.80000000e+01,  7.88325165e-11)
( 4.90000000e+01,  6.09854737e-11)
( 5.00000000e+01,  5.69702746e-11)
( 5.10000000e+01,  3.66213311e-11)
( 5.20000000e+01,  3.28249496e-11)
( 5.30000000e+01,  3.58366511e-11)
( 5.40000000e+01,  2.26699804e-11)
( 5.50000000e+01,  2.06945694e-11)
( 5.60000000e+01,  9.19210090e-12)
( 5.70000000e+01,  1.30357105e-11)
( 5.80000000e+01,  9.75257842e-12)
( 5.90000000e+01,  9.66582661e-12)
( 6.00000000e+01,  5.49338279e-12)
( 6.10000000e+01,  8.45314907e-12)
( 6.20000000e+01,  6.38041790e-12)
( 6.30000000e+01,  6.36596486e-12)
( 6.40000000e+01,  3.97461798e-12)
( 6.50000000e+01,  4.16851211e-12)
( 6.60000000e+01,  4.41656801e-12)
( 6.70000000e+01,  3.03189894e-12)
( 6.80000000e+01,  3.18888644e-12)
( 6.90000000e+01,  3.63026033e-12)
( 7.00000000e+01,  2.37306339e-12)
( 7.10000000e+01,  2.18764491e-12)
( 7.20000000e+01,  1.11056838e-12)
( 7.30000000e+01,  1.50844352e-12)
( 7.40000000e+01,  1.12910030e-12)
( 7.50000000e+01,  1.10989643e-12)
( 7.60000000e+01,  6.68369298e-13)
( 7.70000000e+01,  9.79251216e-13)
( 7.80000000e+01,  7.38337919e-13)
( 7.90000000e+01,  7.31266321e-13)
( 8.00000000e+01,  4.10644412e-13)
( 8.10000000e+01,  6.42735549e-13)
( 8.20000000e+01,  4.87465885e-13)
( 8.30000000e+01,  4.85412129e-13)
( 8.40000000e+01,  3.03114960e-13)
( 8.50000000e+01,  3.14887020e-13)
( 8.60000000e+01,  3.34714918e-13)
( 8.70000000e+01,  2.27903454e-13)
( 8.80000000e+01,  2.41697854e-13)
( 8.90000000e+01,  2.74846364e-13)
( 9.00000000e+01,  1.80434186e-13)
( 9.10000000e+01,  2.13465189e-13)
( 9.20000000e+01,  2.38790536e-13)
( 9.30000000e+01,  1.61864193e-13)
( 9.40000000e+01,  1.46439163e-13)
( 9.50000000e+01,  8.25108691e-14)
( 9.60000000e+01,  1.04810223e-13)
( 9.70000000e+01,  8.12844458e-14)
( 9.80000000e+01,  7.69632458e-14)
( 9.90000000e+01,  5.12479772e-14)};\label{line:advec2d_beta0_c2s2_nel8x4_porder0x2:upwind}

\addplot [solid, thick, color=blue, mark=mark options={solid, thin}, mark=square*, mark size=1, color=blue]
coordinates {
( 1.00000000e+00,  4.42356778e-02)
( 2.00000000e+00,  7.11619923e-03)
( 3.00000000e+00,  2.01594396e-03)
( 4.00000000e+00,  1.01705126e-03)
( 5.00000000e+00,  1.17022641e-03)
( 6.00000000e+00,  1.34504231e-03)
( 7.00000000e+00,  1.43669766e-03)
( 8.00000000e+00,  1.17452348e-03)
( 9.00000000e+00,  1.33900033e-03)
( 1.00000000e+01,  1.00280747e-03)
( 1.10000000e+01,  8.73027794e-04)
( 1.20000000e+01,  6.56909739e-04)
( 1.30000000e+01,  2.73483658e-04)
( 1.40000000e+01,  3.19226582e-04)
( 1.50000000e+01,  3.02389633e-04)
( 1.60000000e+01,  3.41914063e-04)
( 1.70000000e+01,  4.54895289e-04)
( 1.80000000e+01,  6.77821900e-04)
( 1.90000000e+01,  6.41083897e-04)
( 2.00000000e+01,  4.93827750e-04)
( 2.10000000e+01,  5.21676319e-04)
( 2.20000000e+01,  4.92890984e-04)
( 2.30000000e+01,  5.42561188e-04)
( 2.40000000e+01,  4.73898819e-04)
( 2.50000000e+01,  5.21081545e-04)
( 2.60000000e+01,  4.72551529e-04)
( 2.70000000e+01,  4.80571627e-04)
( 2.80000000e+01,  4.79229620e-04)
( 2.90000000e+01,  5.03052339e-04)
( 3.00000000e+01,  4.74172602e-04)
( 3.10000000e+01,  4.85195594e-04)
( 3.20000000e+01,  4.78551270e-04)
( 3.30000000e+01,  4.80279927e-04)
( 3.40000000e+01,  4.76441367e-04)
( 3.50000000e+01,  4.77796915e-04)
( 3.60000000e+01,  4.70601671e-04)
( 3.70000000e+01,  4.72457680e-04)
( 3.80000000e+01,  4.68480106e-04)
( 3.90000000e+01,  4.70252336e-04)
( 4.00000000e+01,  4.67507537e-04)
( 4.10000000e+01,  4.70260438e-04)
( 4.20000000e+01,  4.65528028e-04)
( 4.30000000e+01,  4.68146326e-04)
( 4.40000000e+01,  4.63992100e-04)
( 4.50000000e+01,  4.63857604e-04)
( 4.60000000e+01,  4.63978489e-04)
( 4.70000000e+01,  4.63796245e-04)
( 4.80000000e+01,  4.63759496e-04)
( 4.90000000e+01,  4.63767525e-04)
( 5.00000000e+01,  4.63716174e-04)
( 5.10000000e+01,  4.63762952e-04)
( 5.20000000e+01,  4.63695145e-04)
( 5.30000000e+01,  4.63877933e-04)
( 5.40000000e+01,  4.63540454e-04)
( 5.50000000e+01,  4.63534038e-04)
( 5.60000000e+01,  4.63533480e-04)
( 5.70000000e+01,  4.63534852e-04)
( 5.80000000e+01,  4.63532873e-04)
( 5.90000000e+01,  4.63538110e-04)
( 6.00000000e+01,  4.63528699e-04)
( 6.10000000e+01,  4.63529884e-04)
( 6.20000000e+01,  4.63526635e-04)
( 6.30000000e+01,  4.63528730e-04)
( 6.40000000e+01,  4.63524821e-04)
( 6.50000000e+01,  4.63526823e-04)
( 6.60000000e+01,  4.63522595e-04)
( 6.70000000e+01,  4.63522710e-04)
( 6.80000000e+01,  4.63521621e-04)
( 6.90000000e+01,  4.63524371e-04)
( 7.00000000e+01,  4.63520655e-04)
( 7.10000000e+01,  4.63523859e-04)
( 7.20000000e+01,  4.63519441e-04)
( 7.30000000e+01,  4.63530377e-04)
( 7.40000000e+01,  4.63510861e-04)
( 7.50000000e+01,  4.63513926e-04)
( 7.60000000e+01,  4.63506734e-04)
( 7.70000000e+01,  4.63511622e-04)
( 7.80000000e+01,  4.63503345e-04)
( 7.90000000e+01,  4.63507809e-04)};\label{line:advec2d_beta0_c2s2_nel8x4_porder0x2:upwind_smooth}

\nextgroupplot[width=0.32\textwidth, ymin=1e-09, xlabel={Iteration ($k$)}, ymax=0.1, xmax=80, mark repeat={10}, xmin=0, ymode=log, ylabel={$\norm{\cbm(\zbm_k)}$}]
\addplot [solid, thick, color=black, mark=mark options={solid, thin}, mark=*, mark size=1, color=black, forget plot]
coordinates {
( 1.00000000e+00,  1.85907909e-02)
( 2.00000000e+00,  2.80793189e-03)
( 3.00000000e+00,  1.19510017e-03)
( 4.00000000e+00,  5.06788666e-04)
( 5.00000000e+00,  3.82919151e-04)
( 6.00000000e+00,  3.04933577e-04)
( 7.00000000e+00,  2.52901792e-04)
( 8.00000000e+00,  1.93291278e-04)
( 9.00000000e+00,  1.34525177e-04)
( 1.00000000e+01,  8.59145682e-05)
( 1.10000000e+01,  6.36638295e-05)
( 1.20000000e+01,  5.11402562e-05)
( 1.30000000e+01,  3.82747192e-05)
( 1.40000000e+01,  2.97175754e-05)
( 1.50000000e+01,  2.43844308e-05)
( 1.60000000e+01,  1.87648605e-05)
( 1.70000000e+01,  1.42743443e-05)
( 1.80000000e+01,  1.15227133e-05)
( 1.90000000e+01,  9.59368205e-06)
( 2.00000000e+01,  7.20576367e-06)
( 2.10000000e+01,  5.45128917e-06)
( 2.20000000e+01,  4.38526484e-06)
( 2.30000000e+01,  3.09538414e-06)
( 2.40000000e+01,  2.27477075e-06)
( 2.50000000e+01,  1.73728657e-06)
( 2.60000000e+01,  1.09190511e-06)
( 2.70000000e+01,  6.69671066e-07)
( 2.80000000e+01,  4.14053867e-07)
( 2.90000000e+01,  4.54361726e-07)
( 3.00000000e+01,  6.16580028e-07)
( 3.10000000e+01,  2.72180448e-07)
( 3.20000000e+01,  8.34706780e-07)
( 3.30000000e+01,  3.06745896e-07)
( 3.40000000e+01,  4.63876668e-07)
( 3.50000000e+01,  1.89294781e-07)
( 3.60000000e+01,  2.67328431e-07)
( 3.70000000e+01,  1.26120857e-07)
( 3.80000000e+01,  3.70109068e-07)
( 3.90000000e+01,  1.41418978e-07)
( 4.00000000e+01,  2.11600640e-07)
( 4.10000000e+01,  9.22407483e-08)
( 4.20000000e+01,  1.25792042e-07)
( 4.30000000e+01,  1.96383490e-07)
( 4.40000000e+01,  8.01348130e-08)
( 4.50000000e+01,  1.16178246e-07)
( 4.60000000e+01,  5.56627561e-08)
( 4.70000000e+01,  1.60880587e-07)
( 4.80000000e+01,  6.26189421e-08)
( 4.90000000e+01,  9.40112403e-08)
( 5.00000000e+01,  4.19426090e-08)
( 5.10000000e+01,  5.69897662e-08)
( 5.20000000e+01,  8.71142490e-08)
( 5.30000000e+01,  3.69477068e-08)
( 5.40000000e+01,  5.17034949e-08)
( 5.50000000e+01,  2.63286785e-08)
( 5.60000000e+01,  7.06678480e-08)
( 5.70000000e+01,  2.86108225e-08)
( 5.80000000e+01,  4.14493772e-08)
( 5.90000000e+01,  1.98349607e-08)
( 6.00000000e+01,  5.75182284e-08)
( 6.10000000e+01,  2.25000683e-08)
( 6.20000000e+01,  3.34660293e-08)
( 6.30000000e+01,  1.52279567e-08)
( 6.40000000e+01,  2.03813767e-08)
( 6.50000000e+01,  3.10945063e-08)
( 6.60000000e+01,  1.34601615e-08)
( 6.70000000e+01,  1.86972777e-08)
( 6.80000000e+01,  2.89613044e-08)
( 6.90000000e+01,  1.20668752e-08)
( 7.00000000e+01,  1.71866186e-08)
( 7.10000000e+01,  8.63430984e-09)
( 7.20000000e+01,  2.36434687e-08)
( 7.30000000e+01,  9.60283872e-09)
( 7.40000000e+01,  1.39381169e-08)
( 7.50000000e+01,  6.76604493e-09)
( 7.60000000e+01,  1.93338967e-08)
( 7.70000000e+01,  7.70379859e-09)
( 7.80000000e+01,  1.13429732e-08)
( 7.90000000e+01,  5.35796188e-09)
( 8.00000000e+01,  1.58288323e-08)
( 8.10000000e+01,  6.21567083e-09)
( 8.20000000e+01,  9.25289259e-09)
( 8.30000000e+01,  4.27628408e-09)
( 8.40000000e+01,  5.69700043e-09)
( 8.50000000e+01,  8.61775036e-09)
( 8.60000000e+01,  3.84503496e-09)
( 8.70000000e+01,  5.24584787e-09)
( 8.80000000e+01,  8.04116247e-09)
( 8.90000000e+01,  3.47599942e-09)
( 9.00000000e+01,  4.84322264e-09)
( 9.10000000e+01,  7.51121035e-09)
( 9.20000000e+01,  3.15571103e-09)
( 9.30000000e+01,  4.48509930e-09)
( 9.40000000e+01,  2.30054007e-09)
( 9.50000000e+01,  6.15836516e-09)
( 9.60000000e+01,  2.55939489e-09)
( 9.70000000e+01,  3.66553766e-09)
( 9.80000000e+01,  1.85258677e-09)
( 9.90000000e+01,  5.05130355e-09)};

\addplot [solid, thick, color=blue, mark=mark options={solid, thin}, mark=square*, mark size=1, color=blue, forget plot]
coordinates {
( 1.00000000e+00,  1.97935732e-02)
( 2.00000000e+00,  4.87426655e-03)
( 3.00000000e+00,  3.21521834e-03)
( 4.00000000e+00,  3.03063736e-03)
( 5.00000000e+00,  3.46862424e-03)
( 6.00000000e+00,  2.51107240e-03)
( 7.00000000e+00,  1.55276445e-03)
( 8.00000000e+00,  2.20702322e-03)
( 9.00000000e+00,  1.60785693e-03)
( 1.00000000e+01,  1.31126186e-03)
( 1.10000000e+01,  8.20646434e-04)
( 1.20000000e+01,  7.33465820e-04)
( 1.30000000e+01,  4.14259345e-04)
( 1.40000000e+01,  2.49113483e-04)
( 1.50000000e+01,  5.50160645e-04)
( 1.60000000e+01,  3.97976388e-04)
( 1.70000000e+01,  2.98521107e-04)
( 1.80000000e+01,  4.44359518e-04)
( 1.90000000e+01,  1.75400491e-04)
( 2.00000000e+01,  5.42643966e-04)
( 2.10000000e+01,  1.15187810e-04)
( 2.20000000e+01,  5.31505536e-04)
( 2.30000000e+01,  1.54854950e-04)
( 2.40000000e+01,  4.58849048e-04)
( 2.50000000e+01,  1.03781393e-04)
( 2.60000000e+01,  2.76182104e-04)
( 2.70000000e+01,  3.80903004e-04)
( 2.80000000e+01,  2.44888156e-04)
( 2.90000000e+01,  1.09899507e-04)
( 3.00000000e+01,  5.06944824e-04)
( 3.10000000e+01,  1.03969629e-04)
( 3.20000000e+01,  4.77609719e-04)
( 3.30000000e+01,  1.01848243e-04)
( 3.40000000e+01,  2.89255559e-04)
( 3.50000000e+01,  4.21228143e-04)
( 3.60000000e+01,  2.63288478e-04)
( 3.70000000e+01,  1.49622387e-04)
( 3.80000000e+01,  2.84085504e-04)
( 3.90000000e+01,  9.58096571e-05)
( 4.00000000e+01,  2.60232282e-04)
( 4.10000000e+01,  1.45702179e-04)
( 4.20000000e+01,  2.81403767e-04)
( 4.30000000e+01,  1.43641540e-04)
( 4.40000000e+01,  4.81435698e-04)
( 4.50000000e+01,  4.71767019e-04)
( 4.60000000e+01,  1.24161701e-04)
( 4.70000000e+01,  9.33821757e-05)
( 4.80000000e+01,  2.58579280e-04)
( 4.90000000e+01,  3.65657416e-04)
( 5.00000000e+01,  2.58571444e-04)
( 5.10000000e+01,  9.33684939e-05)
( 5.20000000e+01,  2.58584216e-04)
( 5.30000000e+01,  3.65683150e-04)
( 5.40000000e+01,  1.20808592e-04)
( 5.50000000e+01,  8.89534674e-05)
( 5.60000000e+01,  2.57028466e-04)
( 5.70000000e+01,  8.89538652e-05)
( 5.80000000e+01,  2.57028974e-04)
( 5.90000000e+01,  8.89554562e-05)
( 6.00000000e+01,  2.32059005e-04)
( 6.10000000e+01,  8.89540662e-05)
( 6.20000000e+01,  2.32058684e-04)
( 6.30000000e+01,  8.89540315e-05)
( 6.40000000e+01,  2.57027617e-04)
( 6.50000000e+01,  8.89541494e-05)
( 6.60000000e+01,  2.32057892e-04)
( 6.70000000e+01,  8.89534544e-05)
( 6.80000000e+01,  2.57027052e-04)
( 6.90000000e+01,  1.20808033e-04)
( 7.00000000e+01,  2.57027138e-04)
( 7.10000000e+01,  8.89540070e-05)
( 7.20000000e+01,  2.57028154e-04)
( 7.30000000e+01,  8.89571894e-05)
( 7.40000000e+01,  2.32056661e-04)
( 7.50000000e+01,  8.89544095e-05)
( 7.60000000e+01,  2.32056019e-04)
( 7.70000000e+01,  8.89543403e-05)
( 7.80000000e+01,  2.57025440e-04)
( 7.90000000e+01,  8.89545765e-05)};

\nextgroupplot[width=0.32\textwidth, ymin=0.0001, xlabel={Iteration ($k$)}, ymax=0.1, xmax=80, mark repeat={10}, xmin=0, ymode=log, ylabel={$\norm{\Rbm(\zbm_k)}$}]
\addplot [solid, thick, color=black, mark=mark options={solid, thin}, mark=*, mark size=1, color=black, forget plot]
coordinates {
( 1.00000000e+00,  3.66585183e-02)
( 2.00000000e+00,  8.55144415e-03)
( 3.00000000e+00,  4.40474972e-03)
( 4.00000000e+00,  2.90423450e-03)
( 5.00000000e+00,  2.14157540e-03)
( 6.00000000e+00,  1.73490401e-03)
( 7.00000000e+00,  1.52102115e-03)
( 8.00000000e+00,  1.33074507e-03)
( 9.00000000e+00,  1.15842029e-03)
( 1.00000000e+01,  9.75358366e-04)
( 1.10000000e+01,  8.64126855e-04)
( 1.20000000e+01,  7.83359411e-04)
( 1.30000000e+01,  6.71432161e-04)
( 1.40000000e+01,  5.99506588e-04)
( 1.50000000e+01,  5.48817534e-04)
( 1.60000000e+01,  4.83120358e-04)
( 1.70000000e+01,  4.45344932e-04)
( 1.80000000e+01,  4.21293042e-04)
( 1.90000000e+01,  4.05182475e-04)
( 2.00000000e+01,  3.85822106e-04)
( 2.10000000e+01,  3.75616883e-04)
( 2.20000000e+01,  3.69848032e-04)
( 2.30000000e+01,  3.64453583e-04)
( 2.40000000e+01,  3.62490145e-04)
( 2.50000000e+01,  3.61532692e-04)
( 2.60000000e+01,  3.61036800e-04)
( 2.70000000e+01,  3.60898384e-04)
( 2.80000000e+01,  3.60648967e-04)
( 2.90000000e+01,  3.60763155e-04)
( 3.00000000e+01,  3.60586699e-04)
( 3.10000000e+01,  3.60615149e-04)
( 3.20000000e+01,  3.60544368e-04)
( 3.30000000e+01,  3.60565892e-04)
( 3.40000000e+01,  3.60514714e-04)
( 3.50000000e+01,  3.60549914e-04)
( 3.60000000e+01,  3.60519824e-04)
( 3.70000000e+01,  3.60539687e-04)
( 3.80000000e+01,  3.60491159e-04)
( 3.90000000e+01,  3.60517585e-04)
( 4.00000000e+01,  3.60489037e-04)
( 4.10000000e+01,  3.60507434e-04)
( 4.20000000e+01,  3.60490560e-04)
( 4.30000000e+01,  3.60511331e-04)
( 4.40000000e+01,  3.60486659e-04)
( 4.50000000e+01,  3.60501427e-04)
( 4.60000000e+01,  3.60487577e-04)
( 4.70000000e+01,  3.60502167e-04)
( 4.80000000e+01,  3.60481650e-04)
( 4.90000000e+01,  3.60493967e-04)
( 5.00000000e+01,  3.60483192e-04)
( 5.10000000e+01,  3.60490239e-04)
( 5.20000000e+01,  3.60478382e-04)
( 5.30000000e+01,  3.60487605e-04)
( 5.40000000e+01,  3.60481794e-04)
( 5.50000000e+01,  3.60487116e-04)
( 5.60000000e+01,  3.60480667e-04)
( 5.70000000e+01,  3.60488259e-04)
( 5.80000000e+01,  3.60483588e-04)
( 5.90000000e+01,  3.60488114e-04)
( 6.00000000e+01,  3.60483323e-04)
( 6.10000000e+01,  3.60489697e-04)
( 6.20000000e+01,  3.60486004e-04)
( 6.30000000e+01,  3.60489793e-04)
( 6.40000000e+01,  3.60487933e-04)
( 6.50000000e+01,  3.60492420e-04)
( 6.60000000e+01,  3.60489231e-04)
( 6.70000000e+01,  3.60492037e-04)
( 6.80000000e+01,  3.60488647e-04)
( 6.90000000e+01,  3.60492014e-04)
( 7.00000000e+01,  3.60490361e-04)
( 7.10000000e+01,  3.60492342e-04)
( 7.20000000e+01,  3.60490775e-04)
( 7.30000000e+01,  3.60493527e-04)
( 7.40000000e+01,  3.60492150e-04)
( 7.50000000e+01,  3.60493781e-04)
( 7.60000000e+01,  3.60492498e-04)
( 7.70000000e+01,  3.60494763e-04)
( 7.80000000e+01,  3.60493628e-04)
( 7.90000000e+01,  3.60494966e-04)
( 8.00000000e+01,  3.60493900e-04)
( 8.10000000e+01,  3.60495761e-04)
( 8.20000000e+01,  3.60494819e-04)
( 8.30000000e+01,  3.60495916e-04)
( 8.40000000e+01,  3.60495469e-04)
( 8.50000000e+01,  3.60496782e-04)
( 8.60000000e+01,  3.60495942e-04)
( 8.70000000e+01,  3.60496758e-04)
( 8.80000000e+01,  3.60495881e-04)
( 8.90000000e+01,  3.60496837e-04)
( 9.00000000e+01,  3.60496411e-04)
( 9.10000000e+01,  3.60497543e-04)
( 9.20000000e+01,  3.60496787e-04)
( 9.30000000e+01,  3.60497482e-04)
( 9.40000000e+01,  3.60497081e-04)
( 9.50000000e+01,  3.60497955e-04)
( 9.60000000e+01,  3.60497329e-04)
( 9.70000000e+01,  3.60497893e-04)
( 9.80000000e+01,  3.60497561e-04)
( 9.90000000e+01,  3.60498268e-04)};

\addplot [solid, thick, color=blue, mark=mark options={solid, thin}, mark=square*, mark size=1, color=blue, forget plot]
coordinates {
( 1.00000000e+00,  3.73918382e-02)
( 2.00000000e+00,  1.03651684e-02)
( 3.00000000e+00,  6.65219422e-03)
( 4.00000000e+00,  5.04568122e-03)
( 5.00000000e+00,  4.43887934e-03)
( 6.00000000e+00,  3.89771867e-03)
( 7.00000000e+00,  3.54646227e-03)
( 8.00000000e+00,  2.89386852e-03)
( 9.00000000e+00,  2.53390660e-03)
( 1.00000000e+01,  1.78740884e-03)
( 1.10000000e+01,  1.53086800e-03)
( 1.20000000e+01,  1.20951707e-03)
( 1.30000000e+01,  1.03988285e-03)
( 1.40000000e+01,  1.01420476e-03)
( 1.50000000e+01,  1.01784391e-03)
( 1.60000000e+01,  9.98087915e-04)
( 1.70000000e+01,  9.82330218e-04)
( 1.80000000e+01,  1.05336419e-03)
( 1.90000000e+01,  1.02427913e-03)
( 2.00000000e+01,  9.84869273e-04)
( 2.10000000e+01,  9.53069250e-04)
( 2.20000000e+01,  9.49742563e-04)
( 2.30000000e+01,  9.43792670e-04)
( 2.40000000e+01,  9.34551820e-04)
( 2.50000000e+01,  9.28112197e-04)
( 2.60000000e+01,  9.34232427e-04)
( 2.70000000e+01,  9.31787615e-04)
( 2.80000000e+01,  9.24743818e-04)
( 2.90000000e+01,  9.29764369e-04)
( 3.00000000e+01,  9.28514803e-04)
( 3.10000000e+01,  9.22682546e-04)
( 3.20000000e+01,  9.22433475e-04)
( 3.30000000e+01,  9.19820625e-04)
( 3.40000000e+01,  9.20364155e-04)
( 3.50000000e+01,  9.20104987e-04)
( 3.60000000e+01,  9.15101966e-04)
( 3.70000000e+01,  9.15298625e-04)
( 3.80000000e+01,  9.14902353e-04)
( 3.90000000e+01,  9.15116278e-04)
( 4.00000000e+01,  9.15021955e-04)
( 4.10000000e+01,  9.15577585e-04)
( 4.20000000e+01,  9.14852171e-04)
( 4.30000000e+01,  9.15407684e-04)
( 4.40000000e+01,  9.15095386e-04)
( 4.50000000e+01,  9.14349366e-04)
( 4.60000000e+01,  9.14276992e-04)
( 4.70000000e+01,  9.14236579e-04)
( 4.80000000e+01,  9.14239879e-04)
( 4.90000000e+01,  9.14237420e-04)
( 5.00000000e+01,  9.14199383e-04)
( 5.10000000e+01,  9.14210284e-04)
( 5.20000000e+01,  9.14219287e-04)
( 5.30000000e+01,  9.14275638e-04)
( 5.40000000e+01,  9.14064178e-04)
( 5.50000000e+01,  9.14062962e-04)
( 5.60000000e+01,  9.14063114e-04)
( 5.70000000e+01,  9.14063497e-04)
( 5.80000000e+01,  9.14064212e-04)
( 5.90000000e+01,  9.14065638e-04)
( 6.00000000e+01,  9.14064792e-04)
( 6.10000000e+01,  9.14064821e-04)
( 6.20000000e+01,  9.14065007e-04)
( 6.30000000e+01,  9.14065431e-04)
( 6.40000000e+01,  9.14065335e-04)
( 6.50000000e+01,  9.14065788e-04)
( 6.60000000e+01,  9.14065562e-04)
( 6.70000000e+01,  9.14065380e-04)
( 6.80000000e+01,  9.14065696e-04)
( 6.90000000e+01,  9.14066465e-04)
( 7.00000000e+01,  9.14066405e-04)
( 7.10000000e+01,  9.14067364e-04)
( 7.20000000e+01,  9.14068601e-04)
( 7.30000000e+01,  9.14071646e-04)
( 7.40000000e+01,  9.14069656e-04)
( 7.50000000e+01,  9.14070013e-04)
( 7.60000000e+01,  9.14070085e-04)
( 7.70000000e+01,  9.14071232e-04)
( 7.80000000e+01,  9.14070847e-04)
( 7.90000000e+01,  9.14071947e-04)};

\end{groupplot}\end{tikzpicture}

%% file: tikz/burg1d_sptm_c1s1_nel8x4_porder1x1.tikz
\begin{tikzpicture}
\begin{groupplot}[
  group style={
      group size=2 by 2,
      horizontal sep=1.0cm,
      vertical sep=1.2cm
  },
  width=0.5\textwidth,
  axis equal image,
  xlabel={$x$},
  ylabel={$t$},
  xtick = {-1, -0.5, 0.0, 0.5, 1.0},
  ytick = {0.0, 0.5, 1.0},
  xmin=-1, xmax=1,
  ymin=0, ymax=1
]

\nextgroupplot[title={Iteration 0}, xlabel={}, xtick=\empty]
\addplot graphics [xmin=-1, xmax=1, ymin=0, ymax=1] {{img/burg1d_sptm_c1s1_nel8x4_porder1x1_mesh_iter000}.png};

\nextgroupplot[title={Iteration 20}, xlabel={}, ylabel={}, xtick=\empty, ytick=\empty]
\addplot graphics [xmin=-1, xmax=1, ymin=0, ymax=1] {{img/burg1d_sptm_c1s1_nel8x4_porder1x1_mesh_iter020}.png};

\nextgroupplot[title={Iteration 30}]
\addplot graphics [xmin=-1, xmax=1, ymin=0, ymax=1] {{img/burg1d_sptm_c1s1_nel8x4_porder1x1_mesh_iter030}.png};

\nextgroupplot[title={Iteration 40}, ylabel={}, ytick=\empty]
\addplot graphics [xmin=-1, xmax=1, ymin=0, ymax=1] {{img/burg1d_sptm_c1s1_nel8x4_porder1x1_mesh_iter040}.png};

\end{groupplot}
\end{tikzpicture}

%% file: tikz/burg1d_sptm_c1s1_nel8x4.tikz
\begin{tikzpicture}
\begin{groupplot}[
  group style={
      group size=2 by 3,
      horizontal sep=1.0cm,
      vertical sep=1.2cm
  },
  width=0.5\textwidth,
  axis equal image,
  xlabel={$x$},
  ylabel={$t$},
  xtick = {-1, -0.5, 0.0, 0.5, 1.0},
  ytick = {0.0, 0.5, 1.0},
  xmin=-1, xmax=1,
  ymin=0, ymax=1
]

\nextgroupplot[title={$p = q = 1$}, xlabel={}, xtick=\empty]
\addplot graphics [xmin=-1, xmax=1, ymin=0, ymax=1] {{img/burg1d_sptm_c1s1_nel8x4_porder1x1_mesh_iter040}.png};

\nextgroupplot[title={$p = q = 1$}, xlabel={}, ylabel={}, xtick=\empty, ytick=\empty]
\addplot graphics [xmin=-1, xmax=1, ymin=0, ymax=1] {{img/burg1d_sptm_c1s1_nel8x4_porder1x1_iter040}.png};

\nextgroupplot[title={$p = q = 2$}, xlabel={}, xtick=\empty]
\addplot graphics [xmin=-1, xmax=1, ymin=0, ymax=1] {{img/burg1d_sptm_c1s1_nel8x4_porder2x2_mesh_iter040}.png};

\nextgroupplot[title={$p = q = 2$}, xlabel={}, ylabel={}, xtick=\empty, ytick=\empty]
\addplot graphics [xmin=-1, xmax=1, ymin=0, ymax=1] {{img/burg1d_sptm_c1s1_nel8x4_porder2x2_iter040}.png};

\nextgroupplot[title={$p = q = 3$}]
\addplot graphics [xmin=-1, xmax=1, ymin=0, ymax=1] {{img/burg1d_sptm_c1s1_nel8x4_porder3x3_mesh_iter040}.png};

\nextgroupplot[title={$p = q = 3$}, ylabel={}, ytick=\empty]
\addplot graphics [xmin=-1, xmax=1, ymin=0, ymax=1] {{img/burg1d_sptm_c1s1_nel8x4_porder3x3_iter040}.png};

\end{groupplot}
\end{tikzpicture}

%% file: py/burg1d_sptm_c1s1_nel8x4_porder4x4_slice.tikz
\begin{tikzpicture}
\begin{axis}[
width=0.7\textwidth,
xlabel={$x$},
ymax=2,
xmax=1,
ylabel={$U(x, t)$},
xmin=-1,
ymin=0,
height=0.4\textwidth]
\addplot [solid, thick, color=black]
coordinates {
( 0.00000000e+00,  0.00000000e+00)
(-9.79800000e-01,  8.16240000e-04)
(-9.59600000e-01,  3.26500000e-03)
(-9.39390000e-01,  7.34620000e-03)
(-9.19190000e-01,  1.30600000e-02)
(-8.98990000e-01,  2.04060000e-02)
(-8.78790000e-01,  2.93850000e-02)
(-8.58590000e-01,  3.99960000e-02)
(-8.38380000e-01,  5.22400000e-02)
(-8.18180000e-01,  6.61160000e-02)
(-7.97980000e-01,  8.16240000e-02)
(-7.77780000e-01,  9.87650000e-02)
(-7.57580000e-01,  1.17540000e-01)
(-7.37370000e-01,  1.37950000e-01)
(-7.17170000e-01,  1.59980000e-01)
(-6.96970000e-01,  1.83650000e-01)
(-6.76770000e-01,  2.08960000e-01)
(-6.56570000e-01,  2.35890000e-01)
(-6.36360000e-01,  2.64460000e-01)
(-6.16160000e-01,  2.94660000e-01)
(-5.95960000e-01,  3.26500000e-01)
(-5.75760000e-01,  3.59960000e-01)
(-5.55560000e-01,  3.95060000e-01)
(-5.35350000e-01,  4.31790000e-01)
(-5.15150000e-01,  4.70160000e-01)
(-4.94950000e-01,  5.10150000e-01)
(-4.74750000e-01,  5.51780000e-01)
(-4.54550000e-01,  5.95040000e-01)
(-4.34340000e-01,  6.39930000e-01)
(-4.14140000e-01,  6.86460000e-01)
(-3.93940000e-01,  7.34620000e-01)
(-3.73740000e-01,  7.84410000e-01)
(-3.53540000e-01,  8.35830000e-01)
(-3.33330000e-01,  8.88890000e-01)
(-3.13130000e-01,  9.43580000e-01)
(-2.92930000e-01,  9.99900000e-01)
(-2.72730000e-01,  1.05790000e+00)
(-2.52530000e-01,  1.11740000e+00)
(-2.32320000e-01,  1.17870000e+00)
(-2.12120000e-01,  1.24150000e+00)
(-1.91920000e-01,  1.30600000e+00)
(-1.71720000e-01,  1.37210000e+00)
(-1.51520000e-01,  1.43990000e+00)
(-1.31310000e-01,  1.50920000e+00)
(-1.11110000e-01,  1.58020000e+00)
(-9.09090000e-02,  1.65290000e+00)
(-7.07070000e-02,  1.72720000e+00)
(-5.05050000e-02,  1.80310000e+00)
(-3.03030000e-02,  1.88060000e+00)
(-1.01010000e-02,  1.95980000e+00)
( 1.01010000e-02,  0.00000000e+00)
( 3.03030000e-02,  0.00000000e+00)
( 5.05050000e-02,  0.00000000e+00)
( 7.07070000e-02,  0.00000000e+00)
( 9.09090000e-02,  0.00000000e+00)
( 1.11110000e-01,  0.00000000e+00)
( 1.31310000e-01,  0.00000000e+00)
( 1.51520000e-01,  0.00000000e+00)
( 1.71720000e-01,  0.00000000e+00)
( 1.91920000e-01,  0.00000000e+00)
( 2.12120000e-01,  0.00000000e+00)
( 2.32320000e-01,  0.00000000e+00)
( 2.52530000e-01,  0.00000000e+00)
( 2.72730000e-01,  0.00000000e+00)
( 2.92930000e-01,  0.00000000e+00)
( 3.13130000e-01,  0.00000000e+00)
( 3.33330000e-01,  0.00000000e+00)
( 3.53540000e-01,  0.00000000e+00)
( 3.73740000e-01,  0.00000000e+00)
( 3.93940000e-01,  0.00000000e+00)
( 4.14140000e-01,  0.00000000e+00)
( 4.34340000e-01,  0.00000000e+00)
( 4.54550000e-01,  0.00000000e+00)
( 4.74750000e-01,  0.00000000e+00)
( 4.94950000e-01,  0.00000000e+00)
( 5.15150000e-01,  0.00000000e+00)
( 5.35350000e-01,  0.00000000e+00)
( 5.55560000e-01,  0.00000000e+00)
( 5.75760000e-01,  0.00000000e+00)
( 5.95960000e-01,  0.00000000e+00)
( 6.16160000e-01,  0.00000000e+00)
( 6.36360000e-01,  0.00000000e+00)
( 6.56570000e-01,  0.00000000e+00)
( 6.76770000e-01,  0.00000000e+00)
( 6.96970000e-01,  0.00000000e+00)
( 7.17170000e-01,  0.00000000e+00)
( 7.37370000e-01,  0.00000000e+00)
( 7.57580000e-01,  0.00000000e+00)
( 7.77780000e-01,  0.00000000e+00)
( 7.97980000e-01,  0.00000000e+00)
( 8.18180000e-01,  0.00000000e+00)
( 8.38380000e-01,  0.00000000e+00)
( 8.58590000e-01,  0.00000000e+00)
( 8.78790000e-01,  0.00000000e+00)
( 8.98990000e-01,  0.00000000e+00)
( 9.19190000e-01,  0.00000000e+00)
( 9.39390000e-01,  0.00000000e+00)
( 9.59600000e-01,  0.00000000e+00)
( 9.79800000e-01,  0.00000000e+00)
( 0.00000000e+00,  0.00000000e+00)};\label{line:burg1d_sptm_slice}

\addplot [solid, thick, color=black]
coordinates {
( 0.00000000e+00,  0.00000000e+00)
(-9.79800000e-01,  8.04670000e-04)
(-9.59600000e-01,  3.14470000e-03)
(-9.39390000e-01,  6.93240000e-03)
(-9.19190000e-01,  1.20960000e-02)
(-8.98990000e-01,  1.85680000e-02)
(-8.78790000e-01,  2.62820000e-02)
(-8.58590000e-01,  3.51770000e-02)
(-8.38380000e-01,  4.51950000e-02)
(-8.18180000e-01,  5.62840000e-02)
(-7.97980000e-01,  6.83960000e-02)
(-7.77780000e-01,  8.14870000e-02)
(-7.57580000e-01,  9.55170000e-02)
(-7.37370000e-01,  1.10450000e-01)
(-7.17170000e-01,  1.26270000e-01)
(-6.96970000e-01,  1.42900000e-01)
(-6.76770000e-01,  1.60340000e-01)
(-6.56570000e-01,  1.78560000e-01)
(-6.36360000e-01,  1.97520000e-01)
(-6.16160000e-01,  2.17190000e-01)
(-5.95960000e-01,  2.37560000e-01)
(-5.75760000e-01,  2.58610000e-01)
(-5.55560000e-01,  2.80300000e-01)
(-5.35350000e-01,  3.02630000e-01)
(-5.15150000e-01,  3.25560000e-01)
(-4.94950000e-01,  3.49080000e-01)
(-4.74750000e-01,  3.73180000e-01)
(-4.54550000e-01,  3.97830000e-01)
(-4.34340000e-01,  4.23020000e-01)
(-4.14140000e-01,  4.48740000e-01)
(-3.93940000e-01,  4.74960000e-01)
(-3.73740000e-01,  5.01680000e-01)
(-3.53540000e-01,  5.28890000e-01)
(-3.33330000e-01,  5.56570000e-01)
(-3.13130000e-01,  5.84700000e-01)
(-2.92930000e-01,  6.13280000e-01)
(-2.72730000e-01,  6.42300000e-01)
(-2.52530000e-01,  6.71740000e-01)
(-2.32320000e-01,  7.01590000e-01)
(-2.12120000e-01,  7.31850000e-01)
(-1.91920000e-01,  7.62500000e-01)
(-1.71720000e-01,  7.93540000e-01)
(-1.51520000e-01,  8.24960000e-01)
(-1.31310000e-01,  8.56740000e-01)
(-1.11110000e-01,  8.88890000e-01)
(-9.09090000e-02,  9.21390000e-01)
(-7.07070000e-02,  9.54230000e-01)
(-5.05050000e-02,  9.87410000e-01)
(-3.03030000e-02,  1.02090000e+00)
(-1.01010000e-02,  1.05480000e+00)
( 1.01010000e-02,  1.08890000e+00)
( 3.03030000e-02,  1.12340000e+00)
( 5.05050000e-02,  1.15810000e+00)
( 7.07070000e-02,  1.19320000e+00)
( 9.09090000e-02,  1.22860000e+00)
( 1.11110000e-01,  1.26420000e+00)
( 1.31310000e-01,  1.30020000e+00)
( 1.51520000e-01,  1.33640000e+00)
( 1.71720000e-01,  1.37280000e+00)
( 1.91920000e-01,  1.40960000e+00)
( 2.12120000e-01, -3.29810000e-06)
( 2.32320000e-01, -1.69600000e-06)
( 2.52530000e-01, -5.98330000e-07)
( 2.72730000e-01,  1.04280000e-07)
( 2.92930000e-01,  5.01300000e-07)
( 3.13130000e-01,  6.66340000e-07)
( 3.33330000e-01,  6.60580000e-07)
( 3.53540000e-01,  5.35400000e-07)
( 3.73740000e-01,  3.34450000e-07)
( 3.93940000e-01,  9.52640000e-08)
( 4.14140000e-01, -1.49400000e-07)
( 4.34340000e-01, -3.70680000e-07)
( 4.54550000e-01, -5.42860000e-07)
( 4.74750000e-01, -6.42900000e-07)
( 4.94950000e-01, -6.50230000e-07)
( 5.15150000e-01, -5.46680000e-07)
( 5.35350000e-01, -3.16600000e-07)
( 5.55560000e-01,  5.26690000e-08)
( 5.75760000e-01,  5.70260000e-07)
( 5.95960000e-01,  1.24090000e-06)
( 6.16160000e-01,  2.06350000e-06)
( 6.36360000e-01,  3.03000000e-06)
( 6.56570000e-01,  3.67690000e-06)
( 6.76770000e-01,  3.69450000e-06)
( 6.96970000e-01,  1.91390000e-06)
( 7.17170000e-01, -2.41510000e-06)
( 7.37370000e-01,  4.54230000e-07)
( 7.57580000e-01,  4.81980000e-08)
( 7.77780000e-01, -1.07490000e-07)
( 7.97980000e-01, -1.12950000e-07)
( 8.18180000e-01, -4.58210000e-08)
( 8.38380000e-01,  3.64980000e-08)
( 8.58590000e-01,  9.49950000e-08)
( 8.78790000e-01,  1.07450000e-07)
( 8.98990000e-01,  6.71580000e-08)
( 9.19190000e-01, -1.81620000e-08)
( 9.39390000e-01, -1.27150000e-07)
( 9.59600000e-01, -2.25490000e-07)
( 9.79800000e-01, -1.96700000e-07)
( 0.00000000e+00,  0.00000000e+00)};\label{line:burg1d_sptm_slice}

\addplot [solid, thick, color=black]
coordinates {
( 0.00000000e+00,  0.00000000e+00)
(-9.79800000e-01,  7.90590000e-04)
(-9.59600000e-01,  3.01490000e-03)
(-9.39390000e-01,  6.50670000e-03)
(-9.19190000e-01,  1.11460000e-02)
(-8.98990000e-01,  1.68260000e-02)
(-8.78790000e-01,  2.34540000e-02)
(-8.58590000e-01,  3.09450000e-02)
(-8.38380000e-01,  3.92280000e-02)
(-8.18180000e-01,  4.82390000e-02)
(-7.97980000e-01,  5.79230000e-02)
(-7.77780000e-01,  6.82320000e-02)
(-7.57580000e-01,  7.91260000e-02)
(-7.37370000e-01,  9.05690000e-02)
(-7.17170000e-01,  1.02530000e-01)
(-6.96970000e-01,  1.14990000e-01)
(-6.76770000e-01,  1.27920000e-01)
(-6.56570000e-01,  1.41250000e-01)
(-6.36360000e-01,  1.55000000e-01)
(-6.16160000e-01,  1.69140000e-01)
(-5.95960000e-01,  1.83650000e-01)
(-5.75760000e-01,  1.98520000e-01)
(-5.55560000e-01,  2.13720000e-01)
(-5.35350000e-01,  2.29250000e-01)
(-5.15150000e-01,  2.45080000e-01)
(-4.94950000e-01,  2.61200000e-01)
(-4.74750000e-01,  2.77610000e-01)
(-4.54550000e-01,  2.94290000e-01)
(-4.34340000e-01,  3.11230000e-01)
(-4.14140000e-01,  3.28420000e-01)
(-3.93940000e-01,  3.45850000e-01)
(-3.73740000e-01,  3.63510000e-01)
(-3.53540000e-01,  3.81400000e-01)
(-3.33330000e-01,  3.99510000e-01)
(-3.13130000e-01,  4.17820000e-01)
(-2.92930000e-01,  4.36340000e-01)
(-2.72730000e-01,  4.55050000e-01)
(-2.52530000e-01,  4.73950000e-01)
(-2.32320000e-01,  4.93040000e-01)
(-2.12120000e-01,  5.12300000e-01)
(-1.91920000e-01,  5.31740000e-01)
(-1.71720000e-01,  5.51340000e-01)
(-1.51520000e-01,  5.71110000e-01)
(-1.31310000e-01,  5.91040000e-01)
(-1.11110000e-01,  6.11120000e-01)
(-9.09090000e-02,  6.31350000e-01)
(-7.07070000e-02,  6.51730000e-01)
(-5.05050000e-02,  6.72250000e-01)
(-3.03030000e-02,  6.92900000e-01)
(-1.01010000e-02,  7.13690000e-01)
( 1.01010000e-02,  7.34620000e-01)
( 3.03030000e-02,  7.55670000e-01)
( 5.05050000e-02,  7.76850000e-01)
( 7.07070000e-02,  7.98150000e-01)
( 9.09090000e-02,  8.19570000e-01)
( 1.11110000e-01,  8.41110000e-01)
( 1.31310000e-01,  8.62760000e-01)
( 1.51520000e-01,  8.84520000e-01)
( 1.71720000e-01,  9.06390000e-01)
( 1.91920000e-01,  9.28370000e-01)
( 2.12120000e-01,  9.50460000e-01)
( 2.32320000e-01,  9.72640000e-01)
( 2.52530000e-01,  9.94930000e-01)
( 2.72730000e-01,  1.01730000e+00)
( 2.92930000e-01,  1.03980000e+00)
( 3.13130000e-01,  1.06240000e+00)
( 3.33330000e-01,  1.08510000e+00)
( 3.53540000e-01,  1.10780000e+00)
( 3.73740000e-01,  1.13060000e+00)
( 3.93940000e-01,  1.15360000e+00)
( 4.14140000e-01, -1.13080000e-05)
( 4.34340000e-01,  2.87780000e-06)
( 4.54550000e-01,  1.06640000e-05)
( 4.74750000e-01,  1.32940000e-05)
( 4.94950000e-01,  1.20220000e-05)
( 5.15150000e-01,  8.06750000e-06)
( 5.35350000e-01,  2.58030000e-06)
( 5.55560000e-01, -3.38690000e-06)
( 5.75760000e-01, -8.90340000e-06)
( 5.95960000e-01, -1.31790000e-05)
( 6.16160000e-01, -1.55790000e-05)
( 6.36360000e-01,  1.95620000e-06)
( 6.56570000e-01,  1.12380000e-05)
( 6.76770000e-01,  3.76220000e-05)
( 6.96970000e-01,  8.48380000e-05)
( 7.17170000e-01,  1.57520000e-04)
( 7.37370000e-01,  2.61290000e-04)
( 7.57580000e-01,  9.70210000e-05)
( 7.77780000e-01,  6.57950000e-05)
( 7.97980000e-01,  3.97200000e-05)
( 8.18180000e-01,  1.84180000e-05)
( 8.38380000e-01,  1.54960000e-06)
( 8.58590000e-01, -1.11860000e-05)
( 8.78790000e-01, -2.00500000e-05)
( 8.98990000e-01, -2.52550000e-05)
( 9.19190000e-01, -2.69720000e-05)
( 9.39390000e-01, -2.47110000e-05)
( 9.59600000e-01, -1.63570000e-05)
( 9.79800000e-01,  7.58310000e-06)
( 0.00000000e+00,  0.00000000e+00)};\label{line:burg1d_sptm_slice}

\addplot [solid, thick, color=black]
coordinates {
( 0.00000000e+00,  0.00000000e+00)
(-9.79800000e-01,  7.75910000e-04)
(-9.59600000e-01,  2.89680000e-03)
(-9.39390000e-01,  6.14040000e-03)
(-9.19190000e-01,  1.03600000e-02)
(-8.98990000e-01,  1.54340000e-02)
(-8.78790000e-01,  2.12610000e-02)
(-8.58590000e-01,  2.77560000e-02)
(-8.38380000e-01,  3.48480000e-02)
(-8.18180000e-01,  4.24760000e-02)
(-7.97980000e-01,  5.05900000e-02)
(-7.77780000e-01,  5.91480000e-02)
(-7.57580000e-01,  6.81120000e-02)
(-7.37370000e-01,  7.74540000e-02)
(-7.17170000e-01,  8.71460000e-02)
(-6.96970000e-01,  9.71670000e-02)
(-6.76770000e-01,  1.07500000e-01)
(-6.56570000e-01,  1.18130000e-01)
(-6.36360000e-01,  1.29010000e-01)
(-6.16160000e-01,  1.40150000e-01)
(-5.95960000e-01,  1.51520000e-01)
(-5.75760000e-01,  1.63120000e-01)
(-5.55560000e-01,  1.74930000e-01)
(-5.35350000e-01,  1.86950000e-01)
(-5.15150000e-01,  1.99160000e-01)
(-4.94950000e-01,  2.11550000e-01)
(-4.74750000e-01,  2.24120000e-01)
(-4.54550000e-01,  2.36850000e-01)
(-4.34340000e-01,  2.49750000e-01)
(-4.14140000e-01,  2.62790000e-01)
(-3.93940000e-01,  2.75990000e-01)
(-3.73740000e-01,  2.89320000e-01)
(-3.53540000e-01,  3.02790000e-01)
(-3.33330000e-01,  3.16390000e-01)
(-3.13130000e-01,  3.30110000e-01)
(-2.92930000e-01,  3.43960000e-01)
(-2.72730000e-01,  3.57920000e-01)
(-2.52530000e-01,  3.72000000e-01)
(-2.32320000e-01,  3.86180000e-01)
(-2.12120000e-01,  4.00470000e-01)
(-1.91920000e-01,  4.14860000e-01)
(-1.71720000e-01,  4.29350000e-01)
(-1.51520000e-01,  4.43940000e-01)
(-1.31310000e-01,  4.58620000e-01)
(-1.11110000e-01,  4.73390000e-01)
(-9.09090000e-02,  4.88240000e-01)
(-7.07070000e-02,  5.03180000e-01)
(-5.05050000e-02,  5.18200000e-01)
(-3.03030000e-02,  5.33310000e-01)
(-1.01010000e-02,  5.48490000e-01)
( 1.01010000e-02,  5.63750000e-01)
( 3.03030000e-02,  5.79080000e-01)
( 5.05050000e-02,  5.94480000e-01)
( 7.07070000e-02,  6.09950000e-01)
( 9.09090000e-02,  6.25500000e-01)
( 1.11110000e-01,  6.41110000e-01)
( 1.31310000e-01,  6.56780000e-01)
( 1.51520000e-01,  6.72520000e-01)
( 1.71720000e-01,  6.88320000e-01)
( 1.91920000e-01,  7.04180000e-01)
( 2.12120000e-01,  7.20090000e-01)
( 2.32320000e-01,  7.36070000e-01)
( 2.52530000e-01,  7.52110000e-01)
( 2.72730000e-01,  7.68200000e-01)
( 2.92930000e-01,  7.84340000e-01)
( 3.13130000e-01,  8.00540000e-01)
( 3.33330000e-01,  8.16790000e-01)
( 3.53540000e-01,  8.33090000e-01)
( 3.73740000e-01,  8.49440000e-01)
( 3.93940000e-01,  8.65850000e-01)
( 4.14140000e-01,  8.82270000e-01)
( 4.34340000e-01,  8.98760000e-01)
( 4.54550000e-01,  9.15320000e-01)
( 4.74750000e-01,  9.31950000e-01)
( 4.94950000e-01,  9.48630000e-01)
( 5.15150000e-01,  9.65320000e-01)
( 5.35350000e-01,  9.81980000e-01)
( 5.55560000e-01,  9.98560000e-01)
( 5.75760000e-01, -7.57880000e-05)
( 5.95960000e-01,  9.97890000e-05)
( 6.16160000e-01,  1.78790000e-04)
( 6.36360000e-01,  1.75040000e-04)
( 6.56570000e-01, -1.52460000e-04)
( 6.76770000e-01, -6.47840000e-05)
( 6.96970000e-01, -2.77040000e-05)
( 7.17170000e-01, -2.74430000e-05)
( 7.37370000e-01, -8.04660000e-06)
( 7.57580000e-01,  9.76040000e-06)
( 7.77780000e-01,  2.33240000e-05)
( 7.97980000e-01,  3.09890000e-05)
( 8.18180000e-01,  3.20360000e-05)
( 8.38380000e-01,  2.66370000e-05)
( 8.58590000e-01,  1.58290000e-05)
( 8.78790000e-01,  1.26830000e-06)
( 8.98990000e-01, -1.36460000e-05)
( 9.19190000e-01, -2.54030000e-05)
( 9.39390000e-01, -2.94660000e-05)
( 9.59600000e-01, -2.04870000e-05)
( 9.79800000e-01,  7.67460000e-06)
( 0.00000000e+00,  0.00000000e+00)};\label{line:burg1d_sptm_slice}

\end{axis}
\end{tikzpicture}

%% file: py/wedge0_geom.tikz
\begin{tikzpicture}
\begin{axis}[
axis equal image,
axis line style={gray},
axis x line*=bottom,
axis y line*=left,
xtick={0.0, 0.5, 1.0, 1.5},
ytick={0.0, 0.5, 1.0},
xlabel={$x_1$},
ymax=1.1,
xmax=1.6,
ylabel={$x_2$},
xmin=-0.1,
ymin=-0.1]
\addplot [opacity=0.6, fill=black!30!white, opacity=0.6, forget plot]
coordinates {
( 0.00000000e+00,  0.00000000e+00)
( 5.00000000e-01,  0.00000000e+00)
( 1.50000000e+00,  1.76326981e-01)
( 1.50000000e+00,  1.00000000e+00)
( 0.00000000e+00,  1.00000000e+00)
( 0.00000000e+00,  0.00000000e+00)};

\addplot [thick, color=blue]
coordinates {
( 0.00000000e+00,  0.00000000e+00)
( 0.00000000e+00,  1.00000000e+00)};\label{line:wedge0:supin}

\addplot [thick, color=black]
coordinates {
( 0.00000000e+00,  0.00000000e+00)
( 5.00000000e-01,  0.00000000e+00)};\label{line:wedge0:wall}

\addplot [thick, color=red]
coordinates {
( 1.50000000e+00,  1.76326981e-01)
( 1.50000000e+00,  1.00000000e+00)};\label{line:wedge0:supout}

\addplot [thick, color=black, forget plot]
coordinates {
( 1.50000000e+00,  1.00000000e+00)
( 0.00000000e+00,  1.00000000e+00)};

\addplot [thick, color=black, forget plot]
coordinates {
( 5.00000000e-01,  0.00000000e+00)
( 1.50000000e+00,  1.76326981e-01)};

\end{axis}
\end{tikzpicture}

%% file: tikz/euler2d0_wedge0_c1s1_nel6x4_porder0x1.tikz
\begin{tikzpicture}
\begin{groupplot}[
  group style={
      group size=2 by 2,
      horizontal sep=1.0cm,
      vertical sep=1.2cm
  },
  width=0.5\textwidth,
  axis equal image,
  xlabel={$x_1$},
  ylabel={$x_2$},
  xtick = {0, 0.5, 1.0, 1.5},
  ytick = {0.0, 0.5, 1.0},
  xmin=0, xmax=1.5,
  ymin=0, ymax=1
]

\nextgroupplot[title={Iteration 0}, xlabel={}, xtick=\empty]
\addplot graphics [xmin=0, xmax=1.5, ymin=0, ymax=1] {{img/euler2d0_wedge0_c1s1_nel6x4_porder0x1_iter000}.png};

\nextgroupplot[title={Iteration 4}, xlabel={}, ylabel={}, xtick=\empty, ytick=\empty]
\addplot graphics [xmin=0, xmax=1.5, ymin=0, ymax=1] {{img/euler2d0_wedge0_c1s1_nel6x4_porder0x1_iter004}.png};

\nextgroupplot[title={Iteration 8}]
\addplot graphics [xmin=0, xmax=1.5, ymin=0, ymax=1] {{img/euler2d0_wedge0_c1s1_nel6x4_porder0x1_iter008}.png};

\nextgroupplot[title={Iteration 20}, ylabel={}, ytick=\empty]
\addplot graphics [xmin=0, xmax=1.5, ymin=0, ymax=1] {{img/euler2d0_wedge0_c1s1_nel6x4_porder0x1_iter020}.png};

\end{groupplot}
\end{tikzpicture}

%% file: tikz/cbar_parula_1p65x2.tikz
\begin{tikzpicture}
\begin{axis}[
   hide axis, scale only axis,
   height=0pt, width=0pt,
   colormap={parula}{rgb255=(62,38,168) rgb255=(62,39,172) rgb255=(63,40,175) rgb255=(63,41,178) rgb255=(64,42,180) rgb255=(64,43,183) rgb255=(65,44,186) rgb255=(65,45,189) rgb255=(66,46,191) rgb255=(66,47,194) rgb255=(67,48,197) rgb255=(67,49,200) rgb255=(67,50,202) rgb255=(68,51,205) rgb255=(68,52,208) rgb255=(69,53,210) rgb255=(69,55,213) rgb255=(69,56,215) rgb255=(70,57,217) rgb255=(70,58,220) rgb255=(70,59,222) rgb255=(70,61,224) rgb255=(71,62,225) rgb255=(71,63,227) rgb255=(71,65,229) rgb255=(71,66,230) rgb255=(71,68,232) rgb255=(71,69,233) rgb255=(71,70,235) rgb255=(72,72,236) rgb255=(72,73,237) rgb255=(72,75,238) rgb255=(72,76,240) rgb255=(72,78,241) rgb255=(72,79,242) rgb255=(72,80,243) rgb255=(72,82,244) rgb255=(72,83,245) rgb255=(72,84,246) rgb255=(71,86,247) rgb255=(71,87,247) rgb255=(71,89,248) rgb255=(71,90,249) rgb255=(71,91,250) rgb255=(71,93,250) rgb255=(70,94,251) rgb255=(70,96,251) rgb255=(70,97,252) rgb255=(69,98,252) rgb255=(69,100,253) rgb255=(68,101,253) rgb255=(67,103,253) rgb255=(67,104,254) rgb255=(66,106,254) rgb255=(65,107,254) rgb255=(64,109,254) rgb255=(63,110,255) rgb255=(62,112,255) rgb255=(60,113,255) rgb255=(59,115,255) rgb255=(57,116,255) rgb255=(56,118,254) rgb255=(54,119,254) rgb255=(53,121,253) rgb255=(51,122,253) rgb255=(50,124,252) rgb255=(49,125,252) rgb255=(48,127,251) rgb255=(47,128,250) rgb255=(47,130,250) rgb255=(46,131,249) rgb255=(46,132,248) rgb255=(46,134,248) rgb255=(46,135,247) rgb255=(45,136,246) rgb255=(45,138,245) rgb255=(45,139,244) rgb255=(45,140,243) rgb255=(45,142,242) rgb255=(44,143,241) rgb255=(44,144,240) rgb255=(43,145,239) rgb255=(42,147,238) rgb255=(41,148,237) rgb255=(40,149,236) rgb255=(39,151,235) rgb255=(39,152,234) rgb255=(38,153,233) rgb255=(38,154,232) rgb255=(37,155,232) rgb255=(37,156,231) rgb255=(36,158,230) rgb255=(36,159,229) rgb255=(35,160,229) rgb255=(35,161,228) rgb255=(34,162,228) rgb255=(33,163,227) rgb255=(32,165,227) rgb255=(31,166,226) rgb255=(30,167,225) rgb255=(29,168,225) rgb255=(29,169,224) rgb255=(28,170,223) rgb255=(27,171,222) rgb255=(26,172,221) rgb255=(25,173,220) rgb255=(23,174,218) rgb255=(22,175,217) rgb255=(20,176,216) rgb255=(18,177,214) rgb255=(16,178,213) rgb255=(14,179,212) rgb255=(11,179,210) rgb255=(8,180,209) rgb255=(6,181,207) rgb255=(4,182,206) rgb255=(2,183,204) rgb255=(1,183,202) rgb255=(0,184,201) rgb255=(0,185,199) rgb255=(0,186,198) rgb255=(1,186,196) rgb255=(2,187,194) rgb255=(4,187,193) rgb255=(6,188,191) rgb255=(9,189,189) rgb255=(13,189,188) rgb255=(16,190,186) rgb255=(20,190,184) rgb255=(23,191,182) rgb255=(26,192,181) rgb255=(29,192,179) rgb255=(32,193,177) rgb255=(35,193,175) rgb255=(37,194,174) rgb255=(39,194,172) rgb255=(41,195,170) rgb255=(43,195,168) rgb255=(44,196,166) rgb255=(46,196,165) rgb255=(47,197,163) rgb255=(49,197,161) rgb255=(50,198,159) rgb255=(51,199,157) rgb255=(53,199,155) rgb255=(54,200,153) rgb255=(56,200,150) rgb255=(57,201,148) rgb255=(59,201,146) rgb255=(61,202,144) rgb255=(64,202,141) rgb255=(66,202,139) rgb255=(69,203,137) rgb255=(72,203,134) rgb255=(75,203,132) rgb255=(78,204,129) rgb255=(81,204,127) rgb255=(84,204,124) rgb255=(87,204,122) rgb255=(90,204,119) rgb255=(94,205,116) rgb255=(97,205,114) rgb255=(100,205,111) rgb255=(103,205,108) rgb255=(107,205,105) rgb255=(110,205,102) rgb255=(114,205,100) rgb255=(118,204,97) rgb255=(121,204,94) rgb255=(125,204,91) rgb255=(129,204,89) rgb255=(132,204,86) rgb255=(136,203,83) rgb255=(139,203,81) rgb255=(143,203,78) rgb255=(147,202,75) rgb255=(150,202,72) rgb255=(154,201,70) rgb255=(157,201,67) rgb255=(161,200,64) rgb255=(164,200,62) rgb255=(167,199,59) rgb255=(171,199,57) rgb255=(174,198,55) rgb255=(178,198,53) rgb255=(181,197,51) rgb255=(184,196,49) rgb255=(187,196,47) rgb255=(190,195,45) rgb255=(194,195,44) rgb255=(197,194,42) rgb255=(200,193,41) rgb255=(203,193,40) rgb255=(206,192,39) rgb255=(208,191,39) rgb255=(211,191,39) rgb255=(214,190,39) rgb255=(217,190,40) rgb255=(219,189,40) rgb255=(222,188,41) rgb255=(225,188,42) rgb255=(227,188,43) rgb255=(230,187,45) rgb255=(232,187,46) rgb255=(234,186,48) rgb255=(236,186,50) rgb255=(239,186,53) rgb255=(241,186,55) rgb255=(243,186,57) rgb255=(245,186,59) rgb255=(247,186,61) rgb255=(249,186,62) rgb255=(251,187,62) rgb255=(252,188,62) rgb255=(254,189,61) rgb255=(254,190,60) rgb255=(254,192,59) rgb255=(254,193,58) rgb255=(254,194,57) rgb255=(254,196,56) rgb255=(254,197,55) rgb255=(254,199,53) rgb255=(254,200,52) rgb255=(254,202,51) rgb255=(253,203,50) rgb255=(253,205,49) rgb255=(253,206,49) rgb255=(252,208,48) rgb255=(251,210,47) rgb255=(251,211,46) rgb255=(250,213,46) rgb255=(249,214,45) rgb255=(249,216,44) rgb255=(248,217,43) rgb255=(247,219,42) rgb255=(247,221,42) rgb255=(246,222,41) rgb255=(246,224,40) rgb255=(245,225,40) rgb255=(245,227,39) rgb255=(245,229,38) rgb255=(245,230,38) rgb255=(245,232,37) rgb255=(245,233,36) rgb255=(245,235,35) rgb255=(245,236,34) rgb255=(245,238,33) rgb255=(246,239,32) rgb255=(246,241,31) rgb255=(246,242,30) rgb255=(247,244,28) rgb255=(247,245,27) rgb255=(248,247,26) rgb255=(248,248,24) rgb255=(249,249,22) rgb255=(249,251,21) },
   colorbar horizontal,
   point meta min=1.65, point meta max=2,
   colorbar style={width=10cm, xtick={1.65, 1.74, 1.825, 1.91, 2}}
]
\addplot [draw=none] coordinates {(0,0)};
\end{axis}
\end{tikzpicture}

%% file: py/naca0012cfg1_geom.tikz
\begin{tikzpicture}
\begin{axis}[
axis equal image,
axis line style={gray},
axis x line*=bottom,
axis y line*=left,
xtick={-0.5, 0, 1, 1.5},
ytick={0.0, 1.0, 2.0, 3.0},
xlabel={$x_1$},
ymax=3.3,
xmax=1.8,
ylabel={$x_2$},
xmin=-0.8,
ymin=-0.30126]
\addplot [opacity=0.6, fill=black!30!white, opacity=0.6, forget plot]
coordinates {
(-5.00000000e-01,  0.00000000e+00)
( 0.00000000e+00,  0.00000000e+00)
( 1.01010101e-02,  1.71187581e-02)
( 2.02020202e-02,  2.37077433e-02)
( 3.03030303e-02,  2.85302608e-02)
( 4.04040404e-02,  3.24196415e-02)
( 5.05050505e-02,  3.56992766e-02)
( 6.06060606e-02,  3.85354937e-02)
( 7.07070707e-02,  4.10274840e-02)
( 8.08080808e-02,  4.32401994e-02)
( 9.09090909e-02,  4.52190323e-02)
( 1.01010101e-01,  4.69972732e-02)
( 1.11111111e-01,  4.86002652e-02)
( 1.21212121e-01,  5.00478855e-02)
( 1.31313131e-01,  5.13561092e-02)
( 1.41414141e-01,  5.25380403e-02)
( 1.51515152e-01,  5.36046150e-02)
( 1.61616162e-01,  5.45650974e-02)
( 1.71717172e-01,  5.54274373e-02)
( 1.81818182e-01,  5.61985338e-02)
( 1.91919192e-01,  5.68844349e-02)
( 2.02020202e-01,  5.74904895e-02)
( 2.12121212e-01,  5.80214656e-02)
( 2.22222222e-01,  5.84816438e-02)
( 2.32323232e-01,  5.88748916e-02)
( 2.42424242e-01,  5.92047234e-02)
( 2.52525253e-01,  5.94743493e-02)
( 2.62626263e-01,  5.96867153e-02)
( 2.72727273e-01,  5.98445362e-02)
( 2.82828283e-01,  5.99503238e-02)
( 2.92929293e-01,  6.00064098e-02)
( 3.03030303e-01,  6.00149658e-02)
( 3.13131313e-01,  5.99780201e-02)
( 3.23232323e-01,  5.98974716e-02)
( 3.33333333e-01,  5.97751029e-02)
( 3.43434343e-01,  5.96125903e-02)
( 3.53535354e-01,  5.94115137e-02)
( 3.63636364e-01,  5.91733641e-02)
( 3.73737374e-01,  5.88995511e-02)
( 3.83838384e-01,  5.85914089e-02)
( 3.93939394e-01,  5.82502020e-02)
( 4.04040404e-01,  5.78771301e-02)
( 4.14141414e-01,  5.74733319e-02)
( 4.24242424e-01,  5.70398898e-02)
( 4.34343434e-01,  5.65778326e-02)
( 4.44444444e-01,  5.60881390e-02)
( 4.54545455e-01,  5.55717402e-02)
( 4.64646465e-01,  5.50295223e-02)
( 4.74747475e-01,  5.44623291e-02)
( 4.84848485e-01,  5.38709632e-02)
( 4.94949495e-01,  5.32561888e-02)
( 5.05050505e-01,  5.26187328e-02)
( 5.15151515e-01,  5.19592864e-02)
( 5.25252525e-01,  5.12785067e-02)
( 5.35353535e-01,  5.05770177e-02)
( 5.45454545e-01,  4.98554118e-02)
( 5.55555556e-01,  4.91142506e-02)
( 5.65656566e-01,  4.83540659e-02)
( 5.75757576e-01,  4.75753607e-02)
( 5.85858586e-01,  4.67786102e-02)
( 5.95959596e-01,  4.59642620e-02)
( 6.06060606e-01,  4.51327375e-02)
( 6.16161616e-01,  4.42844321e-02)
( 6.26262626e-01,  4.34197159e-02)
( 6.36363636e-01,  4.25389345e-02)
( 6.46464646e-01,  4.16424090e-02)
( 6.56565657e-01,  4.07304373e-02)
( 6.66666667e-01,  3.98032935e-02)
( 6.76767677e-01,  3.88612294e-02)
( 6.86868687e-01,  3.79044742e-02)
( 6.96969697e-01,  3.69332350e-02)
( 7.07070707e-01,  3.59476972e-02)
( 7.17171717e-01,  3.49480250e-02)
( 7.27272727e-01,  3.39343615e-02)
( 7.37373737e-01,  3.29068287e-02)
( 7.47474747e-01,  3.18655285e-02)
( 7.57575758e-01,  3.08105424e-02)
( 7.67676768e-01,  2.97419317e-02)
( 7.77777778e-01,  2.86597380e-02)
( 7.87878788e-01,  2.75639834e-02)
( 7.97979798e-01,  2.64546705e-02)
( 8.08080808e-01,  2.53317827e-02)
( 8.18181818e-01,  2.41952842e-02)
( 8.28282828e-01,  2.30451206e-02)
( 8.38383838e-01,  2.18812185e-02)
( 8.48484848e-01,  2.07034861e-02)
( 8.58585859e-01,  1.95118132e-02)
( 8.68686869e-01,  1.83060711e-02)
( 8.78787879e-01,  1.70861130e-02)
( 8.88888889e-01,  1.58517740e-02)
( 8.98989899e-01,  1.46028715e-02)
( 9.09090909e-01,  1.33392048e-02)
( 9.19191919e-01,  1.20605554e-02)
( 9.29292929e-01,  1.07666875e-02)
( 9.39393939e-01,  9.45734727e-03)
( 9.49494949e-01,  8.13226386e-03)
( 9.59595960e-01,  6.79114881e-03)
( 9.69696970e-01,  5.43369641e-03)
( 9.79797980e-01,  4.05958371e-03)
( 9.89898990e-01,  2.66847063e-03)
( 1.00000000e+00, -1.26000000e-03)
( 1.50000000e+00,  0.00000000e+00)
( 1.50000000e+00,  3.00000000e+00)
(-5.00000000e-01,  3.00000000e+00)
(-5.00000000e-01,  0.00000000e+00)};

\addplot [thick, color=black]
coordinates {
(-5.00000000e-01,  0.00000000e+00)
( 0.00000000e+00,  0.00000000e+00)
( 1.01010101e-02,  1.71187581e-02)
( 2.02020202e-02,  2.37077433e-02)
( 3.03030303e-02,  2.85302608e-02)
( 4.04040404e-02,  3.24196415e-02)
( 5.05050505e-02,  3.56992766e-02)
( 6.06060606e-02,  3.85354937e-02)
( 7.07070707e-02,  4.10274840e-02)
( 8.08080808e-02,  4.32401994e-02)
( 9.09090909e-02,  4.52190323e-02)
( 1.01010101e-01,  4.69972732e-02)
( 1.11111111e-01,  4.86002652e-02)
( 1.21212121e-01,  5.00478855e-02)
( 1.31313131e-01,  5.13561092e-02)
( 1.41414141e-01,  5.25380403e-02)
( 1.51515152e-01,  5.36046150e-02)
( 1.61616162e-01,  5.45650974e-02)
( 1.71717172e-01,  5.54274373e-02)
( 1.81818182e-01,  5.61985338e-02)
( 1.91919192e-01,  5.68844349e-02)
( 2.02020202e-01,  5.74904895e-02)
( 2.12121212e-01,  5.80214656e-02)
( 2.22222222e-01,  5.84816438e-02)
( 2.32323232e-01,  5.88748916e-02)
( 2.42424242e-01,  5.92047234e-02)
( 2.52525253e-01,  5.94743493e-02)
( 2.62626263e-01,  5.96867153e-02)
( 2.72727273e-01,  5.98445362e-02)
( 2.82828283e-01,  5.99503238e-02)
( 2.92929293e-01,  6.00064098e-02)
( 3.03030303e-01,  6.00149658e-02)
( 3.13131313e-01,  5.99780201e-02)
( 3.23232323e-01,  5.98974716e-02)
( 3.33333333e-01,  5.97751029e-02)
( 3.43434343e-01,  5.96125903e-02)
( 3.53535354e-01,  5.94115137e-02)
( 3.63636364e-01,  5.91733641e-02)
( 3.73737374e-01,  5.88995511e-02)
( 3.83838384e-01,  5.85914089e-02)
( 3.93939394e-01,  5.82502020e-02)
( 4.04040404e-01,  5.78771301e-02)
( 4.14141414e-01,  5.74733319e-02)
( 4.24242424e-01,  5.70398898e-02)
( 4.34343434e-01,  5.65778326e-02)
( 4.44444444e-01,  5.60881390e-02)
( 4.54545455e-01,  5.55717402e-02)
( 4.64646465e-01,  5.50295223e-02)
( 4.74747475e-01,  5.44623291e-02)
( 4.84848485e-01,  5.38709632e-02)
( 4.94949495e-01,  5.32561888e-02)
( 5.05050505e-01,  5.26187328e-02)
( 5.15151515e-01,  5.19592864e-02)
( 5.25252525e-01,  5.12785067e-02)
( 5.35353535e-01,  5.05770177e-02)
( 5.45454545e-01,  4.98554118e-02)
( 5.55555556e-01,  4.91142506e-02)
( 5.65656566e-01,  4.83540659e-02)
( 5.75757576e-01,  4.75753607e-02)
( 5.85858586e-01,  4.67786102e-02)
( 5.95959596e-01,  4.59642620e-02)
( 6.06060606e-01,  4.51327375e-02)
( 6.16161616e-01,  4.42844321e-02)
( 6.26262626e-01,  4.34197159e-02)
( 6.36363636e-01,  4.25389345e-02)
( 6.46464646e-01,  4.16424090e-02)
( 6.56565657e-01,  4.07304373e-02)
( 6.66666667e-01,  3.98032935e-02)
( 6.76767677e-01,  3.88612294e-02)
( 6.86868687e-01,  3.79044742e-02)
( 6.96969697e-01,  3.69332350e-02)
( 7.07070707e-01,  3.59476972e-02)
( 7.17171717e-01,  3.49480250e-02)
( 7.27272727e-01,  3.39343615e-02)
( 7.37373737e-01,  3.29068287e-02)
( 7.47474747e-01,  3.18655285e-02)
( 7.57575758e-01,  3.08105424e-02)
( 7.67676768e-01,  2.97419317e-02)
( 7.77777778e-01,  2.86597380e-02)
( 7.87878788e-01,  2.75639834e-02)
( 7.97979798e-01,  2.64546705e-02)
( 8.08080808e-01,  2.53317827e-02)
( 8.18181818e-01,  2.41952842e-02)
( 8.28282828e-01,  2.30451206e-02)
( 8.38383838e-01,  2.18812185e-02)
( 8.48484848e-01,  2.07034861e-02)
( 8.58585859e-01,  1.95118132e-02)
( 8.68686869e-01,  1.83060711e-02)
( 8.78787879e-01,  1.70861130e-02)
( 8.88888889e-01,  1.58517740e-02)
( 8.98989899e-01,  1.46028715e-02)
( 9.09090909e-01,  1.33392048e-02)
( 9.19191919e-01,  1.20605554e-02)
( 9.29292929e-01,  1.07666875e-02)
( 9.39393939e-01,  9.45734727e-03)
( 9.49494949e-01,  8.13226386e-03)
( 9.59595960e-01,  6.79114881e-03)
( 9.69696970e-01,  5.43369641e-03)
( 9.79797980e-01,  4.05958371e-03)
( 9.89898990e-01,  2.66847063e-03)
( 1.00000000e+00, -1.26000000e-03)
( 1.50000000e+00,  0.00000000e+00)};\label{line:naca0012cfg1:wall}

\addplot [thick, color=blue]
coordinates {
( 1.50000000e+00,  0.00000000e+00)
( 1.50000000e+00,  3.00000000e+00)
(-5.00000000e-01,  3.00000000e+00)
(-5.00000000e-01,  0.00000000e+00)};\label{line:naca0012cfg1:farf}

\end{axis}
\end{tikzpicture}

%% file: tikz/euler2d0_naca0012cfg1_nref1_c1s1.tikz
\begin{tikzpicture}
\begin{groupplot}[
  group style={
      group size=3 by 2,
      horizontal sep=0.9cm,
      vertical sep=0.7cm
  },
  width=0.5\textwidth,
  axis equal image,
  xlabel={$x_1$},
  ylabel={$x_2$},
  xtick = {-0.5, 0, 0.5, 1.0, 1.5},
  ytick = {0.0, 1.0, 2.0, 3.0},
  xmin=-0.5, xmax=1.5,
  ymin=0, ymax=3
]

\nextgroupplot[title={Initialization}, xlabel={}, xtick=\empty]
\addplot graphics [xmin=-0.5, xmax=1.5, ymin=0, ymax=3] {{img/euler2d0_naca0012cfg1_nref1p1_c1s1_mesh_iter000}.png};

\nextgroupplot[title={Converged ($p=1$)}, xlabel={}, ylabel={}, xtick=\empty, ytick=\empty]
\addplot graphics [xmin=-0.5, xmax=1.5, ymin=0, ymax=3] {{img/euler2d0_naca0012cfg1_nref1p1_c1s1_mesh_iter100}.png};

\nextgroupplot[title={Converged ($p=2$)}, xlabel={}, ylabel={}, xtick=\empty, ytick=\empty]
\addplot graphics [xmin=-0.5, xmax=1.5, ymin=0, ymax=3] {{img/euler2d0_naca0012cfg1_nref1p2_c1s1_mesh_iter100}.png};

\nextgroupplot
\addplot graphics [xmin=-0.5, xmax=1.5, ymin=0, ymax=3] {{img/euler2d0_naca0012cfg1_nref1p1_c1s1_iter000}.png};

\nextgroupplot[ylabel={}, ytick=\empty]
\addplot graphics [xmin=-0.5, xmax=1.5, ymin=0, ymax=3] {{img/euler2d0_naca0012cfg1_nref1p1_c1s1_iter100}.png};

\nextgroupplot[ylabel={}, ytick=\empty]
\addplot graphics [xmin=-0.5, xmax=1.5, ymin=0, ymax=3] {{img/euler2d0_naca0012cfg1_nref1p2_c1s1_iter100}.png};

\end{groupplot}
\end{tikzpicture}

%% file: tikz/cbar_parula_0x1p75.tikz
\begin{tikzpicture}
\begin{axis}[
   hide axis, scale only axis,
   height=0pt, width=0pt,
   colormap={parula}{rgb255=(62,38,168) rgb255=(62,39,172) rgb255=(63,40,175) rgb255=(63,41,178) rgb255=(64,42,180) rgb255=(64,43,183) rgb255=(65,44,186) rgb255=(65,45,189) rgb255=(66,46,191) rgb255=(66,47,194) rgb255=(67,48,197) rgb255=(67,49,200) rgb255=(67,50,202) rgb255=(68,51,205) rgb255=(68,52,208) rgb255=(69,53,210) rgb255=(69,55,213) rgb255=(69,56,215) rgb255=(70,57,217) rgb255=(70,58,220) rgb255=(70,59,222) rgb255=(70,61,224) rgb255=(71,62,225) rgb255=(71,63,227) rgb255=(71,65,229) rgb255=(71,66,230) rgb255=(71,68,232) rgb255=(71,69,233) rgb255=(71,70,235) rgb255=(72,72,236) rgb255=(72,73,237) rgb255=(72,75,238) rgb255=(72,76,240) rgb255=(72,78,241) rgb255=(72,79,242) rgb255=(72,80,243) rgb255=(72,82,244) rgb255=(72,83,245) rgb255=(72,84,246) rgb255=(71,86,247) rgb255=(71,87,247) rgb255=(71,89,248) rgb255=(71,90,249) rgb255=(71,91,250) rgb255=(71,93,250) rgb255=(70,94,251) rgb255=(70,96,251) rgb255=(70,97,252) rgb255=(69,98,252) rgb255=(69,100,253) rgb255=(68,101,253) rgb255=(67,103,253) rgb255=(67,104,254) rgb255=(66,106,254) rgb255=(65,107,254) rgb255=(64,109,254) rgb255=(63,110,255) rgb255=(62,112,255) rgb255=(60,113,255) rgb255=(59,115,255) rgb255=(57,116,255) rgb255=(56,118,254) rgb255=(54,119,254) rgb255=(53,121,253) rgb255=(51,122,253) rgb255=(50,124,252) rgb255=(49,125,252) rgb255=(48,127,251) rgb255=(47,128,250) rgb255=(47,130,250) rgb255=(46,131,249) rgb255=(46,132,248) rgb255=(46,134,248) rgb255=(46,135,247) rgb255=(45,136,246) rgb255=(45,138,245) rgb255=(45,139,244) rgb255=(45,140,243) rgb255=(45,142,242) rgb255=(44,143,241) rgb255=(44,144,240) rgb255=(43,145,239) rgb255=(42,147,238) rgb255=(41,148,237) rgb255=(40,149,236) rgb255=(39,151,235) rgb255=(39,152,234) rgb255=(38,153,233) rgb255=(38,154,232) rgb255=(37,155,232) rgb255=(37,156,231) rgb255=(36,158,230) rgb255=(36,159,229) rgb255=(35,160,229) rgb255=(35,161,228) rgb255=(34,162,228) rgb255=(33,163,227) rgb255=(32,165,227) rgb255=(31,166,226) rgb255=(30,167,225) rgb255=(29,168,225) rgb255=(29,169,224) rgb255=(28,170,223) rgb255=(27,171,222) rgb255=(26,172,221) rgb255=(25,173,220) rgb255=(23,174,218) rgb255=(22,175,217) rgb255=(20,176,216) rgb255=(18,177,214) rgb255=(16,178,213) rgb255=(14,179,212) rgb255=(11,179,210) rgb255=(8,180,209) rgb255=(6,181,207) rgb255=(4,182,206) rgb255=(2,183,204) rgb255=(1,183,202) rgb255=(0,184,201) rgb255=(0,185,199) rgb255=(0,186,198) rgb255=(1,186,196) rgb255=(2,187,194) rgb255=(4,187,193) rgb255=(6,188,191) rgb255=(9,189,189) rgb255=(13,189,188) rgb255=(16,190,186) rgb255=(20,190,184) rgb255=(23,191,182) rgb255=(26,192,181) rgb255=(29,192,179) rgb255=(32,193,177) rgb255=(35,193,175) rgb255=(37,194,174) rgb255=(39,194,172) rgb255=(41,195,170) rgb255=(43,195,168) rgb255=(44,196,166) rgb255=(46,196,165) rgb255=(47,197,163) rgb255=(49,197,161) rgb255=(50,198,159) rgb255=(51,199,157) rgb255=(53,199,155) rgb255=(54,200,153) rgb255=(56,200,150) rgb255=(57,201,148) rgb255=(59,201,146) rgb255=(61,202,144) rgb255=(64,202,141) rgb255=(66,202,139) rgb255=(69,203,137) rgb255=(72,203,134) rgb255=(75,203,132) rgb255=(78,204,129) rgb255=(81,204,127) rgb255=(84,204,124) rgb255=(87,204,122) rgb255=(90,204,119) rgb255=(94,205,116) rgb255=(97,205,114) rgb255=(100,205,111) rgb255=(103,205,108) rgb255=(107,205,105) rgb255=(110,205,102) rgb255=(114,205,100) rgb255=(118,204,97) rgb255=(121,204,94) rgb255=(125,204,91) rgb255=(129,204,89) rgb255=(132,204,86) rgb255=(136,203,83) rgb255=(139,203,81) rgb255=(143,203,78) rgb255=(147,202,75) rgb255=(150,202,72) rgb255=(154,201,70) rgb255=(157,201,67) rgb255=(161,200,64) rgb255=(164,200,62) rgb255=(167,199,59) rgb255=(171,199,57) rgb255=(174,198,55) rgb255=(178,198,53) rgb255=(181,197,51) rgb255=(184,196,49) rgb255=(187,196,47) rgb255=(190,195,45) rgb255=(194,195,44) rgb255=(197,194,42) rgb255=(200,193,41) rgb255=(203,193,40) rgb255=(206,192,39) rgb255=(208,191,39) rgb255=(211,191,39) rgb255=(214,190,39) rgb255=(217,190,40) rgb255=(219,189,40) rgb255=(222,188,41) rgb255=(225,188,42) rgb255=(227,188,43) rgb255=(230,187,45) rgb255=(232,187,46) rgb255=(234,186,48) rgb255=(236,186,50) rgb255=(239,186,53) rgb255=(241,186,55) rgb255=(243,186,57) rgb255=(245,186,59) rgb255=(247,186,61) rgb255=(249,186,62) rgb255=(251,187,62) rgb255=(252,188,62) rgb255=(254,189,61) rgb255=(254,190,60) rgb255=(254,192,59) rgb255=(254,193,58) rgb255=(254,194,57) rgb255=(254,196,56) rgb255=(254,197,55) rgb255=(254,199,53) rgb255=(254,200,52) rgb255=(254,202,51) rgb255=(253,203,50) rgb255=(253,205,49) rgb255=(253,206,49) rgb255=(252,208,48) rgb255=(251,210,47) rgb255=(251,211,46) rgb255=(250,213,46) rgb255=(249,214,45) rgb255=(249,216,44) rgb255=(248,217,43) rgb255=(247,219,42) rgb255=(247,221,42) rgb255=(246,222,41) rgb255=(246,224,40) rgb255=(245,225,40) rgb255=(245,227,39) rgb255=(245,229,38) rgb255=(245,230,38) rgb255=(245,232,37) rgb255=(245,233,36) rgb255=(245,235,35) rgb255=(245,236,34) rgb255=(245,238,33) rgb255=(246,239,32) rgb255=(246,241,31) rgb255=(246,242,30) rgb255=(247,244,28) rgb255=(247,245,27) rgb255=(248,247,26) rgb255=(248,248,24) rgb255=(249,249,22) rgb255=(249,251,21) },
   colorbar horizontal,
   point meta min=0, point meta max=1.75,
   colorbar style={width=10cm, xtick={0, 0.4375, 0.875, 1.3125, 1.75}}
]
\addplot [draw=none] coordinates {(0,0)};
\end{axis}
\end{tikzpicture}

%% file: py/euler2d0_naca0012cfg1_nref1p3_c1s1_slice.tikz
\begin{tikzpicture}
\begin{groupplot} [
group style={group size = 3 by 1, horizontal sep = 1.5cm}]
\nextgroupplot[width=0.3\textwidth, xlabel={$x$}, ymax=2.6, xmax=1.5, ylabel={$\rho$}, xmin=-0.5, ymin=0.9, height=0.3\textwidth]
\addplot [solid, thick, color=black, forget plot]
coordinates {
(-5.00000000e-01,  1.40000000e+00)
(-4.93310000e-01,  1.40000000e+00)
(-4.86620000e-01,  1.40000000e+00)
(-4.79930000e-01,  1.40000000e+00)
(-4.73240000e-01,  1.40000000e+00)
(-4.66560000e-01,  1.40000000e+00)
(-4.59870000e-01,  1.40000000e+00)
(-4.53180000e-01,  1.40000000e+00)
(-4.46490000e-01,  1.40000000e+00)
(-4.39800000e-01,  1.40000000e+00)
(-4.33110000e-01,  1.40000000e+00)
(-4.26420000e-01,  1.40000000e+00)
(-4.19730000e-01,  1.40000000e+00)
(-4.13040000e-01,  1.40000000e+00)
(-4.06350000e-01,  1.40000000e+00)
(-3.99670000e-01,  1.40000000e+00)
(-3.92980000e-01,  1.40000000e+00)
(-3.86290000e-01,  1.40000000e+00)
(-3.79600000e-01,  1.40000000e+00)
(-3.72910000e-01,  1.40000000e+00)
(-3.66220000e-01,  1.40000000e+00)
(-3.59530000e-01,  1.40000000e+00)
(-3.52840000e-01,  1.40000000e+00)
(-3.46150000e-01,  1.40000000e+00)
(-3.39460000e-01,  1.40000000e+00)
(-3.32780000e-01,  1.40000000e+00)
(-3.26090000e-01,  1.40000000e+00)
(-3.19400000e-01,  1.40000000e+00)
(-3.12710000e-01,  1.40000000e+00)
(-3.06020000e-01,  1.40000000e+00)
(-2.99330000e-01,  1.40000000e+00)
(-2.92640000e-01,  1.40000000e+00)
(-2.85950000e-01,  1.40000000e+00)
(-2.79260000e-01,  1.40000000e+00)
(-2.72580000e-01,  1.40000000e+00)
(-2.65890000e-01,  1.40000000e+00)
(-2.59200000e-01,  1.40000000e+00)
(-2.52510000e-01,  1.40000000e+00)
(-2.45820000e-01,  1.40000000e+00)
(-2.39130000e-01,  1.40000000e+00)
(-2.32440000e-01,  1.40000000e+00)
(-2.25750000e-01,  1.40000000e+00)
(-2.19060000e-01,  1.40000000e+00)
(-2.12370000e-01,  1.40000000e+00)
(-2.05690000e-01,  1.40000000e+00)
(-1.99000000e-01,  1.40000000e+00)
(-1.92310000e-01,  1.40000000e+00)
(-1.85620000e-01,  1.40000000e+00)
(-1.78930000e-01,  1.40000000e+00)
(-1.72240000e-01,  1.40000000e+00)
(-1.65550000e-01,  1.40000000e+00)
(-1.58860000e-01,  1.40000000e+00)
(-1.52170000e-01,  1.40000000e+00)
(-1.45480000e-01,  1.40000000e+00)
(-1.38800000e-01,  1.40000000e+00)
(-1.32110000e-01,  1.40000000e+00)
(-1.25420000e-01,  1.40000000e+00)
(-1.18730000e-01,  1.40000000e+00)
(-1.12040000e-01,  1.40000000e+00)
(-1.05350000e-01,  1.40000000e+00)
(-9.86620000e-02,  1.40000000e+00)
(-9.19730000e-02,  1.40000000e+00)
(-8.52840000e-02,  1.40000000e+00)
(-7.85950000e-02,  1.40000000e+00)
(-7.19060000e-02,  1.40010000e+00)
(-6.52170000e-02,  2.48530000e+00)
(-5.85280000e-02,  2.48730000e+00)
(-5.18390000e-02,  2.48700000e+00)
(-4.51510000e-02,  2.48430000e+00)
(-3.84620000e-02,  2.47910000e+00)
(-3.17730000e-02,  2.47130000e+00)
(-2.50840000e-02,  2.46070000e+00)
(-1.83950000e-02,  2.44730000e+00)
(-1.17060000e-02,  2.43140000e+00)
(-5.01670000e-03,  2.41270000e+00)
( 1.67220000e-03,  2.39150000e+00)
( 8.36120000e-03,  2.36780000e+00)
( 1.50500000e-02,  2.34180000e+00)
( 2.17390000e-02,  2.31400000e+00)
( 2.84280000e-02,  2.28430000e+00)
( 3.51170000e-02,  2.25390000e+00)
( 4.18060000e-02,  2.22240000e+00)
( 4.84950000e-02,  2.19030000e+00)
( 5.51840000e-02,  2.15770000e+00)
( 6.18730000e-02,  2.12480000e+00)
( 6.85620000e-02,  2.09240000e+00)
( 7.52510000e-02,  2.06010000e+00)
( 8.19400000e-02,  2.02870000e+00)
( 8.86290000e-02,  1.99810000e+00)
( 9.53180000e-02,  1.96840000e+00)
( 1.02010000e-01,  1.93980000e+00)
( 1.08700000e-01,  1.91220000e+00)
( 1.15380000e-01,  1.88570000e+00)
( 1.22070000e-01,  1.86010000e+00)
( 1.28760000e-01,  1.83550000e+00)
( 1.35450000e-01,  1.81170000e+00)
( 1.42140000e-01,  1.78920000e+00)
( 1.48830000e-01,  1.76780000e+00)
( 1.55520000e-01,  1.74720000e+00)
( 1.62210000e-01,  1.72710000e+00)
( 1.68900000e-01,  1.70780000e+00)
( 1.75590000e-01,  1.68920000e+00)
( 1.82270000e-01,  1.67130000e+00)
( 1.88960000e-01,  1.65420000e+00)
( 1.95650000e-01,  1.63800000e+00)
( 2.02340000e-01,  1.62280000e+00)
( 2.09030000e-01,  1.60800000e+00)
( 2.15720000e-01,  1.59350000e+00)
( 2.22410000e-01,  1.57960000e+00)
( 2.29100000e-01,  1.56620000e+00)
( 2.35790000e-01,  1.55320000e+00)
( 2.42470000e-01,  1.54080000e+00)
( 2.49160000e-01,  1.52880000e+00)
( 2.55850000e-01,  1.51710000e+00)
( 2.62540000e-01,  1.50590000e+00)
( 2.69230000e-01,  1.49500000e+00)
( 2.75920000e-01,  1.48440000e+00)
( 2.82610000e-01,  1.47420000e+00)
( 2.89300000e-01,  1.46440000e+00)
( 2.95990000e-01,  1.45480000e+00)
( 3.02680000e-01,  1.44570000e+00)
( 3.09360000e-01,  1.43680000e+00)
( 3.16050000e-01,  1.42830000e+00)
( 3.22740000e-01,  1.42010000e+00)
( 3.29430000e-01,  1.41160000e+00)
( 3.36120000e-01,  1.40390000e+00)
( 3.42810000e-01,  1.39630000e+00)
( 3.49500000e-01,  1.38890000e+00)
( 3.56190000e-01,  1.38160000e+00)
( 3.62880000e-01,  1.37460000e+00)
( 3.69570000e-01,  1.36780000e+00)
( 3.76250000e-01,  1.36110000e+00)
( 3.82940000e-01,  1.35460000e+00)
( 3.89630000e-01,  1.34820000e+00)
( 3.96320000e-01,  1.34210000e+00)
( 4.03010000e-01,  1.33600000e+00)
( 4.09700000e-01,  1.33030000e+00)
( 4.16390000e-01,  1.32460000e+00)
( 4.23080000e-01,  1.31900000e+00)
( 4.29770000e-01,  1.31350000e+00)
( 4.36450000e-01,  1.30820000e+00)
( 4.43140000e-01,  1.30300000e+00)
( 4.49830000e-01,  1.29800000e+00)
( 4.56520000e-01,  1.29300000e+00)
( 4.63210000e-01,  1.28820000e+00)
( 4.69900000e-01,  1.28360000e+00)
( 4.76590000e-01,  1.27900000e+00)
( 4.83280000e-01,  1.27460000e+00)
( 4.89970000e-01,  1.27030000e+00)
( 4.96660000e-01,  1.26610000e+00)
( 5.03340000e-01,  1.26210000e+00)
( 5.10030000e-01,  1.25820000e+00)
( 5.16720000e-01,  1.25440000e+00)
( 5.23410000e-01,  1.24960000e+00)
( 5.30100000e-01,  1.24600000e+00)
( 5.36790000e-01,  1.24240000e+00)
( 5.43480000e-01,  1.23890000e+00)
( 5.50170000e-01,  1.23550000e+00)
( 5.56860000e-01,  1.23210000e+00)
( 5.63550000e-01,  1.22870000e+00)
( 5.70230000e-01,  1.22550000e+00)
( 5.76920000e-01,  1.22220000e+00)
( 5.83610000e-01,  1.21910000e+00)
( 5.90300000e-01,  1.21590000e+00)
( 5.96990000e-01,  1.21290000e+00)
( 6.03680000e-01,  1.20980000e+00)
( 6.10370000e-01,  1.20690000e+00)
( 6.17060000e-01,  1.20400000e+00)
( 6.23750000e-01,  1.20110000e+00)
( 6.30430000e-01,  1.19830000e+00)
( 6.37120000e-01,  1.19550000e+00)
( 6.43810000e-01,  1.19280000e+00)
( 6.50500000e-01,  1.19010000e+00)
( 6.57190000e-01,  1.18750000e+00)
( 6.63880000e-01,  1.18500000e+00)
( 6.70570000e-01,  1.18240000e+00)
( 6.77260000e-01,  1.18000000e+00)
( 6.83950000e-01,  1.17750000e+00)
( 6.90640000e-01,  1.17510000e+00)
( 6.97320000e-01,  1.17270000e+00)
( 7.04010000e-01,  1.17040000e+00)
( 7.10700000e-01,  1.16810000e+00)
( 7.17390000e-01,  1.16590000e+00)
( 7.24080000e-01,  1.16370000e+00)
( 7.30770000e-01,  1.16160000e+00)
( 7.37460000e-01,  1.15950000e+00)
( 7.44150000e-01,  1.15740000e+00)
( 7.50840000e-01,  1.15540000e+00)
( 7.57530000e-01,  1.15330000e+00)
( 7.64210000e-01,  1.15130000e+00)
( 7.70900000e-01,  1.14930000e+00)
( 7.77590000e-01,  1.14740000e+00)
( 7.84280000e-01,  1.14540000e+00)
( 7.90970000e-01,  1.14350000e+00)
( 7.97660000e-01,  1.14160000e+00)
( 8.04350000e-01,  1.13970000e+00)
( 8.11040000e-01,  1.13780000e+00)
( 8.17730000e-01,  1.13600000e+00)
( 8.24410000e-01,  1.13420000e+00)
( 8.31100000e-01,  1.13240000e+00)
( 8.37790000e-01,  1.13060000e+00)
( 8.44480000e-01,  1.12880000e+00)
( 8.51170000e-01,  1.12700000e+00)
( 8.57860000e-01,  1.12530000e+00)
( 8.64550000e-01,  1.12350000e+00)
( 8.71240000e-01,  1.12180000e+00)
( 8.77930000e-01,  1.12010000e+00)
( 8.84620000e-01,  1.11840000e+00)
( 8.91300000e-01,  1.11680000e+00)
( 8.97990000e-01,  1.11510000e+00)
( 9.04680000e-01,  1.11340000e+00)
( 9.11370000e-01,  1.11180000e+00)
( 9.18060000e-01,  1.11020000e+00)
( 9.24750000e-01,  1.10860000e+00)
( 9.31440000e-01,  1.10690000e+00)
( 9.38130000e-01,  1.10530000e+00)
( 9.44820000e-01,  1.10370000e+00)
( 9.51510000e-01,  1.10210000e+00)
( 9.58190000e-01,  1.10050000e+00)
( 9.64880000e-01,  1.09900000e+00)
( 9.71570000e-01,  1.09740000e+00)
( 9.78260000e-01,  1.09580000e+00)
( 9.84950000e-01,  1.09420000e+00)
( 9.91640000e-01,  1.09260000e+00)
( 9.98330000e-01,  1.09110000e+00)
( 1.00500000e+00,  1.08950000e+00)
( 1.01170000e+00,  1.08790000e+00)
( 1.01840000e+00,  1.08640000e+00)
( 1.02510000e+00,  1.08480000e+00)
( 1.03180000e+00,  1.08320000e+00)
( 1.03850000e+00,  1.08170000e+00)
( 1.04520000e+00,  1.08010000e+00)
( 1.05180000e+00,  1.07850000e+00)
( 1.05850000e+00,  1.07700000e+00)
( 1.06520000e+00,  1.07540000e+00)
( 1.07190000e+00,  1.07360000e+00)
( 1.07860000e+00,  1.07210000e+00)
( 1.08530000e+00,  1.07070000e+00)
( 1.09200000e+00,  1.06920000e+00)
( 1.09870000e+00,  1.06760000e+00)
( 1.10540000e+00,  1.06600000e+00)
( 1.11200000e+00,  1.06440000e+00)
( 1.11870000e+00,  1.06280000e+00)
( 1.12540000e+00,  1.06110000e+00)
( 1.13210000e+00,  1.05950000e+00)
( 1.13880000e+00,  1.05780000e+00)
( 1.14550000e+00,  1.05620000e+00)
( 1.15220000e+00,  1.05450000e+00)
( 1.15890000e+00,  1.05290000e+00)
( 1.16560000e+00,  1.05110000e+00)
( 1.17220000e+00,  1.04950000e+00)
( 1.17890000e+00,  1.04800000e+00)
( 1.18560000e+00,  1.04670000e+00)
( 1.19230000e+00,  1.04570000e+00)
( 1.19900000e+00,  1.37480000e+00)
( 1.20570000e+00,  1.37480000e+00)
( 1.21240000e+00,  1.37490000e+00)
( 1.21910000e+00,  1.37490000e+00)
( 1.22580000e+00,  1.37490000e+00)
( 1.23240000e+00,  1.37490000e+00)
( 1.23910000e+00,  1.37490000e+00)
( 1.24580000e+00,  1.37490000e+00)
( 1.25250000e+00,  1.37490000e+00)
( 1.25920000e+00,  1.37500000e+00)
( 1.26590000e+00,  1.37500000e+00)
( 1.27260000e+00,  1.37500000e+00)
( 1.27930000e+00,  1.37500000e+00)
( 1.28600000e+00,  1.37500000e+00)
( 1.29260000e+00,  1.37510000e+00)
( 1.29930000e+00,  1.37510000e+00)
( 1.30600000e+00,  1.37510000e+00)
( 1.31270000e+00,  1.37510000e+00)
( 1.31940000e+00,  1.37520000e+00)
( 1.32610000e+00,  1.37520000e+00)
( 1.33280000e+00,  1.37520000e+00)
( 1.33950000e+00,  1.37520000e+00)
( 1.34620000e+00,  1.37530000e+00)
( 1.35280000e+00,  1.37530000e+00)
( 1.35950000e+00,  1.37530000e+00)
( 1.36620000e+00,  1.37520000e+00)
( 1.37290000e+00,  1.37520000e+00)
( 1.37960000e+00,  1.37530000e+00)
( 1.38630000e+00,  1.37530000e+00)
( 1.39300000e+00,  1.37530000e+00)
( 1.39970000e+00,  1.37530000e+00)
( 1.40640000e+00,  1.37530000e+00)
( 1.41300000e+00,  1.37530000e+00)
( 1.41970000e+00,  1.37540000e+00)
( 1.42640000e+00,  1.37540000e+00)
( 1.43310000e+00,  1.37540000e+00)
( 1.43980000e+00,  1.37540000e+00)
( 1.44650000e+00,  1.37540000e+00)
( 1.45320000e+00,  1.37540000e+00)
( 1.45990000e+00,  1.37540000e+00)
( 1.46660000e+00,  1.37540000e+00)
( 1.47320000e+00,  1.37540000e+00)
( 1.47990000e+00,  1.37540000e+00)
( 1.48660000e+00,  1.37540000e+00)
( 1.49330000e+00,  1.37540000e+00)};

\nextgroupplot[width=0.3\textwidth, xlabel={$x$}, ymax=1.75, xmax=1.5, ylabel={$M$}, xmin=-0.5, ymin=0.75, height=0.3\textwidth]
\addplot [solid, thick, color=black, forget plot]
coordinates {
(-5.00000000e-01,  1.50000000e+00)
(-4.93310000e-01,  1.50000000e+00)
(-4.86620000e-01,  1.50000000e+00)
(-4.79930000e-01,  1.50000000e+00)
(-4.73240000e-01,  1.50000000e+00)
(-4.66560000e-01,  1.50000000e+00)
(-4.59870000e-01,  1.50000000e+00)
(-4.53180000e-01,  1.50000000e+00)
(-4.46490000e-01,  1.50000000e+00)
(-4.39800000e-01,  1.50000000e+00)
(-4.33110000e-01,  1.50000000e+00)
(-4.26420000e-01,  1.50000000e+00)
(-4.19730000e-01,  1.50000000e+00)
(-4.13040000e-01,  1.50000000e+00)
(-4.06350000e-01,  1.50000000e+00)
(-3.99670000e-01,  1.50000000e+00)
(-3.92980000e-01,  1.50000000e+00)
(-3.86290000e-01,  1.50000000e+00)
(-3.79600000e-01,  1.50000000e+00)
(-3.72910000e-01,  1.50000000e+00)
(-3.66220000e-01,  1.50000000e+00)
(-3.59530000e-01,  1.50000000e+00)
(-3.52840000e-01,  1.50000000e+00)
(-3.46150000e-01,  1.50000000e+00)
(-3.39460000e-01,  1.50000000e+00)
(-3.32780000e-01,  1.50000000e+00)
(-3.26090000e-01,  1.50000000e+00)
(-3.19400000e-01,  1.50000000e+00)
(-3.12710000e-01,  1.50000000e+00)
(-3.06020000e-01,  1.50000000e+00)
(-2.99330000e-01,  1.50000000e+00)
(-2.92640000e-01,  1.50000000e+00)
(-2.85950000e-01,  1.50000000e+00)
(-2.79260000e-01,  1.50000000e+00)
(-2.72580000e-01,  1.50000000e+00)
(-2.65890000e-01,  1.50000000e+00)
(-2.59200000e-01,  1.50000000e+00)
(-2.52510000e-01,  1.50000000e+00)
(-2.45820000e-01,  1.50000000e+00)
(-2.39130000e-01,  1.50000000e+00)
(-2.32440000e-01,  1.50000000e+00)
(-2.25750000e-01,  1.50000000e+00)
(-2.19060000e-01,  1.50000000e+00)
(-2.12370000e-01,  1.50000000e+00)
(-2.05690000e-01,  1.50000000e+00)
(-1.99000000e-01,  1.50000000e+00)
(-1.92310000e-01,  1.50000000e+00)
(-1.85620000e-01,  1.50000000e+00)
(-1.78930000e-01,  1.50000000e+00)
(-1.72240000e-01,  1.50000000e+00)
(-1.65550000e-01,  1.50000000e+00)
(-1.58860000e-01,  1.50000000e+00)
(-1.52170000e-01,  1.50000000e+00)
(-1.45480000e-01,  1.50000000e+00)
(-1.38800000e-01,  1.50000000e+00)
(-1.32110000e-01,  1.50000000e+00)
(-1.25420000e-01,  1.50000000e+00)
(-1.18730000e-01,  1.50000000e+00)
(-1.12040000e-01,  1.50000000e+00)
(-1.05350000e-01,  1.50000000e+00)
(-9.86620000e-02,  1.50000000e+00)
(-9.19730000e-02,  1.50000000e+00)
(-8.52840000e-02,  1.50000000e+00)
(-7.85950000e-02,  1.49997000e+00)
(-7.19060000e-02,  1.49988269e+00)
(-6.52170000e-02,  7.96309594e-01)
(-5.85280000e-02,  7.94893158e-01)
(-5.18390000e-02,  7.94772972e-01)
(-4.51510000e-02,  7.95968648e-01)
(-3.84620000e-02,  7.98598469e-01)
(-3.17730000e-02,  8.02736448e-01)
(-2.50840000e-02,  8.08375554e-01)
(-1.83950000e-02,  8.15622808e-01)
(-1.17060000e-02,  8.24292581e-01)
(-5.01670000e-03,  8.34485418e-01)
( 1.67220000e-03,  8.46081665e-01)
( 8.36120000e-03,  8.59046240e-01)
( 1.50500000e-02,  8.73330017e-01)
( 2.17390000e-02,  8.88646945e-01)
( 2.84280000e-02,  9.04938131e-01)
( 3.51170000e-02,  9.21750421e-01)
( 4.18060000e-02,  9.39101277e-01)
( 4.84950000e-02,  9.56873031e-01)
( 5.51840000e-02,  9.74940201e-01)
( 6.18730000e-02,  9.93176592e-01)
( 6.85620000e-02,  1.01133251e+00)
( 7.52510000e-02,  1.02939165e+00)
( 8.19400000e-02,  1.04714242e+00)
( 8.86290000e-02,  1.06445218e+00)
( 9.53180000e-02,  1.08133166e+00)
( 1.02010000e-01,  1.09772854e+00)
( 1.08700000e-01,  1.11360127e+00)
( 1.15380000e-01,  1.12893423e+00)
( 1.22070000e-01,  1.14387094e+00)
( 1.28760000e-01,  1.15828106e+00)
( 1.35450000e-01,  1.17231884e+00)
( 1.42140000e-01,  1.18572846e+00)
( 1.48830000e-01,  1.19846795e+00)
( 1.55520000e-01,  1.21082763e+00)
( 1.62210000e-01,  1.22297067e+00)
( 1.68900000e-01,  1.23483081e+00)
( 1.75590000e-01,  1.24620985e+00)
( 1.82270000e-01,  1.25730718e+00)
( 1.88960000e-01,  1.26793152e+00)
( 1.95650000e-01,  1.27800894e+00)
( 2.02340000e-01,  1.28756948e+00)
( 2.09030000e-01,  1.29692778e+00)
( 2.15720000e-01,  1.30621312e+00)
( 2.22410000e-01,  1.31508146e+00)
( 2.29100000e-01,  1.32369837e+00)
( 2.35790000e-01,  1.33214776e+00)
( 2.42470000e-01,  1.34019708e+00)
( 2.49160000e-01,  1.34806847e+00)
( 2.55850000e-01,  1.35574369e+00)
( 2.62540000e-01,  1.36314246e+00)
( 2.69230000e-01,  1.37033195e+00)
( 2.75920000e-01,  1.37744848e+00)
( 2.82610000e-01,  1.38435695e+00)
( 2.89300000e-01,  1.39091301e+00)
( 2.95990000e-01,  1.39747851e+00)
( 3.02680000e-01,  1.40368989e+00)
( 3.09360000e-01,  1.40978246e+00)
( 3.16050000e-01,  1.41556475e+00)
( 3.22740000e-01,  1.42123380e+00)
( 3.29430000e-01,  1.42716462e+00)
( 3.36120000e-01,  1.43249587e+00)
( 3.42810000e-01,  1.43781831e+00)
( 3.49500000e-01,  1.44307438e+00)
( 3.56190000e-01,  1.44819918e+00)
( 3.62880000e-01,  1.45324698e+00)
( 3.69570000e-01,  1.45809445e+00)
( 3.76250000e-01,  1.46290589e+00)
( 3.82940000e-01,  1.46750297e+00)
( 3.89630000e-01,  1.47217171e+00)
( 3.96320000e-01,  1.47655605e+00)
( 4.03010000e-01,  1.48097299e+00)
( 4.09700000e-01,  1.48524165e+00)
( 4.16390000e-01,  1.48932502e+00)
( 4.23080000e-01,  1.49349299e+00)
( 4.29770000e-01,  1.49752947e+00)
( 4.36450000e-01,  1.50146080e+00)
( 4.43140000e-01,  1.50534448e+00)
( 4.49830000e-01,  1.50898862e+00)
( 4.56520000e-01,  1.51283263e+00)
( 4.63210000e-01,  1.51639643e+00)
( 4.69900000e-01,  1.51983455e+00)
( 4.76590000e-01,  1.52330546e+00)
( 4.83280000e-01,  1.52667824e+00)
( 4.89970000e-01,  1.52998069e+00)
( 4.96660000e-01,  1.53311238e+00)
( 5.03340000e-01,  1.53623287e+00)
( 5.10030000e-01,  1.53914052e+00)
( 5.16720000e-01,  1.54213309e+00)
( 5.23410000e-01,  1.54593450e+00)
( 5.30100000e-01,  1.54860302e+00)
( 5.36790000e-01,  1.55139482e+00)
( 5.43480000e-01,  1.55420616e+00)
( 5.50170000e-01,  1.55686315e+00)
( 5.56860000e-01,  1.55950545e+00)
( 5.63550000e-01,  1.56223433e+00)
( 5.70230000e-01,  1.56473326e+00)
( 5.76920000e-01,  1.56732292e+00)
( 5.83610000e-01,  1.56978058e+00)
( 5.90300000e-01,  1.57229044e+00)
( 5.96990000e-01,  1.57474143e+00)
( 6.03680000e-01,  1.57724424e+00)
( 6.10370000e-01,  1.57960961e+00)
( 6.17060000e-01,  1.58184751e+00)
( 6.23750000e-01,  1.58431156e+00)
( 6.30430000e-01,  1.58645758e+00)
( 6.37120000e-01,  1.58876190e+00)
( 6.43810000e-01,  1.59103560e+00)
( 6.50500000e-01,  1.59324488e+00)
( 6.57190000e-01,  1.59535187e+00)
( 6.63880000e-01,  1.59742461e+00)
( 6.70570000e-01,  1.59958701e+00)
( 6.77260000e-01,  1.60148810e+00)
( 6.83950000e-01,  1.60362555e+00)
( 6.90640000e-01,  1.60561613e+00)
( 6.97320000e-01,  1.60757702e+00)
( 7.04010000e-01,  1.60938905e+00)
( 7.10700000e-01,  1.61131900e+00)
( 7.17390000e-01,  1.61320631e+00)
( 7.24080000e-01,  1.61506137e+00)
( 7.30770000e-01,  1.61676313e+00)
( 7.37460000e-01,  1.61858469e+00)
( 7.44150000e-01,  1.62037377e+00)
( 7.50840000e-01,  1.62204936e+00)
( 7.57530000e-01,  1.62381197e+00)
( 7.64210000e-01,  1.62557407e+00)
( 7.70900000e-01,  1.62734258e+00)
( 7.77590000e-01,  1.62895449e+00)
( 7.84280000e-01,  1.63054120e+00)
( 7.90970000e-01,  1.63232394e+00)
( 7.97660000e-01,  1.63387251e+00)
( 8.04350000e-01,  1.63546954e+00)
( 8.11040000e-01,  1.63714652e+00)
( 8.17730000e-01,  1.63874496e+00)
( 8.24410000e-01,  1.64026632e+00)
( 8.31100000e-01,  1.64183609e+00)
( 8.37790000e-01,  1.64336816e+00)
( 8.44480000e-01,  1.64502344e+00)
( 8.51170000e-01,  1.64652588e+00)
( 8.57860000e-01,  1.64806656e+00)
( 8.64550000e-01,  1.64969774e+00)
( 8.71240000e-01,  1.65120623e+00)
( 8.77930000e-01,  1.65267855e+00)
( 8.84620000e-01,  1.65419943e+00)
( 8.91300000e-01,  1.65550525e+00)
( 8.97990000e-01,  1.65715492e+00)
( 9.04680000e-01,  1.65860328e+00)
( 9.11370000e-01,  1.66009072e+00)
( 9.18060000e-01,  1.66137644e+00)
( 9.24750000e-01,  1.66287656e+00)
( 9.31440000e-01,  1.66447082e+00)
( 9.38130000e-01,  1.66593750e+00)
( 9.44820000e-01,  1.66741148e+00)
( 9.51510000e-01,  1.66884548e+00)
( 9.58190000e-01,  1.67016547e+00)
( 9.64880000e-01,  1.67147605e+00)
( 9.71570000e-01,  1.67297483e+00)
( 9.78260000e-01,  1.67443505e+00)
( 9.84950000e-01,  1.67590064e+00)
( 9.91640000e-01,  1.67741983e+00)
( 9.98330000e-01,  1.67876059e+00)
( 1.00500000e+00,  1.68024592e+00)
( 1.01170000e+00,  1.68173864e+00)
( 1.01840000e+00,  1.68309989e+00)
( 1.02510000e+00,  1.68460713e+00)
( 1.03180000e+00,  1.68611991e+00)
( 1.03850000e+00,  1.68732838e+00)
( 1.04520000e+00,  1.68885548e+00)
( 1.05180000e+00,  1.69039025e+00)
( 1.05850000e+00,  1.69179125e+00)
( 1.06520000e+00,  1.69316561e+00)
( 1.07190000e+00,  1.69485601e+00)
( 1.07860000e+00,  1.69641374e+00)
( 1.08530000e+00,  1.69758915e+00)
( 1.09200000e+00,  1.69892024e+00)
( 1.09870000e+00,  1.70050115e+00)
( 1.10540000e+00,  1.70209668e+00)
( 1.11200000e+00,  1.70362391e+00)
( 1.11870000e+00,  1.70493552e+00)
( 1.12540000e+00,  1.70667781e+00)
( 1.13210000e+00,  1.70805598e+00)
( 1.13880000e+00,  1.70981609e+00)
( 1.14550000e+00,  1.71121000e+00)
( 1.15220000e+00,  1.71280901e+00)
( 1.15890000e+00,  1.71439570e+00)
( 1.16560000e+00,  1.71609062e+00)
( 1.17220000e+00,  1.71764219e+00)
( 1.17890000e+00,  1.71906010e+00)
( 1.18560000e+00,  1.72032073e+00)
( 1.19230000e+00,  1.72110971e+00)
( 1.19900000e+00,  1.45619331e+00)
( 1.20570000e+00,  1.45616419e+00)
( 1.21240000e+00,  1.45607979e+00)
( 1.21910000e+00,  1.45605067e+00)
( 1.22580000e+00,  1.45605067e+00)
( 1.23240000e+00,  1.45605067e+00)
( 1.23910000e+00,  1.45613642e+00)
( 1.24580000e+00,  1.45613642e+00)
( 1.25250000e+00,  1.45610729e+00)
( 1.25920000e+00,  1.45602291e+00)
( 1.26590000e+00,  1.45599379e+00)
( 1.27260000e+00,  1.45599379e+00)
( 1.27930000e+00,  1.45596467e+00)
( 1.28600000e+00,  1.45607953e+00)
( 1.29260000e+00,  1.45596604e+00)
( 1.29930000e+00,  1.45596604e+00)
( 1.30600000e+00,  1.45593692e+00)
( 1.31270000e+00,  1.45590781e+00)
( 1.31940000e+00,  1.45582346e+00)
( 1.32610000e+00,  1.45590918e+00)
( 1.33280000e+00,  1.45590918e+00)
( 1.33950000e+00,  1.45588007e+00)
( 1.34620000e+00,  1.45576662e+00)
( 1.35280000e+00,  1.45573752e+00)
( 1.35950000e+00,  1.45585233e+00)
( 1.36620000e+00,  1.45588007e+00)
( 1.37290000e+00,  1.45588006e+00)
( 1.37960000e+00,  1.45576662e+00)
( 1.38630000e+00,  1.45576662e+00)
( 1.39300000e+00,  1.45573751e+00)
( 1.39970000e+00,  1.45573751e+00)
( 1.40640000e+00,  1.45573751e+00)
( 1.41300000e+00,  1.45570841e+00)
( 1.41970000e+00,  1.45562409e+00)
( 1.42640000e+00,  1.45562409e+00)
( 1.43310000e+00,  1.45562409e+00)
( 1.43980000e+00,  1.45562409e+00)
( 1.44650000e+00,  1.45562409e+00)
( 1.45320000e+00,  1.45562409e+00)
( 1.45990000e+00,  1.45562409e+00)
( 1.46660000e+00,  1.45562409e+00)
( 1.47320000e+00,  1.45550931e+00)
( 1.47990000e+00,  1.45550931e+00)
( 1.48660000e+00,  1.45550931e+00)
( 1.49330000e+00,  1.45553841e+00)};

\nextgroupplot[width=0.3\textwidth, xlabel={$x$}, ymax=2.4, xmax=1.5, ylabel={$P$}, xmin=-0.5, ymin=0.6, height=0.3\textwidth]
\addplot [solid, thick, color=black, forget plot]
coordinates {
(-5.00000000e-01,  1.00000000e+00)
(-4.93310000e-01,  1.00000000e+00)
(-4.86620000e-01,  1.00000000e+00)
(-4.79930000e-01,  1.00000000e+00)
(-4.73240000e-01,  1.00000000e+00)
(-4.66560000e-01,  1.00000000e+00)
(-4.59870000e-01,  1.00000000e+00)
(-4.53180000e-01,  1.00000000e+00)
(-4.46490000e-01,  1.00000000e+00)
(-4.39800000e-01,  1.00000000e+00)
(-4.33110000e-01,  1.00000000e+00)
(-4.26420000e-01,  1.00000000e+00)
(-4.19730000e-01,  1.00000000e+00)
(-4.13040000e-01,  1.00000000e+00)
(-4.06350000e-01,  1.00000000e+00)
(-3.99670000e-01,  1.00000000e+00)
(-3.92980000e-01,  1.00000000e+00)
(-3.86290000e-01,  1.00000000e+00)
(-3.79600000e-01,  1.00000000e+00)
(-3.72910000e-01,  1.00000000e+00)
(-3.66220000e-01,  1.00000000e+00)
(-3.59530000e-01,  1.00000000e+00)
(-3.52840000e-01,  1.00000000e+00)
(-3.46150000e-01,  1.00000000e+00)
(-3.39460000e-01,  1.00000000e+00)
(-3.32780000e-01,  1.00000000e+00)
(-3.26090000e-01,  1.00000000e+00)
(-3.19400000e-01,  1.00000000e+00)
(-3.12710000e-01,  1.00000000e+00)
(-3.06020000e-01,  1.00000000e+00)
(-2.99330000e-01,  1.00000000e+00)
(-2.92640000e-01,  1.00000000e+00)
(-2.85950000e-01,  1.00000000e+00)
(-2.79260000e-01,  1.00000000e+00)
(-2.72580000e-01,  1.00000000e+00)
(-2.65890000e-01,  1.00000000e+00)
(-2.59200000e-01,  1.00000000e+00)
(-2.52510000e-01,  1.00000000e+00)
(-2.45820000e-01,  1.00000000e+00)
(-2.39130000e-01,  1.00000000e+00)
(-2.32440000e-01,  1.00000000e+00)
(-2.25750000e-01,  1.00000000e+00)
(-2.19060000e-01,  1.00000000e+00)
(-2.12370000e-01,  1.00000000e+00)
(-2.05690000e-01,  1.00000000e+00)
(-1.99000000e-01,  1.00000000e+00)
(-1.92310000e-01,  1.00000000e+00)
(-1.85620000e-01,  1.00000000e+00)
(-1.78930000e-01,  1.00000000e+00)
(-1.72240000e-01,  1.00000000e+00)
(-1.65550000e-01,  1.00000000e+00)
(-1.58860000e-01,  1.00000000e+00)
(-1.52170000e-01,  1.00000000e+00)
(-1.45480000e-01,  1.00000000e+00)
(-1.38800000e-01,  1.00000000e+00)
(-1.32110000e-01,  1.00000000e+00)
(-1.25420000e-01,  1.00000000e+00)
(-1.18730000e-01,  1.00000000e+00)
(-1.12040000e-01,  1.00000000e+00)
(-1.05350000e-01,  1.00000000e+00)
(-9.86620000e-02,  1.00000000e+00)
(-9.19730000e-02,  1.00000000e+00)
(-8.52840000e-02,  1.00000000e+00)
(-7.85950000e-02,  1.00004000e+00)
(-7.19060000e-02,  1.00008500e+00)
(-6.52170000e-02,  2.28430116e+00)
(-5.85280000e-02,  2.28705544e+00)
(-5.18390000e-02,  2.28688749e+00)
(-4.51510000e-02,  2.28367880e+00)
(-3.84620000e-02,  2.27719519e+00)
(-3.17730000e-02,  2.26731316e+00)
(-2.50840000e-02,  2.25405187e+00)
(-1.83950000e-02,  2.23710109e+00)
(-1.17060000e-02,  2.21694916e+00)
(-5.01670000e-03,  2.19340459e+00)
( 1.67220000e-03,  2.16670704e+00)
( 8.36120000e-03,  2.13696187e+00)
( 1.50500000e-02,  2.10439145e+00)
( 2.17390000e-02,  2.06971641e+00)
( 2.84280000e-02,  2.03290397e+00)
( 3.51170000e-02,  1.99529092e+00)
( 4.18060000e-02,  1.95666954e+00)
( 4.84950000e-02,  1.91736560e+00)
( 5.51840000e-02,  1.87773515e+00)
( 6.18730000e-02,  1.83809343e+00)
( 6.85620000e-02,  1.79906101e+00)
( 7.52510000e-02,  1.76055979e+00)
( 8.19400000e-02,  1.72318442e+00)
( 8.86290000e-02,  1.68708238e+00)
( 9.53180000e-02,  1.65231543e+00)
( 1.02010000e-01,  1.61889965e+00)
( 1.08700000e-01,  1.58690032e+00)
( 1.15380000e-01,  1.55627204e+00)
( 1.22070000e-01,  1.52691443e+00)
( 1.28760000e-01,  1.49883815e+00)
( 1.35450000e-01,  1.47182347e+00)
( 1.42140000e-01,  1.44634309e+00)
( 1.48830000e-01,  1.42228004e+00)
( 1.55520000e-01,  1.39926793e+00)
( 1.62210000e-01,  1.37695375e+00)
( 1.68900000e-01,  1.35542621e+00)
( 1.75590000e-01,  1.33488433e+00)
( 1.82270000e-01,  1.31518132e+00)
( 1.88960000e-01,  1.29645227e+00)
( 1.95650000e-01,  1.27880797e+00)
( 2.02340000e-01,  1.26226456e+00)
( 2.09030000e-01,  1.24620142e+00)
( 2.15720000e-01,  1.23053273e+00)
( 2.22410000e-01,  1.21554221e+00)
( 2.29100000e-01,  1.20117294e+00)
( 2.35790000e-01,  1.18731270e+00)
( 2.42470000e-01,  1.17406437e+00)
( 2.49160000e-01,  1.16128885e+00)
( 2.55850000e-01,  1.14893759e+00)
( 2.62540000e-01,  1.13707697e+00)
( 2.69230000e-01,  1.12560950e+00)
( 2.75920000e-01,  1.11450554e+00)
( 2.82610000e-01,  1.10378445e+00)
( 2.89300000e-01,  1.09354685e+00)
( 2.95990000e-01,  1.08356153e+00)
( 3.02680000e-01,  1.07405011e+00)
( 3.09360000e-01,  1.06485319e+00)
( 3.16050000e-01,  1.05609573e+00)
( 3.22740000e-01,  1.04764162e+00)
( 3.29430000e-01,  1.03885638e+00)
( 3.36120000e-01,  1.03090843e+00)
( 3.42810000e-01,  1.02312441e+00)
( 3.49500000e-01,  1.01550776e+00)
( 3.56190000e-01,  1.00811558e+00)
( 3.62880000e-01,  1.00091778e+00)
( 3.69570000e-01,  9.93950293e-01)
( 3.76250000e-01,  9.87152400e-01)
( 3.82940000e-01,  9.80587505e-01)
( 3.89630000e-01,  9.74134918e-01)
( 3.96320000e-01,  9.67921623e-01)
( 4.03010000e-01,  9.61837575e-01)
( 4.09700000e-01,  9.56020548e-01)
( 4.16390000e-01,  9.50313574e-01)
( 4.23080000e-01,  9.44682411e-01)
( 4.29770000e-01,  9.39207044e-01)
( 4.36450000e-01,  9.33873320e-01)
( 4.43140000e-01,  9.28676934e-01)
( 4.49830000e-01,  9.23683976e-01)
( 4.56520000e-01,  9.18683908e-01)
( 4.63210000e-01,  9.13928626e-01)
( 4.69900000e-01,  9.09318806e-01)
( 4.76590000e-01,  9.04803334e-01)
( 4.83280000e-01,  9.00394462e-01)
( 4.89970000e-01,  8.96126754e-01)
( 4.96660000e-01,  8.92021539e-01)
( 5.03340000e-01,  8.88004067e-01)
( 5.10030000e-01,  8.84190356e-01)
( 5.16720000e-01,  8.80418884e-01)
( 5.23410000e-01,  8.75665718e-01)
( 5.30100000e-01,  8.72180971e-01)
( 5.36790000e-01,  8.68672101e-01)
( 5.43480000e-01,  8.65186206e-01)
( 5.50170000e-01,  8.61824926e-01)
( 5.56860000e-01,  8.58501237e-01)
( 5.63550000e-01,  8.55193709e-01)
( 5.70230000e-01,  8.52019006e-01)
( 5.76920000e-01,  8.48834432e-01)
( 5.83610000e-01,  8.45761941e-01)
( 5.90300000e-01,  8.42719990e-01)
( 5.96990000e-01,  8.39709948e-01)
( 6.03680000e-01,  8.36730495e-01)
( 6.10370000e-01,  8.33863839e-01)
( 6.17060000e-01,  8.31057884e-01)
( 6.23750000e-01,  8.28206482e-01)
( 6.30430000e-01,  8.25544485e-01)
( 6.37120000e-01,  8.22819286e-01)
( 6.43810000e-01,  8.20159104e-01)
( 6.50500000e-01,  8.17578061e-01)
( 6.57190000e-01,  8.15045145e-01)
( 6.63880000e-01,  8.12577791e-01)
( 6.70570000e-01,  8.10100030e-01)
( 6.77260000e-01,  8.07760065e-01)
( 6.83950000e-01,  8.05347663e-01)
( 6.90640000e-01,  8.03024258e-01)
( 6.97320000e-01,  8.00740269e-01)
( 7.04010000e-01,  7.98545232e-01)
( 7.10700000e-01,  7.96327596e-01)
( 7.17390000e-01,  7.94177443e-01)
( 7.24080000e-01,  7.92067075e-01)
( 7.30770000e-01,  7.90046324e-01)
( 7.37460000e-01,  7.88002141e-01)
( 7.44150000e-01,  7.85997406e-01)
( 7.50840000e-01,  7.84042145e-01)
( 7.57530000e-01,  7.82077009e-01)
( 7.64210000e-01,  7.80137914e-01)
( 7.70900000e-01,  7.78198998e-01)
( 7.77590000e-01,  7.76349461e-01)
( 7.84280000e-01,  7.74512608e-01)
( 7.90970000e-01,  7.72598887e-01)
( 7.97660000e-01,  7.70828233e-01)
( 8.04350000e-01,  7.69017466e-01)
( 8.11040000e-01,  7.67222673e-01)
( 8.17730000e-01,  7.65439000e-01)
( 8.24410000e-01,  7.63734154e-01)
( 8.31100000e-01,  7.61989340e-01)
( 8.37790000e-01,  7.60284583e-01)
( 8.44480000e-01,  7.58555772e-01)
( 8.51170000e-01,  7.56889931e-01)
( 8.57860000e-01,  7.55211028e-01)
( 8.64550000e-01,  7.53521211e-01)
( 8.71240000e-01,  7.51882016e-01)
( 8.77930000e-01,  7.50281768e-01)
( 8.84620000e-01,  7.48641756e-01)
( 8.91300000e-01,  7.47133350e-01)
( 8.97990000e-01,  7.45469505e-01)
( 9.04680000e-01,  7.43909057e-01)
( 9.11370000e-01,  7.42335263e-01)
( 9.18060000e-01,  7.40864388e-01)
( 9.24750000e-01,  7.39289370e-01)
( 9.31440000e-01,  7.37702377e-01)
( 9.38130000e-01,  7.36166751e-01)
( 9.44820000e-01,  7.34630174e-01)
( 9.51510000e-01,  7.33133338e-01)
( 9.58190000e-01,  7.31659624e-01)
( 9.64880000e-01,  7.30213606e-01)
( 9.71570000e-01,  7.28674717e-01)
( 9.78260000e-01,  7.27174877e-01)
( 9.84950000e-01,  7.25674803e-01)
( 9.91640000e-01,  7.24133778e-01)
( 9.98330000e-01,  7.22684802e-01)
( 1.00500000e+00,  7.21182668e-01)
( 1.01170000e+00,  7.19679580e-01)
( 1.01840000e+00,  7.18227980e-01)
( 1.02510000e+00,  7.16723036e-01)
( 1.03180000e+00,  7.15217884e-01)
( 1.03850000e+00,  7.13828817e-01)
( 1.04520000e+00,  7.12321837e-01)
( 1.05180000e+00,  7.10813890e-01)
( 1.05850000e+00,  7.09357759e-01)
( 1.06520000e+00,  7.07913256e-01)
( 1.07190000e+00,  7.06197809e-01)
( 1.07860000e+00,  7.04777612e-01)
( 1.08530000e+00,  7.03456433e-01)
( 1.09200000e+00,  7.02079075e-01)
( 1.09870000e+00,  7.00566945e-01)
( 1.10540000e+00,  6.99051407e-01)
( 1.11200000e+00,  6.97518805e-01)
( 1.11870000e+00,  6.96089208e-01)
( 1.12540000e+00,  6.94459684e-01)
( 1.13210000e+00,  6.92987332e-01)
( 1.13880000e+00,  6.91355675e-01)
( 1.14550000e+00,  6.89881367e-01)
( 1.15220000e+00,  6.88312501e-01)
( 1.15890000e+00,  6.86772130e-01)
( 1.16560000e+00,  6.85128976e-01)
( 1.17220000e+00,  6.83627378e-01)
( 1.17890000e+00,  6.82241426e-01)
( 1.18560000e+00,  6.81001035e-01)
( 1.19230000e+00,  6.80083946e-01)
( 1.19900000e+00,  9.99862141e-01)
( 1.20570000e+00,  9.99902141e-01)
( 1.21240000e+00,  9.99945319e-01)
( 1.21910000e+00,  9.99985319e-01)
( 1.22580000e+00,  9.99985319e-01)
( 1.23240000e+00,  9.99985319e-01)
( 1.23910000e+00,  9.99966547e-01)
( 1.24580000e+00,  9.99966547e-01)
( 1.25250000e+00,  1.00000655e+00)
( 1.25920000e+00,  1.00004972e+00)
( 1.26590000e+00,  1.00008972e+00)
( 1.27260000e+00,  1.00008972e+00)
( 1.27930000e+00,  1.00012972e+00)
( 1.28600000e+00,  1.00007095e+00)
( 1.29260000e+00,  1.00015413e+00)
( 1.29930000e+00,  1.00015413e+00)
( 1.30600000e+00,  1.00019413e+00)
( 1.31270000e+00,  1.00023413e+00)
( 1.31940000e+00,  1.00027730e+00)
( 1.32610000e+00,  1.00025853e+00)
( 1.33280000e+00,  1.00025853e+00)
( 1.33950000e+00,  1.00029853e+00)
( 1.34620000e+00,  1.00038170e+00)
( 1.35280000e+00,  1.00042170e+00)
( 1.35950000e+00,  1.00036293e+00)
( 1.36620000e+00,  1.00029853e+00)
( 1.37290000e+00,  1.00029853e+00)
( 1.37960000e+00,  1.00038170e+00)
( 1.38630000e+00,  1.00038170e+00)
( 1.39300000e+00,  1.00042170e+00)
( 1.39970000e+00,  1.00042170e+00)
( 1.40640000e+00,  1.00042170e+00)
( 1.41300000e+00,  1.00046170e+00)
( 1.41970000e+00,  1.00050486e+00)
( 1.42640000e+00,  1.00050486e+00)
( 1.43310000e+00,  1.00050486e+00)
( 1.43980000e+00,  1.00050486e+00)
( 1.44650000e+00,  1.00050486e+00)
( 1.45320000e+00,  1.00050486e+00)
( 1.45990000e+00,  1.00050486e+00)
( 1.46660000e+00,  1.00050486e+00)
( 1.47320000e+00,  1.00056362e+00)
( 1.47990000e+00,  1.00056362e+00)
( 1.48660000e+00,  1.00056361e+00)
( 1.49330000e+00,  1.00052361e+00)};

\end{groupplot}\end{tikzpicture}

%% file: tikz/euler2d0_naca0012cfg1_nref0_c2s1.tikz
\begin{tikzpicture}
\begin{groupplot}[
  group style={
      group size=3 by 2,
      horizontal sep=0.9cm,
      vertical sep=0.7cm
  },
  width=0.5\textwidth,
  axis equal image,
  xlabel={$x_1$},
  ylabel={$x_2$},
  xtick = {-0.5, 0, 0.5, 1.0, 1.5},
  ytick = {0.0, 1.0, 2.0, 3.0},
  xmin=-0.5, xmax=1.5,
  ymin=0, ymax=3
]

\nextgroupplot[title={Initialization}, xlabel={}, xtick=\empty]
\addplot graphics [xmin=-0.5, xmax=1.5, ymin=0, ymax=3] {{img/euler2d0_naca0012cfg1_nref0p1_c2s1_mesh_iter000}.png};

\nextgroupplot[title={Converged ($p=1$)}, xlabel={}, ylabel={}, xtick=\empty, ytick=\empty]
\addplot graphics [xmin=-0.5, xmax=1.5, ymin=0, ymax=3] {{img/euler2d0_naca0012cfg1_nref0p1_c2s1_mesh_iter030}.png};

\nextgroupplot[title={Converged ($p=2$)}, xlabel={}, ylabel={}, xtick=\empty, ytick=\empty]
\addplot graphics [xmin=-0.5, xmax=1.5, ymin=0, ymax=3] {{img/euler2d0_naca0012cfg1_nref0p2_c2s1_mesh_iter030}.png};

\nextgroupplot
\addplot graphics [xmin=-0.5, xmax=1.5, ymin=0, ymax=3] {{img/euler2d0_naca0012cfg1_nref0p1_c2s1_iter000}.png};

\nextgroupplot[ylabel={}, ytick=\empty]
\addplot graphics [xmin=-0.5, xmax=1.5, ymin=0, ymax=3] {{img/euler2d0_naca0012cfg1_nref0p1_c2s1_iter030}.png};

\nextgroupplot[ylabel={}, ytick=\empty]
\addplot graphics [xmin=-0.5, xmax=1.5, ymin=0, ymax=3] {{img/euler2d0_naca0012cfg1_nref0p2_c2s1_iter030}.png};

\end{groupplot}
\end{tikzpicture}

%% file: tikz/cbar_parula_0x1p4.tikz
\begin{tikzpicture}
\begin{axis}[
   hide axis, scale only axis,
   height=0pt, width=0pt,
   colormap={parula}{rgb255=(62,38,168) rgb255=(62,39,172) rgb255=(63,40,175) rgb255=(63,41,178) rgb255=(64,42,180) rgb255=(64,43,183) rgb255=(65,44,186) rgb255=(65,45,189) rgb255=(66,46,191) rgb255=(66,47,194) rgb255=(67,48,197) rgb255=(67,49,200) rgb255=(67,50,202) rgb255=(68,51,205) rgb255=(68,52,208) rgb255=(69,53,210) rgb255=(69,55,213) rgb255=(69,56,215) rgb255=(70,57,217) rgb255=(70,58,220) rgb255=(70,59,222) rgb255=(70,61,224) rgb255=(71,62,225) rgb255=(71,63,227) rgb255=(71,65,229) rgb255=(71,66,230) rgb255=(71,68,232) rgb255=(71,69,233) rgb255=(71,70,235) rgb255=(72,72,236) rgb255=(72,73,237) rgb255=(72,75,238) rgb255=(72,76,240) rgb255=(72,78,241) rgb255=(72,79,242) rgb255=(72,80,243) rgb255=(72,82,244) rgb255=(72,83,245) rgb255=(72,84,246) rgb255=(71,86,247) rgb255=(71,87,247) rgb255=(71,89,248) rgb255=(71,90,249) rgb255=(71,91,250) rgb255=(71,93,250) rgb255=(70,94,251) rgb255=(70,96,251) rgb255=(70,97,252) rgb255=(69,98,252) rgb255=(69,100,253) rgb255=(68,101,253) rgb255=(67,103,253) rgb255=(67,104,254) rgb255=(66,106,254) rgb255=(65,107,254) rgb255=(64,109,254) rgb255=(63,110,255) rgb255=(62,112,255) rgb255=(60,113,255) rgb255=(59,115,255) rgb255=(57,116,255) rgb255=(56,118,254) rgb255=(54,119,254) rgb255=(53,121,253) rgb255=(51,122,253) rgb255=(50,124,252) rgb255=(49,125,252) rgb255=(48,127,251) rgb255=(47,128,250) rgb255=(47,130,250) rgb255=(46,131,249) rgb255=(46,132,248) rgb255=(46,134,248) rgb255=(46,135,247) rgb255=(45,136,246) rgb255=(45,138,245) rgb255=(45,139,244) rgb255=(45,140,243) rgb255=(45,142,242) rgb255=(44,143,241) rgb255=(44,144,240) rgb255=(43,145,239) rgb255=(42,147,238) rgb255=(41,148,237) rgb255=(40,149,236) rgb255=(39,151,235) rgb255=(39,152,234) rgb255=(38,153,233) rgb255=(38,154,232) rgb255=(37,155,232) rgb255=(37,156,231) rgb255=(36,158,230) rgb255=(36,159,229) rgb255=(35,160,229) rgb255=(35,161,228) rgb255=(34,162,228) rgb255=(33,163,227) rgb255=(32,165,227) rgb255=(31,166,226) rgb255=(30,167,225) rgb255=(29,168,225) rgb255=(29,169,224) rgb255=(28,170,223) rgb255=(27,171,222) rgb255=(26,172,221) rgb255=(25,173,220) rgb255=(23,174,218) rgb255=(22,175,217) rgb255=(20,176,216) rgb255=(18,177,214) rgb255=(16,178,213) rgb255=(14,179,212) rgb255=(11,179,210) rgb255=(8,180,209) rgb255=(6,181,207) rgb255=(4,182,206) rgb255=(2,183,204) rgb255=(1,183,202) rgb255=(0,184,201) rgb255=(0,185,199) rgb255=(0,186,198) rgb255=(1,186,196) rgb255=(2,187,194) rgb255=(4,187,193) rgb255=(6,188,191) rgb255=(9,189,189) rgb255=(13,189,188) rgb255=(16,190,186) rgb255=(20,190,184) rgb255=(23,191,182) rgb255=(26,192,181) rgb255=(29,192,179) rgb255=(32,193,177) rgb255=(35,193,175) rgb255=(37,194,174) rgb255=(39,194,172) rgb255=(41,195,170) rgb255=(43,195,168) rgb255=(44,196,166) rgb255=(46,196,165) rgb255=(47,197,163) rgb255=(49,197,161) rgb255=(50,198,159) rgb255=(51,199,157) rgb255=(53,199,155) rgb255=(54,200,153) rgb255=(56,200,150) rgb255=(57,201,148) rgb255=(59,201,146) rgb255=(61,202,144) rgb255=(64,202,141) rgb255=(66,202,139) rgb255=(69,203,137) rgb255=(72,203,134) rgb255=(75,203,132) rgb255=(78,204,129) rgb255=(81,204,127) rgb255=(84,204,124) rgb255=(87,204,122) rgb255=(90,204,119) rgb255=(94,205,116) rgb255=(97,205,114) rgb255=(100,205,111) rgb255=(103,205,108) rgb255=(107,205,105) rgb255=(110,205,102) rgb255=(114,205,100) rgb255=(118,204,97) rgb255=(121,204,94) rgb255=(125,204,91) rgb255=(129,204,89) rgb255=(132,204,86) rgb255=(136,203,83) rgb255=(139,203,81) rgb255=(143,203,78) rgb255=(147,202,75) rgb255=(150,202,72) rgb255=(154,201,70) rgb255=(157,201,67) rgb255=(161,200,64) rgb255=(164,200,62) rgb255=(167,199,59) rgb255=(171,199,57) rgb255=(174,198,55) rgb255=(178,198,53) rgb255=(181,197,51) rgb255=(184,196,49) rgb255=(187,196,47) rgb255=(190,195,45) rgb255=(194,195,44) rgb255=(197,194,42) rgb255=(200,193,41) rgb255=(203,193,40) rgb255=(206,192,39) rgb255=(208,191,39) rgb255=(211,191,39) rgb255=(214,190,39) rgb255=(217,190,40) rgb255=(219,189,40) rgb255=(222,188,41) rgb255=(225,188,42) rgb255=(227,188,43) rgb255=(230,187,45) rgb255=(232,187,46) rgb255=(234,186,48) rgb255=(236,186,50) rgb255=(239,186,53) rgb255=(241,186,55) rgb255=(243,186,57) rgb255=(245,186,59) rgb255=(247,186,61) rgb255=(249,186,62) rgb255=(251,187,62) rgb255=(252,188,62) rgb255=(254,189,61) rgb255=(254,190,60) rgb255=(254,192,59) rgb255=(254,193,58) rgb255=(254,194,57) rgb255=(254,196,56) rgb255=(254,197,55) rgb255=(254,199,53) rgb255=(254,200,52) rgb255=(254,202,51) rgb255=(253,203,50) rgb255=(253,205,49) rgb255=(253,206,49) rgb255=(252,208,48) rgb255=(251,210,47) rgb255=(251,211,46) rgb255=(250,213,46) rgb255=(249,214,45) rgb255=(249,216,44) rgb255=(248,217,43) rgb255=(247,219,42) rgb255=(247,221,42) rgb255=(246,222,41) rgb255=(246,224,40) rgb255=(245,225,40) rgb255=(245,227,39) rgb255=(245,229,38) rgb255=(245,230,38) rgb255=(245,232,37) rgb255=(245,233,36) rgb255=(245,235,35) rgb255=(245,236,34) rgb255=(245,238,33) rgb255=(246,239,32) rgb255=(246,241,31) rgb255=(246,242,30) rgb255=(247,244,28) rgb255=(247,245,27) rgb255=(248,247,26) rgb255=(248,248,24) rgb255=(249,249,22) rgb255=(249,251,21) },
   colorbar horizontal,
   point meta min=0, point meta max=1.4,
   colorbar style={width=10cm, xtick={0, 0.35, 0.70, 1.05, 1.4}}
]
\addplot [draw=none] coordinates {(0,0)};
\end{axis}
\end{tikzpicture}

%% file: py/euler2d0_naca0012cfg1_nref0p3_c2s1_slice.tikz
\begin{tikzpicture}
\begin{groupplot} [
group style={group size = 3 by 1, horizontal sep = 1.5cm}]
\nextgroupplot[width=0.3\textwidth, xlabel={$x$}, ymax=1.65, xmax=1.5, ylabel={$\rho$}, xmin=-0.5, ymin=0.9, height=0.3\textwidth]
\addplot [solid, thick, color=black, forget plot]
coordinates {
(-5.00000000e-01,          nan)
(-4.93310000e-01,  1.52900000e+00)
(-4.86620000e-01,  1.52910000e+00)
(-4.79930000e-01,  1.52920000e+00)
(-4.73240000e-01,  1.52930000e+00)
(-4.66560000e-01,  1.52940000e+00)
(-4.59870000e-01,  1.52960000e+00)
(-4.53180000e-01,  1.52980000e+00)
(-4.46490000e-01,  1.53010000e+00)
(-4.39800000e-01,  1.53040000e+00)
(-4.33110000e-01,  1.53070000e+00)
(-4.26420000e-01,  1.53110000e+00)
(-4.19730000e-01,  1.53150000e+00)
(-4.13040000e-01,  1.53190000e+00)
(-4.06350000e-01,  1.53240000e+00)
(-3.99670000e-01,  1.53290000e+00)
(-3.92980000e-01,  1.53340000e+00)
(-3.86290000e-01,  1.53400000e+00)
(-3.79600000e-01,  1.53470000e+00)
(-3.72910000e-01,  1.53530000e+00)
(-3.66220000e-01,  1.53600000e+00)
(-3.59530000e-01,  1.53680000e+00)
(-3.52840000e-01,  1.53760000e+00)
(-3.46150000e-01,  1.53840000e+00)
(-3.39460000e-01,  1.53930000e+00)
(-3.32780000e-01,  1.54020000e+00)
(-3.26090000e-01,  1.54110000e+00)
(-3.19400000e-01,  1.54200000e+00)
(-3.12710000e-01,  1.54310000e+00)
(-3.06020000e-01,  1.54410000e+00)
(-2.99330000e-01,  1.54520000e+00)
(-2.92640000e-01,  1.54640000e+00)
(-2.85950000e-01,  1.54760000e+00)
(-2.79260000e-01,  1.54890000e+00)
(-2.72580000e-01,  1.55020000e+00)
(-2.65890000e-01,  1.55150000e+00)
(-2.59200000e-01,  1.55290000e+00)
(-2.52510000e-01,  1.55430000e+00)
(-2.45820000e-01,  1.55580000e+00)
(-2.39130000e-01,  1.55730000e+00)
(-2.32440000e-01,  1.55890000e+00)
(-2.25750000e-01,  1.56050000e+00)
(-2.19060000e-01,  1.56210000e+00)
(-2.12370000e-01,  1.56380000e+00)
(-2.05690000e-01,  1.56550000e+00)
(-1.99000000e-01,  1.56730000e+00)
(-1.92310000e-01,  1.56900000e+00)
(-1.85620000e-01,  1.57090000e+00)
(-1.78930000e-01,  1.57270000e+00)
(-1.72240000e-01,  1.57440000e+00)
(-1.65550000e-01,  1.57620000e+00)
(-1.58860000e-01,  1.57800000e+00)
(-1.52170000e-01,  1.57980000e+00)
(-1.45480000e-01,  1.58160000e+00)
(-1.38800000e-01,  1.58320000e+00)
(-1.32110000e-01,  1.58480000e+00)
(-1.25420000e-01,  1.58620000e+00)
(-1.18730000e-01,  1.58750000e+00)
(-1.12040000e-01,  1.58870000e+00)
(-1.05350000e-01,  1.58960000e+00)
(-9.86620000e-02,  1.59020000e+00)
(-9.19730000e-02,  1.59060000e+00)
(-8.52840000e-02,  1.59050000e+00)
(-7.85950000e-02,  1.59000000e+00)
(-7.19060000e-02,  1.58900000e+00)
(-6.52170000e-02,  1.58730000e+00)
(-5.85280000e-02,  1.58500000e+00)
(-5.18390000e-02,  1.58190000e+00)
(-4.51510000e-02,  1.57820000e+00)
(-3.84620000e-02,  1.57340000e+00)
(-3.17730000e-02,  1.56760000e+00)
(-2.50840000e-02,  1.56070000e+00)
(-1.83950000e-02,  1.55280000e+00)
(-1.17060000e-02,  1.54390000e+00)
(-5.01670000e-03,  1.53370000e+00)
( 1.67220000e-03,  1.52250000e+00)
( 8.36120000e-03,  1.51030000e+00)
( 1.50500000e-02,  1.49730000e+00)
( 2.17390000e-02,  1.48350000e+00)
( 2.84280000e-02,  1.46900000e+00)
( 3.51170000e-02,  1.45400000e+00)
( 4.18060000e-02,  1.43850000e+00)
( 4.84950000e-02,  1.42270000e+00)
( 5.51840000e-02,  1.40670000e+00)
( 6.18730000e-02,  1.39120000e+00)
( 6.85620000e-02,  1.37580000e+00)
( 7.52510000e-02,  1.36040000e+00)
( 8.19400000e-02,  1.34540000e+00)
( 8.86290000e-02,  1.33090000e+00)
( 9.53180000e-02,  1.31680000e+00)
( 1.02010000e-01,  1.30300000e+00)
( 1.08700000e-01,  1.28960000e+00)
( 1.15380000e-01,  1.27670000e+00)
( 1.22070000e-01,  1.26420000e+00)
( 1.28760000e-01,  1.25230000e+00)
( 1.35450000e-01,  1.24090000e+00)
( 1.42140000e-01,  1.23000000e+00)
( 1.48830000e-01,  1.21980000e+00)
( 1.55520000e-01,  1.21020000e+00)
( 1.62210000e-01,  1.20130000e+00)
( 1.68900000e-01,  1.19210000e+00)
( 1.75590000e-01,  1.18340000e+00)
( 1.82270000e-01,  1.17510000e+00)
( 1.88960000e-01,  1.16700000e+00)
( 1.95650000e-01,  1.15920000e+00)
( 2.02340000e-01,  1.15170000e+00)
( 2.09030000e-01,  1.14450000e+00)
( 2.15720000e-01,  1.13750000e+00)
( 2.22410000e-01,  1.13070000e+00)
( 2.29100000e-01,  1.12430000e+00)
( 2.35790000e-01,  1.11810000e+00)
( 2.42470000e-01,  1.11210000e+00)
( 2.49160000e-01,  1.10640000e+00)
( 2.55850000e-01,  1.10090000e+00)
( 2.62540000e-01,  1.09560000e+00)
( 2.69230000e-01,  1.09050000e+00)
( 2.75920000e-01,  1.08550000e+00)
( 2.82610000e-01,  1.08070000e+00)
( 2.89300000e-01,  1.07600000e+00)
( 2.95990000e-01,  1.07150000e+00)
( 3.02680000e-01,  1.06720000e+00)
( 3.09360000e-01,  1.06290000e+00)
( 3.16050000e-01,  1.05890000e+00)
( 3.22740000e-01,  1.05490000e+00)
( 3.29430000e-01,  1.05110000e+00)
( 3.36120000e-01,  1.04740000e+00)
( 3.42810000e-01,  1.04380000e+00)
( 3.49500000e-01,  1.04030000e+00)
( 3.56190000e-01,  1.03690000e+00)
( 3.62880000e-01,  1.03360000e+00)
( 3.69570000e-01,  1.03040000e+00)
( 3.76250000e-01,  1.02730000e+00)
( 3.82940000e-01,  1.02420000e+00)
( 3.89630000e-01,  1.02130000e+00)
( 3.96320000e-01,  1.01840000e+00)
( 4.03010000e-01,  1.01560000e+00)
( 4.09700000e-01,  1.01290000e+00)
( 4.16390000e-01,  1.01020000e+00)
( 4.23080000e-01,  1.00770000e+00)
( 4.29770000e-01,  1.00520000e+00)
( 4.36450000e-01,  1.00280000e+00)
( 4.43140000e-01,  1.00050000e+00)
( 4.49830000e-01,  9.98210000e-01)
( 4.56520000e-01,  9.96030000e-01)
( 4.63210000e-01,  9.93910000e-01)
( 4.69900000e-01,  9.91850000e-01)
( 4.76590000e-01,  9.89870000e-01)
( 4.83280000e-01,  9.87950000e-01)
( 4.89970000e-01,  9.86090000e-01)
( 4.96660000e-01,  9.84300000e-01)
( 5.03340000e-01,  9.82570000e-01)
( 5.10030000e-01,  9.80010000e-01)
( 5.16720000e-01,  9.78410000e-01)
( 5.23410000e-01,  9.76850000e-01)
( 5.30100000e-01,  9.75320000e-01)
( 5.36790000e-01,  9.73820000e-01)
( 5.43480000e-01,  9.72350000e-01)
( 5.50170000e-01,  9.70910000e-01)
( 5.56860000e-01,  9.69490000e-01)
( 5.63550000e-01,  9.68110000e-01)
( 5.70230000e-01,  9.66750000e-01)
( 5.76920000e-01,  9.65420000e-01)
( 5.83610000e-01,  9.64120000e-01)
( 5.90300000e-01,  9.62840000e-01)
( 5.96990000e-01,  9.61590000e-01)
( 6.03680000e-01,  9.60360000e-01)
( 6.10370000e-01,  9.59160000e-01)
( 6.17060000e-01,  9.57980000e-01)
( 6.23750000e-01,  9.56820000e-01)
( 6.30430000e-01,  9.55690000e-01)
( 6.37120000e-01,  9.54580000e-01)
( 6.43810000e-01,  9.53490000e-01)
( 6.50500000e-01,  9.52420000e-01)
( 6.57190000e-01,  9.51370000e-01)
( 6.63880000e-01,  9.50350000e-01)
( 6.70570000e-01,  9.49340000e-01)
( 6.77260000e-01,  9.48360000e-01)
( 6.83950000e-01,  9.47390000e-01)
( 6.90640000e-01,  9.46440000e-01)
( 6.97320000e-01,  9.45510000e-01)
( 7.04010000e-01,  9.44720000e-01)
( 7.10700000e-01,  9.43840000e-01)
( 7.17390000e-01,  9.42970000e-01)
( 7.24080000e-01,  9.42110000e-01)
( 7.30770000e-01,  9.41260000e-01)
( 7.37460000e-01,  9.40410000e-01)
( 7.44150000e-01,  9.39580000e-01)
( 7.50840000e-01,  9.38760000e-01)
( 7.57530000e-01,  9.37950000e-01)
( 7.64210000e-01,  9.37160000e-01)
( 7.70900000e-01,  9.36380000e-01)
( 7.77590000e-01,  9.35620000e-01)
( 7.84280000e-01,  9.34880000e-01)
( 7.90970000e-01,  9.34150000e-01)
( 7.97660000e-01,  9.33450000e-01)
( 8.04350000e-01,  9.32760000e-01)
( 8.11040000e-01,  9.32100000e-01)
( 8.17730000e-01,  1.46690000e+00)
( 8.24410000e-01,  1.46800000e+00)
( 8.31100000e-01,  1.46900000e+00)
( 8.37790000e-01,  1.47020000e+00)
( 8.44480000e-01,  1.47140000e+00)
( 8.51170000e-01,  1.47260000e+00)
( 8.57860000e-01,  1.47390000e+00)
( 8.64550000e-01,  1.47520000e+00)
( 8.71240000e-01,  1.47670000e+00)
( 8.77930000e-01,  1.47810000e+00)
( 8.84620000e-01,  1.47970000e+00)
( 8.91300000e-01,  1.48140000e+00)
( 8.97990000e-01,  1.48300000e+00)
( 9.04680000e-01,  1.48470000e+00)
( 9.11370000e-01,  1.48650000e+00)
( 9.18060000e-01,  1.48830000e+00)
( 9.24750000e-01,  1.48990000e+00)
( 9.31440000e-01,  1.49170000e+00)
( 9.38130000e-01,  1.49350000e+00)
( 9.44820000e-01,  1.49520000e+00)
( 9.51510000e-01,  1.49680000e+00)
( 9.58190000e-01,  1.49830000e+00)
( 9.64880000e-01,  1.49970000e+00)
( 9.71570000e-01,  1.50090000e+00)
( 9.78260000e-01,  1.50200000e+00)
( 9.84950000e-01,  1.50290000e+00)
( 9.91640000e-01,  1.50360000e+00)
( 9.98330000e-01,  1.50420000e+00)
( 1.00500000e+00,  1.50450000e+00)
( 1.01170000e+00,  1.50430000e+00)
( 1.01840000e+00,  1.50400000e+00)
( 1.02510000e+00,  1.50370000e+00)
( 1.03180000e+00,  1.50320000e+00)
( 1.03850000e+00,  1.50200000e+00)
( 1.04520000e+00,  1.50090000e+00)
( 1.05180000e+00,  1.49970000e+00)
( 1.05850000e+00,  1.49840000e+00)
( 1.06520000e+00,  1.49690000e+00)
( 1.07190000e+00,  1.49530000e+00)
( 1.07860000e+00,  1.49350000e+00)
( 1.08530000e+00,  1.49170000e+00)
( 1.09200000e+00,  1.48990000e+00)
( 1.09870000e+00,  1.48810000e+00)
( 1.10540000e+00,  1.48620000e+00)
( 1.11200000e+00,  1.48430000e+00)
( 1.11870000e+00,  1.48240000e+00)
( 1.12540000e+00,  1.48040000e+00)
( 1.13210000e+00,  1.47850000e+00)
( 1.13880000e+00,  1.47660000e+00)
( 1.14550000e+00,  1.47470000e+00)
( 1.15220000e+00,  1.47280000e+00)
( 1.15890000e+00,  1.47090000e+00)
( 1.16560000e+00,  1.46900000e+00)
( 1.17220000e+00,  1.46720000e+00)
( 1.17890000e+00,  1.46540000e+00)
( 1.18560000e+00,  1.46350000e+00)
( 1.19230000e+00,  1.46170000e+00)
( 1.19900000e+00,  1.45990000e+00)
( 1.20570000e+00,  1.45820000e+00)
( 1.21240000e+00,  1.45640000e+00)
( 1.21910000e+00,  1.45480000e+00)
( 1.22580000e+00,  1.45310000e+00)
( 1.23240000e+00,  1.45140000e+00)
( 1.23910000e+00,  1.44970000e+00)
( 1.24580000e+00,  1.44800000e+00)
( 1.25250000e+00,  1.44640000e+00)
( 1.25920000e+00,  1.44470000e+00)
( 1.26590000e+00,  1.44310000e+00)
( 1.27260000e+00,  1.44150000e+00)
( 1.27930000e+00,  1.43990000e+00)
( 1.28600000e+00,  1.43830000e+00)
( 1.29260000e+00,  1.43680000e+00)
( 1.29930000e+00,  1.43520000e+00)
( 1.30600000e+00,  1.43370000e+00)
( 1.31270000e+00,  1.43210000e+00)
( 1.31940000e+00,  1.43060000e+00)
( 1.32610000e+00,  1.42910000e+00)
( 1.33280000e+00,  1.42760000e+00)
( 1.33950000e+00,  1.42610000e+00)
( 1.34620000e+00,  1.42460000e+00)
( 1.35280000e+00,  1.42310000e+00)
( 1.35950000e+00,  1.42160000e+00)
( 1.36620000e+00,  1.42020000e+00)
( 1.37290000e+00,  1.41870000e+00)
( 1.37960000e+00,  1.41720000e+00)
( 1.38630000e+00,  1.41570000e+00)
( 1.39300000e+00,  1.41420000e+00)
( 1.39970000e+00,  1.41270000e+00)
( 1.40640000e+00,  1.41120000e+00)
( 1.41300000e+00,  1.40970000e+00)
( 1.41970000e+00,  1.40820000e+00)
( 1.42640000e+00,  1.40670000e+00)
( 1.43310000e+00,  1.40520000e+00)
( 1.43980000e+00,  1.40360000e+00)
( 1.44650000e+00,  1.40210000e+00)
( 1.45320000e+00,  1.40050000e+00)
( 1.45990000e+00,  1.39900000e+00)
( 1.46660000e+00,  1.39740000e+00)
( 1.47320000e+00,  1.39580000e+00)
( 1.47990000e+00,  1.39420000e+00)
( 1.48660000e+00,  1.39260000e+00)
( 1.49330000e+00,  1.39100000e+00)};

\nextgroupplot[width=0.3\textwidth, xlabel={$x$}, ymax=1.4, xmax=1.5, ylabel={$M$}, xmin=-0.5, ymin=0.6, height=0.3\textwidth]
\addplot [solid, thick, color=black, forget plot]
coordinates {
(-5.00000000e-01,          nan)
(-4.93310000e-01,  7.47695831e-01)
(-4.86620000e-01,  7.47580266e-01)
(-4.79930000e-01,  7.47525970e-01)
(-4.73240000e-01,  7.47397553e-01)
(-4.66560000e-01,  7.47256122e-01)
(-4.59870000e-01,  7.47073454e-01)
(-4.53180000e-01,  7.46891025e-01)
(-4.46490000e-01,  7.46606406e-01)
(-4.39800000e-01,  7.46295762e-01)
(-4.33110000e-01,  7.45985487e-01)
(-4.26420000e-01,  7.45647453e-01)
(-4.19730000e-01,  7.45283568e-01)
(-4.13040000e-01,  7.44846113e-01)
(-4.06350000e-01,  7.44381187e-01)
(-3.99670000e-01,  7.43890719e-01)
(-3.92980000e-01,  7.43400996e-01)
(-3.86290000e-01,  7.42797122e-01)
(-3.79600000e-01,  7.42166332e-01)
(-3.72910000e-01,  7.41538449e-01)
(-3.66220000e-01,  7.40810153e-01)
(-3.59530000e-01,  7.40129104e-01)
(-3.52840000e-01,  7.39349895e-01)
(-3.46150000e-01,  7.38572323e-01)
(-3.39460000e-01,  7.37682437e-01)
(-3.32780000e-01,  7.36781679e-01)
(-3.26090000e-01,  7.35883041e-01)
(-3.19400000e-01,  7.34986898e-01)
(-3.12710000e-01,  7.33939650e-01)
(-3.06020000e-01,  7.32922672e-01)
(-2.99330000e-01,  7.31856233e-01)
(-2.92640000e-01,  7.30692377e-01)
(-2.85950000e-01,  7.29641447e-01)
(-2.79260000e-01,  7.28360371e-01)
(-2.72580000e-01,  7.27071226e-01)
(-2.65890000e-01,  7.25714485e-01)
(-2.59200000e-01,  7.24383563e-01)
(-2.52510000e-01,  7.22986146e-01)
(-2.45820000e-01,  7.21483924e-01)
(-2.39130000e-01,  7.19976419e-01)
(-2.32440000e-01,  7.18450183e-01)
(-2.25750000e-01,  7.16907595e-01)
(-2.19060000e-01,  7.15302396e-01)
(-2.12370000e-01,  7.13609513e-01)
(-2.05690000e-01,  7.11903187e-01)
(-1.99000000e-01,  7.10182782e-01)
(-1.92310000e-01,  7.08488792e-01)
(-1.85620000e-01,  7.06686189e-01)
(-1.78930000e-01,  7.04841609e-01)
(-1.72240000e-01,  7.03129951e-01)
(-1.65550000e-01,  7.01347161e-01)
(-1.58860000e-01,  6.99559931e-01)
(-1.52170000e-01,  6.97747720e-01)
(-1.45480000e-01,  6.96029704e-01)
(-1.38800000e-01,  6.94399893e-01)
(-1.32110000e-01,  6.92834629e-01)
(-1.25420000e-01,  6.91444502e-01)
(-1.18730000e-01,  6.90150137e-01)
(-1.12040000e-01,  6.88997417e-01)
(-1.05350000e-01,  6.88091946e-01)
(-9.86620000e-02,  6.87441618e-01)
(-9.19730000e-02,  6.87070930e-01)
(-8.52840000e-02,  6.87179893e-01)
(-7.85950000e-02,  6.87668937e-01)
(-7.19060000e-02,  6.88657482e-01)
(-6.52170000e-02,  6.90368451e-01)
(-5.85280000e-02,  6.92721069e-01)
(-5.18390000e-02,  6.95837128e-01)
(-4.51510000e-02,  6.99812847e-01)
(-3.84620000e-02,  7.04418084e-01)
(-3.17730000e-02,  7.10076487e-01)
(-2.50840000e-02,  7.16858652e-01)
(-1.83950000e-02,  7.24626605e-01)
(-1.17060000e-02,  7.33357950e-01)
(-5.01670000e-03,  7.43306770e-01)
( 1.67220000e-03,  7.54317172e-01)
( 8.36120000e-03,  7.66089109e-01)
( 1.50500000e-02,  7.78667002e-01)
( 2.17390000e-02,  7.91917053e-01)
( 2.84280000e-02,  8.05750761e-01)
( 3.51170000e-02,  8.20037449e-01)
( 4.18060000e-02,  8.34784755e-01)
( 4.84950000e-02,  8.49764774e-01)
( 5.51840000e-02,  8.64784854e-01)
( 6.18730000e-02,  8.79602443e-01)
( 6.85620000e-02,  8.94110767e-01)
( 7.52510000e-02,  9.08686517e-01)
( 8.19400000e-02,  9.22815909e-01)
( 8.86290000e-02,  9.36532746e-01)
( 9.53180000e-02,  9.49922348e-01)
( 1.02010000e-01,  9.62856614e-01)
( 1.08700000e-01,  9.75477182e-01)
( 1.15380000e-01,  9.87751414e-01)
( 1.22070000e-01,  9.99576636e-01)
( 1.28760000e-01,  1.01089877e+00)
( 1.35450000e-01,  1.02173595e+00)
( 1.42140000e-01,  1.03216096e+00)
( 1.48830000e-01,  1.04187550e+00)
( 1.55520000e-01,  1.05111126e+00)
( 1.62210000e-01,  1.05966524e+00)
( 1.68900000e-01,  1.06856095e+00)
( 1.75590000e-01,  1.07695413e+00)
( 1.82270000e-01,  1.08492005e+00)
( 1.88960000e-01,  1.09281988e+00)
( 1.95650000e-01,  1.10034125e+00)
( 2.02340000e-01,  1.10770402e+00)
( 2.09030000e-01,  1.11477214e+00)
( 2.15720000e-01,  1.12171809e+00)
( 2.22410000e-01,  1.12840612e+00)
( 2.29100000e-01,  1.13468404e+00)
( 2.35790000e-01,  1.14076899e+00)
( 2.42470000e-01,  1.14671064e+00)
( 2.49160000e-01,  1.15239462e+00)
( 2.55850000e-01,  1.15781325e+00)
( 2.62540000e-01,  1.16315328e+00)
( 2.69230000e-01,  1.16824047e+00)
( 2.75920000e-01,  1.17320651e+00)
( 2.82610000e-01,  1.17800287e+00)
( 2.89300000e-01,  1.18269806e+00)
( 2.95990000e-01,  1.18727833e+00)
( 3.02680000e-01,  1.19159957e+00)
( 3.09360000e-01,  1.19587924e+00)
( 3.16050000e-01,  1.20005147e+00)
( 3.22740000e-01,  1.20403329e+00)
( 3.29430000e-01,  1.20783530e+00)
( 3.36120000e-01,  1.21167580e+00)
( 3.42810000e-01,  1.21537316e+00)
( 3.49500000e-01,  1.21899647e+00)
( 3.56190000e-01,  1.22235968e+00)
( 3.62880000e-01,  1.22571163e+00)
( 3.69570000e-01,  1.22901665e+00)
( 3.76250000e-01,  1.23230630e+00)
( 3.82940000e-01,  1.23544446e+00)
( 3.89630000e-01,  1.23851420e+00)
( 3.96320000e-01,  1.24142617e+00)
( 4.03010000e-01,  1.24431250e+00)
( 4.09700000e-01,  1.24713770e+00)
( 4.16390000e-01,  1.24987043e+00)
( 4.23080000e-01,  1.25260109e+00)
( 4.29770000e-01,  1.25508213e+00)
( 4.36450000e-01,  1.25764680e+00)
( 4.43140000e-01,  1.26001999e+00)
( 4.49830000e-01,  1.26240188e+00)
( 4.56520000e-01,  1.26477427e+00)
( 4.63210000e-01,  1.26702600e+00)
( 4.69900000e-01,  1.26903941e+00)
( 4.76590000e-01,  1.27126386e+00)
( 4.83280000e-01,  1.27324669e+00)
( 4.89970000e-01,  1.27510115e+00)
( 4.96660000e-01,  1.27701637e+00)
( 5.03340000e-01,  1.27884255e+00)
( 5.10030000e-01,  1.28159177e+00)
( 5.16720000e-01,  1.28334970e+00)
( 5.23410000e-01,  1.28487493e+00)
( 5.30100000e-01,  1.28653972e+00)
( 5.36790000e-01,  1.28814052e+00)
( 5.43480000e-01,  1.28971971e+00)
( 5.50170000e-01,  1.29123385e+00)
( 5.56860000e-01,  1.29273503e+00)
( 5.63550000e-01,  1.29416001e+00)
( 5.70230000e-01,  1.29569126e+00)
( 5.76920000e-01,  1.29703530e+00)
( 5.83610000e-01,  1.29847537e+00)
( 5.90300000e-01,  1.29973664e+00)
( 5.96990000e-01,  1.30109357e+00)
( 6.03680000e-01,  1.30243551e+00)
( 6.10370000e-01,  1.30370786e+00)
( 6.17060000e-01,  1.30496463e+00)
( 6.23750000e-01,  1.30620555e+00)
( 6.30430000e-01,  1.30754209e+00)
( 6.37120000e-01,  1.30865139e+00)
( 6.43810000e-01,  1.30974373e+00)
( 6.50500000e-01,  1.31098641e+00)
( 6.57190000e-01,  1.31216749e+00)
( 6.63880000e-01,  1.31315320e+00)
( 6.70570000e-01,  1.31425463e+00)
( 6.77260000e-01,  1.31532887e+00)
( 6.83950000e-01,  1.31639391e+00)
( 6.90640000e-01,  1.31739770e+00)
( 6.97320000e-01,  1.31838244e+00)
( 7.04010000e-01,  1.31927420e+00)
( 7.10700000e-01,  1.32016540e+00)
( 7.17390000e-01,  1.32121879e+00)
( 7.24080000e-01,  1.32204666e+00)
( 7.30770000e-01,  1.32303681e+00)
( 7.37460000e-01,  1.32390322e+00)
( 7.44150000e-01,  1.32487723e+00)
( 7.50840000e-01,  1.32562200e+00)
( 7.57530000e-01,  1.32657754e+00)
( 7.64210000e-01,  1.32746756e+00)
( 7.70900000e-01,  1.32817683e+00)
( 7.77590000e-01,  1.32899176e+00)
( 7.84280000e-01,  1.32978696e+00)
( 7.90970000e-01,  1.33057291e+00)
( 7.97660000e-01,  1.33145592e+00)
( 8.04350000e-01,  1.33215662e+00)
( 8.11040000e-01,  1.33282804e+00)
( 8.17730000e-01,  7.69042273e-01)
( 8.24410000e-01,  7.68044786e-01)
( 8.31100000e-01,  7.67079995e-01)
( 8.37790000e-01,  7.65950825e-01)
( 8.44480000e-01,  7.64824413e-01)
( 8.51170000e-01,  7.63685014e-01)
( 8.57860000e-01,  7.62423815e-01)
( 8.64550000e-01,  7.61149508e-01)
( 8.71240000e-01,  7.59800056e-01)
( 8.77930000e-01,  7.58387974e-01)
( 8.84620000e-01,  7.56809913e-01)
( 8.91300000e-01,  7.55248568e-01)
( 8.97990000e-01,  7.53650117e-01)
( 9.04680000e-01,  7.52010393e-01)
( 9.11370000e-01,  7.50305786e-01)
( 9.18060000e-01,  7.48585474e-01)
( 9.24750000e-01,  7.46893976e-01)
( 9.31440000e-01,  7.45109153e-01)
( 9.38130000e-01,  7.43333689e-01)
( 9.44820000e-01,  7.41672195e-01)
( 9.51510000e-01,  7.40125472e-01)
( 9.58190000e-01,  7.38645124e-01)
( 9.64880000e-01,  7.37293893e-01)
( 9.71570000e-01,  7.36100536e-01)
( 9.78260000e-01,  7.34975276e-01)
( 9.84950000e-01,  7.34096830e-01)
( 9.91640000e-01,  7.33465380e-01)
( 9.98330000e-01,  7.32916692e-01)
( 1.00500000e+00,  7.32644030e-01)
( 1.01170000e+00,  7.32719704e-01)
( 1.01840000e+00,  7.33126149e-01)
( 1.02510000e+00,  7.33648122e-01)
( 1.03180000e+00,  7.34208158e-01)
( 1.03850000e+00,  7.35345411e-01)
( 1.04520000e+00,  7.36465761e-01)
( 1.05180000e+00,  7.37669694e-01)
( 1.05850000e+00,  7.39093748e-01)
( 1.06520000e+00,  7.40628658e-01)
( 1.07190000e+00,  7.42232218e-01)
( 1.07860000e+00,  7.44044693e-01)
( 1.08530000e+00,  7.45768827e-01)
( 1.09200000e+00,  7.47607867e-01)
( 1.09870000e+00,  7.49458669e-01)
( 1.10540000e+00,  7.51350082e-01)
( 1.11200000e+00,  7.53280229e-01)
( 1.11870000e+00,  7.55130867e-01)
( 1.12540000e+00,  7.57112496e-01)
( 1.13210000e+00,  7.58984057e-01)
( 1.13880000e+00,  7.60955371e-01)
( 1.14550000e+00,  7.62859959e-01)
( 1.15220000e+00,  7.64759096e-01)
( 1.15890000e+00,  7.66574540e-01)
( 1.16560000e+00,  7.68490193e-01)
( 1.17220000e+00,  7.70290867e-01)
( 1.17890000e+00,  7.72098928e-01)
( 1.18560000e+00,  7.73945172e-01)
( 1.19230000e+00,  7.75753472e-01)
( 1.19900000e+00,  7.77489850e-01)
( 1.20570000e+00,  7.79280800e-01)
( 1.21240000e+00,  7.81015860e-01)
( 1.21910000e+00,  7.82822779e-01)
( 1.22580000e+00,  7.84459415e-01)
( 1.23240000e+00,  7.86087004e-01)
( 1.23910000e+00,  7.87735063e-01)
( 1.24580000e+00,  7.89373909e-01)
( 1.25250000e+00,  7.90890930e-01)
( 1.25920000e+00,  7.92540507e-01)
( 1.26590000e+00,  7.94148153e-01)
( 1.27260000e+00,  7.95680153e-01)
( 1.27930000e+00,  7.97201789e-01)
( 1.28600000e+00,  7.98824387e-01)
( 1.29260000e+00,  8.00307706e-01)
( 1.29930000e+00,  8.01843532e-01)
( 1.30600000e+00,  8.03351040e-01)
( 1.31270000e+00,  8.04880866e-01)
( 1.31940000e+00,  8.06397534e-01)
( 1.32610000e+00,  8.07903231e-01)
( 1.33280000e+00,  8.09331257e-01)
( 1.33950000e+00,  8.10830660e-01)
( 1.34620000e+00,  8.12267145e-01)
( 1.35280000e+00,  8.13790402e-01)
( 1.35950000e+00,  8.15235431e-01)
( 1.36620000e+00,  8.16634719e-01)
( 1.37290000e+00,  8.18088044e-01)
( 1.37960000e+00,  8.19545623e-01)
( 1.38630000e+00,  8.21007482e-01)
( 1.39300000e+00,  8.22457577e-01)
( 1.39970000e+00,  8.23927974e-01)
( 1.40640000e+00,  8.25402702e-01)
( 1.41300000e+00,  8.26881769e-01)
( 1.41970000e+00,  8.28381467e-01)
( 1.42640000e+00,  8.29869335e-01)
( 1.43310000e+00,  8.31377991e-01)
( 1.43980000e+00,  8.32910181e-01)
( 1.44650000e+00,  8.34343240e-01)
( 1.45320000e+00,  8.35901154e-01)
( 1.45990000e+00,  8.37428120e-01)
( 1.46660000e+00,  8.39012310e-01)
( 1.47320000e+00,  8.40516100e-01)
( 1.47990000e+00,  8.42110205e-01)
( 1.48660000e+00,  8.43709396e-01)
( 1.49330000e+00,  8.45227768e-01)};

\nextgroupplot[width=0.3\textwidth, xlabel={$x$}, ymax=1.3, xmax=1.5, ylabel={$P$}, xmin=-0.5, ymin=0.5, height=0.3\textwidth]
\addplot [solid, thick, color=black, forget plot]
coordinates {
(-5.00000000e-01,          nan)
(-4.93310000e-01,  1.13141533e+00)
(-4.86620000e-01,  1.13149725e+00)
(-4.79930000e-01,  1.13158867e+00)
(-4.73240000e-01,  1.13171045e+00)
(-4.66560000e-01,  1.13187215e+00)
(-4.59870000e-01,  1.13208534e+00)
(-4.53180000e-01,  1.13229845e+00)
(-4.46490000e-01,  1.13255342e+00)
(-4.39800000e-01,  1.13288827e+00)
(-4.33110000e-01,  1.13322299e+00)
(-4.26420000e-01,  1.13356910e+00)
(-4.19730000e-01,  1.13399505e+00)
(-4.13040000e-01,  1.13445117e+00)
(-4.06350000e-01,  1.13491859e+00)
(-3.99670000e-01,  1.13546577e+00)
(-3.92980000e-01,  1.13601270e+00)
(-3.86290000e-01,  1.13664109e+00)
(-3.79600000e-01,  1.13728058e+00)
(-3.72910000e-01,  1.13798827e+00)
(-3.66220000e-01,  1.13873719e+00)
(-3.59530000e-01,  1.13946686e+00)
(-3.52840000e-01,  1.14030620e+00)
(-3.46150000e-01,  1.14114500e+00)
(-3.39460000e-01,  1.14206463e+00)
(-3.32780000e-01,  1.14302358e+00)
(-3.26090000e-01,  1.14398183e+00)
(-3.19400000e-01,  1.14493923e+00)
(-3.12710000e-01,  1.14602816e+00)
(-3.06020000e-01,  1.14710489e+00)
(-2.99330000e-01,  1.14827164e+00)
(-2.92640000e-01,  1.14947851e+00)
(-2.85950000e-01,  1.15074384e+00)
(-2.79260000e-01,  1.15206885e+00)
(-2.72580000e-01,  1.15343228e+00)
(-2.65890000e-01,  1.15482367e+00)
(-2.59200000e-01,  1.15627455e+00)
(-2.52510000e-01,  1.15775300e+00)
(-2.45820000e-01,  1.15930955e+00)
(-2.39130000e-01,  1.16090356e+00)
(-2.32440000e-01,  1.16250563e+00)
(-2.25750000e-01,  1.16418472e+00)
(-2.19060000e-01,  1.16588984e+00)
(-2.12370000e-01,  1.16763121e+00)
(-2.05690000e-01,  1.16944856e+00)
(-1.99000000e-01,  1.17127215e+00)
(-1.92310000e-01,  1.17312051e+00)
(-1.85620000e-01,  1.17501372e+00)
(-1.78930000e-01,  1.17695965e+00)
(-1.72240000e-01,  1.17878163e+00)
(-1.65550000e-01,  1.18063307e+00)
(-1.58860000e-01,  1.18255710e+00)
(-1.52170000e-01,  1.18442142e+00)
(-1.45480000e-01,  1.18624779e+00)
(-1.38800000e-01,  1.18800345e+00)
(-1.32110000e-01,  1.18962747e+00)
(-1.25420000e-01,  1.19110999e+00)
(-1.18730000e-01,  1.19244796e+00)
(-1.12040000e-01,  1.19369324e+00)
(-1.05350000e-01,  1.19466162e+00)
(-9.86620000e-02,  1.19531447e+00)
(-9.19730000e-02,  1.19571240e+00)
(-8.52840000e-02,  1.19563279e+00)
(-7.85950000e-02,  1.19511611e+00)
(-7.19060000e-02,  1.19407860e+00)
(-6.52170000e-02,  1.19232389e+00)
(-5.85280000e-02,  1.18984998e+00)
(-5.18390000e-02,  1.18660799e+00)
(-4.51510000e-02,  1.18266510e+00)
(-3.84620000e-02,  1.17772939e+00)
(-3.17730000e-02,  1.17166605e+00)
(-2.50840000e-02,  1.16437955e+00)
(-1.83950000e-02,  1.15614028e+00)
(-1.17060000e-02,  1.14692673e+00)
(-5.01670000e-03,  1.13629378e+00)
( 1.67220000e-03,  1.12466081e+00)
( 8.36120000e-03,  1.11212466e+00)
( 1.50500000e-02,  1.09871148e+00)
( 2.17390000e-02,  1.08455521e+00)
( 2.84280000e-02,  1.06977103e+00)
( 3.51170000e-02,  1.05449078e+00)
( 4.18060000e-02,  1.03878903e+00)
( 4.84950000e-02,  1.02283499e+00)
( 5.51840000e-02,  1.00674815e+00)
( 6.18730000e-02,  9.91290635e-01)
( 6.85620000e-02,  9.75976026e-01)
( 7.52510000e-02,  9.60656947e-01)
( 8.19400000e-02,  9.45895879e-01)
( 8.86290000e-02,  9.31609548e-01)
( 9.53180000e-02,  9.17791900e-01)
( 1.02010000e-01,  9.04424152e-01)
( 1.08700000e-01,  8.91446826e-01)
( 1.15380000e-01,  8.78931133e-01)
( 1.22070000e-01,  8.66941792e-01)
( 1.28760000e-01,  8.55523455e-01)
( 1.35450000e-01,  8.44661665e-01)
( 1.42140000e-01,  8.34321968e-01)
( 1.48830000e-01,  8.24669041e-01)
( 1.55520000e-01,  8.15608454e-01)
( 1.62210000e-01,  8.07206692e-01)
( 1.68900000e-01,  7.98508581e-01)
( 1.75590000e-01,  7.90396724e-01)
( 1.82270000e-01,  7.82656450e-01)
( 1.88960000e-01,  7.75124853e-01)
( 1.95650000e-01,  7.67910483e-01)
( 2.02340000e-01,  7.60932339e-01)
( 2.09030000e-01,  7.54235616e-01)
( 2.15720000e-01,  7.47757313e-01)
( 2.22410000e-01,  7.41542297e-01)
( 2.29100000e-01,  7.35640117e-01)
( 2.35790000e-01,  7.30002394e-01)
( 2.42470000e-01,  7.24551307e-01)
( 2.49160000e-01,  7.19353065e-01)
( 2.55850000e-01,  7.14349684e-01)
( 2.62540000e-01,  7.09534675e-01)
( 2.69230000e-01,  7.04877930e-01)
( 2.75920000e-01,  7.00391910e-01)
( 2.82610000e-01,  6.96062528e-01)
( 2.89300000e-01,  6.91865925e-01)
( 2.95990000e-01,  6.87748267e-01)
( 3.02680000e-01,  6.83870436e-01)
( 3.09360000e-01,  6.80102675e-01)
( 3.16050000e-01,  6.76437147e-01)
( 3.22740000e-01,  6.72928341e-01)
( 3.29430000e-01,  6.69543279e-01)
( 3.36120000e-01,  6.66211078e-01)
( 3.42810000e-01,  6.63017848e-01)
( 3.49500000e-01,  6.59884183e-01)
( 3.56190000e-01,  6.56896881e-01)
( 3.62880000e-01,  6.54004356e-01)
( 3.69570000e-01,  6.51133074e-01)
( 3.76250000e-01,  6.48357184e-01)
( 3.82940000e-01,  6.45662801e-01)
( 3.89630000e-01,  6.43051421e-01)
( 3.96320000e-01,  6.40522453e-01)
( 4.03010000e-01,  6.38090167e-01)
( 4.09700000e-01,  6.35681531e-01)
( 4.16390000e-01,  6.33389527e-01)
( 4.23080000e-01,  6.31102413e-01)
( 4.29770000e-01,  6.28979382e-01)
( 4.36450000e-01,  6.26834395e-01)
( 4.43140000e-01,  6.24834677e-01)
( 4.49830000e-01,  6.22835580e-01)
( 4.56520000e-01,  6.20897902e-01)
( 4.63210000e-01,  6.19055118e-01)
( 4.69900000e-01,  6.17314800e-01)
( 4.76590000e-01,  6.15540819e-01)
( 4.83280000e-01,  6.13869777e-01)
( 4.89970000e-01,  6.12294589e-01)
( 4.96660000e-01,  6.10730697e-01)
( 5.03340000e-01,  6.09222966e-01)
( 5.10030000e-01,  6.06975703e-01)
( 5.16720000e-01,  6.05575227e-01)
( 5.23410000e-01,  6.04272993e-01)
( 5.30100000e-01,  6.02930751e-01)
( 5.36790000e-01,  6.01636340e-01)
( 5.43480000e-01,  6.00349771e-01)
( 5.50170000e-01,  5.99111111e-01)
( 5.56860000e-01,  5.97877598e-01)
( 5.63550000e-01,  5.96695058e-01)
( 5.70230000e-01,  5.95509690e-01)
( 5.76920000e-01,  5.94380525e-01)
( 5.83610000e-01,  5.93251480e-01)
( 5.90300000e-01,  5.92175893e-01)
( 5.96990000e-01,  5.91100471e-01)
( 6.03680000e-01,  5.90030455e-01)
( 6.10370000e-01,  5.89008829e-01)
( 6.17060000e-01,  5.87992662e-01)
( 6.23750000e-01,  5.86982011e-01)
( 6.30430000e-01,  5.85971562e-01)
( 6.37120000e-01,  5.85054941e-01)
( 6.43810000e-01,  5.84143961e-01)
( 6.50500000e-01,  5.83190248e-01)
( 6.57190000e-01,  5.82282140e-01)
( 6.63880000e-01,  5.81431119e-01)
( 6.70570000e-01,  5.80574350e-01)
( 6.77260000e-01,  5.79726074e-01)
( 6.83950000e-01,  5.78881092e-01)
( 6.90640000e-01,  5.78081254e-01)
( 6.97320000e-01,  5.77287455e-01)
( 7.04010000e-01,  5.76601526e-01)
( 7.10700000e-01,  5.75862433e-01)
( 7.17390000e-01,  5.75077465e-01)
( 7.24080000e-01,  5.74384320e-01)
( 7.30770000e-01,  5.73645311e-01)
( 7.37460000e-01,  5.72915242e-01)
( 7.44150000e-01,  5.72181834e-01)
( 7.50840000e-01,  5.71540928e-01)
( 7.57530000e-01,  5.70813700e-01)
( 7.64210000e-01,  5.70132348e-01)
( 7.70900000e-01,  5.69502676e-01)
( 7.77590000e-01,  5.68870175e-01)
( 7.84280000e-01,  5.68243633e-01)
( 7.90970000e-01,  5.67620065e-01)
( 7.97660000e-01,  5.66996166e-01)
( 8.04350000e-01,  5.66424113e-01)
( 8.11040000e-01,  5.65860613e-01)
( 8.17730000e-01,  1.07927314e+00)
( 8.24410000e-01,  1.08035739e+00)
( 8.31100000e-01,  1.08142888e+00)
( 8.37790000e-01,  1.08263493e+00)
( 8.44480000e-01,  1.08384021e+00)
( 8.51170000e-01,  1.08508540e+00)
( 8.57860000e-01,  1.08641394e+00)
( 8.64550000e-01,  1.08778253e+00)
( 8.71240000e-01,  1.08921596e+00)
( 8.77930000e-01,  1.09070929e+00)
( 8.84620000e-01,  1.09233827e+00)
( 8.91300000e-01,  1.09403014e+00)
( 8.97990000e-01,  1.09570293e+00)
( 9.04680000e-01,  1.09742733e+00)
( 9.11370000e-01,  1.09924734e+00)
( 9.18060000e-01,  1.10110890e+00)
( 9.24750000e-01,  1.10282133e+00)
( 9.31440000e-01,  1.10467540e+00)
( 9.38130000e-01,  1.10649179e+00)
( 9.44820000e-01,  1.10822869e+00)
( 9.51510000e-01,  1.10988573e+00)
( 9.58190000e-01,  1.11141275e+00)
( 9.64880000e-01,  1.11281929e+00)
( 9.71570000e-01,  1.11409379e+00)
( 9.78260000e-01,  1.11523722e+00)
( 9.84950000e-01,  1.11617821e+00)
( 9.91640000e-01,  1.11691661e+00)
( 9.98330000e-01,  1.11748304e+00)
( 1.00500000e+00,  1.11783515e+00)
( 1.01170000e+00,  1.11756157e+00)
( 1.01840000e+00,  1.11729522e+00)
( 1.02510000e+00,  1.11691326e+00)
( 1.03180000e+00,  1.11641186e+00)
( 1.03850000e+00,  1.11512437e+00)
( 1.04520000e+00,  1.11401759e+00)
( 1.05180000e+00,  1.11277355e+00)
( 1.05850000e+00,  1.11137259e+00)
( 1.06520000e+00,  1.10982370e+00)
( 1.07190000e+00,  1.10817855e+00)
( 1.07860000e+00,  1.10627738e+00)
( 1.08530000e+00,  1.10448141e+00)
( 1.09200000e+00,  1.10257104e+00)
( 1.09870000e+00,  1.10065672e+00)
( 1.10540000e+00,  1.09872693e+00)
( 1.11200000e+00,  1.09671340e+00)
( 1.11870000e+00,  1.09476695e+00)
( 1.12540000e+00,  1.09273456e+00)
( 1.13210000e+00,  1.09078142e+00)
( 1.13880000e+00,  1.08875488e+00)
( 1.14550000e+00,  1.08675586e+00)
( 1.15220000e+00,  1.08479429e+00)
( 1.15890000e+00,  1.08290133e+00)
( 1.16560000e+00,  1.08093487e+00)
( 1.17220000e+00,  1.07904948e+00)
( 1.17890000e+00,  1.07716197e+00)
( 1.18560000e+00,  1.07526005e+00)
( 1.19230000e+00,  1.07340837e+00)
( 1.19900000e+00,  1.07158618e+00)
( 1.20570000e+00,  1.06974304e+00)
( 1.21240000e+00,  1.06795717e+00)
( 1.21910000e+00,  1.06621170e+00)
( 1.22580000e+00,  1.06446661e+00)
( 1.23240000e+00,  1.06276001e+00)
( 1.23910000e+00,  1.06101191e+00)
( 1.24580000e+00,  1.05930233e+00)
( 1.25250000e+00,  1.05767610e+00)
( 1.25920000e+00,  1.05596372e+00)
( 1.26590000e+00,  1.05430282e+00)
( 1.27260000e+00,  1.05267277e+00)
( 1.27930000e+00,  1.05108152e+00)
( 1.28600000e+00,  1.04941680e+00)
( 1.29260000e+00,  1.04787612e+00)
( 1.29930000e+00,  1.04628127e+00)
( 1.30600000e+00,  1.04469839e+00)
( 1.31270000e+00,  1.04314118e+00)
( 1.31940000e+00,  1.04155603e+00)
( 1.32610000e+00,  1.04000976e+00)
( 1.33280000e+00,  1.03849504e+00)
( 1.33950000e+00,  1.03698638e+00)
( 1.34620000e+00,  1.03546961e+00)
( 1.35280000e+00,  1.03391900e+00)
( 1.35950000e+00,  1.03240017e+00)
( 1.36620000e+00,  1.03093389e+00)
( 1.37290000e+00,  1.02941308e+00)
( 1.37960000e+00,  1.02789127e+00)
( 1.38630000e+00,  1.02636844e+00)
( 1.39300000e+00,  1.02488459e+00)
( 1.39970000e+00,  1.02335973e+00)
( 1.40640000e+00,  1.02183386e+00)
( 1.41300000e+00,  1.02030696e+00)
( 1.41970000e+00,  1.01873903e+00)
( 1.42640000e+00,  1.01721009e+00)
( 1.43310000e+00,  1.01564010e+00)
( 1.43980000e+00,  1.01409505e+00)
( 1.44650000e+00,  1.01255651e+00)
( 1.45320000e+00,  1.01096931e+00)
( 1.45990000e+00,  1.00939513e+00)
( 1.46660000e+00,  1.00776569e+00)
( 1.47320000e+00,  1.00616888e+00)
( 1.47990000e+00,  1.00453720e+00)
( 1.48660000e+00,  1.00290436e+00)
( 1.49330000e+00,  1.00130430e+00)};

\end{groupplot}\end{tikzpicture}